\documentclass[openany,reqno,xxxnumber,index]{amsbook}
\usepackage{math,people,amssymb,amsthm,array,backref,fancyhdr,ifthen,
            pb-diagram,graphics,curves,eepic,
           lamsarrow,pb-lams}

\makeindex
\def\thedate{}
\def\date#1{\global\def\thedate{#1}}
\font\cyr=wncyr8                

\newif\iffast\fastfalse         

\renewcommand\emdef[2][!]{\emph{#2}
  \ifthenelse{\equal{#1}{!}}{\index{#2|textbf}}{\index{#1|textbf}}}
\catcode`"=\active\def"#1"{#1\index{#1}}

\theoremstyle{plain}
\newtheorem{theorem}{Theorem}[chapter]
\newtheorem{lemma}[theorem]{Lemma}
\newtheorem{proposition}[theorem]{Proposition}
\newtheorem{definition}[theorem]{Definition}
\newtheorem{corollary}[theorem]{Corollary}

\newcolumntype{C}{>{$}c<{$}}

\newcommand\Y{{\overline{Y}}}
\newcommand\m{{\overline{m}}}
\newcommand\treeY{{\tree^{(\Y)}}}

\newcommand\LL{{\mathcal L}}

\newcommand\om{{\omega}}
\newcommand\omover{{\overline{\om}}}
\newcommand\Omhat{{\widehat\Omega}}
\newcommand\Omr{{\Omega^{(r)}}}

\newcommand\Stab{\operatorname{\mathsf{St}}}
\newcommand\Rist{\operatorname{\mathsf{Rst}}}
\newcommand\Aut{\operatorname{\mathsf{Aut}}}

\newcommand\Isom{\operatorname{\mathsf{Isom}}}
\newcommand\Sym{\operatorname{\mathsf{Sym}}}
\newcommand\PSL{\operatorname{\mathsf{PSL}}}
\newcommand\SL{\operatorname{\mathsf{SL}}}
\newcommand\GL{\operatorname{\mathsf{GL}}}
\newcommand\Alt{\operatorname{\mathsf{Alt}}}

\newcommand\SPRES[2]{\langle#1|\,#2\rangle}

\newcommand\PRES[2]{\big\langle#1\big|\,#2\big\rangle}
\newcommand\LPRES[4]{\PRES{#1}{#2\big|\,#3\big|\,#4}}
\newcommand\ELPRES[3]{\PRES{#1\big\|\,#2}{#3}}

\def\tarrow#1{\arrow{#1,blue}}
\def\aarrow#1{\arrow{#1,red,..}}
\makeatletter
\def\redvector(#1,#2)#3{\red{\lamsvector(#1,#2){#3}}}
\@namedef{dgo@red}{\let\dg@VECTOR=\redvector}
\def\bluevector(#1,#2)#3{\blue{\lamsvector(#1,#2){#3}}}
\@namedef{dgo@blue}{\let\dg@VECTOR=\bluevector}
\makeatother

\newcommand\FGg{{\Gamma}}       
\newcommand\BGg{{\overline\Gamma}} 
\newcommand\GSg{{\doverline\Gamma}} 
\newcommand\Sg{{\tilde\Gg}}     
\newcommand\GGS{{\textsf{GGS}}} 
\newcommand\GG{{\textsf{G}}}   

\newcommand\gr{{\operatorname{\mathsf{growth}}}}
\newcommand\lcm{\operatorname{\text{lcm}}}
\newcommand\charact{{\operatorname{char}}} 

\begin{document}
\title{Branch Groups}
\author{Laurent Bartholdi}
\address{Section de Math\'ematiques\\
    Universit\'e de Gen\`eve\\
    CP 240, 1211 Gen\`eve 24\\
    Switzerland}
  \email{laurent@math.berkeley.edu}
\author{Rostislav I. Grigorchuk}
\address{Department of Ordinary Differential Equations\\
    Steklov Mathematical Institute\\
    Gubkina Street 8\\
    Moscow 119991\\
    Russia}
  \email{grigorch@mi.ras.ru}
\author{Zoran \v{S}uni\'{k}}
\address{Department of Mathematics\\
    Cornell University\\
    431 Malott Hall\\
    Ithaca, NY 14853\\
    USA}
  \email{sunik@math.cornell.edu}
\date{October 27, 2002}
\maketitle

\tableofcontents
\date{October 27, 2002}
\chapter*{Introduction}
Branch groups were defined only recently although they make
non-explicit appearances in the literature in the past, starting
with the article of \JWilson~\cite{wilson:jig}. Moreover, such
examples as infinitely iterated wreath products or the group of
tree automorphisms $\Aut(\tree)$, where $\tree$ is a regular
rooted tree, go back to the work of \LKaloujnine, \BNeumann\,
\PHall\, and others.

Branch groups were explicitly defined for the first time at the
1997 St-Andrews conference in Bath in a talk by the second author.
Immediately, this sparked a great interest among group theorists,
who started investigating numerous properties of branch groups
(see~\cite{grigorchuk:jibg,bartholdi-g:spectrum,bartholdi-g:parabolic,grigorchuk-w:conj})
as well as \JWilson's classification of just-infinite groups.

There are two new approaches to the definition of a branch group,
given in~\cite{grigorchuk:jibg}. The first one is purely
algebraic, defining branch groups as groups whose lattice of
subnormal subgroups is similar to the structure of a spherically
homogeneous rooted tree. The second one is based on a geometric
point of view according to which branch groups are groups acting
spherically transitively on a spherically homogeneous rooted tree
and having structure of subnormal subgroups similar to the
corresponding structure in the full group $\Aut(\tree)$ of
automorphisms of the tree.

Until 1980 no examples of finitely generated branch groups were
known and the first such examples were constructed
in~\cite{grigorchuk:burnside}. These examples are usually referred
to as the first and the second Grigorchuk groups, following Pierre
de la Harpe (\cite{harpe:ggt}). Other examples soon appeared
in~\cite{grigorchuk:growth,grigorchuk:gdegree,grigorchuk:pgps,gupta-s:burnside,gupta-s:infinitep,gupta-s:ext,neumann:pride}
and these examples are the basic examples of branch groups, the
study of which continues at the present time. Let us mention that
the  examples of \SAlioshin~\cite{aleshin:burnside} and
\VSushchansky~\cite{sushchansky:pgps} that appeared earlier also
belong to the class of finitely generated branch groups, but the
methods used in the study the groups from~\cite{aleshin:burnside}
and~\cite{sushchansky:pgps} did not allow the discovery of the
branch structure and this was done much later.

Already in~\cite{grigorchuk:burnside} the main features of a general
method which works for almost any finitely generated branch group had
appeared: one considers the stabilizer of a vertex on the first level
and projects it on the corresponding subtree. Then either this
projection is equal to the initial group and thus one gets a
self-similarity property, or otherwise one gets a finite or infinite
chain of branch groups related by some homomorphisms. One of the
essential properties of this chain is that these homomorphisms satisfy
a ``Lipschitz'' property of norm reduction, which lends itself to
arguments using direct induction on length in the case of a
self-similar group, or simultaneous induction on length for all groups
in a chain.

Before we give more information of a historical character and briefly
describe the main directions of investigation and the main results in the
area, let us explain why the class of branch groups is important. There is
a lot of evidence that this is indeed the case. For example
\begin{enumerate}
\item The class of branch groups is one of the classes into which
  the class of just-infinite groups naturally splits (just-infinite
  groups are groups whose proper quotients are finite).
\item The class contains groups with many extraordinary properties,
  like infinite finitely generated torsion groups, groups of intermediate
  growth, amenable but not elementary amenable groups, groups of finite
  width, etc.
\item Branch groups have many applications and are related to
  analysis, geometry, combinatorics, probability, computer science, etc.
\item They are relatively easy to handle and usually the proofs even
  of deep theorems are short and do not require special
  techniques. Therefore the branch groups constitute an easy-to-study
  class of groups, whose basic examples have already appeared in many
  textbooks and lecture notes, for example
  in~\cite{kargapolov-merzlyakov:otg,baumslag:cgt,robinson:ctg,harpe:ggt}.
\end{enumerate}

This survey article deals almost exclusively with abstract branch
groups. The theory of profinite branch groups is also being
actively developed at present
(see~\cite{grigorchuk-h-z:profinite,grigorchuk:jibg,wilson:segal}),
but we hardly touch on this subject.

The survey does not pretend to be complete. There are several
topics that we did not include in the text due to the lack of
space and time. Among them, we mention the results of \SSidki\
from~\cite{sidki:primitive} on thin algebras associated to
Gupta-Sidki groups, the results on automorphisms of branch groups
of \SSidki\ from~\cite{sidki:subgroups} and the recent results of
Lavreniuk and Nekrashevich from~\cite{lavreniuk-n:rigidity}, the
results of \CRover\ on embeddings of Grigorchuk groups into
finitely presented simple groups and on abstract commensurators
from~\cite{rover:fpsg,rover:commensurator}, the results of
\VSushchansky\ on factorizations based on the use of torsion
branch
groups~\cite{sushchansky:pfactorizable,sushchansky:factorizable},
the results of B. Fine, A. Gaglione, A. Myasnikov and D. Spellman
from~\cite{fine-g-m-s:discriminating} on discriminating groups,
etc.

\section{Just-infinite groups}
Let $\mathcal P$ be any property which is preserved under homomorphic
images (we call such a property an \emdef[$H$-property]{$\mathcal
  H$-property}).  Any infinite finitely generated group can be mapped
onto a just-infinite group (see~\cite{grigorchuk:jibg,harpe:ggt}),
so if there is an infinite finitely generated group with the
$\mathcal H$-property $\mathcal P$ then there is a just-infinite
finitely generated group with the same property. Among the
$\mathcal H$-properties let us mention the property of being a
torsion group, not containing the free group $F_2$ on two elements
as a subgroup, having subexponential growth, being amenable,
satisfying a given identity, having bounded generation, finite
width, trivial space of pseudocharacters (for a relation to
bounded cohomology see~\cite{grigorchuk:bounded}), only
finite-index maximal subgroups, $T$-property of Kazhdan etc.

The branch just-infinite groups are precisely the just-infinite
groups whose structure lattice of subnormal subgroups (with some
identifications) is isomorphic to the lattice of closed and open
subsets of a Cantor set. This is the approach of \JWilson\
from~\cite{wilson:jig}.

In that paper, \JWilson\ split the class of just-infinite groups
into two subclasses -- the groups with finite and the groups with
infinite structure lattice. The dichotomy of \JWilson\ can be
reformulated (see~\cite{grigorchuk:jibg}) in the form of a
trichotomy according to which any finitely generated just-infinite
group is either a branch group or can easily be constructed from a
simple group or from a hereditarily just-infinite group (i.e., a
residually finite group all of whose subgroups of finite index are
just-infinite).

Therefore the study of finitely generated just-infinite groups
naturally splits into the study of branch groups, infinite simple
groups and hereditarily just-infinite groups. Unfortunately, at
the moment, none of these classes of groups are well understood,
but we have several (classes of) examples.

There are several examples and constructions of finitely generated
infinite simple groups, probably starting with the example of
\GHigman~in~\cite{higman:fgisg}, followed by the finitely presented
example of \RThompson, generalized by \GHigman~in~\cite{higman:fpisg} (see
also the survey~\cite{cannon-f-p:thompson} and~\cite{brown:finiteness}),
the constructions of different monsters by \AOlshanskii\
(see~\cite{olshansky:gdr}), as well as by \SAdyan\ and \ILysionok\
in~\cite{adyan-l:monsters}, and more recently some finitely presented
examples by \CRover\ in~\cite{rover:fpsg}. The $\mathcal H$-properties
that can be satisfied by such groups are, for instance, the Burnside
identity $x^p$, for large prime $p$, and triviality of the space of
pseudocharacters. The latter holds for the simple groups $T$ and $V$ of
\RThompson~(this follows from the results on finiteness of commutator
length, see~\cite{ghys-s:thompson,brown:finiteness}).

All known hereditarily just-infinite groups (like the projective groups
$\PSL(n,\Z)$ for $n>2$) are linear (in the profinite case there are extra
examples like Nottingham group), so by the alternative of Tits they
contain $F_2$ as a subgroup and therefore cannot be amenable, of
intermediate growth, torsion etc. However, they can have bounded
generation: it is shown in~\cite{carter-k:boundedg} that this holds for
$\SL(n,\Z)$, $n>2$, and therefore also for $\PSL(n,\Z)$.

It seems that there are fewer constraints in the class of branch
groups and that they can have various $\mathcal H$-properties,
some of which are listed below. It is conjectured that many of
these properties do not hold for groups from the other two
classes. On the other hand, branch groups cannot satisfy
nontrivial identities (see~\cite{leonov:identities}
and~\cite{wilson:segal} where the proof is given for the
just-infinite case).

\section{Algorithmic aspects}
Branch groups have good algorithmic properties. In the branch groups of
\GG\ or \GGS\ type (or more generally spinal type groups) the word problem
is solvable by an universal branch algorithm described
in~\cite{grigorchuk:gdegree}. This algorithm is very fast and requires a
minimal amount of memory.

The conjugacy problem was unsettled for a long time, and it was solved for
the basic examples of branch groups just recently. The
article~\cite{wilson-z:conj} solves the problem for regular branch
$p$-groups, where $p$ is an odd prime, and the argument uses the property
of ``conjugacy separability'' as well as profinite group machinery.
In~\cite{leonov:conj} and~\cite{rozhkov:conj} a different approach was
used, which also works in case $p=2$. This ideas were developed
in~\cite{grigorchuk-w:conj} in different directions. For instance, it was
shown that, under certain conditions, the conjugacy problem is solvable
for all subgroups of finite index in a given branch group (we mention here
that the property of solvability of conjugacy problem, in contrast with
the word problem, is not preserved when one passes to subgroups of finite
index). Still, we are far from understanding if the conjugacy problem is
solvable in all branch groups with solvable word problem.

The isomorphism problem was also considered in~\cite{grigorchuk:gdegree}
where it is proven that each of the uncountably many constructed groups
$G_\om$ is isomorphic to at most countably many of them, thus showing that
the construction gives uncountably many non-isomorphic examples. It would
be very interesting to distinguish all these examples.

Branch groups are related to groups of finite automata. A brief
account is given in Section~\ref{sec:automata} (see
also~\cite{grigorchuk-n-s:automata}). Every group generated by
finite automata has a solvable word problem. It is unclear if
every such group has solvable conjugacy problem. On the other
hand, it seems that the isomorphism problem cannot be solved in
this particular case. Indeed, according to the results
in~\cite{klarner-b-s:undecidability}, the freeness of a matrix
group with integer entries cannot be determined, and the general
linear group $\GL(n,\Z)$ can be embedded in the group of automata
defined over an alphabet on $2^n$ letters as shown by A. Brunner
and \SSidki\ (see~\cite{brunner-s:glnz}).

\section{Group presentations}
In Chapter~\ref{chapter:presentation} we study presentations of branch
groups by generators and relations. It seems probable that no branch group
is finitely presented. However, the regular branch groups have nice
recursive presentations called $L$-presentations. The first such
presentation was found for the first Grigorchuk group by \ILysionok\
in~\cite{lysionok:pres}. Shortly afterwards, \SSidki\ devised a general
method yielding recursive definitions of such groups, and applied it to
the Gupta-Sidki group~\cite{sidki:pres}.
In~\cite{grigorchuk:notEG,grigorchuk:bath} the idea and the result of
\ILysionok\ were developed in different directions.
In~\cite{grigorchuk:bath} it was proven that the \ILysionok\ system of
relations is minimal and the Schur multiplier of the group was computed:
it is $(\Z/2)^\infty$. Thus the second homology group of the first
Grigorchuk group is infinite dimensional. In~\cite{grigorchuk:bath} it was
indicated that the Gupta-Sidki $p$-groups also have finite
$L$-presentations. On the other hand, it was shown
in~\cite{grigorchuk:notEG} how $L$-presentations can be used to embed some
branch groups into finitely presented groups using just one \textsf{HNN}
extension. The important feature of this embedding is that it preserves
the amenability. The first examples of finitely presented, amenable but
not elementary amenable groups were constructed this way, thus providing
new examples of good fundamental groups (in the terminology of \Freedmann\
and \Teichler~\cite{freedman-t:subexponential}).

In~\cite{bartholdi:lpres} the notion of an $L$-presentation was
slightly extended to the notion of an endomorphic presentation in
a way that allowed to show that a finitely generated, fractal,
regular branch group satisfying some natural extra conditions has
a finite endomorphic presentation. A number of concrete $L$ and
endomorphic presentations of branch groups appear in the article
along with general facts on such presentations.

As was stated, no known branch group has finite presentation. For
the first Grigorchuk group this was already mentioned
in~\cite{grigorchuk:burnside} with a sketch of a proof that was
given completely in~\cite{grigorchuk:bath}.
In~\cite{grigorchuk:gdegree} two other proofs were presented. Yet
another approach in proving the absence of finite presentations is
used by \NGupta\ in~\cite{gupta:recursive}. More on the history of
the presentation problem and related methods appears in
Section~\ref{sec:presentation:hist}.

\section{Burnside groups}
The third part of the survey is devoted to the algebraic properties of
branch groups in general and of the most important examples.

The first examples of branch groups appeared
in~\cite{grigorchuk:burnside} as examples of infinite finitely
generated torsion groups. Thus the branch groups are related to
the Burnside Problem on torsion groups. This difficult problem has
three branches: the Unbounded Burnside Problem, the Bounded
Burnside Problem and the Restricted Burnside problem
(see~\cite{adyan:burnside,kostrikin:around,grigorchuk-l:burnside})
and, in one way or another, all of them are solved. However, there
is still a series of unsolved problems in the neighborhood of the
Burnside problem and which are very important to the theory of
groups --- to name one, ``is there a finitely presented infinite
torsion group?''. The first example of an infinite finitely
generated torsion group, which provided a negative answer to the
General Burnside Problem, was constructed by \Golod\
in~\cite{golod:nil} and it was based on the Golod-Shafarevich
Theorem.  The actual problem of constructing simple examples which
do not require the use of such deep results as Golod-Shafarevich
Theorem remained open until such examples appeared
in~\cite{grigorchuk:burnside}. Soon, more examples appeared
in~\cite{grigorchuk:growth,grigorchuk:gdegree,grigorchuk:pgps,gupta-s:burnside,gupta-s:infinitep,gupta-s:ext}
and more recently
in~\cite{bartholdi-s:wpg,grigorchuk:jibg,bartholdi:ggs,sunik:phd}.
The early examples are finitely generated infinite $p$-groups, for
$p$ a prime, and the latter papers contain interesting examples
that are not $p$-groups.

We already mentioned the idea of using induction on word length,
based on the fact that the projections on coordinates decrease the
length. In conjunction, the idea of fixing larger and larger
layers of the tree under taking powers was developing. The
stabilization occurs in the first Grigorchuk group
from~\cite{grigorchuk:burnside} after three steps, and for the
second example from the same article it occurs after the second
step of taking $p$-th powers. Using a slightly modified metric on
the group~\cite{bartholdi:upperbd}, the stabilization can be made
to appear after just one step; this is extensively developed
in~\cite{bartholdi-s:wpg}. Examples with strong stabilization
properties for the standard word metric are constructed
in~\cite{gupta-s:burnside}, where stabilization takes place after
the first step. The notion of \emdef{depth} of an element, i.e.,
the number of decompositions one must perform to decrease the
length down to $1$ was introduced by \SSidki\
in~\cite{sidki:subgroups}, and is very useful in some situations.

One of the important principles of modern group theory is to try to
develop asymptotic methods related to growth, amenability and other
asymptotic notions. In~\cite{grigorchuk:gdegree} the torsion growth
functions were introduced for finitely generated torsion groups and it was
shown that some examples
from~\cite{grigorchuk:burnside,grigorchuk:growth,grigorchuk:gdegree} have
polynomial growth in this sense. These results were improved in many
directions
in~\cite{lysionok:order,leonov:periodbd,leonov:estimation,bartholdi-s:wpg}
and some of them are described in Chapter~\ref{chapter:torsion}.

Among the main consequences of the theory of \EZelmanov\
(see~\cite{zelmanov:odd,zelmanov:2}) is that if a finitely
generated torsion residually finite group has finite exponent
(i.e., there exists $n\neq 0$ such that $g^n=1$ for every element
$g$ in the group) then the group is finite. Although the results
of \EZelmanov\ do not depend on the classification of finite
simple groups, the above mentioned consequence does and it would
be nice to produce a proof that is independent of the
classification. A simple proof that finitely generated torsion
branch groups always have infinite exponent is provided
(Theorem~\ref{theorem:branch-ie}).

The profinite completion of a group of finite exponent is a torsion
profinite group. If it has a just-infinite quotient (and this is the case
if the group $G$ is a virtually pro-$p$-group) then one gets a profinite
just-infinite torsion group of bounded exponent. By Wilson's alternative
such a group is either just-infinite branch, or hereditarily
just-infinite; but by the results from~\cite{grigorchuk-h-z:profinite} if
it were branch it would have unbounded exponent, so the search for
profinite groups of bounded exponent can be narrowed to hereditarily
just-infinite groups. We believe such groups do not exist.

In Chapter~\ref{chapter:torsion} we give a simple proof of the fact that a
finitely generated torsion branch group has infinite exponent. An
interesting question is to describe the type of torsion growth that
distinguishes the finite and infinite groups (it follows from the results
of \EZelmanov\ that there exists a recursive, unbounded function $z_k(n)$,
depending on the number of generators $k$, such that every torsion group
on $k$ generators whose torsion growth is bounded above by $z_k(n)$ is
finite). It seems likely that the problem can be reduced to the case of
branch torsion groups.

In Chapter~\ref{chapter:torsion} we also analyze carefully the
idea of the first example in~\cite{grigorchuk:burnside} which is
not related to the stabilization, but rather uses a covering of a
group by kernels of homomorphisms. The torsion groups (called \GG\
groups) that generalize the constructions
of~\cite{grigorchuk:gdegree,grigorchuk:pgps} are investigated in
greater detail in~\cite{bartholdi-s:wpg,sunik:phd}. We exhibit the
construction and some interesting examples, based on existence of
finite groups with certain required properties (\DHolt's example,
for instance).

Until recently, it was not known whether there exist torsion-free
just-infinite branch groups. Such an example was constructed
in~\cite{bartholdi-g:parabolic}.

\section{Subgroups of branch groups}
The study of the subgroup structure of any class of groups is an
important part of the investigation. Branch groups have a rich and
nice subgroup structure which has not yet been completely
investigated. In the early works attention was paid to some
particular subgroups of small index such as the stabilizers of the
first few levels, the initial members of the lower central series,
derived series, etc. A fundamental observation made by \NGupta\
and \SSidki\ in~\cite{gupta-s:ext} is that many \GGS\ groups
contain a normal subgroup K of finite index with the property that
$K$ contains \emdef{geometrically} $K^m$ ($m$ is the degree of the
tree and ``geometrically'' means that the product $K^m$ acts on
the subtrees on the first level) as a subgroup of finite index. It
so happens that all the main examples of branch groups have such a
subgroup and this fact lies at the base of the definition of
regular branch groups.

Rigid stabilizers and stabilizers of the first Grigorchuk group are
described in~\cite{rozhkov:stab,bartholdi-g:parabolic}. For the
Gupta-Sidki $p$-groups this is done in~\cite{sidki:subgroups}. The
structure of normal subgroups of \PNeumann's groups~\cite{neumann:pride}
is very simple, since the normal subgroups coincide with the (rigid)
stabilizers. For branch $p$-groups it is more difficult to obtain the
structure of the lattice of normal subgroups (this is related to the fact
that pro-$p$-groups usually have rich structure of subgroups of finite
index). The lattice of normal subgroups of the first Grigorchuk group was
recently described by the first author in~\cite{bartholdi:lcs} (see
also~\cite{tcs-s-t:normal} where the normal subgroups are described up to
the fourth level).

In the study of infinite finitely generated groups an important
role is played by the maximal and weakly maximal subgroups (i.e.,
subgroups of infinite index maximal with respect to this
property). It is strange that little attention was paid to the
latter until recently. The main result of
\EPervova~\cite{pervova:edsubgroups} claims that in the basic
examples of branch groups every proper maximal subgroup has finite
index. This is in contrast with lattices in semisimple Lie groups,
as follows from results of \Margulis\ and
\Soifer~\cite{margulis-s:maximal}.

Important examples of weakly maximal subgroups are the parabolic
subgroups, i.e., the stabilizers of infinite paths in the tree.
The structure of parabolic subgroups is described
in~\cite{bartholdi-g:parabolic} for some particular examples. They
are not finitely generated and have a tree-like structure. It
would be interesting to obtain a complete description of weakly
maximal subgroups in branch groups.

\section{Lie algebras}
Chapter~\ref{chapter:lie} deals with central series and associated
Lie algebras of branch groups. There is a canonical way, due to
\WMagnus\, in which a central series corresponds to a graded Lie
ring or Lie algebra. The most interesting central series are the
lower central series and the series of dimension subgroups. It was
proved in~\cite{grigorchuk:hp} that the Cesar\`o averages of the
ranks of the factors in the lower central series of the first
Grigorchuk group (which are elementary $2$-groups) are uniformly
bounded and it was conjectured that the ranks themselves were
uniformly bounded, i.e., that the first Grigorchuk group has
finite width. An important step in proving this conjecture was
made in~\cite{rozhkov:lcs}, and a complete proof appears
in~\cite{bartholdi-g:lie}, using ideas of
\LKaloujnine~\cite{kaloujnine:struct} and the notion of uniserial
module. Moreover, a negative answer to a problem of \EZelmanov\ on
the classification of just-infinite profinite groups of finite
width is provided in~\cite{bartholdi-g:lie}, a new example of a
group of finite width was constructed and the structure of the
Cayley graph of the associated Lie algebras was described. This is
one of the few cases of a nontrivial computation of a Cayley (or
Lie) graph of a graded Lie algebra.

The question of the finiteness of width of other basic branch groups
(first of all the Gupta-Sidki groups) was open for a long time and
recently answered negatively by the first author
in~\cite{bartholdi:lcs}.  These results are also presented in
Chapter~\ref{chapter:lie}.

An important role in the study of profinite completions of branch
group is played by the ``congruence subgroup property'' with
respect to the sequence of stabilizers, meaning that every
finite-index subgroup contains a level stabilizer $\Stab_G(n)$ for
some $n$, and which holds for many branch groups.  Nevertheless,
there are branch groups without this property and the complete
solution of the congruence subgroup property problem for the class
of all branch groups is not completely resolved.

\section{Growth}
The fourth part of the paper deals with some geometric and analytic
properties of branch groups. The main notion in asymptotic group theory is
the notion of \emph{growth} of a finitely generated group. The growth
function $\gamma(n)$ of a finitely generated group $G$ with respect to a
system of generators $S$ counts the number of group elements of length at
most $n$.  The group's type of growth --- exponential, intermediate,
polynomial
--- does not depend on the choice of $S$. One can easyly construct an
example of a group of polynomial growth of any given degree $d$ (for
instance $\Z^d$) or a group of exponential growth (like $F_2$, the free
group on two generators) but it is a highly non-trivial task to construct
a group of intermediate growth. The question of existence of such groups
of intermediate growth was posed by \JMilnor~\cite{milnor:5603} and solved
fifteen years later
in~\cite{grigorchuk:growth,grigorchuk:gdegree,grigorchuk:pgps}, were the
second author shows that the first group in~\cite{grigorchuk:burnside} and
all $p$-groups $G_\omega$ in~\cite{grigorchuk:gdegree,grigorchuk:pgps}
have intermediate growth; the estimates are of the form
\begin{equation}\label{eq:intro:interm}
  e^{\sqrt n}\precsim \gamma(n)\precsim e^{n^\beta},
\end{equation}
for some $\beta<1$, where for two functions $f,g:\N\to\N$ we write
$f\precsim g$ to mean that there exists a constant $C$ with $f(n)
\leq g(Cn)$ for all $n\in\N$.

Milnor's problem was therefore solved using branch groups. Up to the
present time all known groups of intermediate growth are either branch
groups or groups constructed using branch groups and we believe that all
just-infinite groups of intermediate growth are branch groups.

By using branch groups, the second author showed
in~\cite{grigorchuk:gdegree} that there exist uncountable chains
and anti-chains of intermediate growth functions of groups acting
on trees.

The upper bound in~(\ref{eq:intro:interm}) was improved
in~\cite{bartholdi:upperbd}, and a general improvement of the
upper bounds for all groups $G_\omega$ was given
in~\cite{bartholdi-s:wpg} and~\cite{muchnik-p:growth}.

One of the main remaining question on growth is whether there exists a
group with growth precisely $e^{\sqrt n}$. It is known that if a group is
residually nilpotent and its growth is strictly less than $e^{\sqrt
  n}$, then the group is virtually nilpotent and therefore has
polynomial growth~\cite{grigorchuk:hp,lubotzky-m:resfinite}. Using
arguments given in~\cite{grigorchuk:hp} (see
also~\cite{bartholdi-g:lie}), it follows that if a group of growth
$e^{\sqrt n}$ exists in the class of residually-$p$ groups, then
it must have finite width. For some time, among all known examples
of groups of intermediate growth only the first Grigorchuk group
was known to have finite width, and the second author conjectured
that this group has precisely this type of growth.  However, his
conjecture was infirmed by \YLeonov~\cite{leonov:lowerbd} and the
first author~\cite{bartholdi:lowerbd}; indeed the growth of the
first Grigorchuk group is bounded below by $e^{n^\alpha}$ for some
$\alpha>\frac12$.

The notion of growth can be defined for other algebraic and
geometric objects as well: algebras, graphs, etc. A very
interesting topic is the study of the growth of graded Lie
algebras $\Lie(G)$ associated to groups. In case of \GGS\ groups
some progress is achieved in~\cite{bartholdi:lcs}, where the
growth of $\Lie(G)$ for the Gupta-Sidki $3$-group and some other
groups is computed; in particular, it is shown that the
Gupta-Sidki group does not have finite width. A connection between
the Lie algebra structure and the tree structure is used in the
majoration of the growth of the associated Lie algebra by the
growth of any homogeneous space $G/P$, where $P$ is a parabolic
subgroup, i.e., the stabilizer of an infinite path in a
tree~\cite{bartholdi:lcs}. As metric spaces, these homogeneous
spaces are equivalent to Schreier graphs. These graphs have an
interesting structure: they are substitutional graphs, and have a
fractal behavior in the case of many fractal branch groups. They
have polynomial growth, usually of non-integral degree. These
results are presented in Section~\ref{sec:pspace}.

One of the promising directions of research is the study of
spectral properties of the above graphs. This question is linked
to several famous problems of operator $K$-theory and theory of
$C^*$-algebras. One of the first works in this direction
is~\cite{bartholdi-g:spectrum}, where it is shown that the
spectrum of the discrete Laplace operator on such graphs can be a
Cantor set, optionally with extra isolated points.  The
computation of these spectra is related to operator recursions
that hold for the Laplace or Hecke-type operators associated to
the dynamical system $(G,\partial\tree,\mu)$, where
$(\partial\tree,\mu)$ is the boundary of the tree endowed with the
uniform measure.  The main results are presented in
Chapter~\ref{chapter:spectrum}.

Finally, there are a great number of open questions on branch groups. Some
of them are listed in the final part of the paper; we hope that they will
stimulate the development of the subject.

\section {Acknowledgments}
Thanks to \DSegal\ and \PdlHarpe\ for valuable suggestions and to
\CRover\ for careful reading of the text and numerous corrections
and clarifications.

The second author would like to thank the Royal Society of the
United Kingdom, University of Birmingham and \JWilson, and also
acknowledge the support from a Russian grant ``Leading Schools of
the Russian Federation'', project N~00-15-96107.

The third author thanks to Fernando Guzm\'an, Laurent Saloff-Coste
and Christophe Pittet, for their advice, help and valuable
conversations.
\section {Some notation}
We include zero in the set of natural numbers
$\N=\{0,1,2,\dots\}$. The set of positive integers is denoted by
$\N_+=\{1,2,\dots\}$.

Expressions as $((1,2)(3,4,5))$ listing non-trivial cycles are used to
describe permutations in the group of symmetries $\Sym(n)$ of $n$ elements

We want all the group actions to be on the right. Thus we conjugate as
follows,
\[ g^h=h^{-1}gh, \]
and we denote
\[[g,h]=g^{-1}h^{-1}gh=g^{-1}g^h. \]
The commutator subgroup of $G$ is denoted by $[G,G]$ and the
abelianization $G/[G,G]$ by $G^{ab}$.

Let the group $A$ act on the right on the group $H$ through $\alpha: A
\rightarrow \Aut(H)$. We define the \emdef{semidirect product} $G=H
\rtimes_\alpha A$ as the group whose elements are the ordered pairs from
the set $H \times A$ and the operation is given by
\[ (h,a)(g,b) = (hg^{(a^{-1})\alpha},ab). \]
After the identification $(h,1)=h$ and $(1,a)=a$ we see that $G$
is a group containing $H$ as a normal subgroup and $A$ as its
complement, i.e., $HA=G$ and $H \cap A=1$. Moreover, the
conjugation of the elements in $H$ by the elements in $A$ is given
by the action $\alpha$, i.e., $h^a = h^{(a)\alpha}$.

If we start with a group $G$ that has a normal subgroup $H$ with a
complement $A$ in $G$, we say that $G$ is the \emdef{internal semidirect
product} of $H$ and $A$. Indeed, $G=H \rtimes A$ where the action of $A$
on $H$ is through conjugation (note that $hagb=hg^{a^{-1}}ab$ for $h,g\in
H$ and $a,b\in A$).

Let $G$ and $A$ be groups acting on the set $X$ and the finite set $Y$,
respectively. We define the \emdef{permutational wreath product} $G \wr_Y
A$ that acts on the set $Y \times X$ (note the change in the order) as
follows: let $A$ act on the direct power $G^Y$ on the right by permuting
the coordinates of $G^Y$ by
\begin{equation}\label{h^a}
(h^a)_y = h_{y^{a^{-1}}},
\end{equation}
for $h\in G^Y$, $a \in A$, $y \in Y$; then define $G \wr_Y A$ as the
semidirect product
\[ G \wr_Y A = G^Y \rtimes A \]
obtained through the action of $A$ on $G^Y$; finally let the wreath
product act on the right on the set $Y \times X$ by
\[ (y,x)^{ha} = (y^a,x^{h_y}), \]
for $y \in Y$, $x \in X$, $h \in G^Y$ and $a \in A$. Note that the
equality (\ref{h^a}) which represents the action of $A$ on $G^Y$, also
represents conjugation in the wreath product, exactly as we want, and that
this wreath product is associative, modulo the necessary natural
identifications.

All actions defined by now were right actions. However, we
achieved this by introducing inversion at several crucial places,
thus introducing left actions through the back door. Another
possibility was to let the semidirect product of $A$ and $H$ as
above be the group whose elements the ordered pairs in $A \times
H$ and define $(a,h)(b,g)=(ab,h^{(b)\alpha}g)$. This works well,
but we choose not to do it.

We introduce here the basic notation of growth series. Growth series will
be used in Chapters~\ref{chapter:lie} and~\ref{chapter:growth}.

Let $X$ be a set on which the group $G$ acts, and fix a base point $*\in
X$ and a set $S$ that generates $G$ as a monoid. The \emdef{growth
function} of $X$ is
\[\gamma_{*,G}(n)=|\setsuch{x\in X}{x=*^{s_1\dots s_n}\text{ for some
    }s_i\in S}|.\]
The \emdef{growth series} of $X$ is
\[\gr(X)=\sum_{n\ge0}\gamma_{*,G}(n)\hbar^n.\]

Let $V=\bigoplus_{n\ge0}V_n$ be a graded vector space. The
\emdef{Hilbert-Poincar\'e} series of $V$ is the formal power series
\[\gr(V)=\sum_{n\ge0}\dim V_n\hbar^n.\]

A preorder $\precsim$ is defined on the set of non-decreasing functions
$\R_{\geq 0} \rightarrow \R_{\geq 0}$ by $f \precsim g$ if there exists a
positive constant $C$ such that $f(n)\leq g(Cn)$, for all $n$ in $\R_{\geq
0}$.  An equivalence relation $\sim$ is defined by $f \sim g$ if $f
\precsim g$ and $g \precsim f$.

Several branch groups are distinguished enough to be given
separate notation. They are the first Grigorhuk group $\Gg$, the
Grigorchuk supergroup $\Sg$, the Gupta-Sidki 3-group $\GSg$, the
Fabrykowski-Gupta group $\FGg$ and the Bartholdi-Grigorchuk group
$\BGg$. See Section~\ref{sec:branch:examples} for the definitions.

\part{Basic Definitions and Examples}
\date{Nov 3, 2002}
\chapter{Branch Groups and Spherically Homogeneous Trees}
\section{Algebraic definition of a branch group}
We start with the main definition of the survey, namely the definition
of a branch group. The definition is given in purely algebraic terms,
emphasizing the subgroup structure\index{subgroup!structure} of the
groups. We give a geometric version of the definition in
Section~\ref{sec:gdef} in terms of actions on rooted trees. The two
definitions are not equivalent and we will say something about the
difference later.

\begin{definition}\label{defn:branch}
  Let $G$ be a group. We say that $G$ is a \emdef[group!branch]{branch
    group} if there exist two decreasing sequences of subgroups
  $(L_i)_{i\in\N}$ and $(H_i)_{i\in\N}$ and a sequence of integers
  $(k_i)_{i\in\N}$ such that $L_0=H_0=G$, $k_0=1$,
  \[\bigcap_{i\in\N} H_i = 1\]
  and, for each $i$,
  \begin{enumerate}
  \item\label{defn:branch:1} $H_i$ is a normal subgroup of $G$ of
    finite index.
  \item $H_i$ is a direct product of $k_i$ copies of the subgroup
    $L_i$, i.e., there are subgroups $L_i^{(1)}, \dots L_i^{(k_i)}$ of
    $G$ such that
    \begin{equation} \label{eq:H_i=}
      H_i = L_i^{(1)} \times \dots \times L_i^{(k_i)}
    \end{equation}
    and each of the factors is isomorphic to $L_i$.
  \item $k_i$ properly divides $k_{i+1}$, i.e., $m_{i+1}=k_{i+1}/k_i
    \geq 2$, and the product decomposition~(\ref{eq:H_i=}) of $H_{i+1}$
    refines the product decomposition~(\ref{eq:H_i=}) of $H_i$ in the
    sense that each factor $L_i^{(j)}$ of $H_i$ contains $m_{i+1}$ of
    the factors of $H_{i+1}$, namely the factors $L_{i+1}^{(\ell)}$
    for $\ell=(j-1)m_{i+1}+1, \dots, jm_{i+1}$.
  \item\label{defn:branch:4} conjugations by the elements in $G$
    transitively permute the factors in the product
    decomposition~(\ref{eq:H_i=}).
  \end{enumerate}
\end{definition}

The definition implies that branch groups are infinite, but residually
finite groups. Note that the subgroups $L_i$ are not normal, but they are
subnormal of defect $2$.

\begin{definition}
  Let $G$ be a branch group. Keeping the notation from the previous
  definition, we call the sequence of pairs $(L_i,H_i)_{i\in\N}$ a
  \emdef{branch structure} on $G$.
\end{definition}

\begin{figure}[!ht]
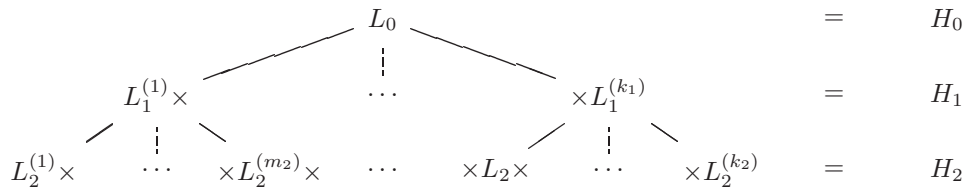

\[\begin{diagram}
  \dgARROWLENGTH=0em
  \node[4]{L_0}\arrow{wsw,-}\arrow{s,-,..}\arrow{ese,-} \node[4]{=} \node{H_0} \\
  \node[2]{L_1^{(1)}\times}\arrow{sw,-}\arrow{s,-,..}\arrow{se,-}
        \node[2]{\dots}
        \node[2]{\times L_1^{(k_{1})}}\arrow{sw,-}\arrow{s,-,..}\arrow{se,-}
        \node[2]{=} \node{H_1} \\
  \node{L_2^{(1)}\times} \node{\dots} \node{\times L_2^{(m_{2})}\times} \node{\dots}
        \node{\times L_2 \times} \node{\dots} \node{\times L_2^{(k_{2})}}
        \node{=} \node{H_2}
\end{diagram}\]
\caption{Branch structure of a branch group} \label{figure:Ls}
\end{figure}

The branch\index{structure!branch} structure of a branch group is
depicted in Figure~\ref{figure:Ls}.
The branch structure on a group $G$ is not unique, since any subsequence
of pairs $(L_{i_j},H_{i_j})_{j=1}^\infty$ is also a branch structure on
$G$.

One can quickly construct examples of branch groups, using infinitely
iterated wreath products. For example, let $p$ be a prime and $\Z/p\Z$ act
on the set $Y=\{1,\dots,p\}$ by cyclic permutations. Define the
permutational wreath product
\[ G_n = \underbrace{((\Z/p\Z \wr_Y \dots) \wr_Y \Z/p\Z )\wr_Y \Z/p\Z}_n, \]
and let $G$ be the inverse limit $\varprojlim G_n$, where the
projections from $G_n$ to $G_{n-1}$ are just the natural
restrictions. Since $G = G \wr_Y \Z/p\Z$, $G$ is a branch group.

Similarly, for $m \geq 2$, let $\Sym(m)$ be the group of permutations of
$Y=\{1,\dots,m\}$, define $G_n$ as the permutational wreath product
\[ G_n = \underbrace{((\Sym(m) \wr_Y \dots ) \wr_Y \Sym(m)) \wr_Y \Sym(m)}_n, \]
and $G$ as the inverse limit $\varprojlim G_n$. Since $G = G \wr_Y
\Sym(m)$, $G$ is a branch group.

In the next section we will look at the last group from a geometric point
of view. We will also develop some terminology for groups acting on rooted
trees that will be used for the second, more geometric,  definition of
branch groups.

\section{Spherically homogeneous rooted trees}
We will define the notion of a spherically homogeneous tree as a set of
words ordered by the prefix relation and then make a connection to the
graph-theoretical version of the same notion. We find it useful to live in
both worlds and use their terminology and notation.

\subsection{The trees}
Let
\[\m = m_1, m_1, m_3, \dots\]
be a sequence of integers with $m_i \geq 2$ and let
\[\Y = Y_1, Y_2, Y_3, \dots\]
be a sequence of alphabets with $|Y_i|=m_i$. A \emdef[word]{word of
  length} $n$ over $\Y$ is any sequence of letters of the form
$w=y_1y_2\dots y_n$ where $y_i \in Y_i$ for all $i$. The unique
word of length 0, the \emph{empty word}, is denoted by
$\emptyset$. The length of the word $u$ is denoted by $|u|$.
Denote the set of words over $\Y$ by $\Y^*$. We introduce a
partial order on the set of all words over $\Y$ by the
\emph{prefix relation} $\leq$. Namely, $u \leq v$ if $u$ is an
initial segment of the sequence $v$, i.e., if $u=u_1\dots u_n$,
$v=v_1\dots v_k$, $n \leq k$, and $u_i=v_i$, for
$i\in\{1,\dots,n\}$. The partially ordered set of words over $\Y$,
denoted by $\treeY$, is called the \emdef[tree!spherically
homogeneous]{spherically homogeneous tree} over $\Y$. The sequence
$\m$ is the \emdef[branch!branching indices]{sequence of branching
  indices} of the tree $\treeY$. If there is no room for confusion we
denote $\treeY$ by $\tree$. For the remainder of the section (and
later on) we think of $\Y$ as being fixed, and let $Y_i = \{
y_{i,1},\dots,y_{i,m_i} \}$, for $i\in\N_+$. In case all the sets
$Y_i$ are equal, say to $Y$, the tree $\treeY$ is said to be
\emdef[tree!regular]{regular}, and is denoted by $\tree^{(Y)}$.

\begin{figure}[!ht]
  \begin{center}
    \includegraphics{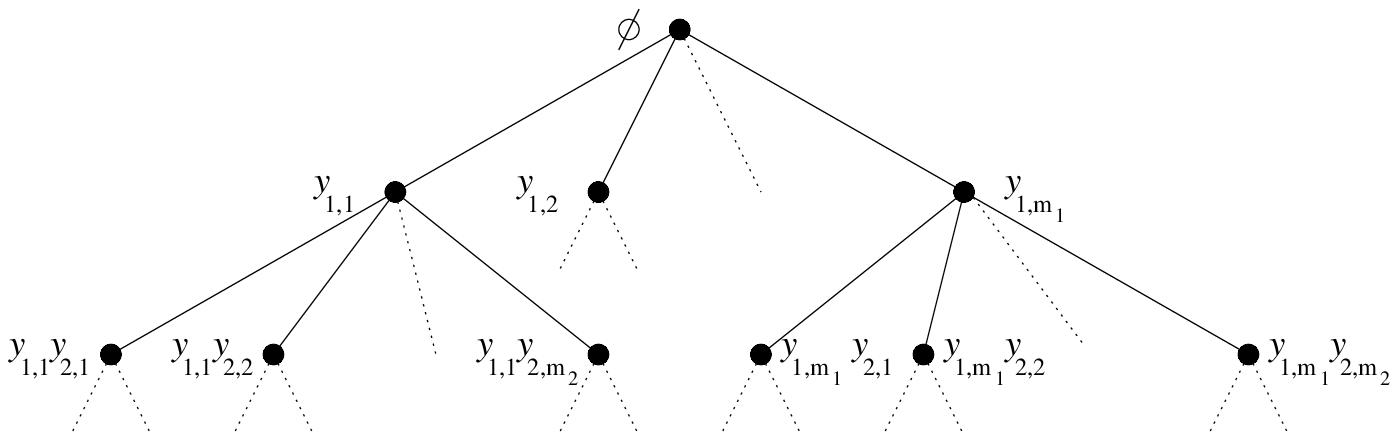}
  \end{center}
  \caption{The tree $\tree$} \label{figure:treeY}
\end{figure}

Let us give now the graph-theoretical interpretation of $\tree$
and thus justify our terminology. Every word over $\Y$ represents
a vertex in a rooted tree. Namely, the empty word $\emptyset$
represents the \emdef{root}\index{tree!root}, the $m_1$ one-letter
words $y_{1,1}$, $\dots$, $y_{1,m_1}$ represent the $m_1$ children
of the root, the $m_2$ two-letter words $y_{1,1}y_{2,1}$, $\dots$,
$y_{1,1}y_{2,m_2}$ represent the $m_2$ children of the vertex
$y_{1,1}$, etc. More generally, if $u$ is a word over $\Y$, then
the words $uy$, for $y$ in $Y_{|u|+1}$, of length $|u|+1$
represent the $m_{|u|+1}$ \emdef[vertex!children]{children} (or
\emdef[vertex!successors]{successors}) of $u$ (see
Figure~\ref{figure:treeY}).

The graph structure of $\tree$ induces a distance function on the
set of words by
\[d(u,v) = |u| + |v| - 2|u \wedge v|,\]
where $u \wedge v$ is the longest common prefix of $u$ and $v$. In
particular, the words of length $n$ represent the vertices that are at
distance $n$ to the root. Such vertices constitute the
\emdef[tree!level]{level} $n$ of the tree, denoted by $\LL_n^{(\Y)}$
or, when $\Y$ is assumed, just by $\LL_n$. In the terminology of
metric spaces, the vertices on the level $n$ are precisely the
elements of the \emph{sphere of radius $n$} with center at the root.
In the sequel, we will rarely make any distinction between a word $u$
over $\Y$, the vertex represented by $u$ and the unique path from the
root to the vertex $u$.

\subsection{Tree automorphisms}
A permutation of the words over $\Y$ that preserves the prefix
relation is an \emdef[tree!automorphism]{automorphism} of the tree
$\tree$. From the graph-theoretical point of view an automorphism of
$\tree$ is just a graph automorphism that fixes the root. We denote
the group of automorphisms of $\tree$ by $\Aut(\tree)$. Clearly, the
orbits of the action of $\Aut(\tree)$ on $\tree$ are precisely the
levels of the tree. The fact that the automorphism group acts
transitively on the spheres centered at the root is precisely the
reason for which these trees are called spherically homogeneous.

Consider an automorphism $f$ of $\tree$ and a word $u$ over $\Y$. The
image of $u$ under $f$ is denoted by $u^f$. For a letter $y$ in
$Y_{|u|+1}$ we have $(uy)^f=u^fy'$ where $y'$ is a uniquely determined
letter in $Y_{|u|+1}$. Clearly, the induced map $y \mapsto y'$ is a
permutation of $Y_{|u|+1}$, we denote this permutation by $(u)f$ and
we call it the \emdef[vertex!permutation]{vertex
  permutation}\index{permutation!vertex} of $f$ at $u$. If we denote
the image of $y$ under $(u)f$ by $y^{(u)f}$, we have
\begin{equation} \label{eq:(uy)^f=}
  (uy)^f=u^fy^{(u)f},
\end{equation}
and this easily extends to
\begin{equation} \label{eq:u^f=}
  (y_1y_2 \dots y_n)^f = y_1^{(\emptyset)f}y_2^{(y_1)f} \dots y_n^{(y_1y_2
    \dots y_{n-1})f} .
\end{equation}

\begin{figure}[!ht]
  \begin{center}
    \includegraphics{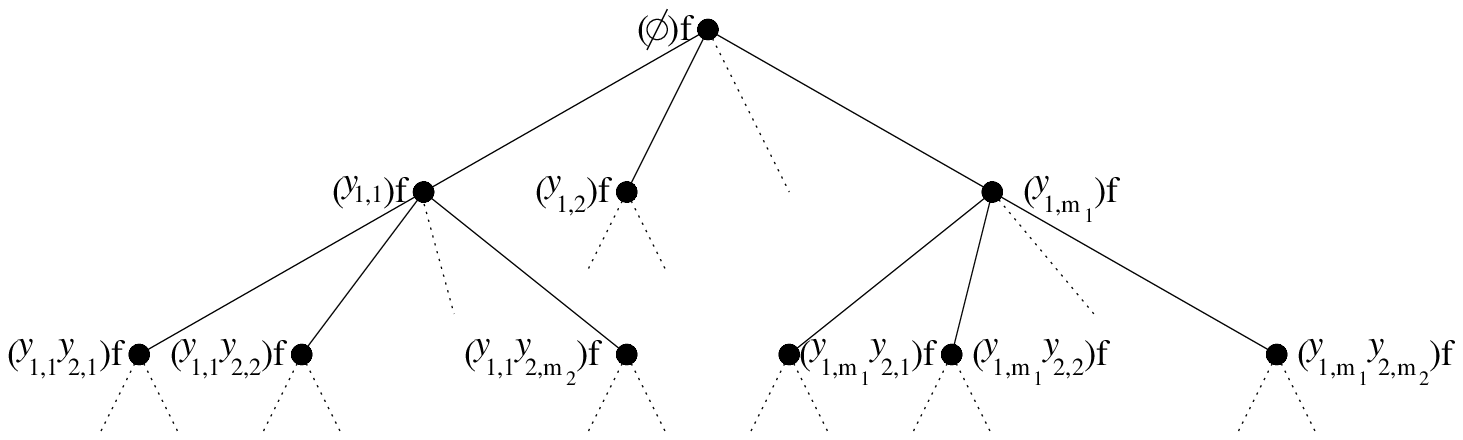}
  \end{center}
  \caption{The automorphism $f$ of $\tree$} \label{figure:f}
\end{figure}

Any tuple $((u)g)_{u \in \Y^*}$, indexed by the words $u$ over $\Y$,
where the entry $(u)g$ is a permutation of the alphabet $Y_{|u|+1}$,
determines an automorphism $g$ of $\tree$ given by
\[(y_1y_2 \dots y_n)^g = y_1^{(\emptyset)g}y_2^{(y_1)g} \dots
y_n^{(y_1y_2 \dots y_{n-1})g} .\] Therefore, we can think of an
automorphism $f$ of $\tree$ as the tuple of vertex permutations $((u)f)_{u
\in \Y^*}$ and we can represent the automorphism $f$ on the tree $\tree$
by decorating each vertex $u$ in $\tree$ by its permutation $(u)f$ (see
Figure~\ref{figure:f}). The decorated tree that represents $f$ is called
the \emdef{portrait}\index{element!portrait} of $f$.

The portrait of $f$ gives an intuitively clear picture of the action of
$f$ on $\tree$, if we can imagine what happens when we perform all the
vertex permutations at once. If only finitely many vertex permutations are
non-trivial this is not difficult to do.

By using~(\ref{eq:(uy)^f=}), we can easily see that
\begin{equation} \label{eq:(u)fg=}
  (u)fg = (u)f \circ (u^f)g \qquad \text{and} \qquad (u)f^{-1}
  = [ (u^{f^{-1}})f]^{-1},
\end{equation}
for all words $u$ over $\Y$ and automorphisms $f$ and $g$ of $\tree$.

\begin{figure}[!ht]
  \begin{center}
    \includegraphics{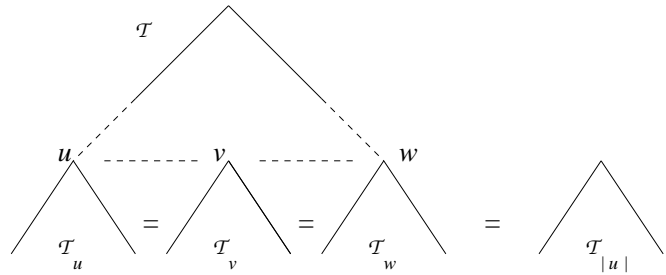}
  \end{center}
  \caption{Canonically isomorphic subtrees and shifted trees}
  \label{figure:t_u}
\end{figure}

We introduce the \emdef{shift operator} $\sigma$ that acts on
sequences as follows:
\[\sigma(s_1,s_2,s_3,\dots) = s_2,s_3,s_4,\dots.\]
Acting $n$ times on the sequence of alphabets $\Y$ gives the shifted
sequence of alphabets
\[\sigma^n\Y = Y_{n+1}, Y_{n+2}, Y_{n+3}, \dots.\]
which has the following shifted sequence of branching indices
\[\sigma^n\m = m_{n+1}, m_{n+2}, m_{n+3}, \dots,\]
and this new sequence of alphabets defines the spherically homogeneous
tree $\tree^{(\sigma^n\Y)}$. Let $u$ be a word over $\Y$ of length $n$
and denote by $\tree_u^{(\Y)}$ the spherically homogeneous tree that
consists of all words over $\Y$ with prefix $u$ ordered by the prefix
relation. It is the subtree of $\treeY$ that is hanging below the
vertex $u$. Clearly, the trees $\tree_u^{(\Y)}$ and
$\tree^{(\sigma^n\Y)}$ are canonically isomorphic under the
isomorphism $\delta_u$ that deletes the prefix $u$ from the words in
$\tree_u^{(\Y)}$, and any two trees $\tree_u^{(\Y)}$ and
$\tree_v^{(\Y)}$, where $u$ and $v$ are words over $\Y$ of the same
length, are canonically isomorphic under the isomorphism that deletes
the prefix $u$ and replaces it by the prefix $v$ (this isomorphism is
just the composition $\delta_u\delta_v^{-1}$). In order to avoid
cumbersome notation we denote the tree $\tree^{(\sigma^n\Y)}$ by
$\tree_n$ and the tree $\tree_u^{(\Y)}$ by $\tree_u$ when $\Y$ is
assumed to be fixed. The previous observations then say that
$\tree_{|u|}$ and $\tree_u$ are canonically isomorphic (see
Figure~\ref{figure:t_u}).

Let $f$ be an automorphism of $\tree$ and $u$ a word over $Y$. The
\emdef{section}\index{element!section} of $f$ at $u$ (other words in
use are \emdef{component}\index{element!component/projection/slice},
\emdef{projection} and \emdef{slice}), is the the automorphism $f_u$
of $\tree_{|u|}$ defined by the vertex permutations
\begin{equation} \label{eq:(w)f_u=}
  (w)f_u=(uw)f,
\end{equation}
for all words $w$ over $\sigma^{|u|}\Y$. Therefore, $f_u$ uses the vertex
permutations of $f$ at and below the vertex $u$ and assigns them to
words over $\sigma^{|u|}\Y$ in a natural way.

The set $G_u = \setsuch{g_u}{g\in G}$ of sections at $u$ of the elements
in $G$ is called the \emdef{section} of $G$ at $u$. We mention that the
section $G_u$ is not necessarily a subgroup of $G$ even if the tree is
regular.

It is easy to show, using~(\ref{eq:u^f=}) and~(\ref{eq:(w)f_u=}), that
\begin{equation} \label{eq:(uv)^f=}
  (uv)^f = u^fv^{f_u}
\end{equation}
for all automorphisms $f$ of $\tree$, words $u$ over $\Y$ and words $v$
over $\sigma^{|u|}\Y$. Using~(\ref{eq:(u)fg=}), (\ref{eq:(w)f_u=})
and~(\ref{eq:(uv)^f=}) we obtain the equalities
\begin{equation} \label{eq:(fg)_u=}
  (fg)_u = f_ug_{u^f} \text{ and } (f^{-1})_u = (f_{u^{f^{-1}}})^{-1},
\end{equation}
that hold for all automorphisms $f$ and $g$ of $\tree$ and words $u$
over $\Y$.

Before we move on, let us look at trees and automorphisms from another
point of view. The set of infinite paths (rays) from the root
$\emptyset$ in $\tree$ is called the \emdef[tree!boundary]{boundary}
of $\tree$ and is denoted by $\partial\tree$. We define a metric $d$
on $\partial\tree$ by
\[d(r,s) = \begin{cases} \frac{1}{2^{|r \wedge s|}}&\text{ if } r \neq s\\
  0&\text{ if } r=s,\end{cases}
\]
for all infinite rays $r$ and $s$ in $\partial\tree$. Any automorphism
$f$ of $\tree$ defines an isometry $\overline{f}$ of the space
$\partial\tree$ given by
\[(y_1y_2y_3\dots)^{\overline{f}}=y_1^{(\emptyset)f}y_2^{(y_1)f}y_3^{(y_1y_2)f}\dots.\]
Conversely, any isometry $\overline{g}$ of $\partial\tree$ defines an
automorphism $g$ of $\tree$ as follows: $u^g$ is the prefix of
$r^{\overline{g}}$ of length $|u|$, where $r$ is any infinite path in
$\partial\tree$ with prefix $u$. Therefore,
$\Aut(\tree)=\Isom(\partial\tree)$.

Note that the definition of the metric $d$ above was very arbitrary. Given
a strictly decreasing sequence $\overline{d}=(d_i)_{i\in\N}$ of positive
numbers with limit 0, we could define a metric on $\partial\tree$ by
$d(r,s)=d_{|r \wedge s|}$ if $r \neq s$, and it can be shown that the
topology of the metric space $\partial\tree$ is independent of the choice
of the sequence $\overline{d}$.

The metric space $(\partial\tree,d)$ is a universal model for
"ultrametric" "homogeneous spaces" as is mentioned
in~\cite{grigorchuk:jibg} and explained in more details in
Proposition~6.2 in~\cite{grigorchuk-n-s:automata}.

\subsection{Level and rigid stabilizers}
We introduce the notions of (rigid) vertex and level stabilizers, as well
as the congruence subgroup property.

\begin{definition}
  Let $G$ be a group of automorphisms of $\tree$. The subgroup
  $\Stab_G(u)$ of $G$, called the \emdef[vertex!stabilizer]{vertex
    stabilizer}\index{stabilizer!vertex} of $u$ in $G$, consists of
  those automorphisms in $G$ that fix the vertex $u$.
\end{definition}

For any two automorphisms $f$ and $g$ in $\Stab_G(u)$, by
using~(\ref{eq:(fg)_u=}), we have
\[(fg)_u = f_ug_{u^f} = f_ug_u,\]
so that the map
\[\varphi_u^G: \Stab_G(u) \to \Aut(\tree_{|u|})\]
given by
\[(f)\varphi_u^G = f_u\]
is a homomorphism. We call this homomorphism the
\emdef[homomorphism!section]{section
  homomorphism}\index{section!homomorphism} at $u$, and we usually
avoid the superscript. We denote the image of the section homomorphism
$\varphi_u$ by $U_u^G$, or just by $U_u$ when $G$ is assumed, and call
it the \emdef[companion!upper]{upper companion} of $G$ at $u$. Note
that the upper companion of $G$ at $u$ is a subgroup of
$\Aut(\tree_{|u|})$, and is not necessarily a subgroup of $G$, even in
case of a regular tree.  Nevertheless, in many important cases in
which the tree $\tree$ is regular the upper companion groups are equal
to $G$ after the canonical identification of the original tree $\tree$
with its subtrees.

\begin{definition}
  Let $G$ be a group of automorphisms of a regular tree $\tree$.  The
  group $G$ is \emdef{fractal}\index{group!fractal} if for every
  vertex $u$ the upper companion group $U_u$ is equal to $G$ (after
  the tree identifications $\tree=\tree_{|u|}$).
\end{definition}

The vertex stabilizers lead to the notion of level stabilizers as follows:
\begin{definition}
  Let $G$ be a group of automorphisms of $\tree$ and let
  $\Stab_G(\LL_n)$, called the $n$-th \emdef[stabilizer!level]{level
    stabilizer}\index{level!stabilizer} in $G$, denote the subgroup of
  $G$ consisting of the automorphisms of $\tree$ that fix all the
  vertices on the level $n$ (and up of course), i.e.,
\[\Stab_G(\LL_n) = \bigcap_{u \in \LL_n} \Stab_G(u).\]
\end{definition}

The homomorphism
\[ \psi_n^G:\Stab_G(\LL_n) \rightarrowtail \prod_{u\in\LL_n} U_u \leq
    \prod_{u\in \LL_n}\Aut(\tree_{|u|})\]
given by
\[(f)\psi_n^G = ((f)\varphi_u^G)_{u \in \LL_n} = (f_u)_{u \in \LL_n}\]
is an embedding, since the only automorphism that fixes all the
vertices at level $n$ and acts trivially on all the subtrees
hanging below the level $n$ is the trivial one. In case $n=1$ we
almost always omit the index 1 in $\psi_1$, and we omit the
superscript $G$, for all $n$. We will see in a moment that the
level stabilizers of $G$ have finite index in $G$. It follows that
the same is true for the vertex stabilizers.

We note that the current literature contains several versions of
definitions of fractal branch groups. In some of them the
sufficient condition from Lemma~\ref{lemma:fractal} below is used
as a definition. One can impose even stronger conditions.

\begin{definition}
  Let $G$ be a group of automorphisms of a regular tree $\tree$.  The
  group $G$ is \emdef{strongly fractal}\index{group!strongly fractal}
  if it is fractal and the embedding
  \[ \psi:Stab_G(\LL_1) \rightarrow \prod_{i=1}^m G \]
  is subdirect, i.e. surjective on each factor.
\end{definition}

\begin{lemma}\label{lemma:fractal}
Let $G$ be a group of automorphisms of $\tree^{(p)}$, $p$ a prime,
and let all vertex permutations of the automorphisms in $G$ are
powers of a fixed cyclic permutation of order $p$. Then, $G$ is
fractal if and only if $G$ is strongly fractal.
\end{lemma}

\begin{definition}
A group $G$ of tree automorphisms satisfies the \emdef{congruence subgroup
property} if every finite index subgroup of $G$ contains a level
stabilizer $\Stab_G(\LL_n)$, for some $n$.
\end{definition}

We now move to the rigid version of stabilizers:

\begin{definition}
  The \emdef[stabilizer!rigid vertex]{rigid vertex
    stabilizer}\index{rigid vertex stabilizer} of $u$ in $G$, denoted
  by $\Rist_G(u)$, is the subgroup of $G$ that consists of those
  automorphisms of $\tree$ that fix all vertices not having $u$ as a
  prefix.
\end{definition}

\begin{figure}[!ht]
  \begin{center}
    \includegraphics{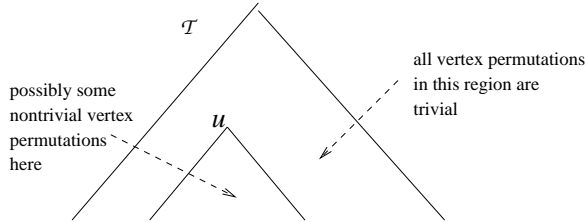}
  \end{center}
  \caption{An automorphism in the rigid stabilizer of $u$}
  \label{figure:rist}
\end{figure}

The automorphisms in $\Rist_G(u)$ must also fix $u$, and the only vertex
permutations that are possibly non-trivial are those corresponding to the
vertices in $\tree_u$ (see Figure~\ref{figure:rist}).

The rigid stabilizer $\Rist_G(u)$ is also known as the
\emdef[companion!lower]{lower companion} of $G$ at $u$, denoted by
$L_u^G$, or by $L_u$ when $G$ is assumed. Clearly, the lower companion
group at $u$ can be embedded in the upper companion group which is
contained in the corresponding section, i.e.
\[L_u \rightarrowtail U_u \subseteq G_u.\]

\begin{definition}
  The subgroup of $G$ generated by all the rigid stabilizers of
  vertices on the level $n$ is the \emdef[stabilizer!rigid $n$-th
  level]{rigid $n$-th level stabilizer} and it is denoted by
  $\Rist_G(\LL_n)$.
\end{definition}

Clearly, automorphisms in different rigid vertex stabilizers on the same
level commute and
\[\Rist_G(\LL_n) = \prod_{u \in \LL_n} \Rist_G(u).\]

The level stabilizer $\Stab_G(\LL_n)$ and the rigid level stabilizer
$\Rist_G(\LL_n)$ are normal in $\Aut(\tree)$. Further, the following
relations hold:
\[\prod_{u\in \LL_n} L_u = \Rist_G(\LL_n) \leq
\Stab_G(\LL_n)  \overset{\psi}{\rightarrowtail} \prod_{u\in \LL_n} U_u.
\]

In contrast to the level stabilizers, the rigid level stabilizers may have
infinite index, and may even be trivial.

Let us restrict our attention, for a moment, to the case when $G$
is the full automorphism group $\Aut(\tree)$. Clearly, every
automorphism of $\tree_{|u|}$ is a section of an automorphism of
$\tree$, since any choice of vertex permutations at and below $u$
is possible for automorphisms of the tree $\tree$ that fix $u$.
Therefore, $\Aut(\tree)_u = \Aut(\tree_{|u|})=\Aut(\tree_u)$,
i.e., the section is equal to the full automorphism group of the
corresponding subtree. Moreover, the section groups are equal to
the corresponding upper companion groups. It is also clear that
the rigid stabilizer $\Rist_{\Aut(\tree)}(u)$ is canonically
isomorphic to $\Aut(\tree_u)$, that the rigid and the level
stabilizer of the same level are equal, and $\psi$ is an
isomorphism.

Consider the subgroup $\Aut_f(\tree)$ of automorphisms that have only
finitely many non-trivial vertex permutations. The automorphisms in
this group are called
\emdef[automorphism!finitary]{finitary}\index{finitary automorphism}.
The group of finitary automorphisms is the union of the chain of
subgroups $\Aut_{[n]}(\tree)$ for $n\in\N$, where $\Aut_{[n]}(\tree)$
denotes the group of tree automorphisms whose vertex permutations at
level $n$ and below are trivial. The group $\Aut_{[n]}(\tree)$ is
canonically isomorphic to the automorphism group $\Aut(\tree_{[n]})$
of the finite tree $\tree_{[n]}$ that consists of the vertices of
$\tree$ represented by words no longer than $n$ (level $n$ and above).
The group $\Aut(\tree_{[n]})$ is isomorphic to the iterated
permutational wreath product
\[\Aut(\tree_{[n]})\cong((\dots(\Sym(Y_n)\wr\Sym(Y_{n-1}))\wr\dots)\wr\Sym(Y_1),\]
and its cardinality is $m_1! (m_2!)^{m_1} (m_3!)^{m_1m_2} \dots
(m_n!)^{m_1m_2 \dots m_{n-1}}$. Also, the equality
\[\Aut(\tree) = \Stab_{\Aut(\tree)}(\LL_n) \rtimes \Aut_{[n]}(\tree),\]
holds. As the intersection of all level stabilizers is trivial we see that
$\Aut(\tree)$ is residually finite and, as a corollary, every subgroup of
$\Aut(\tree)$ is residually finite.

We organize some of the remarks we already made in the following
diagram:
\[\begin{diagram}
  \dgARROWLENGTH=1em
  \node{\Aut(\tree)}\arrow{s,=} \node{\Aut(\tree)}\arrow{s,=}
  \node{\Aut(\tree)}\arrow{s,=} \\
  \node{\Stab_{\Aut(\tree)}(\LL_0)}
  \node{\Stab_{\Aut(\tree)}(\LL_1)}\arrow{w,L}
  \node{\Stab_{\Aut(\tree)}(\LL_2)}\arrow{w,L}
  \node{\dots}\arrow{w,L}\\
  \node{\rtimes}   \node{\rtimes}   \node{\rtimes} \\
  \node{\Aut_{[0]}(\tree)} \node{\Aut_{[1]}(\tree)}\arrow{w,A}
  \node{\Aut_{[2]}(\tree)}\arrow{w,A} \node{\dots}\arrow{w,A}
\end{diagram}\]
The homomorphisms in the bottom row are the natural restrictions, and
$\Aut(\tree)$ is the inverse limit of the inverse system represented by
this row. Thus $\Aut(\tree)$ is a profinite group with topology that
coincides with the Tychonoff product topology.

In this topological setting, we recall the \emph{Hausdorff dimension}
of a subgroup of $\Aut(\tree)$:
\begin{definition}[\cite{barnea-s:hausdorff}]
  Let $G<\Aut(\tree)$ be a closed subgroup. Its
  \emdef[dimension!Hausdorff]{Hausdorff dimension} is
  \[\limsup_{n\to\infty}\frac{\log|G/\Stab_G(\LL_n)|}{\log|\Aut(\tree)/\Stab_{\Aut(\tree)}(\LL_n)|},\]
  a real number in $[0,1]$.
\end{definition}

Note that, according to our agreements, the iterated permutational wreath
product
\[ \prod_{i=1}^n(\wr) \Sym(m_i) = ((\dots(\Sym(m_n)\wr\Sym(m_{n-1}))\wr\dots)\wr\Sym(m_1),\]
naturally acts on $Y_1 \times Y_2 \times \dots Y_n$ which is exactly the
set of words of length $n$. The action is by permutations $f$ that respect
prefixes in the sense that
\[ | u \wedge v | = | u^f \wedge v^f |, \]
for all words $u$ and $v$ of length $n$. This allows us to define the
action on the set of all words of length at most $n$, which is exactly why
we may think of $\prod_{i=1}^n(\wr) \Sym(m_i)$ as being the automorphism
group $\Aut(\tree_{[n]})$ of the finite tree $\tree_{[n]}$.

Since $\prod_{i=1}^n(\wr) \Sym(m_i)$ acts on the words of length $n$, the
inverse limit $\varprojlim_n \prod_{i=1}^n(\wr) \Sym(m_i)$ acts on the set
of infinite words by isometries, which is one of the interpretations of
$\Aut(\tree)$ we already mentioned.

We agree on a simplified notation concerning the word
$u=y_{1,j_1}y_{2,j_2}\dots y_{n,j_n}$ over $\Y$, the section $f_u$ of the
automorphism $f$ and the homomorphism $\varphi_u$. We will write sometimes
just $u=j_1 j_2\dots j_n$ since the sequence of indices $j_1 j_2\dots j_n$
uniquely determines and is uniquely determined by the word $u$. Also, we
will write $f_{j_1j_2\dots j_n}$ and $\varphi_{j_1j_2\dots j_n}$ for the
appropriate section $f_u$ and section homomorphism $\varphi_u$. Actually,
we could agree that $Y_i=\{1,2,\dots,m_i\}$, for $i\in\N_+$, in which case
the original and the simplified notation are the same.

\section{Geometric definition of a branch group}\label{sec:gdef}
For the length of this section we make an important assumption that $G$ is
a group of automorphisms of $\tree$ that acts transitively on each level
of the tree. In this case we say that $G$ acts \emph{spherically
transitively}\index{action!spherically transitive}.

It follows easily from~(\ref{eq:(fg)_u=}), that all vertex
stabilizers $\Stab_G(u)$ corresponding to vertices $u$ on the same
level are conjugate in $G$. Indeed, if $h$ fixes $u$ then $h^g$
fixes $u^g$ and $(h_u)^{g_u} = (h^g)_{u^g}$. This also shows that,
in case of a spherically transitive action, the upper companion
groups $U_u$, $u$ a vertex in the level $\LL_n$, are conjugate in
$\Aut(\tree_n)$, we denote by $U_n^G$, or just by $U_n$, their
isomorphism type, and we call it the \emdef{upper companion group}
of $G$ at level $n$ or the $n$-th upper companion:
\[\begin{diagram}
  \node{L_u=\Rist_G(u)}\arrow{e,t,J}{\trianglelefteq}\arrow[2]{s,l,A,V}{(.)^g}
  \node{\Stab_G(u)}\arrow[2]{e,t,A}{\varphi_u}\arrow[2]{s,l,A,V}{(.)^g}
  \node[2]{U_u}\arrow[2]{s,r,A,V}{(.)^{{g_u}}} \\
  \node[3]{(h_u)^{g_u}=(h^g)_{u^g}}\\
  \node{L_{u^g}=\Rist_G(u^g)}\arrow{e,t,J}{\trianglelefteq}
  \node{\Stab_G(u^g)}\arrow[2]{e,b,A}{\varphi_{u^g}} \node[2]{U_{u^g}}
\end{diagram}\]

Moreover, we note that if the section $g_u$ is trivial then the upper
companion groups $U_u$ and $U_{u^g}$ are not only conjugate, but they are
equal.

Similarly, the rigid vertex stabilizers of vertices on the same level
are also conjugate in $G$, we denote by $L_n^G$, or just by $L_n$,
their isomorphism type, and we call it the
\emdef[companion!lower]{lower companion group} of $G$ at level $n$ or
the $n$-th lower companion (see the above diagram).  We note that the
rigid vertex stabilizer $\Rist_G(u)$ is a normal subgroup in the
corresponding vertex stabilizer $\Stab_G(u)$. Moreover, the lower
companion group $L_u$ naturally embeds via the section homomorphisms
$\varphi_u$ in the upper companion group $U_u$ as a normal subgroup.
In case of the full automorphism group, we already remarked that this
embedding is an isomorphism. In general, this is not true, and we will
study more closely the ``next best case'' when the embedded subgroup
has a finite index.

\begin{proposition}\label{prop:rist}
  Let $G$ be a group of automorphisms of $\tree$ acting spherically
  transitively. If $\Rist_G(\LL_n)$ has finite index in $G$ for all
  $n$, then $G$ is branch group with branch structure
  $(L_n^G,\Rist_G(\LL_n))_{i=1}^\infty$.
\end{proposition}

\begin{definition}\label{defn:gbranch}
  Let $G$ be a group of automorphisms of $\tree$ acting spherically
  transitively. We say that $G$ is a
\begin{enumerate}
\item \emdef[group!branch]{branch group acting on a
    tree}\index{branch!on a tree} if all rigid stabilizers of $G$ have
  finite index in $G$.
\item \emdef[group!weakly branch]{weakly branch group acting on a
    tree}\index{weakly branch!on a tree} if all rigid stabilizers
    of $G$ are non-trivial (which implies that they are infinite).
\item \emdef[group!rough]{rough group acting on a tree}\index{rough!on
    a tree} if all rigid stabilizers of $G$ are trivial.
\end{enumerate}
\end{definition}

The branch structure from the previous proposition is not unique, as
usual, and we see that for $G$ to be branch group it is enough if we
require that each rigid vertex stabilizer group $\Rist_G(u)$, $u$ a vertex
in $\tree$, has a subgroup $L(u)$ such that
$H_n=\prod\setsuch{L(u)}{u\in\LL_n}$ is normal of finite index in $G$, for
all $n$.

Particularly important type of branch groups is introduced by the
following definitions.

\begin{definition}
  A fractal branch group $G$ acting on the regular tree $\tree^{(m)}$
  is a \emdef[group!regular branch]{regular branch group} if there
  exists a finite index subgroup $K$ of $G$ such that $K^m$ is
  contained in $(K)\psi$ as a subgroup of finite index. In such a
  case, we say that $G$ is \emph{branching over} $K$. We also say that
  $K$ \emph{geometrically contains}\index{geometric inclusion} $K^m$.
  In case $K$ contains $K^m$ but the index is infinite we say that $G$
  is \emdef[branch!weakly regular]{weakly regular branch} over $K$.
\end{definition}

\begin{definition}\label{defn:contract}
  Let $G$ be a regular branch group generated as a monoid by a finite
  set $S$, and consider the induced word metric on $G$. We say $G$ is
  \emdef[group!contracting]{contracting}\index{contracting
    homomorphism} if there exist positive constants $\lambda<1$ and $C$ such
    that for every word
  $w\in S^*$ representing an element of $\Stab_G(\LL_1)$, writing
  $(w)\psi=(w_1,\dots,w_m)$, we have
  \begin{equation}\label{eq:contract}
    |w_i|<\lambda|w|\text{ for all }i\in Y,\text{ as soon as }|w|>C.
  \end{equation}
  The constant $\lambda$ is called a \emph{contracting constant}.
\end{definition}

In a loose sense, the abstract branch groups are groups that
remind us of the full automorphism groups of the spherically
homogeneous rooted trees. Any branch group has a natural action on
a rooted tree. Indeed, let $G$ be a branch group with branch
structure $(L_i,H_i)_{i\in\N}$. The set of subgroups
$\setsuch{L_i^{(j)}}{i\in\N,\;j=1,\dots,k_i}$ ordered by inclusion
forms a spherically homogeneous tree with branching sequence $\m =
m_1,m_2,\dots$, where $m_i=k_i/k_{i-1}$.  The group $G$ acts on
this set by conjugation and, because of the refinement conditions,
the resulting permutation is a tree automorphism (see
Figure~\ref{figure:Ls})

The action of the branch group $G$ on the tree determined by the branch
structure is not faithful in general. Indeed, it is known that a branch
group that satisfies the conditions in Proposition~\ref{prop:rist} is
centerless (see~\cite{grigorchuk:jibg}). On the other hand, a direct
product of a branch group $G$ in the sense of our algebraic definition
with a finite group $H$ is still a branch group. If $H$ has non-trivial
center, then $G \times H$ is a branch group with non-trivial center.

It would be interesting to understand the nature of the kernel of the
action in the passing from an abstract branch group to a group that acts
on a tree. In particular, is it correct that this kernel is always in the
center (Question~\ref{question:kernel-action})? We note that the kernel is
trivial in case of a just-infinite branch group (see
Chapter~\ref{chapter:just-infinite}).

\section{Portraits and branch portraits}\label{sec:bp}
A tree automorphism can be described by its \emph{portrait}, already
defined before, and repeated here in the following form:
\begin{definition}\label{defn:label}
  Let $f$ be an automorphism of $\tree$. The
  \emdef[portrait!branch]{portrait}\index{branch!portrait} of $g$ is a
  decoration of the tree $\tree$, where the decoration of the vertex
  $u$ belongs to $\Sym(Y_{|u|+1})$, and is defined inductively as
  follows: first, there is $\pi_\emptyset\in\Sym(Y_1)$ such that
  $g=h\pi_\emptyset$ and $h$ stabilizes the first level. This
  $\pi_\emptyset$ is the label of the root vertex. Then, for all $y\in
  Y_1$, label the tree below $y$ with the portrait of the section
  $g_y$.
\end{definition}

The following notion of a branch portrait based on the branch structure of
the group in question is useful in some considerations:
\begin{definition}\label{defn:bp}
  Let $G$ be a branch group, with branch structure
  $(L_i,H_i)_{i\in\N}$. The \emdef[portrait!branch]{branch portrait}
  of $g\in G$ is a decoration of the tree $\treeY$, where the
  decoration of the root vertex belongs to $G/H_1$ and the decoration
  of the vertex $y_1\dots y_i$ belongs to $L^{(y_1\dots
    y_i)}/\prod_{y\in Y_{i+1}}L^{(y_1\dots y_iy)}$. Fix once and for
  all transversals for the above coset spaces. The branch portrait of
  $g$ is defined inductively as follows: the decoration of the root
  vertex is $H_1g$, and the choice of the transversal gives us an
  element $(g_{y_1})_{y_1\in Y_1}$ of $H_1$.  Decorate then $y_1\in
  Y_1$ by $\prod_{y_2\in Y_2}L^{(y_1y_2)}g_{y_1}$; again the choice of
  transversals gives us elements $g_{y_1y_2}\in L^{(y_1y_2)}$; etc.
\end{definition}

There are uncountably many possible branch portraits that use the chosen
transversals, even when $G$ is a countable branch group. We therefore
introduce the following notion:
\begin{definition}\label{defn:brcompletion}
  Let $G$ be a branch group. Its \emdef[completion!tree]{tree
  completion}\index{tree!completion} $\overline G$ is
  the inverse limit
  \[\varprojlim G/\Stab_G(n).\]
  This is also the closure in $\Aut\tree$ of $G$ in the "topology" given
  by its action on the tree $\tree$.
\end{definition}
Note that since $\overline G$ is closed in $\Aut(\tree)$ it is a
profinite\index{group!profinite} group, and thus is compact, and
totally disconnected. If $G$ has the "congruence subgroup
property"~\cite{grigorchuk:jibg}, then $\overline G$ is also the
profinite completion of $G$.

\begin{lemma}\label{lemma:portraits}
  Let $G$ be a branch group and $\overline G$ its tree completion. Then
  Definition~\ref{defn:bp} yields a bijection between the set of
  branch portraits over $G$ and $\overline G$.
\end{lemma}

Branch portraits are very useful to express, for instance, the lower
central series. They appear also, in more or less hidden manner, in most
results on growth and torsion.

\section{Groups of finite automata and recursively defined automorphisms}\label{sec:automata}
We introduce two more ways to think about tree automorphisms in the case
of a regular tree. It is not impossible to extend the definitions to more
general cases, but we choose not to do so. Thus for the length of this
section we set $Y=\{1,\dots,m\}$ and we work with the regular tree
$\tree=\tree^{(Y)}=\tree^{(m)}$.

\subsection{Recursively defined automorphisms}
Let $X=\{x^{(1)},\dots,x^{(n)}\}$ be a set of symbols and $F$ a finite set
of finitary automorphisms of $\tree$. The following equations
\begin{equation}\label{eq:x^(i)}
x^{(i)} = (w_{i,1},\dots,w_{i,m})f^{(i)}, \qquad i \in \{1,\dots,n\},
\end{equation}
where $w_{i,j}$ are words over $X \cup F$ and $f^{(i)}$ are elements in
$F$, define an automorphism of $\tree$, still written $x^{(i)}$, for each
symbol $x^{(i)}$ in $X$. The way in which the equations~(\ref{eq:x^(i)})
define automorphisms recursively is as follows: we interpret the $m$-tuple
$(w_{i,1},\dots,w_{i,m})$ as an automorphism fixing the first level and
$w_{i,j}$ are just the sections at $j$, $j \in\{1,\dots,m\}$; the
equations~(\ref{eq:x^(i)}) clearly define the vertex permutation at the
root for all $x^{(i)}$; if the vertex permutations of all the $x^{(i)}$
are defined for the first $k$ levels then the vertex permutations of all
the $w_{i,j}$ are defined for the first $k$ levels, which in turn defines
all the vertex permutations of all the $x^{(i)}$ for the first $k+1$
levels.

Every automorphism of $\tree$ that can be defined as a member of some
set $X$ of recursively defined automorphisms as above is called a
\emph{recursively definable}\index{automorphism!recursively definable}
automorphism of $\tree$. The set of recursively definable
automorphisms of $\tree$ forms a subgroup $\Aut_r(\tree)$ of
$\Aut(\tree)$. The group $\Aut_r(\tree)$ is a regular branch group
which properly contains $\Aut_f(\tree)$.

\begin{figure}[!ht]
  \begin{center}
    \includegraphics{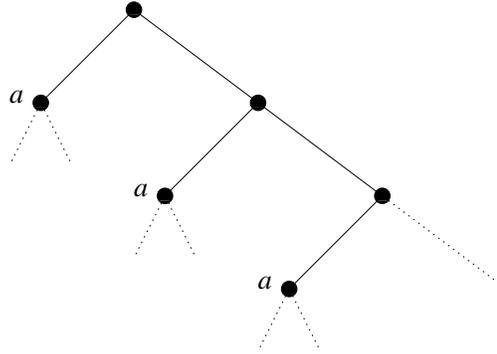}
  \end{center}
  \caption{The recursively defined automorphism $b$}
  \label{figure:dih}
\end{figure}

When one defines automorphisms recursively it is customary to choose all
finitary automorphisms $f^{(i)}$ to be rooted automorphisms (see
Definition~\ref{defn:rooted}). The advantage in that case is that
$w_{i,j}$ is exactly the section of $x^{(i)}$ at $j$. As an example of a
recursively defined automorphisms consider
\[ b = (a,b) \]
acting on the binary tree $\tree^{(2)}$, where $a=((1,2))$ is the
non-trivial rooted automorphism of $\tree^{(2)}$. Clearly, the diagram in
Figure~\ref{figure:dih} represents $b$ through its vertex permutations.

It is easy to see that all tree automorphisms are recursively definable if
we extend our definition to allow infinite sets $X$. Indeed
\[g_u = (g_{u1},\dots,g_{um})(u)g, \qquad u \in Y^* \]
defines recursively $g$ and all of its sections.

\subsection{Groups of finite automata}
Since we want to define automata that behave like tree
automorphisms we need automata that transform words rather then
recognize them, i.e., we will be working with transducers. The
fact that we want our automata to preserve lengths and permute
words while preserving prefixes strongly suggests the choices made
in the following definition.

\begin{definition}
  A \emdef[transducer!synchronous invertible finite]{synchronous
    invertible finite transducer} is a quadruple
  $T=(Q,Y,\tau,\lambda)$ where
\begin{enumerate}
\item $Q$ is a finite set (set of \emph{states} of $T$),
\item $Y$ is a finite set (the \emph{alphabet} of $T$),
\item $\tau$ is a map $\tau:Q \times Y \rightarrow Q $ (the
  \emph{transition function} of $T$), and
\item $\lambda$ is a map $\lambda:Q \times Y \rightarrow Y$ (the
  \emph{output function} of $T$) such that the induced map
  $\lambda_q:Y \rightarrow Y$ obtained by fixing a state $q$ is a
  permutation of $Y$, for all states $q \in Q$.
\end{enumerate}

If $T$ is a synchronous invertible finite transducer and $q$ a state
in $Q$ we sometimes give $q$ a distinguished status, call it the
\emph{initial state} and define the \emph{initial synchronous
  invertible finite transducer} $T_q$ as the transducer $T$ with
initial state $q$.
\end{definition}

We say just (initial) transducer in the sequel rather then (initial)
synchronous invertible finite transducer.

\begin{figure}[!ht]
  \begin{center}
    \includegraphics{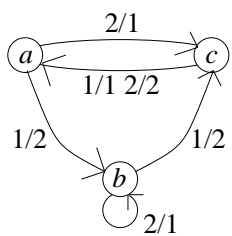}
  \end{center}
  \caption{An example of a transducer}
  \label{figure:free}
\end{figure}

It is customary to represent transducers with directed labelled
graphs where $Q$ is the set of vertices and there exists an edge
from $q_0$ to $q_1$ if and only if $\tau(q_0,y)=q_1$, for some $y
\in Y$, in which case the edge is labelled by $y|\lambda(q_0,y)$.
The diagram in Figure~\ref{figure:free} gives an example.

Informally speaking, given an initial transducer $T_q$ and an
input word $w$ over $Y$ we start at the vertex $q$ and we travel
through the graph by reading $w$ one letter at a time and
following the values of the transition function. Thus if we find
ourselves at the state $q'$ and we read the letter $y$ we move to
the state $\tau(q',y)$ by following the edge labelled by
$y|\lambda(q',y)$. In the same time, we write down an output word,
one letter at a time, simply by writing down the letters after the
vertical bar in the labels of the edges we used in out journey.

More formally, given an initial transducer $T_q$ we define recursively the
maps $\tau_q: Y^* \rightarrow Q$ and $\lambda_q: Y^* \rightarrow Y^*$ as
follows:
\begin{align*}
    \tau_q(\emptyset) &= q \\
    \tau_q(wy) &= \tau(\tau_q(w),y) \qquad\text{for}\; w \in Y^* \\
    \lambda_q(\emptyset) &= \emptyset \\
    \lambda_q(wy) &= \lambda_q(w)\lambda(\tau_q(w),y) \qquad\text{for}\; w \in Y^*
\end{align*}
It is not difficult to see that $\lambda_q$, the \emph{output
  function} of the initial transducer $T_q$, represents an
automorphism of $\tree$.  The set of all tree automorphisms
$\Aut_{ft}(\tree)$ that can be realized as output functions of some
initial transducer forms a subgroup of $\Aut(\tree)$ and this subgroup
is a regular branch group sitting properly between the group of
finitary and the group of recursively definable automorphisms of
$\tree$.

If we allow infinitely many states, then every automorphism $g$ of $\tree$
can be realized by an initial transducer. Indeed, we may define the set of
states $Q$ to be the set of sections of $g$, i.e. $Q=\setsuch{g_u}{u \in
Y^*}$ and
\[ \tau(g_u,y) = g_{uy}  \qquad\text{and}\qquad \lambda(g_u,y) = y^{(u)g}. \]
We could index the states by the vertices in $\tree$, but by indexing them
by the sections of $g$ we see that an automorphism $g$ can be defined by a
finite initial transducer if it has only finitely many distinct sections.
The converse is also true.

\begin{proposition}
  An automorphism of $\tree$ is the output function of some initial
  transducer if and only if it has finitely many distinct sections.
\end{proposition}

Every transducer $T$ defines a group $G_T$ of tree automorphisms
generated by the initial transducers of $T$ (one for each state).
Groups that are defined by transducers are known as \emdef{groups
of automata}.

Note that the notion of a group of automata is different from the
notion of automatic group in the sense of \JCannon. For more
information on groups of automata we refer the reader
to~\cite{grigorchuk-n-s:automata}.
\section{Examples of branch groups}\label{sec:branch:examples}
We have already seen a couple of examples of branch groups acting on a
regular tree $\tree$. Namely, the groups $\Aut_f(\tree)$,
$\Aut_{ft}(\tree)$, $\Aut_r(\tree)$ and the full automorphism group
$\Aut(\tree)$.

\begin{proposition}
Let $\tree$ be regular tree and $G$ be any of the groups $\Aut_f(\tree)$,
$\Aut_{ft}(\tree)$, $\Aut_r(\tree)$ or $\Aut(\tree)$. Then $G$ is a
regular branch group with $G = G \wr_{Y_1} \Sym(Y_1)$.
\end{proposition}

None of the groups in the previous proposition is finitely
generated, but the first three are countable. Another example of a
regular branch group is, for a permutation group $A$ of $Y$, the
group $\Aut_A(\tree^{(Y)})$ that consists of those automorphisms
of the regular tree $\tree^{(Y)}$ whose vertex permutations come
from $A$. A special case of the last example was mentioned before
as the infinitely iterated wreath product of copies of the cyclic
group $\Z/p\Z$ and the full automorphism group is another special
case.

In the sequel we give some examples of finitely generated branch groups.
We make use of rooted and directed automorphisms.

\begin{definition}\label{defn:rooted}
  An automorphism of $\tree$ is \emdef[automorphism!rooted]{rooted} if
  all of its vertex permutations that correspond to non-empty words
  are trivial.
\end{definition}

Clearly, the rooted automorphisms are precisely the finitary automorphisms
from $\Aut_{[1]}$. A rooted automorphism $f$ just permutes rigidly the
$m_1$ trees $\tree_1$,$\dots$,$\tree_{m_1}$ as prescribed by the root
permutation $(\emptyset)f$. It is convenient not to make too much
difference between the root vertex permutation $(\emptyset)f$ and the
rooted automorphism $f$ defined by it. Therefore, if $a$ is a permutation
of $Y_1$ we also say that $a$ is a rooted automorphism of $\tree$. More
generally, if $a$ is the vertex permutation of $f$ at $u$ and all the
vertex permutations below $u$ are trivial, then we do not distinguish $a$
from the section $f_u$ defined by it, i.e., we write $(f)\varphi_u = f_u =
a = (u)f$.

\begin{definition}\label{defn:directed}
  Let $\ell = y_1y_2y_3\dots$ be an infinite ray in $\tree$. We say
  that the automorphism $f$ of $\tree$ is
  \emdef[automorphism!directed]{directed} along $\ell$ and we call
  $\ell$ the \emdef[automorphism!spine]{spine}\index{spine} of $f$ if
  all vertex permutations along the ray $\ell$ and all vertex
  permutations corresponding to vertices whose distance to the ray
  $\ell$ is at least 2 are trivial.
\end{definition}

In the sequel, we define many directed automorphisms that use the
rightmost infinite ray in $\tree$ as a spine, i.e., the spine is
$m_1m_2m_3\dots$ . Therefore, the only vertices that can have a
nontrivial permutation are the vertices of the form $m_1m_2\dots
m_nj$ where $j \neq m_{n+1}$. Note that directed automorphisms fix
the first level, i.e., their root vertex permutation is trivial.

\subsection{The first Grigorchuk group $\Gg$}\label{subs:firstgg}
A description of the first Grigorchuk group, denoted by $\Gg$,
appeared for the first time in 1980 in~\cite{grigorchuk:burnside}.
Since then, the group $\Gg$ has been used as an example or
counter-example in many non-trivial situations.

\begin{figure}[!ht]
  \begin{center}
    \includegraphics{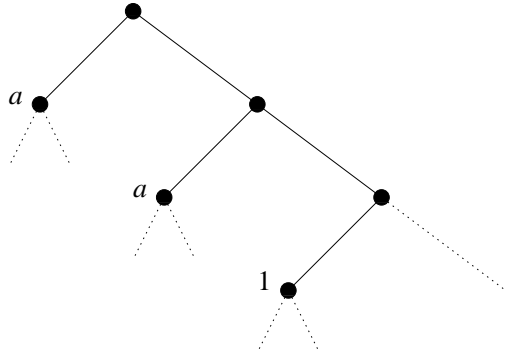}
  \end{center}
  \caption{The directed automorphism $b$} \label{figure:b}
\end{figure}

\begin{figure}[!ht]
  \begin{center}
    \includegraphics{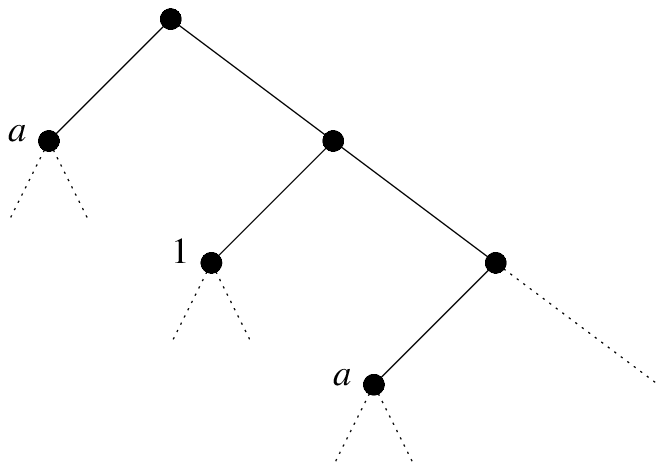}
  \end{center}
  \caption{The directed automorphism $c$} \label{figure:c}
\end{figure}

\begin{figure}[!ht]
  \begin{center}
    \includegraphics{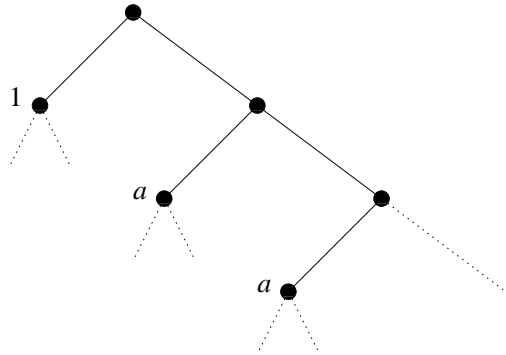}
  \end{center}
  \caption{The directed automorphism $d$} \label{figure:d}
\end{figure}

The group $\Gg$ acts on the rooted binary tree $\tree^{(2)}$ and
it is generated by the four automorphisms $a$, $b$, $c$ and $d$
defined below. The automorphism $a$ is the only possible rooted
automorphism $a=((1,2))$ that permutes rigidly the two subtrees
below the root. Parts of the portraits along the spine of the
generators $b$, $c$, and $d$ are depicted in
Figure~\ref{figure:b}, Figure~\ref{figure:c} and
Figure~\ref{figure:d}. We implicitly assume that the patterns that
are visible in the diagrams repeat indefinitely along the spine,
i.e., along the rightmost ray.

Another way to define the directed generators of $\Gg$ is by the following
recursive definition:
\[ b = (a,c) \qquad c=(a,d) \qquad d=(1,b). \]
It is clear from this recursive definition that $\Gg$ can also be defined
as a group of automata.

The group $\Gg$ is a $2$-group, has a solvable "word problem" and
intermediate growth\index{growth!intermediate}
(see~\cite{grigorchuk:gdegree}). The best known estimates of the
growth of the first Grigorchuk group are given by the first author
in~\cite{bartholdi:upperbd} and~\cite{bartholdi:lowerbd} (see
also~\cite{leonov:pont,leonov:lowerbd}). The subgroup
structure\index{subgroup!structure} of $\Gg$ is a subject of many
articles (see~\cite{rozhkov:lcs,bartholdi-g:parabolic} and
Chapter~\ref{ch:subgroup}) and it turns out that $\Gg$ has finite
width.  An infinite set of defining relations is given by
\ILysionok\ in~\cite{lysionok:pres}, and the second author shows
that this system is minimal in~\cite{grigorchuk:bath}. The
conjugacy problem is solved by \YLeonov\ in~\cite{leonov:conj} and
\ARozhkov\ in~\cite{rozhkov:conj}. A detailed exposition of many
of the known properties of $\Gg$ is included in the
book~\cite{harpe:ggt} by \PdlHarpe. Another exposition (in
Italian) appears in~\cite{ceccherini-m-s:grigorchuk}.

Most properties of $\Gg$, in one way or another, follow from the
following
\begin{proposition}
  The group $\Gg$ is a regular branch group over the subgroup $K =
  \langle[a,b]\rangle^\Gg$.
\end{proposition}
\begin{proof}
  The first step is to prove that $K$ has finite index in $\Gg$. We
  check that $a^2,b^2,c^2,d^2,bcd,(ad)^4$ are relators in $\Gg$. It
  follows that $B=\langle b\rangle^\Gg$ has index $8$,since $\Gg/B$
  is $\langle aB,dB\rangle$, a dihedral group of order
  $8$. Consequently $B=\langle b,b^a,b^{ad},b^{ada}\rangle$. Now $B/K$
  is $\langle bK\rangle$, of order $2$, so $K$ has index $2$ in $B$
  and thus $16$ in $\Gg$.

  Then we consider $L=\langle[b,d^a]\rangle^\Gg$ in $K$. A simple
  computation gives $([b,d^a])\psi=([a,b],1)$, so $(L)\psi=K\times 1$, and
  we get $(K)\psi \geq (L\times L^a)\psi = K \times K$.
\end{proof}

The index of $K$ in $\Gg$ is $16$ and the index of $K \times K$ in
$(K)\psi$ is $4$.

\subsection{The second Grigorchuk group}
\begin{figure}[!ht]
  \begin{center}
    \includegraphics{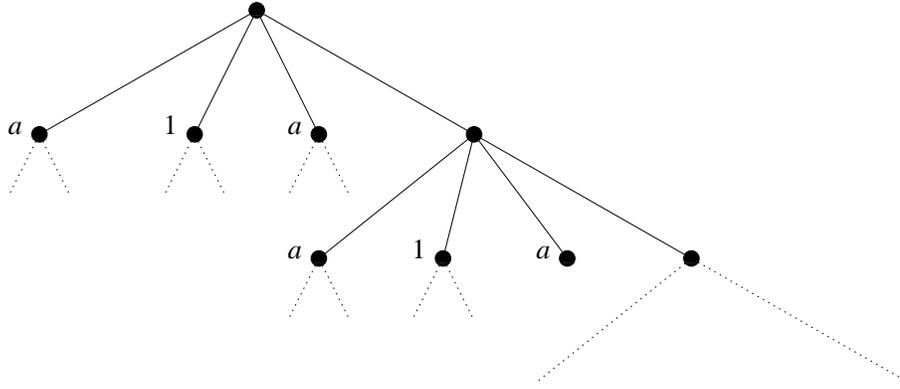}
  \end{center}
  \caption{The directed automorphism $b$ in the second Grigorchuk group}
\label{figure:b_second}
\end{figure}

The second Grigorchuk group was introduced in the same paper as the first
one (\cite{grigorchuk:burnside}). It acts on the $4$-regular tree
$\tree^{(4)}$ and it is generated by the cyclic rooted automorphism
$a=((1,2,3,4))$ and the directed automorphism $b$ whose portrait is given
in Figure~\ref{figure:b_second}.

A recursive definition of $b$ is
\[ b = (a,1,a,b). \]

The second Grigorchuk group is not investigated nearly as thoroughly as
the first one. It is finitely generated, infinite, residually finite and
centerless. In addition, it is not finitely presented, and is a torsion
group with a solvable word problem. More on this group can be found
in~\cite{vovkivsky:ggs}.

\subsection{Gupta-Sidki $p$-groups}\label{subs:gspg}
\begin{figure}[!ht]
  \begin{center}
    \includegraphics{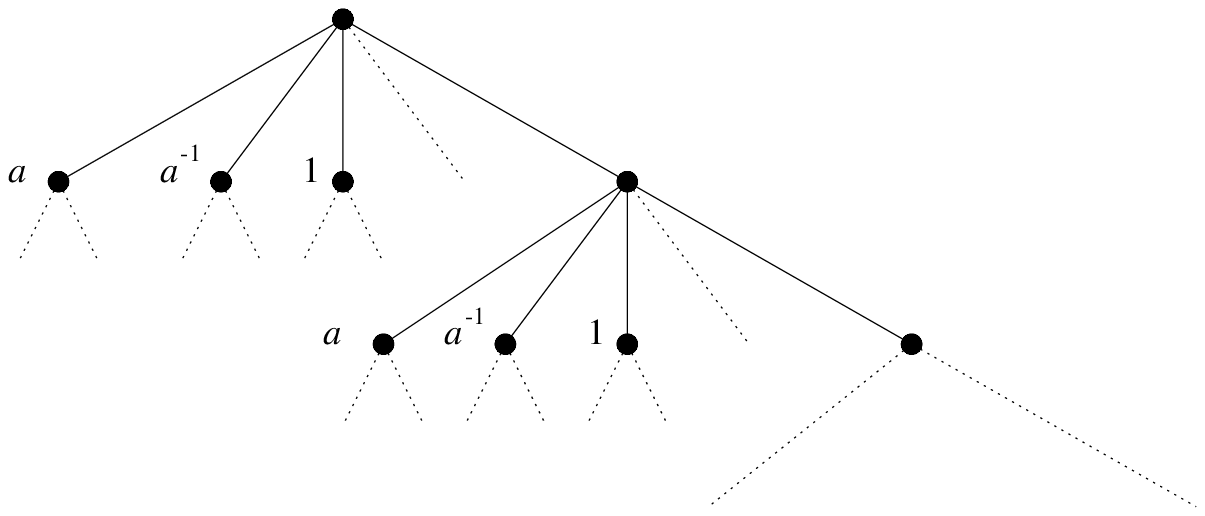}
  \end{center}
  \caption{The directed automorphism $b$ in the Gupta-Sidki groups}
  \label{figure:b_gs}
\end{figure}

The first Gupta-Sidki $p$-groups were introduced
in~\cite{gupta-s:burnside}. For odd prime $p$, we define the rooted
automorphism $a=((1,2,\dots,p))$ and the directed automorphism $b$ of
$\tree^{(p)}$ whose portrait is given in Figure~\ref{figure:b_gs}.

The group $G=\langle a, b \rangle$ is a $2$-generated $p$-group
with solvable word problem and no finite presentation
(consider~\cite{sidki:pres}). In case $p=3$ the automorphism
group, centralizers and derived groups were calculated by \SSidki\
in~\cite{sidki:subgroups}.

More examples of $2$-generated $p$-groups along the same lines were
constructed by \NGupta\ and \SSidki\ in~\cite{gupta-s:infinitep}. In
this paper the directed automorphism $b$ is defined recursively by
\[ b = (a,a^{-1},a,a^{-1},\dots,a,a^{-1},b). \]
It is shown in~\cite{gupta-s:infinitep} that these groups are
just-infinite\index{group!just-infinite} $p$-groups. Also, every
finite $p$-group is contained in the corresponding Gupta-Sidki
infinite $p$-group.

\subsection{Three groups acting on the ternary tree}\label{subsec:three}
We define three groups acting on the ternary tree\index{tree!ternary}
$\tree^{(3)}$. Each of them is $2$-generated, with generators $a$ and
$b$, where $a$ is the rooted automorphism $a=((1,2,3))$ and $b$ is one
of the following three directed automorphisms
\[ b= (a,1,b) \quad\text{or}\quad b=(a,a,b) \quad\text{or}\quad b=(a,a^2,b). \]
The corresponding group $G=\langle a,b \rangle$ is denoted by $\FGg$,
$\BGg$ and $\GSg$, respectively. The group $\FGg$ is called the
Fabrykowski-Gupta group and is the first example of a group of
intermediate growth that is not constructed by the second author. The
construction appears in~\cite{fabrykowski-g:growth1}, with an incorrect
proof, and in~\cite{fabrykowski-g:growth2}.  The group $\BGg$ is called
the Bartholdi-Grigorchuk group and is studied
in~\cite{bartholdi-g:parabolic}.  The article shows that both $\FGg$ and
$\BGg$ are virtually "torsion-free" with a torsion-free subgroup of index 3.
The group $\GSg$ is known as the Gupta-Sidki group, it is the first one of
the three to appear in print in~\cite{gupta-s:burnside}. All three groups
have intermediate growth.  The first author has calculated the central
lower series for $\FGg$ and $\GSg$
(see~\cite{bartholdi:ggs,bartholdi:phd}). We note here that $\BGg$ is not
branch group but only weakly branch group, and the other two are regular
branch groups over their commutator subgroups.

\subsection{Generalization of the Fabrykowsky-Gupta example}
\begin{figure}[!ht]
  \begin{center}
    \includegraphics{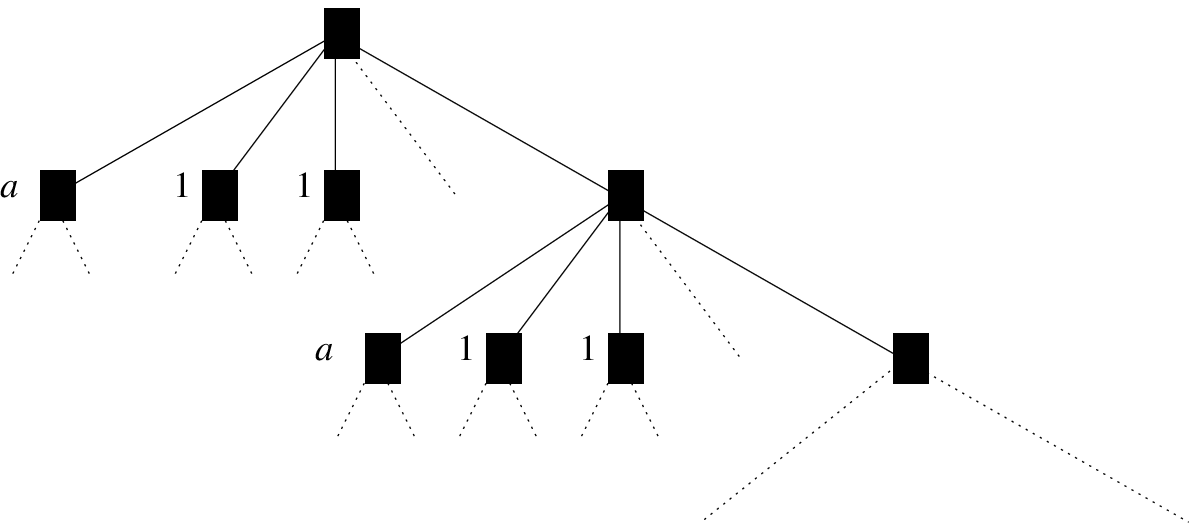}
  \end{center}
  \caption{The directed automorphism $b$}
  \label{figure:b_fg}
\end{figure}

The following examples of branch groups acting on the regular tree
$\tree^{(m)}$, for $m \geq 3$, are studied
in~\cite{grigorchuk:jibg}. The group $G=\langle a,b \rangle$ is
generated by the rooted automorphism $a=((1,2,\dots,m))$ and the
directed automorphism $b$ whose portrait is given in
Figure~\ref{figure:b_fg}.

The group $G$ is a regular branch group over its commutator. Moreover, the
rigid vertex stabilizers are isomorphic to the commutator subgroup.
Clearly, for $m=3$ we obtain the Fabrykowsky-Gupta group $\FGg$. For
$m\geq 5$, $G$ is just-infinite, and for a prime $m \geq 7$ the group $G$
has the congruence subgroup property. The last two results can probably be
extended to other branching indices.

\subsection{Examples of \PNeumann}
The following example is
constructed in~\cite{neumann:pride}. Let $A=\Alt(6)$ be the alternating
group acting on the alphabet $Y=\{1,\dots,6\}$. For each pair $(a,y)$ with
$y \in Y$ and $a \in \Stab_A(y)$, define an automorphism $b_{(a,y)}$ of
the regular tree $\tree^{(6)}$ recursively by
\[b_{(a,y)} = (1,\dots,1,b_{(a,y)},1,\dots,1)a,\]
where the only nontrivial section appears at the vertex $y$. Let $G$ be
the group generated by all these tree automorphisms, i.e.
\[G = \big\langle \setsuch{b_{(a,y)}}{ y\in Y,\;a\in\Stab_A(y)}\big\rangle.\]
Since $G$ is generated by 6 perfect\label{group!perfect} subgroups
(one copy of $\Alt(5)$ for each $y \in Y$) we see that $G$ is perfect.
It is also easy to see that if $a$ and $a'$ fix both $x$ and $y$, with
$x \neq y$, then $[b_{(a,x)},b_{(a',y)}] = [a,a']$. Since
\[\Big\langle \bigcup_{x\neq y\in Y}
[\Stab_A(\{x,y\}),\Stab_A(\{x,y\})] \Big\rangle = A,\]
we see that G contains $A$ as a rooted subgroup. It follows that $G \cong
G \wr_Y A$.

\begin{theorem}[\PNeumann~\cite{neumann:pride}]
  Let $A$ be a non-abelian finite simple\index{group!simple} group
  acting faithfully and transitively on the set $Y$. If $G$ is a
  perfect, residually finite group such that $G \cong G \wr_Y A$ then
  \begin{enumerate}
  \item All non-trivial normal subgroups subgroups of $G$ have finite
    index, i.e., $G$ is just-infinite\index{group!just-infinite} (see
    Chapter~\ref{chapter:just-infinite}).
  \item Every subnormal subgroup of $G$ is isomorphic to a finite
    direct power of $G$, but $G$ does not satisfy the ascending chain
    condition on subnormal subgroups.
    \item
    $G$ is minimal (see Chapter~\ref{chapter:just-infinite} again).
  \end{enumerate}
\end{theorem}

A group that satisfies the conditions of the previous theorem is a regular
branch group over itself, acting on the regular tree $\tree^{(Y)}$.
Furthermore, the only normal subgroups of $G$ are the (rigid) level
stabilizers $\Stab_G(\LL_n)=\Rist_G(\LL_n) \cong G^{m^n}$.
\subsection{The examples of \DSegal}\label{subs:segal}
For more details on the following examples
check~\cite{segal:fimages}. Precursors can be found
in~\cite{segal:fgfimages}.

\begin{figure}[!ht]
 \begin{center}
   \includegraphics{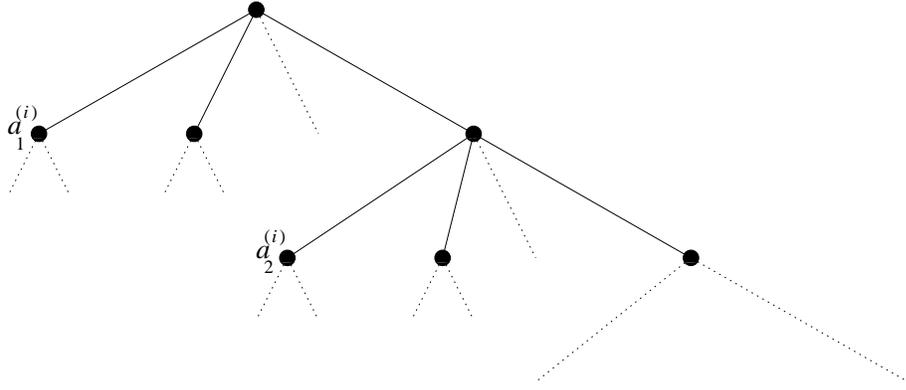}
 \end{center}
 \caption{The directed automorphism in the examples of Segal} \label{figure:segal}
\end{figure}

For $i \in \N$, let $A_i$ be a finite, perfect, transitive
permutation subgroup of $\Sym(Y_{i+1})$, where $Y_{i+1}$ is a set
of $m_{i+1}$ elements. We assume that all stabilizers
$\Stab_{A_i}(y)$ are distinct, for fixed $i$ and $y\in Y_{i+1}$.
Let $A_i=\langle a_i^{(1)}, \dots, a_ i^{(k)} \rangle$. For $j \in
\{1,\dots,k\}$, the diagram in Figure~\ref{figure:segal}
represents the directed automorphism $b_0^{(j)}$ of $\treeY$,
which is recursively defined through
\[ b_i^{(j)} = (a_{i+1}^{j},1,\dots,1,b_{i+1}^{(j)}), \quad
i\in\N, \quad j\in\{1,\dots,k\}. \]

We define a group $G_i = \langle A_i \cup B_i \rangle$ where
$B_i=\langle b_i^{(1)}, \dots, b_i^{(k)} \rangle$, for $i\in\N$.

Let $x$ be an element in $A_0$ such that $x$ fixes 1 but not $m_1$. Then
\[ [ b_0^{(i)}, (b_0^{(j)})^x ] = ([a_1^{(i)},a_1^{(j)}],1,\dots,1) \]
which, by the perfectness of $A_1$, shows that the rigid
stabilizers of the vertices on the first level contain the rooted
subgroup $A_1$. Since $((a_1^{(j)})^{-1},1,\dots,1) b_0^{j} =
(1,\dots,b_1^{(j)})$ we see that the rigid vertex stabilizers of
the vertices on the first level are exactly the upper companion
groups. Similar claims hold for the other levels and we see that
$G=G_0$ is a branch group with
\[\Stab_G(\LL_n)=\Rist_G(\LL_n)\cong \prod^{m_1m_2 \dots m_n}G_n.\]

In a similar fashion, in case $k=2$ we can construct a branch
group on only three generators as follows. For each $i\in\N$
choose a permutation $\mu_i\in A_i$ of $Y_{i+1}$ that fixes
$m_{i+1}$ but $\mu_i^2$ does not fix the symbol $r_{i+1}\in
Y_{i+1}$. Then define the directed automorphism $b_0$ recursively
through
\[ b_i =
 (1,\dots,1,a_{i+1}^{(1)},1,\dots,1,a_{i+1}^{(2)},1,\dots,1,b_{i+1}),
 \quad i\in\N, \]
where $a_{i+1}^{(1)}$ is on position $r_{i+1}$ and $a_{i+1}^{(2)}$
is on position $r_{i+1}^{\mu_i}$, for $i\in\N$. Define
$G_i=\langle A_i \cup \{b_i\} \rangle$, for $i\in\N$. Then we also
get a branch group $G=G_0$ with
\[\Stab_G(\LL_n)=\Rist_G(\LL_n)\cong \prod^{m_1m_2 \dots m_n}G_n.\]

Different choices of groups $A_i$ together with appropriate
actions give various groups with various interesting properties.
We list two interesting results from~\cite{segal:fimages}, one
below and another in Theorem~\ref{thm:segalsbgp} that use the
above examples:

\begin{theorem}
  For every collection $\mathcal S$ of finite non-abelian simple
  groups, there exists a $63$-generated just-infinite group $G$ whose
  upper composition factors (composition factors of the finite
  quotients) are precisely the members of $\mathcal S$. In
  addition, there exists a $3$-generated just-infinite group
  $\bar{G}$ whose non-cyclic upper composition factors are
  precisely the members of $\mathcal S$.
\end{theorem}

\date{October 27, 2002}
\chapter{Spinal Groups}\label{chapter:spinal}
In this chapter we introduce the class of spinal groups of tree
automorphisms. This class is rich in examples of finitely generated branch
groups with various exceptional properties, constructed by the second
author in~\cite{grigorchuk:burnside,grigorchuk:gdegree,grigorchuk:pgps},
\NGupta\ and \SSidki\ in~\cite{gupta-s:burnside,gupta-s:infinitep},
\ARozhkov\ in~\cite{rozhkov:aleshin}, \JFabrykowski\ and \NGupta\
in~\cite{fabrykowski-g:growth1,fabrykowski-g:growth2}, and more recently
the first and the third author
in~\cite{bartholdi-s:wpg,bartholdi:ggs,bartholdi:phd,sunik:phd} and
\DSegal\ in~\cite{segal:fgfimages}.  We will discuss some of these
examples in the sequel, after we give a definition that covers all of
them. Many of these examples were already presented in
Section~\ref{sec:branch:examples}.

All examples of finitely generated branch groups that we mentioned
by now are spinal groups, except for the examples of \PNeumann\
from~\cite{neumann:pride} (see Section\ref{sec:branch:examples}).
It is not known at present if the groups of \PNeumann\ are
conjugate to spinal groups
(Question~\ref{question:neumann-spinal}).

\section{Construction, basic tools and properties}

\subsection{Definition of spinal groups}
\begin{figure}[!ht]
\begin{center}
\includegraphics{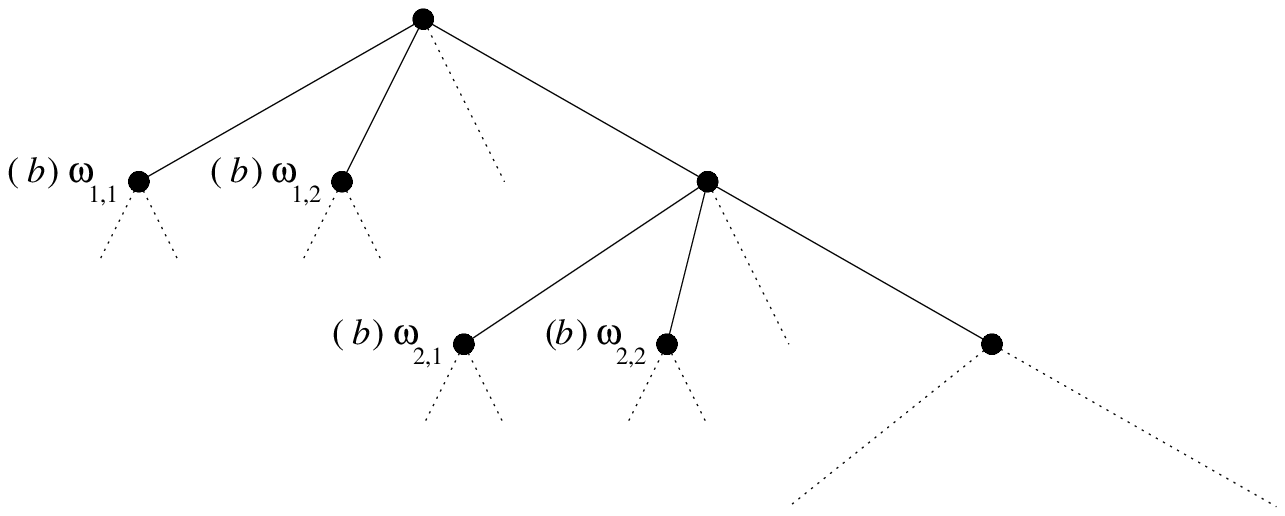}
\end{center}
\caption{The directed automorphism $b_\om$} \label{figure:b_om}
\end{figure}

Let $\om$ be a triple consisting of a group of rooted automorphisms
$A_\om$, a group $B$ and a doubly indexed family $\omover$ of
homomorphisms:
\[\om_{ij}: B \to\Sym(Y_{i+1}),\qquad i\in\N,\;j\in\{1,\dots,m_i-1\}.\]
Such a triple is called a \emdef{defining triple}. Each $b \in B$ defines
a directed automorphism $b_\om$ whose portrait is depicted in the diagram
in Figure~\ref{figure:b_om}.

Therefore, $B_\om=\setsuch{b_\om}{b\in B}$ is a set of directed
automorphisms. We can think of $B$ as acting on the tree $\tree$
by automorphisms. We define now the group $G_\om$, where $\om$ is
a defining triple, as the group of tree automorphisms generated by
$A_\om$ and $B_\om$, namely $G_\om = \langle A_\om \cup B_\om
\rangle$. We call $A_\om$ the \emdef{rooted part}, or the
\emdef{root group}, and $B$ the \emdef{directed part} of $G_\om$.
The family $\omover$ is sometimes referred to as the
\emdef{defining family of homomorphisms}.

\begin{figure}[!ht]
  \begin{center}
    \includegraphics{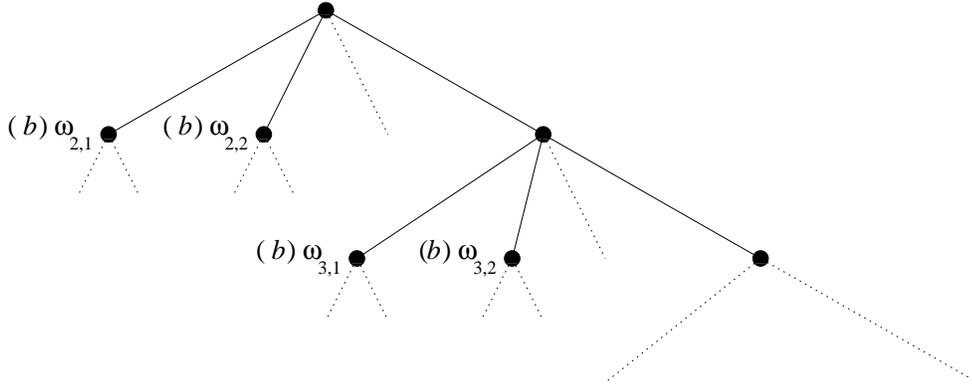}
  \end{center}
  \caption{The directed automorphism $b_{\sigma\om}$} \label{figure:b_sigma}
\end{figure}

Let us define the \emdef{shifted triple} $\sigma^r\om$, for $r\in\N_+$.
The triple $\sigma^r\om$ consists of the group $A_{\sigma^r\om}$ of rooted
automorphisms of $\tree_r$ defined by
\[A_{\sigma^r\om} = \left\langle \bigcup_{j=1}^{m_r-1} (B)\om_{rj} \right\rangle,\]
the same group $B$ as in $\om$, and the shifted family $\sigma^r\omover$
of homomorphisms
\[\om_{i+r,j}: B\to \Sym(Y_{i+r+1}), \qquad i\in\N,\qquad j\in\{1,\dots,m_{i+r}-1\}.\]
With the natural agreement that $\sigma^0\om=\om$, we see that
$\sigma^r\om$ defines a group $G_{\sigma^r\om}$ of tree automorphisms of
$\tree_r$ for each $r\in\N$ . Note that the diagram in
Figure~\ref{figure:b_sigma} describes the action of $b_{\sigma\om}$ on the
shifted tree $\tree^{(\sigma\Y)}=\tree_1$, and that $b_{\sigma\om}$ is
just the section of $b_\om$ at $m_1$.

We are ready now for the
\begin{definition}\label{defn:spinal}
  Let $\om=(A_\om,B,\omover)$ be a defining triple. The group
  \[G_\om= \langle A_\om \cup B_\om \rangle\]
  is called the \emdef{spinal group} defined by $\om$ if the
  following two conditions are satisfied:
  \begin{enumerate}
  \item \emdef{spherical transitivity condition}: $A_{\sigma^r\om}$ acts
    transitively on the corresponding alphabet $Y_{r+1}$, for all
    $r\in\N$.
  \item \emdef{strong kernel intersection condition}:
    \[\bigcap_{i\ge r}\bigcap_{j=1}^{m_i-1} \Ker(\om_{ij}) =
    1,\text{ for all } r.\]
  \end{enumerate}
\end{definition}

The spherical transitivity condition guarantees that $G_\om$ acts
spherically transitively on $\tree$, as well as that the same is true for
the actions of the shifted groups $G_{\sigma^r\om}$ on the corresponding
shifted trees. Similarly, the strong intersection condition guarantees
that the action of $B$ on $\tree$ is faithful, and that the same is true
for the actions of the shifted groups $G_{\sigma^r\om}$ on the
corresponding shifted trees.

The class of defining triples $\om$ that satisfy the above two
conditions will be denoted by $\Omega$. The above considerations
indicate that $\Omega$ is closed under the shift, i.e., if $\om$
is in $\Omega$ then so is any shift $\sigma^r\om$. This fact is
crucial in many arguments involving spinal groups, but we will
rarely mention it explicitly.

In the following subsections we introduce the tools and constructions we
use in the investigation of the spinal groups along with some basic
properties that follow quickly from the given considerations.

\subsection{Simple reductions}\label{subs:simplered}
The abstract group $B$ is canonically isomorphic to the group of tree
automorphisms $B_{\sigma^r\om}$, for any $r$, so that we will not make too
much difference between them and will frequently omit the index in the
notation. Letters like $b,b_1,b',\dots$ are exclusively reserved for the
nontrivial elements in $B$ and are called $B$-letters. Letters like
$a,a_1,a',\dots$ are exclusively reserved for the nontrivial elements in
$A_{\sigma^r\om}$, $r\in\N$, and are called $A$-letters. Note that the
groups $A_{\sigma^r\om}$, $r\in\N$ are not necessarily isomorphic but we
omit the index sometimes anyway.

The set $S_\om = (A_\om \cup B_\om)\setminus\{1\}$ is the \emdef{canonical
generating set} of $G_\om$. The generators in $A_\om\setminus\{1\}$ are
called \emdef{$A$-generators} and the generators in $B_\om\setminus\{1\}$
are called \emdef{$B$-generators}. Note that $S_\om$ does not contain the
identity and generates $G_\om$ as a monoid, since it is closed under
inversion.

We will not be very careful to distinguish the group elements in $G_\om$
from the words in the canonical generators that represent them. At first,
it is a sacrifice to the clarity of presentation, but our opinion is that
in the long run we only gain by avoiding useless distinctions.

Define the \emdef{length of an element $g$ of $G_\om$} to be the shortest
length of a word over $S$ that represents $g$, and denote this length by
$|g|$. There may be more than one word of shortest length representing the
same element.

Clearly, $G_\om$ is a homomorphic image of the free product $A_\om*B_\om$.
Therefore, every $g$ in $G_\om$ can be written in the form
\begin{equation} \label{eq:nform}
  [a_0]b_1a_1b_2a_2 \dots a_{k-1}b_k[a_k]
\end{equation}
where the appearances of $a_0$ and $a_k$ are optional. Relations of the
following $4$ types:
\[a_1a_2 \to 1, \quad a_3a_4 \to a_5, \quad b_1b_2 \to 1,\quad b_3b_4\to b_5,\]
that follow from the corresponding relations in $A$ and $B$ are called
\emdef{simple relations}. A \emdef{simple reduction} is any single
application of a simple relation from left to right (indicated above by
the arrows).  Any word of the form~(\ref{eq:nform}) is called a
\emdef{reduced word} and any word can be rewritten in unique reduced form
using simple reductions.  Among all the words that represent an element
$g$, the ones of shortest length are necessarily reduced, but those that
are reduced do not necessarily have the shortest length.

Note that the system of reductions described above is complete,
i.e., it always terminates with a word in reduced form and the
order in which we apply the reductions does not change the final
reduced word obtained by the reduction.

In some considerations one needs to perform cyclic reductions. A reduced
word $F$ over $S$ of the form $F=s_1us_2$ for some $s_1,s_2 \in S$, $u \in
S^*$, is \emdef{cyclically reduced} if the word $us_2s_1$ obtained from
$F$ by a cyclic shift is also reduced. If $F=s_1us_2$ is reduced but not
cyclically reduced, the word obtained form $us_2s_1$ after one application
of a simple reduction is said to be obtained from $F$ by one \emdef{cyclic
reduction}. On the group level the simple cyclic reduction described above
corresponds to conjugation by $s_1$.

\subsection{Level stabilizers}
In order to simplify the notation we denote the level stabilizer
$\Stab_{G_\om}(\LL_n)$ by $\Stab_\om(\LL_n)$. In the sequel we
often simplify notation by replacing $G_\om$ as a superscript or
subscript just by $\om$, and we do this without warning. Since
each element in $B$ fixes the first level, a word $u$ over $S$
represents an element in $\Stab_\om(\LL_1)$ if and only if the
word in $A$-letters obtained after deleting all $B$-letters in $u$
represents the identity element.

Further, $\Stab_\om(\LL_1)$ is the normal closure of $B_\om$ in $G_\om$,
$G_\om = \Stab_\om(\LL_1) \rtimes A_\om$, and $\Stab_\om(\LL_1)$ is
generated by the elements $b_\om^g=g^{-1}b_\om g$, for $b_\om$ in $B_\om$
and $g$ in $A_\om$.

Clearly,
$(b_\om)\psi=((b)\om_{1,1},(b)\om_{1,2},\dots,(b)\om_{1,m_1-1},b_{\sigma\om})$.
For any $a$ in $A_\om$, $(b_\om^a)\psi$ has the same components as
$(b_\om)\psi$ does, but in different positions depending on $a$. More
precisely:

\begin{lemma}
  For any $h$ in $\Stab_\om(\LL_1)$, $g$ in $A_\om$, $b$ in $B$ and
  $i\in\{1,\dots,m_1\}$, we have
  \begin{enumerate}
  \item $(h^g)\varphi_i = (h)\varphi_{i^{g^{-1}}}$.
  \item The coordinates of $(b^g)\psi$ are: $(b)\om_{1,j}$ at the
    coordinate $j^g$, for $j\in\{1,\dots,m_1-1\}$, and $b$ at $m_1^g$.
  \end{enumerate}
\end{lemma}

For example, if $a$ is the cyclic permutation $a=((m_1,\dots,2,1))$, the
images of $b_\om^{a^j}$ under various $\varphi_i^\om$ are given in
Table~\ref{table:phi}.
\begin{table}[!ht]
  \[\begin{tabular}{C|C|C|C|C|C|C}
    &\varphi_1  &\varphi_2    &\varphi_3    &\cdots &\varphi_{m_1-1} &\varphi_{m_1} \\[2pt]\hline
    b_\om       &(b)\om_{1,1} &(b)\om_{1,2} &(b)\om_{1,3} &\cdots &(b)\om_{1,m_1-1} &b_{\sigma\om}\\
    b_\om^{a}   &(b)\om_{1,2} &(b)\om_{1,3} &(b)\om_{1,4} &\cdots &b_{\sigma\om} &(b)\om_{1,1} \\
    b_\om^{a^2} &(b)\om_{1,3} &(b)\om_{1,4} &(b)\om_{1,5} &\cdots &(b)\om_{1,1} &(b)\om_{1,2} \\
    \vdots      &\vdots       &\vdots       &\vdots       &\ddots &\vdots &\vdots\\
    b_\om^{a^{m_1-2}} &(b)\om_{1,m_1-1} &b_{\sigma\om} &(b)\om_{1,1} &\cdots
    &(b)\om_{1,m_1-3} &(b)\om_{1,m_1-2} \\
    b_\om^{a^{m_1-1}} &b_{\sigma\om} &(b)\om_{1,1} &(b)\om_{1,2} &\cdots
    &(b)\om_{1,m_1-2} &(b)\om_{1,m_1-1}
  \end{tabular}\]
  \caption{The maps $\varphi_i$ associated with the permutation
    $a=((m_1,\dots,2,1))$} \label{table:phi}
\end{table}

Since the root group $A_\om$ acts transitively on $Y_1$ we get all
the elements from $A_{\sigma\om}$ and all the elements from
$B_{\sigma\om}$ in the image of $\Stab_\om(\LL_1)$ under
$\varphi_i$, for any $i$. Therefore, the shifted group
$G_{\sigma\om}$ is precisely the image of the first level
stabilizer $\Stab_\om(\LL_1)$ under any of the section
homomorphisms $\varphi_i$, $i\in\{1,\dots,m_1\}$.

Let $g$ be an element in $G_\om$. There are unique elements $h$ in
$\Stab_\om(\LL_n)$ and $a$ in $A_\om$ such that $g=ha$. Clearly,
$a$ is the vertex permutation of $g$ at the root, i.e.,
$a=(\emptyset)g$, and $(h)\psi=(g_1,\dots,g_{m_1})$ where $g_i$ is
the section automorphisms of $g$ at the vertex $y_i$, $i \in
\{1,\dots,m_1\}$. Therefore, the section $G_{y_i}$ is contained in
the image of $\Stab_\om(\LL_1)$ under $\varphi_i$, and we already
established that this image is the shifted group $G_{\sigma\om}$.
Therefore, for all $i=1,\dots,m_1$,
\[G_{y_i} = U_{y_i} = G_{\sigma\om},  \]
i.e., the sections on the first level, upper companion group and
the shifted $G_{\sigma\om}$ group coincide for spinal groups.

Further, the $n$-th upper companion group of $G_\om$ is the
shifted group $G_{\sigma^n\om}$ and the homomorphism
\[\psi_n: \Stab_\om(\LL_n) \to\prod_{i=1}^{m_1m_2\dots m_n} G_{\sigma^{n}\om}\]
given by
\[(g)\psi_n=((g)\varphi_{1\dots1},\dots,(g)\varphi_{m_1 \dots m_n}) =(g_{1\dots1},\dots,g_{m_1\dots m_n})\]
is an embedding.

We end this subsection with an easy proposition:
\begin{proposition}
  Every spinal group $G_\om$ is infinite.
\end{proposition}
\begin{proof}
  The proper subgroup $\Stab_\om(\LL_1)$ of $G_\om$ maps under $\varphi_1$
  onto $G_{\sigma\om}$, which is also spinal.
\end{proof}

The above proof is very simple, so let us use the opportunity here to
point out an important feature that is shared by many of the proofs
involving spinal groups. The proof does not work with a fixed sequence
$\om$, but rather involves arguments and facts about all the shifts of
$\om$ (note the importance of the fact that $\Omega$ is closed for
shifts). In other words, the group $G_\om$ is always considered together
with all of its companions, and the only way we extract some information
about some spinal group $G_\om$ is through such a synergic cooperation
between the group and its companions. Unavoidably, we also make
observations about the companion groups.

\subsection{Portraits}\label{subs:portraits}
The following construction corresponds to the constructions
exhibited
in~\cite{grigorchuk:gdegree,grigorchuk:pgps,bartholdi:upperbd},
but it is more general and allows different modifications and
applications.

A \emdef{profile} is a sequence $\overline{P} = (P_t)_{t\in\N}$ of
sets of automorphisms, called \emdef{profile sets}, where $P_t
\subseteq \Aut(\tree_t)$. We define now a \emdef{portrait} of an
element $f$ of $G_\om$ with respect to the profile $\mathcal P$.
The portrait is defined inductively as follows: if $f$ belongs to
the profile set $P_0$ then the portrait of $f$ is the tree that
consists of one vertex decorated by $f$; otherwise the portrait of
$f$ is the tree that is decorated by $a=(\emptyset)f$ at the root
and has the portraits of the sections $f_1$, $\dots$, $f_{m_1}$,
with respect to the shifted sequence $\sigma\overline{P}$ of
profile sets, hanging on the $m_1$ labelled vertices below the
root. Therefore, the portrait of $f$ is a subtree (finite or
infinite) of the tree on which $f$ acts, its interior vertices are
decorated by the corresponding vertex permutations of $f$, and its
leaves are decorated by elements in the chosen profile sets and
are equal to the corresponding sections.

\begin{figure}[!ht]
  \begin{center}
    \includegraphics{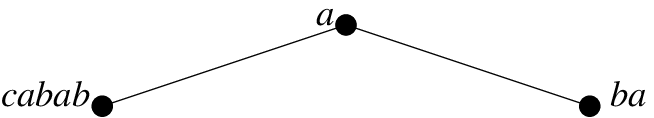}
  \end{center}
  \caption{Depth-1 decomposition of $g$}\label{figure:depth1}
\end{figure}

\begin{figure}[!ht]
  \begin{center}
    \includegraphics{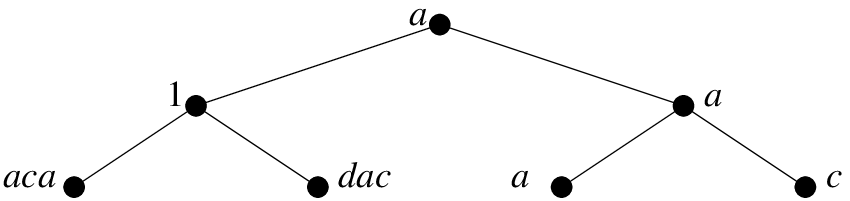}
  \end{center}
  \caption{Depth-2 decomposition of $g$}\label{figure:depth2}
\end{figure}

\begin{figure}[!ht]
  \begin{center}
    \includegraphics{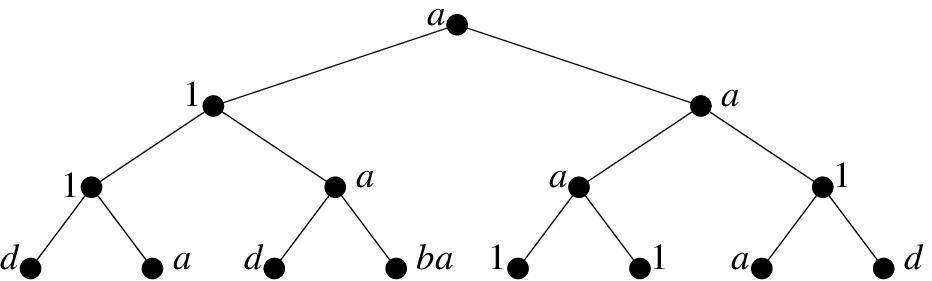}
  \end{center}
  \caption{Depth-3 decomposition of $g$}\label{figure:depth3}
\end{figure}

For example, if all the profile sets are empty, we obtain the
portrait representation of $f$ through its vertex permutations
that we already defined (see Figure~\ref{figure:f}). We sometimes
refer to this portrait as the \emdef{full portrait} of $f$. If
$P_t$ is empty for $t=0,\dots,r-1$ and equal to $G_{\sigma^r\om}$
for $t=r$, then the portrait of $f$ is the subtree of $\tree$ that
consists of the first $r$ levels, the vertices at the levels $0$
through $r-1$ are decorated by their vertex permutations and the
vertices at level $r$ are decorated by their corresponding
sections. Such a portrait is called the \emdef{depth-$r$
decomposition} of $f$. The depth-1, depth-2 and depth-3
decompositions of the element $g=abacadacabadac$ in $\Gg$ are
given in the Figure~\ref{figure:depth1},
Figure~\ref{figure:depth2} and Figure~\ref{figure:depth3}.
The decompositions can be easily calculated from Table~\ref{table:firstgg}
that describes $\varphi_1$ and $\varphi_2$ (recall the definition of $\Gg$
from Subsection~\ref{subs:firstgg}.
\begin{table}[!ht]
  \[\begin{tabular}{C|C|C}
        &\varphi_1  &\varphi_2 \\[2pt]\hline
    b   & a         & c \\
    c   & a         & d \\
    d   & 1         & b \\ \hline
    b^a & c         & a \\
    c^a & d         & a \\
    d^a & b         & 1
  \end{tabular}\]
  \caption{The table for $\varphi_1$ and $\varphi_2$ in $\Gg$}\label{table:firstgg}
\end{table}
The calculations follow, in which we identified the elements and their
images under the decomposition map $\psi$:
\begin{gather*}
 g = abacadacabadac = b^acd^acb^adc^aa = (cabac1d,ad1daba)a = (cabab,ba)a\\
 cabab = cb^ab = (aca,dac) \\
 ba = (a,c)a\\
 aca = c^a = (d,a)\\
 dac = dc^aa = (d,ba)a\\
 a = (1,1)a.
\end{gather*}

Let us note that, in general, the leaves of a portrait do not have to be
all at the same level, and it is possible that some paths from the root
end with a leaf and some are infinite. The various types of portraits
carry important information about the elements they represent, which is
not surprising, since different elements have different portraits.

Consider the construction of a portrait more closely. Let $f$ be
represented by a reduced word
\[F = [a_0]b_1a_1b_2a_2 \dots a_{k-1}b_k[a_k]\]
over $S$. We rewrite $f=F$ in the form
\begin{align}
  f &= b_1^{[a_0^{-1}]}b_2^{([a_0]a_1)^{-1}} \dots b_k^{([a_0]a_1\dots a_{k-1})^{-1}}
  [a_0]a_1\dots a_{k-1}[a_k] \notag \\
  &= b_1^{g_1}b_2^{g_2} \dots b_k^{g_k}g, \label{eq:bform}
\end{align}
where $g_i=([a_0]a_1\dots a_{i-1})^{-1}$ and $[a_0]a_1\dots
a_{k-1}[a_k]=g$ are in $A_\om$. Next, using the homomorphisms
$\om_{1,j}$, $j\in\{1,\dots,m_1-1\}$, and a table similar to
Table~\ref{table:phi} (but for all possible $a$) we compute the
(not necessarily reduced) words $\overline{F_1}, \dots,
\overline{F_{m_1}}$ representing the first level sections
$f_1,\dots,f_{m_1}$, respectively. Then we reduce these $m_1$
words using simple reductions and obtain the reduced words
$F_1,\dots,F_{m_1}$. Note that, on the word level, the order in
which we perform the reductions is unimportant since the system of
simple reductions is complete.

Let us consider any of the words
$\overline{F_1},\dots,\overline{F_{m_1}}$, say $\overline{F_j}$, and study
its possible content. Each factor $b^g$ from~(\ref{eq:bform}) contributes
to the word $\overline{F_j}$ either
\begin{itemize}
\item
one appearance of the $B$-letter $b$, in case $j=m_1^g$ or
\item
one appearance of the $A$-letter $(b)\om_{1,j^{g^{-1}}}$, in case $j \neq
m_1^g$ and $b \not\in \Ker(\om_{1,j^{g^{-1}}})$ or
\item
the empty word, in case $j \neq m_1^g$ and $b \in
\Ker(\om_{1,j^{g^{-1}}})$.
\end{itemize}
Therefore, the length of any of the reduced words $F_1,\dots,F_{m_1}$ does
not exceed $k$, i.e\ does not exceed $(n+1)/2$ where $n$ is the length of
$F$.  As a consequence we obtain the following:

\begin{lemma}
  For any $f \in G_\om$,
  \[|f_i|\leq (|f|+1)/2.\]
\end{lemma}

\begin{figure}[!ht]
  \begin{center}
    \includegraphics{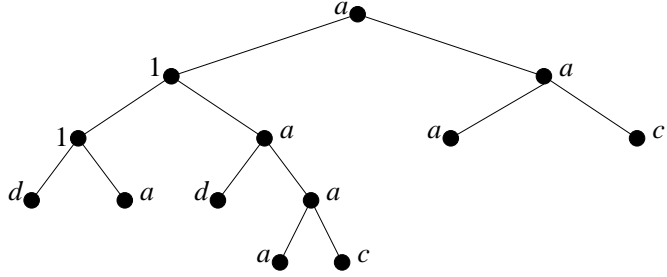}
  \end{center}
  \caption{The canonical portrait of $g$}\label{figure:can}
\end{figure}

The \emdef{canonical profile} is the profile
$\overline{P}=(S_{\sigma^i\om} \cup \{1\})_{i \in \N}$, whose profile set,
on each level, consists of the canonical  generators together with the
identity. For example the canonical portrait of $g = abacadacabadac$ is
given in Figure~\ref{figure:can}.

\begin{corollary}
  The depth of the canonical portrait of a word $w$ over $S_\om$ is
  no larger than $\lceil \log_2(|w|) \rceil + 1$.
\end{corollary}

\section{\GG\ groups} \label{sec:grigorchukexamples}
The \GG\ groups are natural generalizations of the first
Grigorchuk group $\Gg$ from~\cite{grigorchuk:burnside} and of the
groups $G_\om$ from~\cite{grigorchuk:gdegree,grigorchuk:pgps}. The
idea of these examples is based on the strong covering property
(see Definition~\ref{defn:GG}). The section presents uncountable
family of spinal groups of \GG\ type. All of them are finitely
generated, just-infinite, residually finite, centerless, amenable,
not elementary amenable, recursively but not finitely presented
torsion groups of intermediate growth (a definition of
intermediate growth will come later in
Chapter~\ref{chapter:growth}). For proofs
see~\cite{grigorchuk:gdegree,grigorchuk:pgps}, as well as the
later chapters in this text.

These examples were generalized in~\cite{bartholdi-s:wpg,sunik:phd}. The
generalized examples share many of the properties of Grigorchuk groups
mentioned above.

In all Grigorchuk examples the homomorphisms $\om_{i,j}$ are trivial for
all $j \neq 1$, and we denote $\om_i=\om_{i,1}$, for all $i \geq 1$.
Therefore, the only homomorphisms in the defining family of homomorphisms
$\omover$ that we need to specify are the homomorphisms in the sequence
$\om_1\om_2\om_3\dots$, which we call the \emdef{defining sequence of the
triple $\om$}, and we avoid complications in our notation by simply
writing $\omover = \om_1\om_2\om_3\dots$ .

\subsection{Grigorchuk $p$-groups}\label{subs:gpgps}
The Grigorchuk 2-groups, which are a natural generalization of $\Gg$, were
introduced in~\cite{grigorchuk:gdegree}. They act on the rooted binary
tree $\tree^{(2)}$. The rooted group $A=\{1,a\}$ and the group
$B=\{1,b,c,d\}=\Z/2\Z \times \Z/2\Z$ are still the same as for $\Gg$.
There are three non-trivial homomorphisms from $B$ to $\Sym(2)=\{1,a\}$,
and we denote them as follows:
\begin{align*}
  0 &= \left(\begin{matrix}1 & b & c & d \\ 1 & a & a & 1\end{matrix}\right),\\
  1 &= \left(\begin{matrix}1 & b & c & d \\ 1 & a & 1 & a\end{matrix}\right),\\
  2 &= \left(\begin{matrix}1 & b & c & d \\ 1 & 1 & a & a\end{matrix}\right).
\end{align*}
Note that the vectors representing the elements of $P$ also
correspond to the lines in the projective 2-dimensional space over
the finite Galois field ${\mathbb F}_p$. A Grigorchuk 2-group is
defined by any infinite sequence of homomorphisms
$\omover=\om_1\om_2\om_3\dots$, where the homomorphisms $\om_1$,
$\om_2$, $\om_3$, $\dots$ come from the set of homomorphisms
$H=\{0,1,2\}$ and each homomorphism from $H$ occurs infinitely
many times in $\omover$. In this setting, the first Grigorchuk
group $\Gg$ is defined by the periodic sequence of homomorphisms
$\omover=012012\dots$ .

Grigorchuk $p$-groups, introduced in 1985 in~\cite{grigorchuk:pgps}, act
on the rooted $p$-ary tree $\tree^{(p)}$, for $p$ a prime. The rooted
group $A= \langle a \rangle$ is the cyclic group of order $p$ generated by
the cyclic permutation $a=((1,2,\dots,p))$ and the group $B$ is the group
$\Z/p\Z \times \Z/p\Z$. Denote by $\begin{bmatrix}
  u\\v \end{bmatrix}$ the homomorphism from $B$ to $\Sym(p)$ sending
$(x,y)$ to $a^{ux+vy}$ and let $P$ be the set of homomorphisms
\[P = \left\{\;
  \begin{bmatrix}0\\1\end{bmatrix},   \;\begin{bmatrix}1\\1\end{bmatrix},
  \; \begin{bmatrix}2\\1\end{bmatrix}, \;\dots, \;
  \begin{bmatrix}p-1\\1\end{bmatrix},
  \;\begin{bmatrix}1\\0\end{bmatrix} \;\right\}.\]

A Grigorchuk $p$-group is defined by any infinite sequence of
homomorphisms $\omover=\om_1\om_2\om_3\dots$, where the homomorphisms
$\om_1$, $\om_2$, $\om_3$, $\dots$ come from the set $P$ and each
homomorphism from $P$ occurs infinitely many times in $\omover$.

Only those Grigorchuk $p$-groups that are defined by a recursive sequence
$\omover$ have solvable word problem (see~\cite{grigorchuk:gdegree}). The
conjugacy problem is solvable under the same condition (see
Chapter~\ref{chapter:word}).

\subsection{\GG\ groups} \label{subs:gg}
The following examples generalize the class of Grigorchuk groups
from the previous subsection.  A subclass of the class of groups
we are about to define is the subject of~\cite{bartholdi-s:wpg}
and the general case is considered in~\cite{sunik:phd}.

\begin{figure}[!ht]
  \begin{center}
    \includegraphics{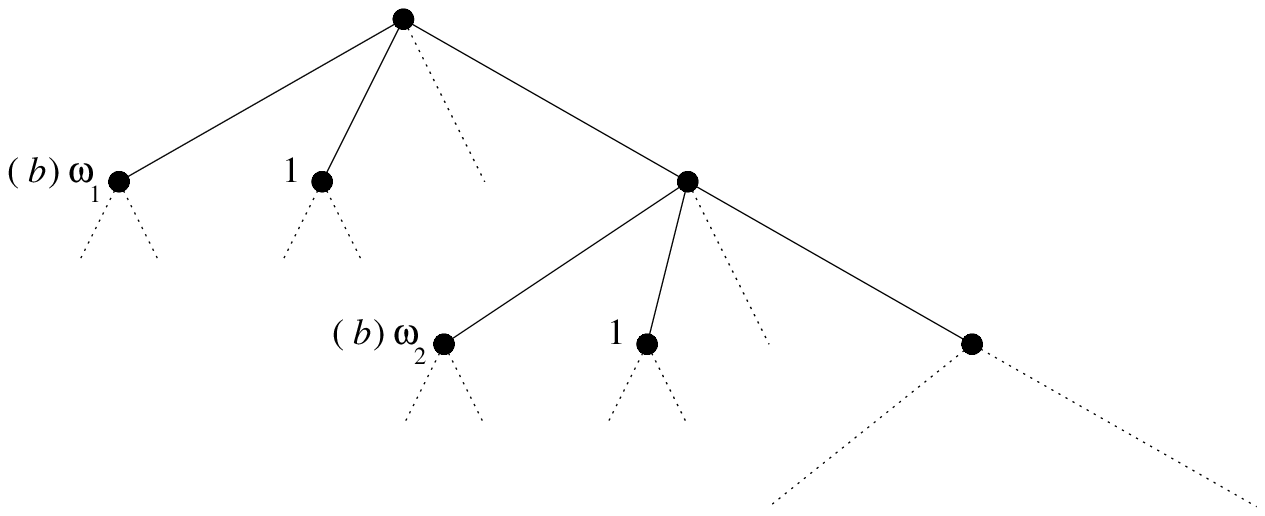}
  \end{center}
  \caption{The directed automorphism $b_\om$ in \GG\ groups}
  \label{figure:b_gg}
\end{figure}

We are now going to specify the triples $\om$ that define the groups in
the class of \GG\ groups. As before, the homomorphisms $\om_{i,j}$ are
trivial for all $j \neq 1$, and we denote $\om_i=\om_{i,1}$, for all $i
\geq 1$. For $b \in B$, the corresponding directed automorphism $b_\om$ is
given by the diagram in Figure~\ref{figure:b_gg}

For all positive $r$, we have
\[A_{\sigma^r\om} = (B)\om_r\]
and we denote
\[K_i = \Ker(\om_i).\]
Therefore, for positive $r$, each of the rooted groups $A_{\sigma^r\om}$
is a homomorphic image of $B$.

\begin{definition}\label{defn:GG}
  Let $G_\om$ be a spinal group defined by the triple
  $\om=(A_\om,B,\omover)$ with $\om_{i,j}=1$ whenever $j\neq 1$.
  The group $G_\om$ is a \GG\ group if
  the following \emdef{strong covering condition} is satisfied:
  \[\bigcup_{i=r}^\infty K_i = B, \qquad \text{for all } r.\]
\end{definition}

The strong kernel intersection condition is equivalent to a statement that
$(b)\om_i$ is nontrivial for infinitely many indices, while the strong
covering conditions is equivalent to a statement that $(b)\om_i$ is
trivial for infinitely many indices. The class of triples $\om$ that
satisfy the above conditions and thus define \GG\ groups will be denoted
by $\Omhat$. Clearly, $\Omhat$ is closed for shifting under $\sigma$, and
the use of this fact is essential but hardly ever emphasized.

Let us show an easy way to build examples of \GG\ groups. Start
with a group $B$ that has a covering by a family of proper normal
subgroups $\setsuch{N_\alpha}{\alpha\in I}$ of finite index and
with trivial intersection, i.e.,
\[\bigcup_{\alpha \in I} N_\alpha=B\text{ and }\bigcap_{\alpha \in I}
N_\alpha = 1.\] Choose a sequence of normal subgroups $N_1N_2N_2\dots$,
where $N_i$ come from the above family of normal subgroups, such that the
strong intersection and covering conditions hold, i.e.,
\[\bigcup_{i=r}^\infty N_i = B, \text{ and }\bigcap_{i=r}^\infty N_i = 1,
    \qquad \text{for all } r.\]
For each $i$, let the factor group $B/N_i$ act transitively and faithfully
as a permutation group on some alphabet $Y_{i+1}$, and let $\om_i$ be the
natural homomorphism from $B$ to the permutation group $B/N_i$ followed by
the embedding of $B/N_i$ in the symmetric group $\Sym(Y_{i+1})$. Choose a
group $A_\om$ and an alphabet $Y_1$ on which $A_\om$ acts transitively and
faithfully. The triple $\om$ that consist of $A_\om$, $B$ and the sequence
of homomorphisms $\omover=\om_1\om_2\om_3\dots$ defines a spinal group
$G_\om$ which is a \GG\ group acting on the tree $\treeY$. Clearly, the
strong intersection and covering conditions are satisfied since $K_i=N_i$,
and the spherical transitivity condition is satisfied since
$A_{\sigma^i\om}= B/N_i$.

It is of special interest to consider the case when $B$ is finite.
The family of all groups $B$ that have a covering by a finite
family of proper normal subgroups can be characterized, according
to a theorem of \MBrodie\, \RChamberlain\ and \LCKappe\
from~\cite{brodie-c-k:coversn}, as the family of those finite
groups that have $\Z/p\Z \times \Z/p\Z$ as a factor group, for
some prime $p$ (actually the theorem holds in general and
characterizes both finite and infinite groups $B$ that have
coverings with a finite family of proper normal subgroups). In
addition, to make sure that there exist a covering with trivial
intersection, we must add the condition that $B$ is not
subdirectly irreducible, i.e., the intersection of all non-trivial
normal subgroups of $B$ must be trivial.

However, for some reason (aesthetics or something deeper) we might
want to restrict our attention only to \GG\ groups that act on
regular trees. For example, such a case appears
in~\cite{bartholdi-s:wpg} where the added restriction is that all
the factors $B/N_i=A_{\sigma^i\om}$ are isomorphic to a fixed
group $A$, and they all act in the same way on the same alphabet
$\{1,2,\dots,m\}$. The \GG\ groups defined in this way act on the
regular tree $\tree^{(m)}$.  It is fairly easy to construct
examples of \GG\ groups with this added restrictions in case of an
abelian group $B$. Any group that is a direct product of proper
powers of cyclic groups can be used as $B$ in the above
construction, i.e., any group of the form
\[B = (\Z/n_1\Z)^{k_1} \times (\Z/n_2\Z)^{k_2} \times \dots (\Z/n_s\Z)^{k_s}\]
with $k_1,\dots,k_s \geq 2$ has a family of normal subgroups of the
required type such that all the factors are isomorphic to
\[A = (\Z/n_1\Z)^{k_1-1} \times (\Z/n_2\Z)^{k_2-1} \times \dots (\Z/n_s\Z)^{k_s-1}.\]
In particular, we see that any finite abelian group can be used in the
role of the rooted group $A$.

Characterizing the family of finite groups $B$ that have a
covering with a family of proper normal subgroups with trivial
intersection and such that all factors are isomorphic is an
interesting problem. The smallest known non-abelian example so far
was communicated to the authors by Derek Holt through the Group
Pub Forum (see \texttt{http://www.bath.ac.uk/\~{}masgcs/gpf.html})

Let $B= \langle b_1,b_2,b_3,b_4,b_5,b_6,x_{12},x_{34}\rangle$ where
$b_1,b_2,b_3,b_4,b_5,b_6$ all have order 3 and commute with each other,
$x_{12}$ and $x_{34}$ have order 2, commute and
\[b_i^{x_{jk}} = \begin{cases} b_i\text{ if }i\in\{j,k\},\\
  b_i^{-1}\text{ otherwise}.\end{cases}\]
In other words, $B$ is the semidirect product $(\Z/3\Z)^6 \rtimes
(\Z/2\Z)^2$ where $(\Z/2\Z)^2=\langle x_{12},x_{34} \rangle$, $x_{12}$
fixes the first two and acts by inversion on the last four coordinates of
$(\Z/3\Z)^6$ , $x_{34}$ fixes the middle two and acts by inversion on the
other four coordinates and, consequently, $x_{56}=x_{12}x_{34}$ fixes the
last two and inverts the first four coordinates.

The following $12$ subgroups are normal in $B$, their intersection is
trivial, their union is $B$, and each factor is isomorphic to the
symmetric group $\Z/3Z \rtimes \Z/2\Z = \Sym(3)$:
\[\begin{matrix}
\langle b_1,b_3,b_4,b_5,b_6,x_{12} \rangle,& \langle
b_1,b_2,b_3,b_5,b_6,x_{34} \rangle,& \langle b_1,b_2,b_3,b_4,b_5,x_{56}
\rangle, \\ \langle b_2,b_3,b_4,b_5,b_6,x_{12} \rangle,& \langle
b_1,b_2,b_4,b_5,b_6,x_{34} \rangle,& \langle b_1,b_2,b_3,b_4,b_6,x_{56}
\rangle, \\ \langle b_1b_2,b_3,b_4,b_5,b_6,x_{12} \rangle,& \langle
b_1,b_2,b_3b_4,b_5,b_6,x_{34} \rangle,& \langle
b_1,b_2,b_3,b_4,b_5b_6,x_{56} \rangle, \\ \langle
b_1b_2^2,b_3,b_4,b_5,b_6,x_{12} \rangle,& \langle
b_1,b_2,b_3b_4^2,b_5,b_6,x_{34} \rangle,& \langle
b_1,b_2,b_3,b_4,b_5b_6^2,x_{56} \rangle.
\end{matrix}\]

The smallest example of a non-abelian group $B$ that has a
covering by normal subgroups with trivial intersection is the
dihedral group $D_6$ on 12 elements. Indeed, $D_6$ is covered by
its 3 normal subgroups, call them $N_1$, $N_2$ and $N_3$, of index
2. The intersection of these 3 groups is not trivial, but if we
include the center $Z(D_{6})$ in the covering, we do get trivial
intersection. The corresponding factors are $\Z/2\Z = D_6/N_1=
D_6/N_2 = D_6/N_3$ and the dihedral group on six elements
$D_3=D_6/Z(D_{6})$. The first three factors can act as the
symmetric group on $\{1,2\}$, and the last one can act either
regularly on $\{1,2,3,4,5,6\}$ by the right regular representation
or as the symmetric group of $\{1,2,3\}$.

We can define even more involved examples of \GG\ groups. For example, we
can let $B$ itself be the first Grigorchuk group $\Gg$. It has $7$
subgroups of index $2$, and $3$ of them cover the group. Since we need
trivial intersection we can use the level stabilizers to accomplish this.
Therefore, we can easily define a \GG\ in which the directed part itself
is a \GG\ group, for example $\Gg$.

It is clear that in the examples when $B$ is infinite the branching
indices of the trees on which the group acts is not bounded. On the other
hand, in case of a finite $B$, the branching indices are bounded by the
order of $B$ and they have to be divisors of the order of $B$, except for
the first branching index $m_1$, which can be arbitrarily large.

\section{\GGS\ groups}\label{sec:ggs}
The \GGS\ groups (Grigorchuk-Gupta-Sidki groups, the terminology
comes from~\cite{baumslag:cgt}) are natural generalizations of the
second Grigorchuk group from~\cite{grigorchuk:burnside} and the
Gupta-Sidki examples from~\cite{gupta-s:burnside}. They act on a
regular tree $\tree^{(m)}$, where in most of the examples we
present $m$ is prime or a prime power, and they have a special
\emdef{stabilization} property, namely if an element $g\in G$ is
not in the level stabilizer $\Stab_G(\LL_1)$ then the first level
sections of the power $g^m$ are either in the stabilizer or are
closer to be in the stabilizer than the original element $g$.

In all examples that we give here, except the general case considered by
Bartholdi in~\cite{bartholdi:phd}, the rooted part $A=\langle a \rangle$
is the cyclic group of order $m$ generated by the permutation
$a=((1,2,\dots,m))$. The group $B=\langle b \rangle$ is also a cyclic
group of order $m$. All homomorphisms $\om_{i,j}$ map the elements from
$B$ to powers of $a$ and $\om_{i,j}=\om_{i',j}$, for all indices.
Therefore, in order to define a spinal group we only need to specify a
vector $E=(\varepsilon_1,\dots,\varepsilon_{m-1})$, where $\varepsilon_j$
are integers, and let $(b)\om_{i,j}=a^{\varepsilon_j}$, for all indices.
The group defined by the vector $E$ is denoted by $G_E$.

\begin{figure}[!ht]
  \begin{center}
    \includegraphics{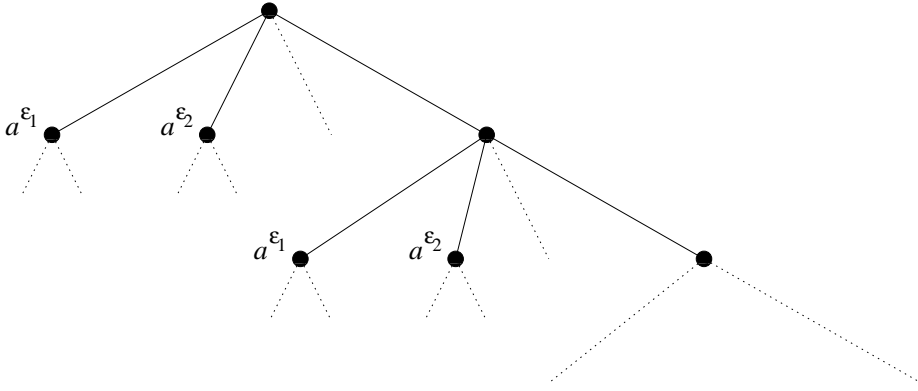}
  \end{center}
  \caption{The directed automorphism $b_E$} \label{figure:b_E}
\end{figure}

Therefore, $G_E=\langle a,b \rangle$ where $a$ is the rooted automorphism
defined by the cyclic permutation $a=((1,2,3,\dots,m))$, and the directed
automorphism $b=b_E$ is defined by the diagram in Figure~\ref{figure:b_E}.

In order for $E$ to define a spinal group it is necessary and sufficient
that $\gcd(\varepsilon_1,\varepsilon_2,\dots,\varepsilon_{m-1},m)=1$. Note
that in this situation $\om=\sigma\om$, $G_\om=G_{\sigma\om}$, etc., and
we have fractal groups.

\subsection{\GGS\ groups with small branching index}
In the case $m=2$, the only possible non-trivial vector $E=(1)$ defines
the infinite dihedral group. This is clear since this spinal group is
infinite and generated by two elements of order 2.

There are only three essentially different vectors in case $m=3$, and they
are $(1,0)$, $(1,1)$ and $(1,2)$. The corresponding groups were already
introduced as the Fabrykowski-Gupta group $\FGg$, Bartholdi-Grigorchuk
group $\BGg$ and Gupta-Sidki group $\GSg$.

We have already mentioned one example of \GGS\ group in case $m=4$ and it
is the second Grigorchuk group. The defining vector for the second
Grigorchuk group is $E=(1,0,1)$.

\subsection{Gupta-Sidki examples}
The Gupta-Sidki groups from~\cite{gupta-s:burnside} act on the regular
tree $\tree^{(p)}$ where $p$ is an odd prime and are defined by the vector
$E=(1,-1,0,0,\dots,0)$. The $p$-groups introduced
in~\cite{gupta-s:infinitep} are defined by the vector
$E=(1,-1,\dots,1,-1)$.
\subsection{More examples of \GGS\ groups}\label{subsec:moreggs}
The defining vector
$E=(\varepsilon_1,\varepsilon_2,\dots,\varepsilon_{m-1})$, where $m=p^n$
is a prime power, defines an infinite 2-generated $p$-group if and only if
\[\sum_{s \in O_k(m)} \varepsilon_s \equiv 0 \pmod{p^{k+1}},\]
for $k=0,\dots,n-1$, where
\[O_k(m) = \{\; p^k,\; 2p^k,\dots, \;(p^{n-k}-1)p^k \;\}.\]
The sufficiency in the above claim (in a more general setting) is
proved by \NGupta\ and \SSidki\ in~\cite{gupta-s:ext} and the
necessity by Vovkivsky in~\cite{vovkivsky:ggs}. The latter article
also shows that in case the defining vector $E$ does satisfy the
condition above, the obtained $p$-group is just-infinite, not
finitely presented branch group.

\subsection{General version of \GGS\ groups}\label{subsec:ggs}
\begin{figure}[!ht]
  \begin{center}
    \includegraphics{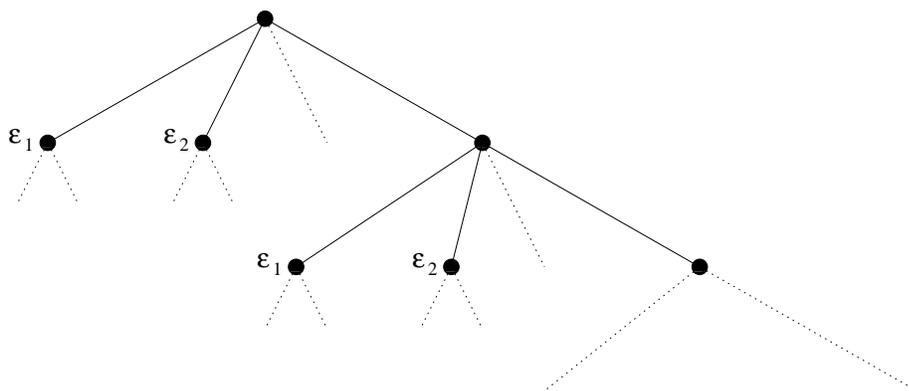}
  \end{center}
  \caption{The directed automorphism $b_E$ in \GGS\ groups}
  \label{figure:b_ggs}
\end{figure}

One chapter of the Ph.D.\ dissertation of Laurent
Bartholdi~\cite{bartholdi:phd} is devoted to a class of groups that comes
as a natural generalization of all of the previous examples of \GGS\
groups. The groups act on the tree $\tree^{(m)}$, where $m$ is arbitrary.
A defining vector
$E=(\varepsilon_1,\varepsilon_2,\dots,\varepsilon_{m-1})$ is a vector of
permutations of the alphabet $Y=\{1,2,\dots,m\}$ of the tree such that the
group of permutations $A=\langle
\varepsilon_1,\varepsilon_2,\dots,\varepsilon_{m-1} \rangle$ acts
transitively on $Y$, and the spinal group $G_E$ is generated by the rooted
automorphisms from $A$ together with the directed automorphism $b_E$
(simply written $b$) defined by the diagram in Figure~\ref{figure:b_ggs}

To see how these examples fit in our general scheme, note that the group
$B=\langle b \rangle$ is a cyclic group of order equal to the least common
multiple of the orders of
$\varepsilon_1,\varepsilon_2,\dots,\varepsilon_{m-1}$, and the
homomorphisms in the triple $\om$ are defined by $(b)\om_{i,j}=
\varepsilon_j$ for all indices.

The first author shows in his dissertation (\cite{bartholdi:phd})
that unless $m=2$ (in which case the group is the infinite
dihedral group) the \GGS\ groups are finitely generated,
residually finite, centerless and not finitely presented groups.
He also presents a sufficient and necessary condition for a \GGS\
group to be torsion group (see Chapter~\ref{chapter:torsion}).

\begin{theorem}\label{theorem:ggswb}
  Let $G$ be a \GGS\ group. Then $G$ is weakly regular branch on the
  non-trivial subgroup
  \[K=\big(\mho_{|A|}(G')\big)',\]
  where $\mho_{|A|}(G')$ denotes the subgroup of $G'$ generated by
  the $|A|$-th powers.
\end{theorem}
\begin{proof}[Sketch of a proof]
  $G'$ contains all elements of the form $[b,b^a]$, which can be
  written through $\psi$ as
  \[ (a_1,\dots,[a_i,b],a_{i+1},\dots,[b,a_m]). \]
  Write $n=|A|$. then
  $\mho_n(G')$ contains all elements of the form
  $(1,\dots,[a_i,b]^n,1,\dots,[b,a_m]^n)$, and its derived subgroup
  contains all elements of the form $(1,\dots,1,[[b,a]^n,[b,a']^n])$.
  We therefore have $1\times\dots\times1\times K \leq \psi(K)$.

  To show that $K$ is not trivial, we first argue that $G'$ is
  non-trivial, since it has finite index in an infinite group. Then,
  since the torsion has unbounded order in $G$, $\mho_n(G')$ is
  non-trivial; now this last group has a subgroup mapping onto $G$
  (namely $\Stab_G(\LL_N)$ for large $N$), and therefore cannot
  be abelian.
\end{proof}

\part{Algorithmic Aspects}
\date{Nov 3 27, 2002}
\chapter{Word and Conjugacy Problem}\label{chapter:word}
Branch groups have good algorithmic properties. A universal
algorithm solving the word problem for \GG\ and \GGS\ groups and
many other branch groups was mentioned in
\cite{grigorchuk:burnside} and described
in~\cite{grigorchuk:gdegree} (see also
\cite{grigorchuk:notEG,grigorchuk:bath}). This algorithm is very
fast and needs a minimal amount of space and we will describe it
below.

However, there are branch groups with unsolvable word problem. For
instance, the following claim is proved in \cite{grigorchuk:gdegree}

\begin{theorem}
The group $G_\om$ has solvable word problem if and only if $\omover$ is a
recursive sequence.
\end{theorem}

Note that the proof in~\cite{grigorchuk:gdegree} is given for the
considered case of $2$-groups (see Section~\ref{subs:gpgps}), but
it can be easily extended to the other cases.

This result inspired the second author to use Kolmogorov complexity to
study word problems in branch groups and other classes of groups in
\cite{grigorchuk:kolmogorov}.

The generalized word problem  (is there an algorithm which for any
element $g \in G$ and a finitely generated subgroup $H \leq G$
given by a generating set decides if $g$ belongs to $H$) was
considered only recently in~\cite{grigorchuk-w:rigidity} for the
first Grigorchuk group $\Gg$ and it was shown that $\Gg$ has a
solvable generalized word problem.

The solution of the conjugacy problem for the basic examples of
branch groups came much later. First, Wilson and Zaleskii solved
the conjugacy problem in \GGS\ $p$-groups, for $p$ an odd prime,
by using the notion of Maltcev's conjugacy separability and
pro-$p$ methods (see \cite{wilson-z:conj}. Slightly later,
simultaneously and independently, Leonov in \cite{leonov:conj} and
Rozhkov in \cite{rozhkov:conj} solved the conjugacy problem for
$p=2$. The paper of Rozhkov deals only with the first Grigorchuk
group $\Gg$, while Leonov considers all 2-groups $G_\om$ from
\cite{grigorchuk:gdegree}. Also, the results of Leonov are
stronger, since upper bounds on the length of conjugating elements
are given in terms of depth.

The idea of Leonov and Rozhkov was developed in \cite{grigorchuk-w:conj}
in different directions. Still, there are some \GGS\ groups with unsettled
conjugacy problem. For example, the results from \cite{grigorchuk-w:conj}
do not apply to the group $\BGg$ (see~\ref{subsec:three}), which is not
branch but only weakly branch.

One of the improvements reached in \cite{grigorchuk-w:conj} is that the
conjugacy problem is solvable not only for the considered groups, but for
their subgroups of finite index (note that there are groups with solvable
conjugacy problem that contain subgroups of index 2 with unsolvable
conjugacy problem, see \cite{collins-m:conj})

Theorem~C from \cite{grigorchuk-w:conj} is the strongest result on the
conjugacy problem.  The corresponding algorithm in Corollary~C in the same
article uses the principle of Dirichlet and is quite different from the
Leonov-Rozhkov algorithm.

There is also a result on the isomorphism problem for the spinal 2-groups
$G_\om$ defined in \cite{grigorchuk:gdegree}. Namely, for every sequence
$\omover$ there are only countably many sequences $\omover'$ with $G_\om
\cong G_{\om'}$. Thus, there are uncountably many finitely generated
branch groups.

\section{The word problem}
We describe an algorithm that solves the word problem in the case of
spinal groups modulo an oracle that has effective knowledge of the
defining triple $\om$ in a sense that we make clear below. In fact, the
way in which we answer the question ``does the word $F$ over $S$ represent
the identity in $G_\om$?'' is by constructing the canonical portrait of
$F$ and we need tables similar to Table~\ref{table:phi}, for all possible
$a$, in order to do that.

Let us assume that we have an oracle $\mathcal O_\om(n)$ that has
\emdef{effective knowledge} of the first $n$ levels of the defining triple
$\om$, meaning that
\begin{enumerate}
\item
the multiplication table of $B$ is known,
\item
the permutation groups $A_\om$, $\dots$, $A_{\sigma^n\om}$ are known, in
the sense that we can actually perform the permutations,
\item
the homomorphisms in the first $n$ levels in $\omover$ are known, i.e. for
all $b \in B$, $i \in \{1,\dots,n\}$, $j \in \{1,\dots,m_i-1\}$, it is
known exactly what permutation in $A_{\sigma^i\om}$ is equal to
$(b)\om_{i,j}$.
\end{enumerate}
Note that if we have an oracle $\mathcal O_\om(n)$, then we also have
oracles $\mathcal O_{\sigma^t\om}(n-t)$, for $t \in \{0,\dots,n\}$.

\begin{proposition}
Let $G_\om$ be a spinal group and let the oracle $\mathcal O_\om(n)$ have
effective knowledge of the first $n$ levels of $\om$. The word problem in
$G_\om$ is solvable, by using the oracle, for all words of length at most
$2^{n+1}+1$.
\end{proposition}
\begin{proof}
If $F$ is not reduced we can reduce it since we have the multiplication
tables for $B$ and $A_\om$ (even in case $n=0$). So assume that $F$ is
reduced and let
\[ F = [a_0]b_1a_1 \dots a_{k-1}b_k[a_k]. \]
Note that if $[a_0]a_1 \dots a_{k-1}[a_k] \neq 1$ in $A_\om$ then $F$ does
not stabilize the first level of $\tree$ and does not represent the
identity.

We prove the claim by induction on $n$, for all spinal groups
simultaneously.

No reduced non-empty word $F$ of length at most 3 represents the identity.
This is true because either $F$ does not fix the first level or it
represents a conjugate of $b$, for some $b \in B$. This completes the base
case $n=0$.

Let $F$ be reduced, $4 \leq |F| \leq 2^{n+1}+1$ and $[a_0]a_1 \dots
a_{k-1}[a_k]=1$ in $A_\om$. Rewrite $F$ as
\[ F = b_1^{g_1}\dots b_k^{g_k}, \]
where $g_i = ([a_0]a_1 \dots a_{i-1})^{-1}$ and then use the oracle's
knowledge of the first level of $\om$ to construct tables similar to
Table~\ref{table:phi} and calculate the possibly unreduced words
$\overline{F_i}$, $i \in \{1,\dots,m_1\}$, that represent the first level
sections of $F$. Each of these words has length no greater than
$\frac{|F|+1}{2} \leq 2^n+1$ and we may apply the inductive hypothesis to
solve the word problem for each of them by using the oracle $\mathcal
O_{\sigma\om}(n-1)$. The word $F$ represents the identity in $G_\om$ if
and only if each of the the words $F_1$, $\dots$, $F_{m_1}$ represents the
identity in $G_{\sigma\om}$ which is clear since $(F)\psi =
(F_1,\dots,F_{m_1})$ and $\psi$ is an embedding.
\end{proof}

\begin{corollary}
If an oracle $\mathcal O_\om$ has effective knowledge of all the
levels of $G_\om$, the word problem in $G_\om$ is solvable.
\end{corollary}

For many spinal groups the converse is also true. Namely, if the word
problem is solvable, one can use the strategy
from~\cite{grigorchuk:gdegree} and construct test-words that would help to
recover the kernels of the homomorphisms used in the definition of $\om$,
which is then enough to recover the actual homomorphisms. For example, if
we know that $G_\om$ is one of the Grigorchuk 2-groups
from~\cite{grigorchuk:gdegree} and we have an oracle that solves the word
problem, we can check the test-words $(ab)^4$, $(ac)^4$ and $(ad)^4$ and
find out which one represents the identity, which then tells us that
$\om_1$ was 2, 1 or 0, respectively. Depending on this result, we
construct longer test-words that give us information about $\om_2$, and
then longer for $\om_3$, etc.

\section{The conjugacy problem in $\Gg$}
Two algorithms which solve the conjugacy problem for regular branch groups
satisfying some natural conditions involving length functions are
described in~\cite{grigorchuk-w:conj}. The first one accumulates the ideas
from~\cite{leonov:conj,rozhkov:conj} and is presented here in the simplest
form for the first Grigorchuk group $\Gg$. The second one, based on the
Dirichlet principle, is quite different and has more potential for
applications.

We describe now the first algorithm for the case of $\Gg$. Recall
that $\Gg$ is fractal regular branch group with $K \times K
\preceq K \trianglelefteq G$, where $K$ is the normal closure of
$[a,b]$. For $g,h \in \Gg$ we define
\[ Q(g,h) = \{ \; Kf \;|\; g^f = h \;\}. \]
Clearly, $Q(g,h) = \emptyset$ if and only if $g$ and $h$ are not conjugate
in $\Gg$.

\begin{theorem} The conjugacy problem is solvable for $\Gg$. \end{theorem}
\begin{proof}

\begin{lemma}\label{Q}
Let $f,g,h \in\Gg$ and let
\begin{equation} \label{g^f=h}
g^f=h.
\end{equation}

\begin{enumerate}
\item Let $g=(g_1,g_2), h=(h_1,h_2) \in\Stab_\Gg(\LL_1)$.

\begin{itemize}
\item[a)]
If $f \in\Stab_\Gg(\LL_1)$ and $f=(f_1,f_2)$ then (\ref{g^f=h}) is
equivalent to
\[ \begin{cases} g_1^{f_1}=h_1 \\ g_2^{f_2}=h_2 \end{cases} \]
\item[b)]
If $f \not\in\Stab_\Gg(\LL_1)$ and $f=(f_1,f_2)a$ then (\ref{g^f=h}) is
equivalent to
\[ \begin{cases} g_1^{f_1}=h_2 \\ g_2^{f_2}=h_1 \end{cases} \]
\end{itemize}

\item
Let $g=(g_1,g_2)a, h=(h_1,h_2)a \not\in\Stab_\Gg(\LL_1)$.
\begin{itemize}
\item[a)]
If $f \in\Stab_\Gg(\LL_1)$ and $f=(f_1,f_2)$ then (\ref{g^f=h}) is
equivalent to
\[ \begin{cases} (g_1g_2)^{f_1}=h_1h_2 \\ f_2 = g_2f_1h_2^{-1} \end{cases} \]
\item[b)]
If $f \not\in\Stab_\Gg(\LL_1)$ and $f=(f_1,f_2)a$ then (\ref{g^f=h}) is
equivalent to
\[ \begin{cases} (g_1g_2)^{f_1}=h_2h_1 \\ f_2 = g_1^{-1}f_1h_2 \end{cases} \]
\end{itemize}

\end{enumerate}
\end{lemma}

\begin{lemma}\label{z}
Let $x=(x_1,y_1)$ and $y=(y_1,y_2)$ be elements of $\Gg$. Let
\begin{gather*}
Q(x,y) = \{Kz_i| i\in I\} \\
Q(x_1,y_1) = \{Kz_j| j\in J\} \\
Q(x_2,y_2) = \{Kz_{j'}| j'\in J'\}
\end{gather*}
For every $i \in I$ there exists $j \in J$ and $j' \in J'$ such
that the element $z_{jj'}=(z_j,z_{j'})$ is in $\Gg$ and
$Kz_i=Kz_{jj'}$.
\end{lemma}

The meaning of Lemma~\ref{z} is that if we have $Q(x_1,y_1)$ and
$Q(x_2,y_2)$ already calculated, then we can calculate $Q(x,y)$.
We just need to check all pairs $(z_j,z_{j'})$ that are in $\Gg$
and chose one in each coset of $K$.

We describe an algorithm that solves the conjugacy problem in $\Gg$ by
calculating $Q(g,h)$, given $g$ and $h$.

Split $Q(g,h)$ as
\[ Q(g,h) = Q_{1a}(g,h) \cup Q_{1b}(g,h) \]
or as \[ Q(g,h) = Q_{2a}(g,h) \cup Q_{2b}(g,h) \] according to the
cases described in Lemma~\ref{Q}.

If we knew how to calculate $Q(g_1,h_1)$ and $Q(g_2,h_2)$ for the
case 1a, we could compute $Q_{1a}(g,h)$. Thus to prove the theorem
we need only to show that there is a reduction in length in each
of the cases. This is obvious for the cases 1a and 1b as
\[ |g_i| + |h_j| < |g|+|h|, \]
for $i,j\in\{1,2\}$, except when $|g|=|h|=1$, which case can be handled
directly. For the cases 2a and 2b we only have
\[ |g_1g_2| + |h_1h_2| \leq |g|+|h|, \]
and equality is possible only if the letter $d$ does not appear in
the word representing $g$ nor in the one representing $f$. In case
of an equality we repeat the argument with the pair
$(g_1g_2,h_1h_2)$ or $(g_1g_2,h_2h_1)$ depending on the case.

After at most three steps there must be some reduction in the length and
we may apply induction on the length.
\end{proof}

\date{October 27, 2002}
\chapter{Presentations and endomorphic presentations of branch groups}
\label{chapter:presentation} Branch groups are probably never
finitely presented and this is established for the known examples.
There are many approaches one can use to prove that a given branch
group cannot have a finite presentation and we will say more about
them below, along with some historical remarks.

Every branch group that has a solvable word problem, as all basic examples
do, has a recursive presentation.

An important discovery was made by \ILysionok\ in 1985, who showed
in~\cite{lysionok:pres} that the first Grigorchuk group can be
given by the presentation~(\ref{eq:lys:pres}). No due attention
was given to this fact for a long time, until the second author
used this presentation in~\cite{grigorchuk:notEG} to construct an
example of finitely presented amenable but not elementary amenable
group, thus answering a question of \MDay\ from~\cite{day:amen}.

The idea of \ILysionok\ was developed in a different direction
in~\cite{grigorchuk:bath}, where the presentations of the
form~(\ref{eq:lys:pres}), which are finite presentations modulo
the iteration of a single endomorphism, were called
$L$-presentations. It was shown in~\cite{grigorchuk:bath} that
some \GGS\ groups (for example the Gupta-Sidki $p$-groups, for $p
> 3$) have such presentations.

The study of $L$-presentations was continued, and some aspects
were completely resolved in~\cite{bartholdi:lpres}. First of all,
the notion of $L$-presentations was extended to the notion of
endomorphic presentation in a way that allows several
endomorphisms to be included in the presentation. On the basis of
this extended definition, a general result was obtained, claiming
that all regular branch groups satisfying certain natural
requirements have finite endomorphic presentations (see
Theorem~\ref{theorem:finitelp} below). It would be interesting to
answer the question what branch groups have $L$-presentations in
the sense of~\cite{grigorchuk:bath} and, in particular, resolve
the status of the Gupta-Sidki group $\GSg$
(Question~\ref{question:L-presentations}).

We finish the chapter by providing several examples of groups with
finite endomorphic presentations, mostly taken
from~\cite{bartholdi:lpres}.

\section{Non-finite presentability}\label{sec:presentation:hist}
In his early work, the second author proved in various ways that $\Gg$ is
not finitely presented. We review briefly these ideas, since they
generalize differently to various other examples.

\begin{theorem}
  The first Grigorchuk group $\Gg$ is not finitely presented.
\end{theorem}
\begin{proof}[First proof, from~\cite{grigorchuk:burnside}]
  More details can be found in~\cite{grigorchuk:bath}.
  Assume for contradiction that $\Gg$ is finitely presented, say as
  $\SPRES SR$. The Reidemeister-Schreier method then gives a
  presentation of $\Stab_\Gg(1)$ with relators $R\cup R^a$; writing for
  each relator $(r)\psi=(r_1,r_2)$ we obtain a presentation of $\Gg$
  with relators $\{r_1,r_2\}_{r\in R}$. Now since $|r_i|\le|r|/2$ for
  all $r\in R$ of length at least $2$ (using cyclic reductions), we
  obtain after enough applications of the above process a presentation
  of $\Gg$ with relations of length $1$, i.e. a free group.

  This contradicts almost every property of $\Gg$: that it is torsion,
  of subexponential growth, just-infinite, or contains elements
  $(x,1)$ and $(1,y)$ that commute.
\end{proof}

\begin{proof}[Second proof, from~\cite{grigorchuk:gdegree}]
  The set of groups with given generator set $\{a_1,\dots,a_m\}$ is a
  topological space, with the ``weak topology'': a sequence $(G_i)$ of
  groups converges to $G$ if for all radii $R$ the sequence of balls $B_{G_i}(R)$
  in the Cayley graphs of the corresponding groups stabilize to the ball $B_G(R)$.

  Assume for contradiction that $\Gg$ is finitely presented, say as
  $\SPRES SR$. Then for any sequence $G_i\to \Gg$ the groups $G_i$ are
  quotients of $\Gg$ for $i$ large enough. However, if one considers all
  Grigorchuk groups $\Gg_\omega$ from~\cite{grigorchuk:gdegree},
  defined through infinite sequences
  $\omover \in \{0,1,2\}^\N$, one notes that the map $\omega\to
  \Gg_\omega$ is continuous, with the Tychonoff topology on $\{0,1,2\}^\N$.
  There are therefore infinite groups converging to $\Gg$,
  which contradicts $\Gg$'s just-infiniteness.
\end{proof}

A third proof involves a complete determination of the presentation, and
of its Schur multiplier. The results are:
\begin{theorem}[Lysionok,~\cite{lysionok:pres}]
  The Grigorchuk group $\Gg$ admits the following presentation:
  \begin{equation}\label{eq:lys:pres}
    \Gg = \PRES{a,c,d}{\phi^i(a^2),\phi^i(ad)^4,
    \phi^i(adacac)^4\;(i\ge0)},
  \end{equation}
  where $\phi:\{a,c,d\}^*\to\{a,c,d\}^*$ is defined by
  $\phi(a)=aca,\phi(c)=cd,\phi(d)=c$.
\end{theorem}

\begin{theorem}[Grigorchuk,~\cite{grigorchuk:bath}]\label{theorem:schurmult}
  The Schur multiplier $H_2(\Gg,\Z)$ of the first Grigorchuk group is
  $(\Z/2)^\infty$, with basis
  $\{\phi^i[d,d^a],\phi^i[d^{ac},d^{aca}]\}_{i\in\N}$.
\end{theorem}

Given a presentation $G=F/R$, where $F$ is a free group, the Schur
multiplier may be computed as $H_2(G,\Z)=(R\cap[F,F])/[F,R]$
(see~\cite{karpilovsky:schur} or~\cite{brown:cohomology} for details).
This implies instantly that $\Gg$ may not be finitely presented, and
moreover that no relation can be omitted.

Another approach that deserves attention is demonstrated
in~\cite{gupta:recursive}. Namely, certain recursively presented groups
are constructed there and the strategy is to build an increasing sequence
of normal subgroups $(R_n)_{n\in\N}$ of the free group $F$ whose union is
the kernel $R$ of the presentation of the constructed group $G$ as $F/R$.

The most general result, at least for spinal groups, is given in
the third author's dissertation~\cite{sunik:phd}, and it follows
the ``topological'' approach from~\cite{grigorchuk:gdegree}, but
without the explicit use of the topology.

\begin{theorem} \label{theorem:presentability}
  Let $\mathcal C$ be a class of groups that is closed under
  homomorphic images and subgroups (of finite index) and
  $\om=(A_\om,B,\omover)$ be a sequence that defines a spinal group in
  $\mathcal C$. Further, assume that, for every $r$, there exists a
  triple $\eta^{(r)}$ of the form $\eta^{(r)}=(A_{\sigma^r\om},B,\overline{\eta})$,
  where $\overline{\eta}^{(r)}$ is a doubly indexed family of homomorphisms
  \[ \eta_{i,j}: B \rightarrow \Sym(Y_{j+1}), \quad i\in\{r+1,r+2,\dots\},
        \quad j\in\{1,\dots,m_i-1\} \]
  defining a group of tree automorphisms (not necessarily spinal)
  $G_{\eta^{(r)}}$ that acts on the shifted tree $\tree^{(\sigma^r Y)}$
  and is not in $\mathcal C$. Then, the spinal group $G_\om$ is not
  finitely presented.
\end{theorem}
\begin{proof}
Assume, on the contrary, that $G_\om$ is finitely presented.

Further, assume that the length of the longest relator in the
finite presentation of $G_\om$ is no greater than $2^{n+1}+1$. If
$\om'$ is any triple (not necessarily defining a spinal group)
that agrees with $\om$ on the first $n$ levels then any word of
length no greater than $2^{n+1}+1$ representing the identity in
$G_\om$ represents the identity in $G_{\om'}$. Thus all relators
from the finite presentation of $G_\om$ represent the identity in
$G_{\om'}$, so that $G_{\om'}$ is a homomorphic image of $G_\om$,
and, therefore, a member of $\mathcal C$.

Define $\om'$ so that it agrees with $\om$ on the first $n$ levels
and it uses the definition of $\eta^{(n)}$ to define the rest of
the levels (just concatenate the definition of $\eta^{(n)}$ to the
definition of the first $n$ levels of $\om$). Since $G_{\om'}$ is
a member of $\mathcal C$, so is each of its upper companion
groups, including $G_{\eta^{(n)}}$, a contradiction.
\end{proof}

Since the class of torsion groups is closed for subgroups and
images and since it is fairly easy to construct triples
$\eta^{(r)}$ that define groups containing elements of infinite
order, we obtain the following:
\begin{corollary} \label{theorem:spinalnotfp}
  No torsion spinal group is finitely presented.
\end{corollary}

\section{Endomorphic presentations of branch groups}
The recursive structure of branch groups appears explicitly in their
presentations by generators and relators, and such presentations have been
described since the mid-80's for the first example, the Grigorchuk group.

In this section, we will mainly consider finitely generated,
regular branch groups, the reason being that the regularity of
presentations becomes much more apparent in these cases. The main
result is best formulated in terms of ``endomorphic
presentations''~\cite{grigorchuk:bath,bartholdi:lpres}:
\begin{definition}
  An \emdef{endomorphic presentation} is an expression of the form
  \[L = \LPRES SQ\Phi R,\]
  where $S$ is an alphabet (i.e., a set of symbols), $Q,R\subset F_S$
  are sets of reduced words (where $F_S$ is the free group on $S$),
  and $\Phi$ is a set of group homomorphisms $\phi:F_S\to F_S$.

  $L$ is \emdef{finite} if $Q,R,S,\Phi$ are finite. It is
  \emdef{ascending} if $Q$ is empty.

  $L$ gives rise to a group $G_L$ defined as
  \[G_L = F_S \Big/\Big\langle Q\cup
  \bigcup_{\phi\in\Phi^*}(R)\phi\Big\rangle^\#,\] where
  $\langle\cdot\rangle^\#$ denotes normal closure and $\Phi^*$ is the
  monoid generated by $\Phi$, i.e., the closure of $\{1\}\cup\Phi$
  under composition.

  As is customary, we identify the endomorphic presentation $L$ and the
  group $G_L$ it defines.

  An endomorphic presentation that has exactly one homomorphism in $\Phi$
  is called \emdef{$L$-presentation}.
\end{definition}

The geometric interpretation of endomorphic presentations in the
context of branch groups is the following: one has a finite
generating set ($S$), a finite collection of relations, some of
which ($R$) are related to the branching and therefore can be
``moved from one tree level down to the next'' by endomorphisms
($\Phi$).

The main result of this chapter is:
\begin{theorem}[\cite{bartholdi:lpres}]\label{theorem:finitelp}
  Let $G$ be a finitely generated, contracting, regular
  branch group. Then $G$ has a finite endomorphic presentation.
  However, $G$ is not finitely presented.
\end{theorem}

The motivations in studying group presentations of branch groups are the
following:
\begin{itemize}
\item They exhibit a regularity that closely parallels the branching
  structure;
\item They allow explicit embeddings of branch groups in
  finitely presented groups (see Theorem~\ref{theorem:answer2md});
\item They may give an explicit basis for the Schur multiplier of
  branch groups (see Theorem~\ref{theorem:schurmult}).
\end{itemize}

This section sums up the proof of Theorem~\ref{theorem:finitelp}. Details
may be found in~\cite{bartholdi:lpres}.

Let $G$ be regular branch on its subgroup $K$, and fix generating
sets $S$ for $G$ and $T$ for $K$. Without loss of generality,
assume $K\le\Stab_G(\LL_1)$, since one may always replace $K$ by
$K\cap\Stab_G(\LL_1)$.

First, there exists a finitely presented group $\Gamma=\SPRES SQ$ with
subgroups $\Delta$ and $\Upsilon=\langle T\rangle$ corresponding to
$\Stab_G(\LL_1)$ and $K$, such that the map $\psi:\Stab_G(\LL_1)\to G^m$
lifts to a map $\Delta\to\Gamma^m$.

The data are summed up in the following diagram:
\[\begin{diagram}
  \node[2]{\Gamma}\arrow{s,-}\arrow{e,A}\node{G}\arrow{s,-}\\
  \node{\Gamma^m}\arrow{s,-}\node{\Delta}\arrow{s,-}
  \arrow{w,t}{\tilde\psi}\arrow{e,A}\node{\Stab_G(\LL_1)}\arrow{s,-}
  \arrow{e,t,J}{\psi}\node{G^m}\arrow{s,-}\\
  \node{\Upsilon^m}\node{\Upsilon}\arrow{e,A}\node{K}\node{K^m}\arrow{w,L}
\end{diagram}\]

Since $\im\tilde\psi$ contains $\Upsilon^m$, it has finite index in
$\Gamma^m$. Since $\Gamma^m$ is finitely presented, $\im\tilde\psi$ too is
finitely presented. Similarly, $\Delta$ is finitely presented, and we may
express $\Ker\tilde\psi$ as the normal closure $\langle R_1\rangle^\#$ in
$\Delta$ of those relators in $\im\tilde\psi$ that are not relators in
$\Delta$. Clearly $R_1$ may be chosen to be finite.

We now use the assumption that $G$ is contracting, with constant $C$. Let
$R_2$ be the set of words over $S$ of length at most $C$ that represent
the identity in $G$. Set $R=R_1\cup R_2$, which clearly is finite.

We consider $T$ as a set distinct from $S$, and not as a subset of
$S^*$. We extend each $\varphi_y$ to a monoid homomorphism
$\tilde\varphi_y:(S\cup T)^*\to(S\cup T)^*$ by defining it arbitrarily
on $S$.

Assume $\Gamma=\SPRES SQ$, and let $w_t\in S^*$ be a
representation of $t\in T$ as a word in $S$. We claim that the
following is an endomorphic presentation of $G$:
\begin{equation}\label{eq:branchpres}
  G = \LPRES{S\cup T}{Q\cup\{t^{-1}w_t\}_{t\in T}}{%
    \{\tilde\varphi_y\}_{y\in Y}}{R_1\cup R_2}.
\end{equation}
For this purpose, consider the following subgroups $\Xi_n$ of $\Gamma$:
first $\Xi_0=\{1\}$, and by induction
\[\Xi_{n+1} = \big\{\gamma\in\Delta\big|\,
(\gamma)\tilde\psi\in\Xi_n^m\big\}.\] We computed $\Xi_1=\langle
R\rangle^\#$. Since $G$ acts transitively on the $n$-th level of the
tree, a set of normal generators for $\Xi_n$ is given by
$\bigcup_{y\in Y^n}(R)\varphi_{y_1}\cdots\varphi_{y_n}$. We also note
that $(\Xi_{n+1})\tilde\psi=\Xi_n^m$.

We will have proven the claim if we show
$G=\Gamma\big/\bigcup_{n\ge0}\Xi_n$. Let $w\in\Gamma$ represent
the identity in $G$. After $\psi$ is applied $|w|$ times, we
obtain $m^{|w|}$ words that are all of length at most $C$, that
is, they belong to $\Xi_1$. Then since
$(\Xi_{n+1})\tilde\psi=\Xi_n^m$, we get $w\in\Xi_{|w|+1}$,
and~(\ref{eq:branchpres}) is a presentation of $G$.

As a bonus, the presentation~(\ref{eq:branchpres}) expresses $K$ as the
subgroup of $G$ generated by $T$.

\section{Examples}\label{sec:l-examples}
We describe here a few examples of branch groups' presentations. As a
first example, let us consider the group $\Aut_f(\tree^{(2)})$ of finitary
automorphisms of the binary tree.
\begin{theorem}
  A presentation of $\Aut_f(\tree^{(2)})$ is given by
  \[T = \PRES{x_0,x_1,\dots}{x_i^2,[x_j,x_k^{x_i}]\quad\forall j,k>i},\]
  and these relators are independent.
\end{theorem}
\begin{proof}
  The generator $x_i$ is interpreted as the element whose portrait has a
  single non-trivial label, at level $i$. The relations are easily
  checked, and they yield a presentation because they are sufficient to
  put words in the $x_i$ in wreath product normal form.

  Finally, if $G_n$ is the quotient of $\Aut_f(\tree)$ acting on the $n$-th level,
  $H_2(G_n,\Z)=(\Z/2)^{\binom{n+1}3}$ by~\cite{blackburn:wreath}
  or~\cite{karpilovsky:schur}.
\end{proof}

We now stick to the $L$-presentation notation, and give presentations for
the following examples:
\begin{description}
\item[The ``first Grigorchuk group'']
  The group $\Gg$ admits the ascending $L$-presentation
  \[\Gg = \ELPRES{a,c,d}{\phi}{a^2,[d,d^a],[d^{ac},(d^{ac})^a]},\]
  where $\phi:\{a,c,d\}^*\to\{a,c,d\}^*$ is defined by
  \[\phi(a)=aca,\quad \phi(c)=cd,\quad \phi(d)=c.\]
  These relators are independent, and $H_2(\Gg,\Z)=(\Z/2)^\infty$.
\item[The ``Grigorchuk supergroup''~\cite{bartholdi-g:parabolic}] The
  group $\Sg=\langle a, \tilde b, \tilde c, \tilde d \rangle$ acting on the binary
  tree, where $a$ is the rooted automorphism $a=((1,2))$ and the other three
  generators are the directed automorphisms defined recursively by
  \[ \tilde b = (a, \tilde c) \quad \tilde c = (1, \tilde d) \quad \tilde
        d = (1, \tilde b), \]
 admits the ascending $L$-presentation
  \[\Sg=\ELPRES{\langle a,\tilde b,\tilde c,\tilde d}{%
  \tilde\phi}{a^2,[\tilde b,\tilde c],
  [\tilde c,\tilde c^a],[\tilde c,\tilde d^a],[\tilde d,\tilde d^a],
  [\tilde c^{a\tilde b},(\tilde c^{a\tilde b})^a],
  [\tilde c^{a\tilde b},(\tilde d^{a\tilde b})^a],
  [\tilde d^{a\tilde b},(\tilde d^{a\tilde b})^a]},\]

  where $\tilde\phi:\{a,\tilde b,\tilde c,\tilde d\}^*\to\{a,\tilde
  b,\tilde c,\tilde d\}^*$ is defined by
  \[a\mapsto a\tilde ba,\quad\tilde b\mapsto\tilde d,\quad\tilde
  c\mapsto\tilde b,\quad\tilde d\mapsto\tilde c.\]
  These relators are independent, and $H_2(\Sg,\Z)=(\Z/2)^\infty$.
\item[The ``Fabrykowski-Gupta
  group''~\cite{fabrykowski-g:growth2,bartholdi-g:parabolic}]
  The group $\FGg$ admits the ascending endomorphic presentation
  \[\ELPRES{a,r}{\phi,\chi_1,\chi_2}{a^3,[r^{1+a^{-1}-1+a+1},a],
    [a^{-1},r^{1+a+a^{-1}}][r^{a+1+a^{-1}},a]},\]
  where $\sigma,\chi_1,\chi_2:\{a,r\}^*\to\{a,r\}^*$ are given by
  \begin{xalignat*}{2}
    \phi(a)&=r^{a^{-1}},&\phi(r)&=r,\\
    \chi_1(a)&=a,&\chi_1(r)&=r^{-1},\\
    \chi_2(a)&=a^{-1},&\chi_2(r)&=r.
  \end{xalignat*}
  These relators are independent, and $H_2(\FGg,\Z)=(\Z/3)^\infty$.
\item[The ``Gupta-Sidki group''~\cite{sidki:pres}]
  The Gupta-Sidki group $\GSg$ admits the endomorphic presentation
  \[\LPRES{a,t,u,v}{a^3,t^3,u^{-1}t^a,v^{-1}t^{a^{-1}}}{\phi,\chi}{%
    (tuv)^3,[v,t][vt,u^{-1}tv^{-1}u],[t,u]^3[u,v]^3[t,v]^3},\]
  where $\phi,\chi:\{t,u,v\}^*\to\{t,u,v\}^*$ are given by
  \begin{gather*}
    \phi:\begin{cases}
      t\mapsto t,\\
      u\mapsto [u^{-1}t^{-1},t^{-1}v^{-1}]t=u^{-1}tv^{-1}tuvt^{-1},\\
      v\mapsto t[tv,ut] = t^{-1}vutv^{-1}tu^{-1},
    \end{cases}\qquad
    \chi:\begin{cases}
      t\mapsto t^{-1},\\ u\mapsto u^{-1},\\ v\mapsto v^{-1}.
    \end{cases}
  \end{gather*}
  These relators are independent, and $H_2(\GSg,\Z)=(\Z/3)^\infty$.
  Note that $\chi$ is induced by the automorphism of $\GSg$
  defined by $a\mapsto a,\;t\mapsto t^{-1}$; however, $\phi$ does not
  extend to an endomorphism of $\GSg$.

  It is precisely for that reason that no ascending endomorphic presentation
  of $\GSg$ is known.
\item[The ``Brunner-Sidki-Vieira group''~\cite{brunner-s-v:nonsolvable}]
  Consider the group $G=\langle\mu,\tau\rangle$ acting on the binary tree,
  where $\mu$ and $\tau$ are defined recursively by
  \[ \mu =(1,\mu^{-1})a, \qquad \tau=(1,\tau)a. \]
  Note that $G$ is not branch, but it is weakly branch.
  Writing $\lambda=\tau\mu^{-1}$, $G$ admits the ascending $L$-presentation
  \[G = \ELPRES{\lambda,\tau}{\phi}{[\lambda,\lambda^\tau],
  [\lambda,\lambda^{\tau^3}]},\]
  where $\phi$ is defined by $\tau\mapsto\tau^2$ and
  $\lambda\mapsto\tau^2\lambda^{-1}\tau^2$.
\end{description}

The above $L$-presentation for $\Gg$ allowed the second author to
answer a question of \MDay~\cite{day:amen} for the class of
finitely presented groups (the question is formulated
in~\cite{cannon-f-p:thompson} in this special setting ):
\begin{theorem}[Grigorchuk,~\cite{grigorchuk:notEG}]\label{theorem:answer2md}
  There exists a finitely presented amenable group that is not
  elementarily amenable.
\end{theorem}
Recall that a group $G$ is amenable if it admits a left-invariant
finitely-additive measure. Examples include the finite groups, the
abelian groups and all groups obtained from previous examples by
short exact sequences and direct limits. The smallest class
containing the finite and abelian groups and closed for the
mentioned basic constructions is known as the class of
\emdef{elementarily amenable} groups.

\begin{proof}
  Consider the presentation of the first Grigorchuk group given above,
  and form the \textsf{HNN} extension $H$ amalgamating $\Gg$ with
  $\phi(\Gg)$. It is an ascending \textsf{HNN} extension, so $H$ is
  amenable; and $H$ admits the (ordinary) finite presentation
  \[H = \PRES{a,c,d,t}{a^2,[d,d^a],[d^{ac},d^{aca}],a^taca,c^tcd,d^tc}.\]
\end{proof}

Another presentation of the group $H$ from the previous proof, due to the
first author, is given in~\cite{grigorchuk-h-s:paradox}
\[ H = \PRES{a,t}{a^2,TaTatataTatataTataT,(Tata)^8,(T^2ataTat^2aTata)^4}, \]
where $T$ denotes the inverse of $t$.

\part{Algebraic Aspects}
\date{October 27, 2002}
\chapter{Just-Infinite Branch Groups}\label{chapter:just-infinite}
\begin{definition}
  A group $G$ is \emdef{just-infinite} if it is infinite but all of
  its proper quotients are finite, i.e., if all of its nontrivial normal
  subgroups have finite index.
\end{definition}

The following simple criterion from~\cite{grigorchuk:jibg} characterizes
the just-infinite branch groups acting on a tree.
\begin{theorem}\label{theorem:jicriterion}
  Let $G$ be a branch group acting on a tree and let $(L_i,H_i)_{i\in\N}$
  be a corresponding branch structure. The following three conditions are equivalent
  \begin{enumerate}
  \item $G$ is just-infinite.
  \item the abelianization $H_i^{ab}$ of $H_i$ is finite, for each $i\in\N$.
  \item the commutator subgroup $H_i'$ of $H_i$ has finite index in $G$,
        for each $i\in\N$.
  \end{enumerate}
\end{theorem}

The statement is a corollary of the fact that $H_i'$, being characteristic
in the normal subgroup $H_n$, is a normal subgroup of $G$, for each
$i\in\N$, and the following useful lemma that  says that weakly branch
groups satisfy the following property:

\begin{lemma}
  Let $G$ be a weakly branch group acting on a tree and let $(L_i,H_i)_{i\in\N}$
  be a corresponding branch structure. Then every non-trivial normal
  subgroup $N$ of $G$ contains the commutator subgroup $H_n'$, for some $n$
  depending on $N$.
\end{lemma}
\begin{proof}
Let $g$ be an element in $G \setminus \Stab_G(\LL_1)$ and let $N=\langle g
\rangle^G$ be its normal closure in G. Then $g=ha$ for some $h \in
\Stab_{\Aut(\tree)}(\LL_1)$ with decomposition $h =(h_1,\dots,h_{m_1})$
and $a$ a nontrivial rooted automorphism of $\tree$. Without loss of
generality we may assume that $1^a=m_1$.

For arbitrary elements $\xi,\nu \in L_1$, we define $f,t\in H_1$ by
$f=(\xi,1,\dots,1)$ and $t=(\nu,1,\dots,1)$ and calculate
\begin{gather*}
    [g,f] = (\xi,1,\dots,1,(\xi^{-1})^{h_1}), \\
    [[g,f],t] = ([\xi,\nu],1,\dots,1).
\end{gather*}
Since $[[g,f],t]$ is always in $N=\langle g \rangle^G$, we obtain
$L_1'\times 1 \times\dots\times 1 \preceq N$ and, by the spherical
transitivity of $G$, it follows that
\[ L_1' \times L_1' \times \dots \times L_1'= H_1' \preceq N. \]
Thus any normal subgroup of $G$ that contains $g$ also contains $H_1'$.

Similarly, if $g$ is an element in $\Stab_G(\LL_n) \setminus
\Stab_G(\LL_{n+1})$ and $N$ is the normal closure $N= \langle g^G
\rangle$, then $N$ contains $H_{n+1}'$.
\end{proof}

In particular, the above results imply that all finitely generated
torsion weakly branch groups are just-infinite branch groups.

The study of just-infinite groups is motivated by their minimality. More
precisely, we have the following
\begin{theorem}\cite{grigorchuk:jibg}
  Every finitely generated infinite group has a just-infinite
  quotient.
\end{theorem}

Therefore, if $\mathcal C$ is a class of groups closed for taking
quotients and it contains a finitely generated infinite group, then it
contains a finitely generated just-infinite group.

Note that there are non-finitely generated groups that do not have
just-infinite quotients, for example, the additive group of rational
numbers $\Q$.

It is known (see~\cite{wilson:jig}) that a just-infinite group with
non-trivial Baer radical is a finite extension of a free abelian group of
finite rank (recall that the \emdef{Baer radical} of the group $G$ is the
subgroup of $G$ generated by the cyclic subnormal subgroups of $G$).
Moreover, the only just-infinite group with non-trivial center is the
infinite cyclic group $\Z$. Therefore, an abelian group has just-infinite
quotient if and only if it can be mapped onto $\Z$. In particular, no
abelian torsion group has just-infinite quotients. The last fact is in a
sharp contrast with the fact that there are large classes of centerless,
torsion, just-infinite, branch groups, for instance \GG\ groups with
finite directed part $B$ (see Chapter~\ref{chapter:torsion}) and many
\GGS\ groups.

\begin{definition}
  A group $G$ is \emdef{hereditarily just-infinite} if it is
  residually finite and all of its non-trivial normal subgroups are
  just-infinite.
\end{definition}

We mention that our definition of hereditarily just-infinite group differs
from the one given in~\cite{wilson:segal} in that we require residual
finiteness. Note that all non-trivial normal subgroups of a group $G$ are
just infinite if and only if all subgroups of finite index in $G$ are just
infinite. This is true since every subgroup of finite index in $G$
contains a normal subgroup of $G$ of finite index.

Examples of hereditarily just-infinite groups are the infinite cyclic
group $\Z$, the infinite dihedral group $D_\infty$ and the projective
special linear groups $\PSL(n,\Z)$, for $n\geq 3$. However, the whole
class is far from well understood and described.

The following result from~\cite{grigorchuk:jibg}, which modifies the
result of \JWilson\ from~\cite{wilson:jig} (see also~\cite{wilson:segal}),
strongly motivates the study of the branch groups.

\begin{theorem}[Trichotomy of just-infinite groups]\label{theorem:tri}
  Let $G$ be a finitely generated just-infinite group. Then exactly one of the
  following holds:
  \begin{enumerate}
  \item\label{theorem:tri:1} $G$ is a branch group.
  \item\label{theorem:tri:2} $G$ has a normal subgroup $H$ of finite index of the form
    \[H = L^{(1)} \times \dots \times L^{(k)}=L^k,\]
    where the factors are copies of a group $L$, the conjugations by
    the elements in $G$ transitively permute the factors of $H$ and $L$ has
    exactly one of the following two properties:
    \begin{enumerate}
    \item\label{theorem:tri:2a} $L$ is hereditary just-infinite (in case
      $G$ is residually finite).
    \item\label{theorem:tri:2b} $L$ is simple (in case $G$ is not
      residually finite).
    \end{enumerate}
  \end{enumerate}
\end{theorem}

The proof of this theorem presented in~\cite{grigorchuk:jibg} uses only
the statement from~\cite{wilson:jig} that every subnormal subgroup in a
just-infinite group with trivial Bear radical has a near complement (but
this is probably one of the most important facts from Wilson's theory).
The proof actually works for any (not necessarily finitely generated)
just-infinite group with trivial Baer radical, for instance just-infinite
groups which are not virtually cyclic.

The results of Wilson in~\cite{wilson:jig} (see also~\cite{wilson:segal})
combined with the above trichotomy result show that the following
characterization of just-infinite branch groups is possible. Define an
equivalence relation on the set of subnormal subgroups of a group $G$ by,
$H \sim K$ if the intersection $H \cap K$ has a finite index booth in $H$
and in $K$. The set of equivalence classes of subnormal subgroups, ordered
by the order induced by inclusion, forms a Boolean lattice, which,
following \JWilson, is called the \emdef{structure lattice} of $G$.

\begin{theorem}
  Let $G$ be a just-infinite group. Then $G$ is branch group if and
  only if it has infinite structure lattice.

  Moreover, in such a case, the structure lattice is isomorphic to the
  lattice of closed and open subsets of the Cantor set.
\end{theorem}

On the intuitive level, the just-infinite groups should be
considered as ``small'' groups in contrast to, say, the free or
non-elementary hyperbolic groups, which are ``large'' groups.

There is a rigorous approach to the concept of largeness in
groups. Namely, following Pride (\cite{pride:large}), we say that
a group $G$ is \emdef{larger} than a group $H$, and we denote $G
\succeq H$, if $H$ has a subgroup of finite index that is a
homomorphic image of a subgroup of $G$ of finite index. The groups
$G$ and $H$ are \emdef{equally large} (or \emdef{Pride
equivalent}) if $G \succeq H$ and $H \succeq G$. The set of
equivalence classes of equally large groups is partially ordered
by $\succeq$ and the class of finite groups is the obvious
smallest element.

We denote the class of groups equally large to $G$ by $[G]$. A group $G$
is called \emdef{minimal} if the only class below $[G]$ is the class of
finite groups $[1]$. The \emdef{height} of a group $G$ is the height of
the class $[G]$ in the ordering, i.e. the length of a maximal chain
between $[1]$ and $[G]$. Therefore, the minimal groups are the groups of
height 1. Such groups are called \emdef{atomic} in~\cite{neumann:pride}

\begin{theorem}[\cite{grigorchuk-w:minimality}]\label{thm:minimal}
The first Grigorchuk group $\Gg$ and the Gupta-Sidki $p$-groups are
minimal.
\end{theorem}

A number of questions about just-infinite groups was asked
in~\cite{pride:large,e-pride:large}. Positive answer to Problem~5
from~\cite{pride:large} (Problem~4' in~\cite{e-pride:large}) that
asks if there exist finitely generated just-infinite groups that
do not satisfy the ascending chain condition on subnormal
subgroups was provided in~\cite{grigorchuk:gdegree}. Later,
\PNeumann\ constructed in~\cite{neumann:pride} more examples of
finitely generated just-infinite regular branch groups answering
the same question (and also some other quastions) raised by
\MEdjvet\ and \SPride\ in~\cite{e-pride:large}. In particular,
\PNeumann\ provided negative answer to the question if every
finitely generated minimal group is finite-by-$D_2$-by-finite
(here $D_2$ denotes the class of groups in which every nontrivial
subnormal subgroup has finite index). Negative answer to this last
question also follows form Theorem~\ref{thm:minimal} above.

The question of possible heights of finitely generated just-infinite
groups is an interesting one. All hereditarily just-infinite and all
infinite simple groups are minimal. It is plausible that the Grigorchuk
2-groups from~\cite{grigorchuk:gdegree} that are defined by non-periodic
sequences have infinite height (see Question~\ref{question:height}).

\date{October 27, 2002}
\chapter{Torsion Branch Groups} \label{chapter:torsion}
As mentioned in the introduction, the most elegant examples of
finitely generated infinite torsion groups were constructed within
the class of branch groups
(\cite{grigorchuk:burnside,gupta-s:burnside}). Therefore, branch
groups play important role in problems of Burnside type. At the
moment, there is a number of constructions of torsion branch
groups. Besides the early works
(\cite{aleshin:burnside,sushchansky:pgps,grigorchuk:burnside,grigorchuk:growth,
grigorchuk:gdegree,gupta-s:burnside,gupta-s:infinitep,gupta-s:ext}
there are more modern and general constructions
(\cite{bartholdi-s:wpg,grigorchuk:jibg,sunik:phd}. Nevertheless,
all these constructions follow the same idea of stabilization and
length reduction.

Recall that finitely generated torsion branch groups acting on a
tree are just-infinite (Theorem~\ref{theorem:jicriterion}).

We provide a proof here that all \GG\ groups with torsion directed
group $B$ are themselves torsion groups, and the argument follows
the mentioned general scheme. This is an improvement over the
result in~\cite{bartholdi-s:wpg} that deals only with regular
trees and all the root actions (the actions of $A_{\sigma^t\om}$,
for all $t$) are isomorphic and regular. Thus whenever $B$ is a
finitely generated torsion group the corresponding \GG\ group is a
Burnside group and there are uncountably many non-isomorphic
examples. We know that torsion spinal groups cannot have finite
presentation (see Corollary~\ref{theorem:spinalnotfp}. The results
in Theorem~\ref{theorem:branch-ie} below show that these groups
cannot have finite exponent.  Therefore, in the finitely generated
case, one of our goals is to give upper bounds on the order of an
element depending on its length. This leads to the notion of
torsion growth.

Let $G$ be finitely generated infinite torsion group and let $S$
be a finite generating set that generates $G$ as a monoid. For any
non-negative real number $n$, the maximal order of an element of
length at most $n$ is finite, and we denote it by $\pi_G^S(n)$.
The function $\pi_G^S$, defined on the non-negative real numbers,
is called the \emdef{torsion growth} function of $G$ with respect
to $S$.

The group order of an element $g$ is denoted by $\pi(g)$. In case $F$ is a
word $\pi(F)$ denotes the order of the element represented by the word
$F$.

We describe a step in a procedure that successfully implements the ideas
and constructions introduced in~\cite{grigorchuk:gdegree}. The final
result of the procedure is a tree that helps us to show that all \GG\
groups whose directed part $B$ is torsion are themselves torsion groups
and also to determine some upper bounds on the torsion growth of the
constructed groups. The construction presented here is slightly more
complicated than the construction from~\cite{grigorchuk:gdegree} because
of the fact that we need to take into account the possibility of
non-cyclic (and even non-regular) actions of the groups $A_{\sigma^t\om}$
on their corresponding alphabets.

Let $G_\om$, be a \GG\ group defined by the triple $\om$ (recall
Definition~\ref{defn:spinal} and Definition~\ref{defn:GG}) and let $F$ be
a reduced word of even length of the form
\begin{equation} \label{eq:crform}
  F= b_1a_1\dots b_ka_k,
\end{equation}
where $a_i$ represent non-trivial rooted automorphisms in $A_\om$ and
$b_j$ represent non-trivial directed automorphisms in $B_\om$. Rewrite $F$
in the form $F= b_1b_2^{g_2}\dots b_k^{g_k}a_1\dots a_k$, where
$g_i=(a_1\dots a_{i-1})^{-1}$, $i=2,\dots,k$. Set $g=a_1\dots a_k\in
A_\om$ and let its order be $s$. Note that $g=1$ corresponds to $F \in
\Stab_\om(\LL_1)$.  Put
\[H=b_1b_2^{g_2}\dots b_k^{g_k}, \]
consider the word $F^s=(Hg)^s\in \Stab_\om(\LL_1)$ and rewrite it in the
form $F^s = HH^{g^{-1}}H^{g^{-2}}\dots H^{g^{-(s-1)}}$. Next, by using
tables similar to Table~\ref{table:phi}, but for all possible $a$, we
calculate the possibly unreduced words
$\overline{F_1},\dots,\overline{F_{m_1}}$ representing the first level
sections $(F^s)\varphi_1$,$\dots$, $(F^s)\varphi_{m_1}$, respectively. We
have
\[\overline{F_i} = H_i H_{i^g} \dots H_{i^{g^{s-1}}},\]
for $i=1,\dots,m_1$, where $H_j$ represents the corresponding section of
$H$. Note that any two words $\overline{F_i}$ and $\overline{F_j}$ that
correspond to two indices in the same orbit of $g$ represent conjugate
elements of $G_\om$. This is clear since
\[\overline{F_{i^g}} = H_{i^g} H_{i^{g^2}} \dots H_{i^{g^{s-1}}} H_i
= \overline{F_i}^{H_i}.\]

For $i=1,\dots,m_1$, let the length of the cycle of $i$ in $g$ be $t_i$.
Then
\[\overline{F_i} = (H_i H_{i^g} \dots H_{i^{g^{t_i-1}}})^{s/t_i}.\]
and let $F_i$ be a cyclically reduced word obtained after applying simple
reductions (including the cyclic ones) to the word $\tilde{F}_i=H_i
H_{i^g} \dots H_{i^{g^{t_i-1}}}$. Clearly
\[(F^s)\psi= (\tilde{F}_1^{s/t_1},\dots,\tilde{F}_{m_1}^{s/t_{m_1}})\]
and $F$ has finite order if and only if $F_1$,$\dots$,$F_{m_1}$
all have finite order. In the case of a finite order, the order
$\pi(F)$ of $F$ is a divisor of $s \cdot
\lcm(\pi(F_1),\dots,\pi(F_{m_1}))$, since the order
$\pi(F_i^{s/t_i})$ divides the order $\pi(F_i)$.

Let us make a couple of simple observations on the structure of the
possibly unreduced words $\tilde{F}_1,\dots,\tilde{F}_{m_1}$ used to
obtain the reduced words $F_1,\dots,F_{m_1}$.  We have
\begin{multline} \label{eq:tildeF}
  \tilde{F}_i = H_i H_{i^g} \dots H_{i^{g^{t_i-1}}} = \\
  = (b_1)\varphi_i (b_2^{g_2})\varphi_i \dots (b_k^{g_k})\varphi_i \;
  (b_1)\varphi_{i^g} (b_2^{g_2})\varphi_{i^g} \dots
  (b_k^{g_k})\varphi_{i^g} \; \dots \\
  (b_1)\varphi_{i^{g^{t_i-1}}} (b_2^{g_2})\varphi_{i^{g^{t_i-1}}}
  \dots (b_k^{g_k})\phi_{i^{g^{t_i-1}}}.
\end{multline}
It is important to note that the indices $i,i^g,\dots, i^{g^{t_i-1}}$ are
all distinct, since the length of the cycle of $i$ in $g$ is $t_i$. This
means that at most one of $(b_1)\varphi_i$,
$(b_1)\varphi_{(i)^g}$,$\dots$, $(b_1)\varphi_{(i)g^{t_i-1}}$ can be equal
to $b_1$, at most one can be equal to $(b_1)\om_1$, and all others are
empty. Thus we conclude that the word $\tilde{F}_i$ consists of some of
the letters $b_1,\dots,b_k$, possibly not in that order, and no more than
$k$ $A$-letters. In particular, it is possible to get $k$ $A$-letters only
if none of the letters $b_1,\dots,b_k$ is in the kernel $K_1$.

We have already observed that any two, possibly unreduced, words
$\tilde{F}_i$ and $\tilde{F}_j$ such that the indices $i$ and $j$
are in the same cycle of $g$ are conjugate. Choose a
representative for each cycle of $g$, and let
$\tilde{F}_{j_1},\dots,\tilde{F}_{j_c}$ be such representatives.
The important feature of our construction is that each letter
$b_i$ from the representation of $F$ as in~(\ref{eq:crform}), as
well as each $(b_i)\om_1$, will appear in exactly one of the
representatives $\tilde{F}_{j_1},\dots,\tilde{F}_{j_c}$ before
reduction.

The tree that has $F$ at its root and the cyclically reduced cycle
representatives $F_{j_1},\dots,F_{j_c}$ at the vertices below the
root is called the \emdef{pruned period decomposition} of $F$.

\begin{theorem}
  A \GG\ group is a torsion group if and only of its directed part $B$
  is a torsion group.
\end{theorem}
\begin{proof}
  We prove that the order of any element $g$ in $G_\om$ is finite in
  case $B$ is a torsion group. The proof is by induction on the length
  $n$ of $g$, for all \GG\ groups simultaneously.

  The statement is clear for $n=0$ and $n=1$. Assume that it is true
  for all words of length less than $n$, where $n\geq 2$, and consider
  an element $g$ of length $n$.

  If $n$ is odd the element $g$ is conjugate to an element of smaller
  length and we are done by the inductive hypothesis. Assume then that
  $n$ is even.  Clearly, $g$ is conjugate to an element that can be
  represented by a word of the form
  \[F= b_1a_1\dots b_ka_k.\]
  If all the cycle representatives $F_i$ from the pruned period
  decomposition of $F$ have length shorter than $n$ we are done by the
  inductive hypothesis.

  Assume that some of the cycle representatives $F_i$ have length $n$.
  This is possible only when $F$ does not have any $B$-letters from
  $K_1$. Also, the words $\tilde{F}_i$ corresponding to the words
  $F_i$ of length $n$ must be reduced, so that the words $F_i$ that
  have length $n$ must have the same $B$-letters as $F$ does. For each
  of these finitely many words we repeat the discussion above. Either
  all of the constructed words $F_{ij}$ are strictly shorter than $n$,
  and we get the result by induction; or some have length $n$, but the
  $B$-letters appearing in them do not come from $K_1\cup K_2$.

  This procedure cannot go on forever since $K_1\cup K_2 \cup
  \dots\cup K_r=B$ holds for some $r\in\N$.  Therefore at some stage we get a
  shortening in all the words and we conclude that the order of $F$ is
  finite.
\end{proof}

In the sequel we just list some estimates on the period growth in case the
directed part $B$ is a finite group. The proofs can be found
in~\cite{bartholdi-s:wpg}.

A finite subsequence $\om_{i+1}\om_{i+2}\dots \om_{i+r}$ of the defining
sequence $\omover = \om_1\om_2\dots$ is \emdef{complete} if each element
of $B$ is sent to the identity by at least one homomorphism from the
sequence $\om_{i+1}\om_{i+2}\dots\om_{i+r}$, i.e., if $\bigcup_{j=1}^r
K_{i+j} = B$. We note that the complete sequence $\om_{i+1}\om_{i+2}\dots
\om_{i+r}$ must have length at least $m+1$, where $m$ is the minimal
branching index in the branching sequence, since each kernel $K_{i+j}$ has
index $|A_{\sigma^{i+j}\om}| \geq m_{i+j+1} \geq m$ in $B$, for all
$j=1,\dots,r$. In particular, the length of a complete sequence is never
shorter than $3$.  By the definition of a \GG\ group, all sequences that
define a \GG\ group can be factored into finite complete subsequences.

A defining sequence $\omover$ is \emdef{$r$-homogeneous}, for $r\geq 3$,
if all of its finite subsequences of length $r$ are complete. A defining
sequence $\omover$ is \emdef{$r$-factorable}, for $r\geq3$, if it can be
factored in complete subsequences of length at most $r$.

\begin{theorem}[Period $\eta$-Estimate] \label{theorem:(r)pgrowth}
  If $\omover$ is an $r$-homogeneous sequence and $B$ has exponent
  $q$, then there exist a positive constant $C$ such that the torsion growth
  function of the group $G_\om$ satisfies
  \[\pi_\om(n) \leq C n^{\log_{1/\eta_r}(q)}\]
  where $\eta_r$ is the positive root of the polynomial $x^r + x^{r-1}
  + x^{r-2} -2$.
\end{theorem}

\begin{theorem}[Period 3/4-Estimate] \label{theorem:[r]pgrowth3/4}
  If $\omover$ is an $r$-factorable sequence and $B$ has exponent $q$,
  then there exists a positive constant $C$ such that the torsion growth function
  of the group $G_\om$ satisfies
  \[\pi_\om(n) \leq C n^{r\log_{4/3}(q)}.\]
\end{theorem}

\begin{theorem}[Period 2/3-Estimate] \label{theorem:[r]pgrowth2/3}
  If $\omover$ is an $r$-factorable sequence such that each factor
  contains three letters whose kernels cover $B$ and $B$ has exponent
  $q$, then there exists a positive constant $C$ such that the torsion growth
  function of the group $G_\om$ satisfies
  \[\pi_\om(n) \leq C n^{r\log_{3/2}(q)}.\]
\end{theorem}

Let us assume now that all the branching indices are prime numbers
and the groups $A_{\sigma^t\om}$ are cyclic of prime order, for
all $t$. There is no loss in generality if we assume that
$A_{\sigma^t\om}$ is generated by the cyclic permutation
$a=(1,2,\dots,m_{t+1})$. Note that our assumptions force $B$ to be
abelian group, since $B$ is always a subdirect product of several
copies of the root groups $A_{\sigma^t\om}$.

\begin{theorem} \label{theorem:r-1/1}
  Let the branching sequence consists only of primes, $B$ have
  exponent $q$ and $\omover$ be an $r$-homogeneous word. There exists a
  positive constant $C$ such that the torsion
  growth function of $G_\om$ satisfies
  \[\pi_\om(n) \leq C n^{(r-1)\log_2(q)}.\]
\end{theorem}

Finally we give a tighter upper bound on the period growth of the
Grigorchuk 2-groups (as defined in~\cite{grigorchuk:gdegree}).
\begin{theorem} \label{theorem:r/2}
  Let $G_\om$ be a Grigorchuk $2$-group. If $\omover$ is an
  $r$-homogeneous word, then there exists a positive constant $C$ such that
  the torsion growth function of the group $G_\om$ satisfies
  \[\pi_\om(n) \leq C n^{r/2}.\]
\end{theorem}

In addition to the above estimates \cite{bartholdi-s:wpg} provides
lower bounds on the torsion growth function in some cases and, in
particular, shows that some Grigorchuk 2-groups have torsion
growth functions $\pi(n)$ that are at least linear in $n$ (for
example the group defined by the sequence
$\omover=01020102\dots$). Previous results of \ILysionok\ and
\YLeonov\ (\cite{lysionok:order,leonov:periodbd} already
established a lower bound of $C_1n^{\frac{1}{2}}$ for the torsion
growth function of the first Grigorchuk group, while
Theorem~\ref{theorem:r/2} establishes an upper bound of
$C_2n^{\frac{3}{2}}$.

The result following the definition below finds its predecessor in
the work of \NGupta\ and \SSidki~\cite{gupta-s:infinitep}, where
they prove that the Gupta-Sidki $p$-groups contain arbitrary long
iterated wreath products of cyclic groups of order $p$ and
therefore contain a copy of each finite $p$-group. The general
argument below follows the approach
from~\cite{grigorchuk-h-z:profinite} and applies well to more
broad settings.

\begin{definition}
Let $\mathcal P$ be a non-empty set of primes. A group $G$ of
automorphisms of $\tree$ has \emdef[omnipresent
torsion]{omnipresent $\mathcal P$-torsion} if, for every vertex $u
\in \tree$, the lower companion group $L_u^G$ (the rigid
stabilizer at $u$) has non-trivial elements of $\mathcal P$-order.
\end{definition}

In case $\mathcal P$ consists of all primes we just say that $G$
has \emdef{omnipresent torsion} and in case $\mathcal P=\{p\}$ we
say that $G$ has omnipresent $p$-torsion. We note that every
spherically transitive group with omnipresent $\mathcal P$-torsion
must be a weakly branch group. We also note that if $G$ has
omnipresent $\mathcal P$-torsion then so does each of its lower
companion groups $L_u$ considered as a group of automorphisms of
$\tree_u$.

\begin{lemma}\label{lemma:tor-ie}
Let $\mathcal P$ be a non-empty set of primes and $G$ a group of
tree automorphisms with omnipresent $\mathcal P$-torsion. Then $G$
contains arbitrary long iterated wreath products of cyclic groups
of the form
\[(( \Z/p_1\Z \wr \dots  )\wr \Z/p_{n-1}\Z )\wr \Z/p_n\Z \]
where $n\in\N$ and each $p_i$ is a $\mathcal P$-prime, for
$i\in\{1,\dots,n\}$.
\end{lemma}
\begin{proof}[Sketch of a proof]
We prove the claim by induction on $n$, simultaneously for all
groups of tree automorphisms with omnipresent $\mathcal
P$-torsion. The claim is obvious for $n=1$. Assume that $n \geq 2$
and that the claim holds for all positive natural numbers $<n$.

Choose an arbitrary non-trivial element $g$ of $G$ of finite
$\mathcal P$-prime order $p_n$. Let $g$ fix the level $\LL_k$ but
not $\LL_{k+1}$ and let $u$ be a vertex on level $k$ with
non-trivial vertex permutation $(u)g=a$. All non-trivial cycles of
$a$ have length $p$ and, without loss of generality, we may assume
that one such cycle is (1,2,\dots,p). Without loss of generality
we may also assume that the sections $g_{u1},\dots,g_{u_p}$ are
trivial (we may accomplish this by conjugation if necessary). By
the inductive hypothesis, the lower companion group $L_{uy}$
contains an iterated wreath product of length $n-1$
\[ Q = (( \Z/p_1\Z \wr \dots ) \wr \Z/p_{n-2}\Z ) \wr \Z/p_{n-1}\Z \]
of the required form. But then
\[ \langle Q, g \rangle \cong Q \wr \langle g \rangle =
        (( \Z/p_1\Z \wr \dots ) \wr \Z/p_{n-1}\Z ) \wr \Z/p_n\Z. \]
\end{proof}

The above lemma has many corollaries, some of which are summed up in the
following theorem:

\begin{theorem}\label{theorem:branch-ie}
Let $G$ be a group of tree automorphisms. If, for each vertex $u$,
the lower companion group $L_u$ of $G$ has an element of finite
order, then $G$ has elements of unbounded finite order. Further,
\begin{enumerate}
\item Every weakly branch torsion group has infinite exponent.
\item
Every weakly branch $p$-group contains a copy of every finite $p$-group.
\item
Every weakly regular branch group is weakly branched over a torsion free
group or contains a copy of every finite $p$-group, for some prime $p$.
\item
A regular branch group is either virtually torsion free or it contains a
copy of every finite $p$-group, for some prime $p$.
\end{enumerate}
\end{theorem}

Finally, we note that the class of \GGS\ groups is also rich with
examples of torsion groups, starting with the second Grigorchuk
group from~\cite{grigorchuk:burnside}, Gupta-Sidki $p$-groups
from~\cite{gupta-s:burnside,gupta-s:infinitep} and certain
Gupta-Sidki extensions from~\cite{gupta-s:ext}. For a wide class
of torsion branch \GGS\ groups see
Subsection~\ref{subsec:moreggs}.

\date{October 27, 2002}
\chapter{Subgroup Structure}\label{ch:subgroup}
We study in this chapter some subgroup series (derived series,
powers series) and general facts about branch groups. We then
describe important small-index subgroups in the examples
$\Gg,\FGg,\BGg,\GSg$ (recall the definitions form
Section~\ref{sec:branch:examples}. The lower central series is
treated in the next chapter. Most of the results come
from~\cite{bartholdi-g:parabolic}.

\section{The derived series}
Let $G$ be a group. The \emdef{derived series}
$(G^{(n)})_{n\in\N}$ of $G$ is defined by $G^{(0)}=G$ and
$G^{(n)}=[G^{(n-1)},G^{(n-1)}]$. A group is \emdef{solvable} if
$G^{(n)}=\{1\}$ for some $n\in\N$. It is \emdef{residually
solvable} if $\bigcap_{n\in\N}G^{(n)}=\{1\}$. Note that if
$(\gamma_n(G))_{n\in\N}$ is the lower central series of $G$, then
a general result states that $G^{(n)} \leq \gamma_{2^n}(G)$ holds
for all $n$.

\subsection{The derived series of $\Gg$}
Since $\Gg$ is regular branch over $K=\langle x \rangle^G$, where
$x=[a,b]$, we consider the finite-index subgroups
$\rist_\Gg(\LL_n)=K\times\dots\times K$ with $2^n$ factors, for $n\ge2$.
\begin{theorem}
  $\Gg^{(n)}=\rist_\Gg(\LL_{2n-3})$ for all $n\ge3$, and
  $K^{(n)}=\rist_\Gg(\LL_{2n})$ for all $n\ge1$.
\end{theorem}
\begin{proof}
  First, one may check by elementary means that
  $\Gg^{(3)}=\rist_\Gg(\LL_3)=K^{\times8}$. Then $K'$ is the normal
  closure in $\Gg$ of $[x^d,x]$, and $[x^d,x]\psi=([(ca)^b,ca],1)$,
  and $[(ca)^b,ca]\psi=(x,1)$. Therefore $K'=K^{\times4}$.
\end{proof}

\subsection{The derived series of $\GSg$}
The result is even slightly simpler for $\GSg$, which is regular
branch over $\GSg'$:
\begin{theorem}
  $\GSg^{(n)}=(\GSg'')^{\times3^{n-2}}$ for all $n\ge2$.
\end{theorem}
\begin{proof}
  The core of the argument is to show that
  $\GSg^{(3)}=\GSg''\times\GSg''\times\GSg''$. This follows from
  $\GSg^{(3)}=\gamma_8(\GSg)$ and $\GSg''=\gamma_5(\GSg)$, but the
  computations are tricky --- see Subsection~\ref{subs:lcsGS} for details.
\end{proof}

\section{The powers series}\label{sec:powers}
Let $G$ be a group and $d$ an integer. The \emdef{powers series}
$(\mho^n_d(G))_{n\in\N}$ is defined by $\mho^0_d(G)=G$ and
\[\mho^n_d(G)=\langle x^d\text{ for all }x\in\mho^{n-1}_d(G)\rangle.\]

\begin{theorem}
  The $2$-powers series of $\Gg$ is as follows $\mho_2(\Gg)=\Gg'$ and
  \[\mho_2^n(\Gg)=\langle\Delta((\mho_2K)^{\times2^{n-2}}),\Delta(K^{\times2^{n-1}})\rangle,\]
  where $\Delta(G^j)=\setsuch{(g,\dots,g)}{g\in G}$ is the diagonal subgroup.
\end{theorem}

\section{Parabolic subgroups}\label{sec:parabolic}
In the context of groups acting on a hyperbolic space, a parabolic
subgroup is the stabilizer of a point on the boundary. We give
here a few general facts concerning parabolic subgroups of branch
groups, and recall some results on growth of groups and sets on
which they act.

More information and uses of parabolic subgroups appear in the context
of representations (Subsection~\ref{subs:quasireg}), Schreier graphs
(Section~\ref{sec:pspace}) and spectrum
(Chapter~\ref{chapter:spectrum}).

\begin{definition}
  A \emdef{ray} $e$ in $\tree$
  is an infinite geodesic starting at the root of $\tree$, or
  equivalently an element of $\partial\tree=Y^\N$.

  Let $G\leq\Aut(\tree)$ be any subgroup acting spherically transitively
  and $e$ be a ray. The associated \emdef{parabolic subgroup} is $P_e=\Stab_G(e)$.
\end{definition}

The following important facts are easy to prove:
\begin{itemize}
\item For any $e\in\partial\tree$, we have
  $\bigcap_{f\in\partial\tree}P_f=\bigcap_{g\in G}P_e^g=1$.
\item Let $e=e_1e_2\dots$ be an infinite ray and define the subgroups
  $P_n=\Stab_G(e_1\dots e_n)$. Then $P_n$ has
  index $m_1m_2 \cdots m_n$ in $G$ (since $G$ acts transitively) and
  \[P_e=\bigcap_{n\in\N}P_n.\]
\item $P$ has infinite index in $G$, and has the same image as $P_n$
  in the quotient $G_n=G/\Stab_G(\LL_n)$.
\end{itemize}

\begin{definition}
  Two infinite sequences $\sigma,\tau:\N\to Y$ are
  \emdef{confinal} if there is an $N\in\N$ such that $\sigma_n=\tau_n$
  for all $n\ge N$.

  Confinality is an equivalence relation, and equivalence classes are
  called \emdef{confinality classes}.
\end{definition}

The following result is due to Volodymyr Nekrashevych and Vitaly
Sushchansky.
\begin{proposition}\label{prop:confinal}
  Let $G$ be a group acting on a regular rooted tree $\tree^{(m)}$, and
  assume that for any generator $g\in G$ and infinite sequence $\tau$,
  the sequences $\tau$ and $\tau^g$ differ only in finitely many
  places. Then the confinality classes are unions of orbits of the
  action of $G$ on $\partial\tree$. If moreover for all $u \in\tree$
  and $v \in\tree\setminus u$ there is some
  $a\in\Stab_G(u)\cap\Stab_G(v)$ transitive on the $m$ subtrees
  below $v$, then the orbits of the action are the confinality classes.
\end{proposition}

\begin{definition}
  The subgroup $H$ of $G$ is \emdef{weakly maximal} if $H$ is of
  infinite index in $G$, but all subgroups of $G$ strictly containing
  $H$ are of finite index in $G$.
\end{definition}
Note that every infinite finitely generated group admits maximal
subgroups, by Zorn's lemma.

However, some branch groups may not contain any infinite-index
maximal subgroups; this is the case for $\Gg$, as was shown by
Ekaterina Pervova (see~\cite{pervova:edsubgroups})

\begin{proposition}
  Let $P$ be a parabolic subgroup of a branch group
  $G$ with branch structure $(L_i,H_i)_{i\in\N}$.
  Then $P$ is weakly maximal.
\end{proposition}
\begin{proof}
  Let $P=Stab_G(e)$ where $e=e_1e_2\dots$ .
  Recall that
  $G$ contains a product of $k_n$ copies of $L_n$ at level $n$, and
  clearly $P$ contains a product of $k_n-1$ copies of $L_n$ at level
  $n$, namely all but the one indexed by the vertex $e_1\dots e_n$.

  Take $g\in G\setminus P$. There is then an $n\in\N$ such that
  $(e_1\dots e_n)^g \neq e_1\dots e_n$, so $\langle P,P^g\rangle$ contains the
  product $H_n=L_n^{1}\times\dots\times L_n^{k_n}$ of $k_n$ copies of $L_n$ at
  level $n$, hence is of finite index in $G$.
\end{proof}

\section{The structure of $\Gg$}\label{sec:Gg}
Recall that $\Gg$, introduced in Subsection~\ref{subs:firstgg}, is the
group acting on the binary tree, generated by the rooted automorphism
$a$ and the directed automorphisms $b,c,d$ satisfying $\psi(b)=(a,c)$,
$\psi(c)=(a,d)$ and $\psi(d)=(1,b)$.

$\Gg$ has $7$ subgroups of index $2$:
\begin{alignat*}{3}
  \langle b,ac\rangle,&&\langle c,ad\rangle,&&\langle d,ab\rangle,\\
  \langle b,a,a^c\rangle,&&\langle c,a,a^d\rangle,&&\langle d,a,a^b\rangle,\\
  &&\Stab_G(1) = \langle b,c,b^a,c^a\rangle.
\end{alignat*}
As can be computed from its presentation~\cite{lysionok:pres} and a
computer algebra system~\cite{gap:manual}, $\Gg$ has the following
subgroup count:

\centerline{\begin{tabular}{c|cc|cc|}
  Index & Subgroups & Normal & In $\Stab_\Gg(\LL_1)$ & Normal\\
  \hline
  1 & 1 & 1 & 0 & 0\\
  2 & 7 & 7 & 1 & 1\\
  4 & 19 & 7 & 9 & 4\\
  8 & 61 & 7 & 41 & 7\\
  16 & 237 & 5 & 169 & 5\\
  32 & 843 & 3 & 609 & 3\\
  \hline
\end{tabular}}

See~\cite{bartholdi:lcs} and~\cite{tcs-s-t:normal} for more
information.

\subsection{Normal closures of generators} They are as follows:
\begin{align*}
  A = \langle a\rangle^\Gg &= \langle a,a^b,a^c,a^d\rangle,&\Gg/A&\cong\Z/2\Z\times\Z/2\Z,\\
  B = \langle b\rangle^\Gg &= \langle b,b^a,b^{ad},b^{ada}\rangle,& \Gg/B &\cong D_8,\\
  C = \langle c\rangle^\Gg &= \langle c,c^a,c^{ad},c^{ada}\rangle,& \Gg/C &\cong D_8,\\
  D = \langle d\rangle^\Gg &= \langle d,d^a,d^{ac},d^{aca}\rangle,& \Gg/D &\cong D_{16}.
\end{align*}

\subsection{Some other subgroups} To complete the picture, we
introduce the following subgroups of $\Gg$:
\begin{align*}
  K &= \langle (ab)^2\rangle^\Gg,&
  L &= \langle (ac)^2\rangle^\Gg,&
  M &= \langle (ad)^2\rangle^\Gg,\\
  \overline B &= \langle B,L\rangle,&
  \overline C &= \langle C,K\rangle,&
  \overline D &= \langle D,K\rangle,\\
  && T = K^2 &= \langle (ab)^4\rangle^\Gg,\\
  T_{(m)} &= \underbrace{T\times\dots\times T}_{2^m},&
  K_{(m)} &= \underbrace{K\times\dots\times K}_{2^m},&
  N_{(m)} &= T_{(m-1)}K_{(m)}.
\end{align*}
\begin{theorem}
  \begin{itemize}
  \item In the Lower Central Series, $\gamma_{2^m+1}(\Gg) = N_{(m)}$ for all
    $m\ge 1$.
  \item In the Derived Series, $K^{(n)}=\rist_\Gg(2n)$ for all $n\ge2$
    and $\Gg^{(n)}=\rist_\Gg(2n-3)$ for all $n\ge3$.
  \item The rigid stabilizers satisfy
    \[\rist_\Gg(n)=\begin{cases} D & \text{ if }n=1,\\
      K_{(n)} & \text{ if }n\ge2.
    \end{cases}\]
  \item The level stabilizers satisfy
    \[\Stab_\Gg(\LL_n)=\begin{cases}\langle b,c,d\rangle^\Gg & \text{ if }n=1,\\
      \langle D,T\rangle & \text{ if }n=2,\\
      \langle N_{(2)},(ab)^4(adabac)^2\rangle & \text{ if }n=3,\\
      \underbrace{\Stab_\Gg(\LL_3)\times\dots\times\Stab_\Gg(\LL_3)}_{2^{n-3}} & \text{ if }n\ge4.
    \end{cases}\]
    Consequently, the index of $\Stab_\Gg(\LL_n)$ is
    \begin{equation}\label{eq:indexStabGg}
      |\Gg/\Stab_\Gg(\LL_n)|=2^{5\cdot2^{n-3}+2}.
    \end{equation}
  \item There is for all $\sigma\in Y^n$ a surjection
    $\cdot_{|\sigma}:\Stab_\Gg(\LL_n)\twoheadrightarrow \Gg$ given by projection
    on the factor indexed by $\sigma$.
  \end{itemize}
  The top of the lattice of normal subgroups of $\Gg$ below
  $\Stab_\Gg(\LL_1)$ is given in Table~\ref{table:Glattice}.
\end{theorem}

\begin{table}
  \[\begin{diagram}\dgARROWLENGTH=1.8em
    \node[3]{\Gg} \arrow{s,-} \node[3]{\text{Index}}\\
    \node[3]{\Stab_\Gg(\LL_1)} \arrow{wsw,-} \arrow{sw,-} \arrow{se,-} \node[3]{2}\\
    \node{\overline B} \arrow{s,-} \arrow{ese,-}
    \node{\overline C} \arrow{s,-} \arrow{se,-}
    \node[2]{\overline D} \arrow{s,-} \arrow{sw,-} \node[2]{4}\\
    \node{B} \arrow{s,-}
    \node{C} \arrow{s,-}
    \node{\Gg'} \arrow{wsw,-} \arrow{sw,-} \arrow{sse,-}
    \node{\Stab_\Gg(\LL_2)} \arrow{s,-} \arrow{ssw,-} \node[2]{8}\\
    \node{K} \arrow{ese,-}
    \node{L} \arrow{se,-}
    \node[2]{D=\rist_\Gg(1)} \arrow{s,-} \node[2]{16}\\
    \node[3]{N_{(1)}=\gamma_3(\Gg)} \arrow{s,-} \arrow{sse,-} \arrow[2]{sw,-}
    \node{M} \arrow{sw,-} \node[2]{32}\\
    \node[3]{K_{(1)}} \arrow{sw,-} \arrow{s,-} \node[3]{64}\\
    \node{} \arrow{s,-}
    \node{\Gg^{(2)}} \arrow{se,-}
    \node{\gamma_4(\Gg)} \arrow{s,-}
    \node{\Stab_\Gg(\LL_3)} \arrow{sw,-} \node[2]{128}\\
    \node{T} \arrow[2]{se,l,..}{64}
    \node[2]{N_{(2)}=\gamma_5(\Gg)} \arrow{wsw,l,..}{64}
    \arrow{s,r,..}{4} \arrow{se,t,..}{16} \node[3]{256}\\
    \node{T_{(1)}}
    \node[2]{K_{(2)}=K'} \arrow{s,r,..}{16}
    \node{\Stab_\Gg(\LL_4)} \arrow{sw,r,..}{4}\\
    \node[3]{N_{(3)}=\gamma_9(\Gg)} \arrow{s,r,..}{16}\\
    \node[3]{K_{(3)}=\Gg^{(3)}}
  \end{diagram}\]
  \caption{The top of the lattice of normal subgroups of $\Gg$ below
    $\Stab_\Gg(\LL_1)$. The index of the inclusions are indicated next to the edges.}
  \label{table:Glattice}
\end{table}

\begin{corollary}
  The closure of $\Gg$ in $\Aut(\tree)$ has Hausdorff
  dimension $5/8$.
\end{corollary}

\subsection{The Subgroup $P$}\label{subs:subP}
Let $e$ be the ray $2^\infty$ and let $P$ be the corresponding
parabolic subgroup. We describe completely its structure as
follows:

\begin{theorem}\label{theorem:decomP}
  $P/P'$ is an infinite elementary $2$-group generated by the images of
  $c$, $d=(1,b)$ and of all elements of the form $(1,\dots,1,(ac)^4)$
  in $\rist_\Gg(n)$ for $n\in\N$. The
  following decomposition holds:
  \[P = \bigg(B\times \Big(\big(K\times ((K\times\dots)\rtimes\langle (ac)^4\rangle)\big)\rtimes\langle b,(ac)^4\rangle\Big)\bigg)\rtimes\langle c,(ac)^4\rangle,\]
  where each factor (of nesting $n$) in the decomposition acts on the
  subtree just below some $e_n$ but not containing $e_{n+1}$.
\end{theorem}
Note that we use the same notation for a subgroup $B$ or $K$ acting on
a subtree, keeping in mind the identification of a subtree with the
original tree. Note also that $\psi$ is omitted when it would make the
notations too heavy.

\begin{proof}
  Define the following subgroups of $\Gg_n$:
  \begin{gather*}
    H_n = \langle b,c\rangle^{\Gg_n};\qquad  B_n = \langle b\rangle^{\Gg_n};\qquad K_{(n)} = \langle (ab)^2\rangle^{\Gg_n};\\
    Q_n = B_n\cap P_n;\qquad R_n = K_{(n)}\cap P_n.
  \end{gather*}
  Then the theorem follows from the following proposition.
\end{proof}

\begin{proposition}
  These subgroups have the following structure:
  \begin{align*}
    P_n &= (B_{n-1}\times Q_{n-1})\rtimes\langle c,(ac)^4\rangle;\\
    Q_n &= (K_{n-1}\times R_{n-1})\rtimes\langle b,(ac)^4\rangle;\\
    R_n &= (K_{n-1}\times R_{n-1})\rtimes\langle (ac)^4\rangle.
  \end{align*}
\end{proposition}
\begin{proof}
  A priori, $P_n$, as a subgroup of $H_n$, maps in $(B_{n-1}\times
  B_{n-1})\rtimes\langle(a,d),(d,a)\rangle$. Restricting to those
  pairs that fix $e_n$ gives the result. Similarly, $Q_n$, as a
  subgroup of $B_n$, maps in $(K_{n-1}\times
  K_{n-1})\rtimes\langle(a,c),(c,a)\rangle$, and $R_n$, as a subgroup
  of $K_n$, maps in $(K_{n-1}\times
  K_{n-1})\rtimes\langle(ac,ca),(ca,ac)\rangle$.
\end{proof}
\begin{corollary} The group $\Gg_n$ and its subgroups
  $H_n,B_n,K_n,P_n,R_n,Q_n$ are arranged in a lattice
  \[\begin{diagram}
    \node[3]{\Gg_n}\arrow{s,l,-}{\langle a\rangle}\\
    \node[3]{H_n}\arrow{sw,l,-}{\langle d,d^a\rangle}\arrow{se,l,-}{2^{n-1}}\\
    \node[2]{B_n}\arrow{sw,l,-}{\langle b\rangle}\arrow{se,l,-}{2^{n-1}}
    \node[2]{P_n}\arrow{sw,r,-}{2^2}\\
    \node{K_n}\arrow{se,r,-}{2^{n-1}}\node[2]{Q_n}\arrow{sw,r,-}{2}\\
    \node[2]{R_n}
  \end{diagram}\]
  where the quotients or the indices are represented next to the edges.
\end{corollary}

\section{The structure of $\FGg$}\label{sec:FGg}
Recall that $\FGg$ is the group acting on the ternary tree,
generated by the rooted automorphism $a=((1,2,3))$ and the
directed automorphism $t$ satisfying $(t)\psi=(a,1,t)$.

Define the elements $x=at$, $y=ta$ of $\FGg$. Let $K$ be the
subgroup of $\FGg$ generated by $x$ and $y$, and let $L$ be the
subgroup of $K$ generated by $K'$ and cubes in $K$. Write
$H=\Stab_\FGg(\LL_1)$.

\begin{proposition}
  We have the following diagram of normal subgroups:
  \[\begingroup\unitlength=0.01pt\begin{diagram}
    \node[2]{\FGg}\arrow{sw,l,-}{\langle a|\,a^3\rangle}
    \arrow{se,l,-}{\langle a|\,a^3\rangle}\\
    \node{K}\arrow{se,r,-}{\langle x|\,x^3\rangle}
    \node[2]{\Stab_\FGg(\LL_1)}\arrow{sw,r,-}{\langle t|\,t^3\rangle}\\
    \node[2]{\FGg'=K\cap H=[K,H]}\arrow{s,-}\\
    \node[2]{L=\langle K',K^3\rangle=\gamma_3(\FGg)}\arrow{sw,-}\arrow{se,-}\\
    \node{K'}\arrow{s,-}
    \node[2]{H'=(\FGg'\times\FGg'\times\FGg')\psi^{-1}=\Stab_\FGg(\LL_2)}\\
    \node{\langle L\times L\times
    L,x^3y^{-3},[x,y^3]\rangle=\gamma_4(\FGg)}\arrow{s,-}\\
    \node{\langle L\times L\times
    L,[x,y^3]\rangle=\gamma_5(\FGg)}
  \end{diagram}\endgroup\]
  where the quotients are represented next to the edges; all edges
  represent normal inclusions of index $3$. Furthermore $L=K\cap
  (K\times K\times K)\psi^{-1}$.
\end{proposition}
\begin{proof}
  First we prove $K$ is normal in $\FGg$, of index $3$, by writing
  $y^t=x^{-1}y^{-1}$, $y^{a^{-1}}=y^{-1}x^{-1}$, $y^{t^{-1}}=y^a=x$;
  similar relations hold for conjugates of $x$.  A transversal of $K$
  in $\FGg$ is $\langle a\rangle$.  All subgroups in the diagram are
  then normal.

  Since $[a,t]=y^{-1}x = t^at^{-1}$, we clearly have $\FGg'<K\cap H$.
  Now as $\FGg'\neq K$ and $\FGg'\neq H$ and $\FGg'$ has index
  $3^2$, we must have $\FGg'=K\cap H$. Finally
  $[a,t]=[x,t]^{t^{-1}}$, so $\FGg'=[K,H]$.

  Next $x^3=[a,t][t,a^{-1}][a^{-1},t^{-1}]$ and similarly for $y$, so
  $K^3<\FGg'$ and $L<\FGg'$. Also,
  $[x,y]\psi=(y^{-1},y^{-1},x^{-1})$ and $(x^3)\psi=(y,x,y)$ both
  belong to $K\times K\times K$, while $[a,t]$ does not; so $L$ is a
  proper subgroup of $\FGg'$, of index $3$ (since $K/L$ is
  the elementary abelian group $(\Z/3\Z)^2$ on $x$ and $y$).

  Consider now $H'$. It is in $\Stab_\FGg(\LL_2)$ since
  $H=\Stab_\FGg(\LL_1)$.  Also, $[t,t^a]=y^3[y^{-1},x]$ and similarly
  for other conjugates of $t$, so $H'<L$, and
  $[t,t^a]\psi=([a,t],1,1)$, so
  $(H')\psi=\FGg'\times\FGg'\times\FGg'$. Finally $H'$ it is of index
  $3$ in $L$ (since $H/H'=(\Z/3\Z)^3$ on $t,t^a,t^{a^{-1}}$), and
  since $\Stab_\FGg(\LL_2)$ is of index $3^4$ in $\FGg$ (with quotient
  $\Z/3\Z\wr\Z/3\Z$) we have all the claimed equalities.
\end{proof}

\begin{proposition}\label{prop:Gammafractal}
  $\FGg$ is a just-infinite fractal group, is regular branch
  over $\FGg'$, and has the congruence property.
\end{proposition}
\begin{proof}
  $\FGg$ is fractal by Lemma~\ref{lemma:fractal} and the nature of the
  map $\psi$.  By direct computation, $[\FGg:\FGg'] =
  [\FGg':(\FGg'\times\FGg'\times\FGg')\psi^{-1}] =
  [(\FGg'\times\FGg'\times\FGg')\psi^{-1}:\FGg''] = 3^2$, so
  $\FGg$ is branched on $\FGg'$. Then $\FGg''=\gamma_5(\FGg)$,
  as is shown in~\cite{bartholdi:lcs}, so $\FGg''$ has finite index
  and $\FGg$ is just-infinite by Theorem~\ref{theorem:jicriterion}.

  $\FGg'\ge\Stab_\FGg(\LL_2)$, so $\FGg$ has the congruence property.
\end{proof}

\begin{proposition}
  Writing $\langle S\rangle$ for the $3$-abelian quotient of $\langle
  S\rangle$, we have exact sequences
  \begin{align*}
    1\to\FGg'\times\FGg'\times\FGg'\to (H)\psi \to\langle t,t^a,t^{a^2}\rangle_{3-ab},\\
    1\to\FGg'\times\FGg'\times\FGg'\to (\FGg')\psi \to\langle [a,t],[a^2,t]\rangle_{3-ab}.
  \end{align*}
\end{proposition}

\begin{theorem}
  The subgroup $K$ of $\FGg$ is torsion-free; thus $\FGg$ is
  virtually torsion-free.
\end{theorem}

\begin{proposition}
  The finite quotients $\FGg_n=\FGg/\Stab_\FGg(\LL_n)$ of $\FGg$ have
  order $3^{3^{n-1}+1}$ for $n\ge2$, and $3$ for $n=1$.
\end{proposition}
\begin{proof}
  Follows immediately from $[\FGg:\FGg']=3^2$ and
  $[\FGg':(\FGg'\times\FGg'\times\FGg')\psi^{-1}]=3^2$.
\end{proof}

\begin{corollary}
  The closure of $\FGg$ in $\Aut(\tree)$ is
  isomorphic to the profinite completion $\widehat\FGg$ and is a
  pro-$3$-group. It has Hausdorff dimension $1/3$.
\end{corollary}

\section{The structure of $\BGg$}\label{sec:BGg}
Recall that $\BGg$ is the group acting on the ternary tree,
generated by the rooted automorphism $a=((1,2,3))$ and the
directed automorphism $t$ defined by $(t)\psi =(a,a,t)$.

Define the elements $x=ta^{-1}$, $y=a^{-1}t$ of $\BGg$, and let $K$
be the subgroup of $\BGg$ generated by $x$ and $y$. Then $K$ is
normal in $\BGg$, because $x^t=y^{-1}x^{-1}$, $x^a=x^{-1}y^{-1}$,
$x^{t^{-1}}=x^{a^{-1}}=y$, and similar relations hold for conjugates of $y$.
Moreover $K$ is of index $3$ in $\BGg$, with transversal $\langle a\rangle$.
Write $H=\Stab_\BGg(\LL_1)$.

\begin{lem}\label{lemma:GammaHK}
  $H$ and $K$ are normal subgroups of index $3$ in $\BGg$, and
  $\BGg'=\Stab_K(\LL_1)=H\cap K$ is of index $9$; furthermore $(H \cap
  K)\psi \triangleleft K\times K\times K$. For any element
  $g=(u,v,w)\in (H\cap K)\psi$ one has $wvu\in H\cap K$.
\end{lem}
\begin{proof}
  First note that $\Stab_K(\LL_1)=\langle x^3,y^3,xy^{-1},y^{-1}x\rangle$,
  for every word in $x$ and $y$ whose number of $a$'s is divisible by
  $3$ can be written in these generators. Then compute
  \begin{alignat*}{2}
    (x^3)\psi &=(y,x^{-1}y^{-1},x),&\qquad (y^3)\psi &=(x^{-1}y^{-1},x,y),\\
    (xy^{-1})\psi &=(1,x^{-1},x),&\qquad (y^{-1}x)\psi &=(y,1,y^{-1}).
  \end{alignat*}
  The last assertion is also checked on this computation.
\end{proof}

\begin{proposition}
  Writing $c=[a,t]=x^{-1}y^{-1}x^{-1}$ and $d=[x,y]$, we have the
  following diagram of normal subgroups:
  \[\begin{diagram}
    \node[2]{\BGg}\arrow{sw,l,-}{\langle a|\,a^3\rangle}
    \arrow{se,l,-}{\langle a|\,a^3\rangle}\\
    \node{K}\arrow{se,r,-}{\langle x,y|\,x^3,y^3,x=y\rangle}
    \node[2]{H}\arrow{sw,r,-}{\langle t_0,t_1,t_2|\,t_0^3,t_0=t_1=t_2\rangle}\\
    \node[2]{\BGg'=\langle c, c^t, c^{a^{-1}}, c^{at} \rangle=K\cap H=[K,H]}\arrow{sw,t,-}{\Z^2}\arrow{se,t,-}{(\Z/3\Z)^2}\\
    \node{K'=\langle d, d^t, d^{a^{-1}}, d^{at} \rangle}\arrow{se,b,-}{\Z^2}
    \node[2]{H'}\arrow{sw,-}\\
    \node[2]{\BGg''=(K'\times K'\times
      K')\psi^{-1}}\arrow{s,r,-}{(\Z/3\Z)^2}\\
    \node[2]{K''}
  \end{diagram}\]
  where the quotients are represented next to the edges; additionally,
  \begin{gather*}
    K/K' = \langle x,y|\,[x,y]\rangle\cong\Z^2,\\
    \BGg'/\BGg'' = \langle
    c,c^t,c^{a^{-1}},c^{at}|\,[c,c^t],\dots\rangle\cong \Z^4,\\
    K'/K'' = \langle d,d^t,d^{a^{-1}},d^{at}|\,[d,d^t],\dots,(d/d^{at})^3,(d^{a^{-1}}/d^t)^3\rangle\cong\Z^2\times(\Z/3\Z)^2.
  \end{gather*}

  Writing each subgroup in the generators of the groups above it, we
  have
  \begin{gather*}
    K = \langle x=at^{-1},y=a^{-1}t\rangle,\\
    H = \langle t,t_1=t^a,t_2=t^{a^{-1}}\rangle,\\
    \BGg' = \langle b_1=xy^{-1},b_2=y^{-1}x,b_3=x^3,b_4=y^3\rangle
    = \langle c_1=tt_1^{-1},c_2=tt_1t,c_3=tt_2^{-1},c_4=tt_2t\rangle.
  \end{gather*}
\end{proposition}

\begin{corollary}
  The congruence property does not hold for $\BGg$; nor is
  it regular branch.
\end{corollary}

\begin{proposition}\label{prop:GammaBfractal}
  $\BGg$ is a fractal group, is weakly branch, and
  just-nonsolvable.
\end{proposition}
\begin{proof}
  $\BGg$ is fractal by Lemma~\ref{lemma:fractal} and the
  nature of the map $\psi$. The subgroup $K$ described above has an
  infinite-index derived subgroup $K'$ (with infinite cyclic
  quotient), from which we conclude that $\BGg$ is not
  just-infinite; indeed $K'$ is normal
  in $\BGg$ and
  $\BGg/K'\cong\Z^2\rtimes\left(\begin{smallmatrix}-1&1\\
      -1&0\end{smallmatrix}\right)$ is infinite.
\end{proof}

\begin{proposition}
  The subgroup $K$ of $\BGg$ is torsion-free; thus
  $\BGg$ is virtually torsion-free.
\end{proposition}
\begin{proof}
  For $1\neq g\in K$, let $|g|_t$, the $t$-length of $g$, denote the minimal
  number of $t^{\pm1}$'s required to write $g$ as a word over the
  alphabet $\{a^{\pm1},t^{\pm1}\}$. We will show by induction on
  $|g|_t$ that $g$ is of infinite order.

  First, if $|g|_t=1$, i.e., $g\in\{x^{\pm1},y^{\pm1}\}$, we
  conclude from $(x^3)\psi=(*,*,x)$ and $(y^3)\psi=(*,*,y)$ that $g$
  is of infinite order.

  Suppose now that $|g|_t>1$, and $g\in
  \Stab_\BGg(\LL_n)\setminus\Stab_\BGg(\LL_{n+1})$. Then there is some
  sequence $\sigma$ of length $n$ that is fixed by $g$ and such that
  $g_{|\sigma}\not\in H$. By Lemma~\ref{lemma:GammaHK},
  $g_{|\sigma}\in K$, so it suffices to show that all $g\in K\setminus
  H$ are of infinite order.

  Such a $g$ can be written as $(u,v,w)\psi^{-1}z$ for some
  $(u,v,w)\in (K\cap H)\psi$ and $z\in\{x^{\pm1},y^{\pm1}\}$; by
  symmetry let us suppose $z=x$. Then
  $g^3=(uavawt,vawtua,wtuava)\psi^{-1}=(g_0,g_1,g_2)\psi^{-1}$, say.
  For any $i$, we have $|g_i|_t\le|g|_t$, because all the components
  of $(x)\psi$ and $(y)\psi$ have $t$-length $\le1$. We distinguish
  three cases:
  \begin{enumerate}
  \item $g_i=1$ for some $i$. Then consider the image $\overline{g_i}$
    of $g_i$ in $\BGg/\BGg'$. By
    Lemma~\ref{lemma:GammaHK}, $wvu\in G'$, so
    $\overline{g_i}=1=\overline{a^2t}$. But this is a contradiction,
    because $\BGg/\BGg'$ is elementary abelian
    of order $9$, generated by the independent images $\overline a$
    and $\overline t$.
  \item $0<|g_i|_t<|g|_t$ for some $i$. Then by induction
    $g_i$ is of infinite order, so $g^3$ too, and $g$ too.
  \item $|g_i|_t=|g|_t$ for all $i$. We repeat the
    argument with $g_i$ substituted for $g$. As there are finitely many
    elements $h$ with $|h|_t=|g|_t$, we will eventually
    reach either an element of shorter length or an element
    already considered. In the latter case we obtain a relation
    of the form $(g^{3^n})\psi^n=(\dots,g,\dots)$ from which $g$ is
    seen to be of infinite order.
  \end{enumerate}
\end{proof}

\begin{proposition}
  The finite quotients $\BGg_n=\BGg/\Stab_\BGg(\LL_n)$ of $\BGg$ have
  order $3^{\frac14(3^n+2n+3)}$ for $n\ge2$, and $3^{\frac12(3^n-1)}$
  for $n\le2$.
\end{proposition}
\begin{proof}
  Define the following family of two-generated finite abelian groups:
  \[A_n = \begin{cases}\langle x,y|\,x^{3^{n/2}},y^{3^{n/2}},[x,y]\rangle&\text{ if }n\equiv0[2],\\
    \langle x,y|\,x^{3^{(n+1)/2}},y^{3^{(n+1)/2}},(xy^{-1})^{3^{(n-1)/2}},[x,y]\rangle&\text{ if }n\equiv1[2].
  \end{cases}\]
  First suppose $n\ge 2$; Consider the diagram of groups described
  above, and quotient all the groups by $\Stab_\BGg(\LL_n)$. Then the
  quotient $K/K'$ is isomorphic to $A_n$, generated by $x$ and $y$,
  and the quotient $K'/\BGg''$ is isomorphic to $A_{n-1}$, generated
  by $[x,y]$ and $[x,y]^t$. As $|A_n|=3^n$, the index of $K_n'$ in
  $\BGg_n$ is $3^{n+1}$ and the index of $\BGg_n''$ is $3^{2n}$.  Then
  as $\BGg''_n\cong K_{n-1}^3$ and $|\BGg_2''|=1$ we deduce by
  induction that $|\BGg''_n|=3^{\frac14(3^n-6n+3)}$ and
  $|K'_n|=3^{\frac14(3^n-2n-1)}$, from which
  $|\BGg_n|=3^{2n}+|\BGg''_n|=3^{\frac14(3^n+2n+3)}$ follows.

  For $n\le2$ we have $\BGg_n=\Aut(\tree)_n=\Z/3\wr\dots\wr\Z/3$.
\end{proof}

\begin{corollary}
  The closure of $\BGg$ in $\Aut(\tree)$ has Hausdorff
  dimension $1/2$.
\end{corollary}

\begin{proposition}
  We have exact sequences
  \begin{align*}
    1\to K'\times K'\times K'\to (H)\psi \to\Z^4\rtimes\Z/3\Z\to 1,\\
    1\to K'\times K'\times K'\to (K')\psi \to\Z^2\to 1.
  \end{align*}
\end{proposition}

\section{The structure of $\GSg$}\label{sec:GSg}
Recall that $\GSg$ is the group acting on the ternary tree,
generated by the rooted automorphism $a=((1,2,3))$ and the
directed automorphism $t$ satisfying $(t)\psi=(a,a^{-1},t)$. Write
$H=\Stab_\GSg(\LL_1)$.

\begin{proposition}
  We have the following diagram of normal subgroups:
  \[\begin{diagram}
    \node{\GSg}\arrow{s,r,-}{\langle a|\,a^3\rangle}\\
    \node{H=\Stab_\GSg(\LL_1)}\arrow{s,r,-}{\langle t|\,t^3\rangle}\\
    \node{\GSg'=[G,H]}\arrow{s,r,-}{[a,t]}\\
    \node{\gamma_3(\GSg)=\GSg^3=
      \Stab_{\GSg}(\LL_2)}\arrow{s,r,-}{(at)^3}\\
    \node{H'=(\GSg'\times\GSg'\times\GSg')\psi^{-1}}
  \end{diagram}\]
  where the quotients are represented next to the arrows; all edges
  represent normal inclusions of index $3$.
\end{proposition}

\begin{proposition}\label{prop:GammaBBfractal}
  $\GSg$ is a just-infinite fractal group, and is a
  regular branch group over $\GSg'$.
\end{proposition}

\begin{proposition}
  $\GSg'\ge\Stab_{\GSg}(\LL_2)$, so
  $\GSg$ has the congruence property.
\end{proposition}

\date{October 27, 2002}
\chapter{Central Series, Finiteness of Width and Associated Lie Algebras}
\label{chapter:lie} In this chapter we study the lower central,
lower $p$-central, and dimension series of basic examples of
branch groups, and describe the associated Lie algebras. This
chapter is very much connected to the previous one.

We exhibit two branch groups of finite width: $\Gg$ and $\FGg$, and
describe the ``Lie graph'' of their associated Lie algebras. We show
that the Gupta-Sidki $3$-group has unbounded width, and its Lie
algebra has growth of degree $n^{\log3/\log(1+\sqrt2)-1}$; we also
describe its Lie graph.

For all regular branch groups, the corresponding Lie algebras have
polynomial growth (usually of non-integral degree), see
Theorem~\ref{theorem:lcs:poly}. The technique used below,
described in~\cite{bartholdi:lcs}, is an extension of the methods
of~\cite{bartholdi-g:lie}.

We start by recalling the famous construction, due to
\WMagnus~\cite{magnus:lie}.

\section{$N$-series}
\begin{definition}
  Let $G$ be a group. An \emdef{$N$-series} is series $\{G_n\}$ of
  normal subgroups with $G_1=G$, $G_{n+1}\le G_n$ and $[G_m,G_n]\leq
  G_{m+n}$ for all $m,n\ge1$.
  The associated Lie ring is
  \[\Lie(G) = \bigoplus_{n=1}^\infty \Lie_n,\]
  with $\Lie_n=G_n/G_{n+1}$ and the bracket operation
  $\Lie_n\otimes\Lie_m\to\Lie_{m+n}$ induced by commutation in $G$.

  For $p$ a prime, an \emdef{$N_p$-series} is an $N$-series $\{G_n\}$
  such that $\mho_p(G_n)\le G_{pn}$, and the associated Lie ring is a
  restricted Lie algebra over $\Fp$~\cite{jacobson:restr},
  \[\Lie_\Fp(G) = \bigoplus_{n=1}^\infty \Lie_n,\]
  with the Frobenius operation $\Lie_n\to\Lie_{pn}$ induced by raising
  to the power $p$ in $G$. (Recall the definition of $\mho_d(G)$ from
  Section~\ref{sec:powers}.)
\end{definition}
The standard example of $N$-series is the \emdef{lower central series},
$\{\gamma_n(G)\}_{n=1}^\infty$, given by $\gamma_1(G)=G$ and

$\gamma_n(G)=[G,\gamma_{n-1}(G)]$, or the \emdef{lower $p$-central
  series} or \emdef{Frattini series} given by $P_1(G)=G$ and
$P_n(G)=[G,P_{n-1}(G)]\mho_p(P_{n-1}(G))$. It differs from the lower
central series in that its successive quotients are all elementary
$p$-groups.

The standard example of $N_p$-series is the \emdef{dimension series}, also
known as the Zassenhaus~\cite{zassenhaus:ordnen},
Jennings~\cite{jennings:gpring}, Lazard~\cite{lazard:nilp} or Brauer
series, given by $G_1=G$ and $G_n=[G,G_{n-1}]\mho_p(G_{\lceil n/p
\rceil})$, where $\lceil n/p \rceil$ is the least integer greater than or
equal to $n/p$. It can alternately be described, by a result of
Lazard~\cite{lazard:nilp}, as
\[G_n = \prod_{i\cdot p^j\ge n}\mho_{p^j}(\gamma_i(G)),\]
or as
\[G_n = \setsuch{g\in G}{g-1\in\Delta^n},\]
where $\Delta$ is the augmentation (or fundamental) ideal of the
group algebra $\Fp G$. Note that this last definition extends to
characteristic $0$, giving a graded Lie algebra $\Lie_\Q(G)$ over
$\Q$. In that case, the subgroup $G_n$ is the isolator of
$\gamma_n(G)$:
\[G_n = \sqrt{\gamma_n(G)} = \setsuch{g\in G}{\langle g\rangle\cap\gamma_n(G)\neq\{1\}}.\]

We mention finally for completeness another $N_p$-series, the \emdef{Lie
dimension series} $L_n(G)$ defined by
\[L_n(G) = \setsuch{g\in G}{g-1\in\Delta^{(n)}},\]
where $\Delta^{(n)}$ is the $n$-th Lie power of $\Delta<\Bbbk G$,
given by $\Delta^{(1)}=\Delta$ and
$\Delta^{(n)}=[\Delta^{(n)},\Delta]=\setsuch{xy-yx}{x\in\Delta^{(n)},\;y\in\Delta}$.
It is then known~\cite{passi-s:liedim} that
\[L_n(G) = \prod_{(i-1)\cdot p^j\ge n}\mho_{p^j}(\gamma_i(G))\]
if $\Bbbk$ is of characteristic $p$, and
\[L_n(G) = \sqrt{\gamma_n(G)}\cap [G,G]\]
if $\Bbbk$ is of characteristic $0$.

\begin{definition}
  An $N$-series $\{G_n\}$ has \emdef{finite width} if there is a
  uniform constant $W$ such that $l_n:=\rank G_n/G_{n+1}\le W$ holds
  for all $n$, where $\rank A$ is the minimal number of generators of
  the abelian group $A$.  A group has \emdef{finite width} if its
  lower central series has finite width --- this definition comes
  from~\cite{klass-lg-p:fw}.
\end{definition}

The following result is well-known, and shows that sometimes the Lie ring
$\Lie(G)$ is actually a Lie algebra over $\Fp$.
\begin{lemma}\label{lemma:liealg}
  Let $G$ be a group generated by a set $S$. Let $\Lie(G)$ be the Lie
  ring associated to the lower central series.
  \begin{enumerate}
  \item If $S$ is finite, then $\Lie_n$ is a finite-rank $\Z$-module
    for all $n$.
  \item If there is a prime $p$ such that all generators $s\in S$ have
    order $p$, then the Lie algebra associated to the lower
    $p$-central series coincides with $\Lie$. As a consequence,
    $\Lie_n$ is a vector space over $\Fp$ for all $n$.
  \end{enumerate}
\end{lemma}

We return to the lower $p$-central series of $G$. Consider the graded
algebra
\[\overline{\Fp G}=\bigoplus_{n\in\N}\Delta^n/\Delta^{n+1}.\]

A fundamental result connecting $\Lie_\Fp(G)$ and $\overline{\Fp G}$ is
the
\begin{proposition}[Quillen~\cite{quillen:ab}]
  $\overline{\Fp G}$ is the enveloping $p$-algebra of $\Lie_\Fp(G)$.
\end{proposition}

The Poincar\'e-Birkhoff-Witt Theorem then gives a basis of $\overline{\Fp
G}$ consisting of monomials over a basis of $\Lie_\Fp(G)$, with exponents
at most $p-1$. As a consequence, we have the
\begin{proposition}[Jennings~\cite{jennings:gpring}]\label{prop:sumprod}
  Let $G$ be a group, and let $\sum_{n\ge1}l_n\hbar^n$ be the
  Hilbert-Poincar\'e series of $\Lie_\Fp(G)$. Then
  \[\gr(\overline{\Fp G})
  =\prod_{n=1}^\infty\left(\frac{1-\hbar^{pn}}{1-\hbar^n}\right)^{l_n}.\]
\end{proposition}

As a consequence, we have the following proposition, firstly
observed by Bereznii (for a proof
see~\cite{petrogradsky:polynilpotent} and~\cite{bartholdi-g:lie}:
\begin{proposition}[]\label{prop:grlgr}
  Let $G$ be a group and expand the power series
  $\gr(\Lie_\Fp(G))=\sum_{n\ge1}l_n\hbar^n$ and
  $\gr(\overline{\Fp G})=\sum_{n\ge0}f_n\hbar^n$. Then
  \begin{enumerate}
  \item $\{f_n\}$ grows exponentially if and only if $\{l_n\}$ does, and
    we have
    \[\limsup_{n\to\infty}\frac{\ln l_n}n
    =\limsup_{n\to\infty}\frac{\ln f_n}n.\]
  \item If $l_n\sim n^d$, then $f_n\sim e^{n^{(d+1)/(d+2)}}$.
  \end{enumerate}
\end{proposition}

Finally, we recall a connection between the growth of $G$ and that of
$\overline{\Fp G}$:
\begin{proposition}[\cite{bartholdi-g:lie}, Lemma~2.5]\label{prop:gpgr}
  Let $G$ be a group generated by a finite set $S$. Then
  \[\frac{\gr(G)}{1-\hbar}\ge\frac{\gr(\overline{\Fp G})}{1-\hbar},\]
  the inequality being valid coefficient-wise.
\end{proposition}

The following result exhibits a ``gap in the spectrum'' of growth, for
residually-$p$ groups:
\begin{corollary}[\cite{grigorchuk:hp,bartholdi-g:lie}]\label{cor:lcslowerbd}
  Let $G$ be a residually-$p$ group for some prime $p$. Then the
  growth of $G$ is either polynomial, in case $G$ is virtually
  nilpotent, or is at least $e^{\sqrt n}$.
\end{corollary}

\section{Lie algebras of branch groups}
Our main purpose, in this section, is to illustrate the following result
by examples:
\begin{theorem}\label{theorem:lcs:poly}
  Let $G$ be a finitely generated regular branch group and $\Lie G$
  the Lie ring associated to its lower central series. Then $\gr(\Lie
  G)$ has polynomial growth (not necessarily of integer degree).
\end{theorem}
Its proof relies on branch portraits, introduced in
Section~\ref{sec:bp}.
\begin{proof}[Sketch of proof]
  Let $G$ be regular branch over $K$. The Lie algebra of $G$ is
  isomorphic to that of $\overline G$, so we consider the latter. For
  each $n\in\N$, consider the set of branch portraits associated to
  $\gamma_n(G)$. Since $K$ has finite index, it suffices to consider
  only the $\gamma_n(G)\le K$. Let $n(\ell)$ be minimal such that the
  branch portraits of $\gamma_{n(\ell)}(G)$ are trivial in their first
  $\ell$ levels. It suffices to show that this function is
  exponential. Consider the portraits of $\gamma_n(G)$, for
  $n(\ell)\le n\le n(\ell+1)$. For $n$ close to $n(\ell)$, there will
  be all portraits that are trivial except in a subtree at level
  $\ell$. Then for larger $n$ there will be, using commutation with a
  generator that is nontrivial at the root vertex, portraits trivial
  except in two subtrees, where they have labels $P$ and $P^{-1}$
  respectively. As $n$ becomes larger and larger, the only remaining portraits will be
  those whose labels in all subtrees at level $\ell$ are identical.
  This passage from one level to the next is exponential.
\end{proof}

We obtain an explicit description of the lower central series in
several cases, and show:
\begin{itemize}
\item For the first Grigorchuk group $\Gg$, the Grigorchuk supergroup
  $\Sg$ and the Fabrykowski-Gupta group $\FGg$, the Lie
  algebras $L$ and $\Lie_\Fp$ have finite width.
\item For the Gupta-Sidki group $\GSg$, the Lie algebras
  $L$ and $\Lie_\Fp$ have polynomial growth of degree
  $d=\log3/\log(1+\sqrt2)-1$.
\end{itemize}

The first result obtained in that direction is due to Alexander Rozhkov.
He proved in~\cite{rozhkov:lcs} that for the first Grigorchuk group $\Gg$
one has
\[\rank \gamma_n(\Gg)/\gamma_{n+1}(\Gg) = \begin{cases}
  3 & \text{ if }n=1,\\
  2 & \text{ if }n=2^m+1+r,\text{ with }0\le r<2^{m-1},\\
  1 & \text{ if }n=2^m+1+r,\text{ with }2^{m-1}\le r<2^m.\end{cases}\]

However, the Lie algebra structure contained in an $N$-series
$\{G_n\}$ is much richer than the series $\{\rank G_n/G_{n+1}\}$,
and we will give a fuller description of the $\gamma_n(\Gg)$
below.

All our examples will satisfy the following conditions:
\begin{enumerate}
\item $G$ is finitely generated by a set $S$;
\item there is a prime $p$ such that all $s\in S$ have order $p$.
\end{enumerate}
It then follows from Lemma~\ref{lemma:liealg} that
$\gamma_n(G)/\gamma_{n+1}(G)$ is a finite-dimensional vector space
over $\F_p$, and therefore that $\Lie(G)$ is a Lie algebra over $\F_p$
that is finite at each dimension. Clearly the same property holds for
the restricted algebra $\Lie_{\F_p}(G)$.

We describe such Lie algebras as oriented labelled graphs, in the
following notation:
\begin{definition}
  Let $\Lie=\bigoplus_{n\ge1}\Lie_n$ be a graded Lie algebra over
  $\F_p$, and choose a basis $B_n$ and a scalar product
  $\langle|\rangle$ of $\Lie_n$ for all $n\ge1$.

  The \emdef{Lie graph} associated to these choices is an abstract
  graph. Its vertex set is $\bigcup_{n\ge1}B_n$, and each vertex $x\in
  B_n$ has a degree, $n=\deg x$. Its edges are labelled as $\alpha x$,
  with $x\in B_1$ and $\alpha\in\F_p$, and may only connect a vertex
  of degree $n$ to a vertex of degree $n+1$. For all $x\in B_1$, $y\in
  B_n$ and $z\in B_{n+1}$, there is an edge labelled $\langle
  [x,y]|z\rangle x$ from $y$ to $z$.

  If $\Lie$ is a restricted algebra of $\F_p$, there are additional
  edges from vertices of degree $n$ to vertices of degree $pn$. For
  all $x\in B_n$ and $y\in B_{pn}$, there is an edge labelled $\langle
  x^p|y\rangle\cdot p$ from $x$ to $y$.

  Edges labelled $0x$ are naturally omitted, and edges labelled $1x$
  are simply written $x$.
\end{definition}

There is some analogy between this definition and that of a Cayley
graph --- this topic will be developed in Section~\ref{sec:pspace}.
The generators (in the Cayley sense) are simply chosen to be the
$\ad(x)$ with $x$ running through $B_1$, a basis of $G/[G,G]$.

As an example of Lie graph, let $G$ be the infinite dihedral group
$D_\infty=\langle a,b|\,a^2,b^2\rangle$. Then $\gamma_n(G)=\langle
(ab)^{2^{n-1}}\rangle$ for all $n\ge2$, and its Lie ring is again a
Lie algebra over $\F_2$, with Lie graph
\begingroup\makeatletter\renewcommand\dggeometry{\dg@XGRID=\thr@@
  \dg@YGRID=\tw@
  \unitlength=0.018pt\relax}
\[\begin{diagram}
  \node{a}\arrow{se,t}{b}\\
  \node[2]{(ab)^2} \arrow{e,t}{a,b} \node{(ab)^4} \arrow{e,t}{a,b}
  \node{(ab)^8} \arrow{e,t,-,..}{a,b} \node{}\\
  \node{b}\arrow{ne,b}{a}
\end{diagram}\]
\endgroup

Note that the lower $2$-central series of $G$ is different: we
have $G_{2^n}=G_{2^n+1}=\dots=G_{2^{n+1}-1}=\gamma_{n+1}(G)$, so
the Lie graph of $\Lie_{\F_2}(G)$ is
\begingroup\makeatletter\renewcommand\dggeometry{\dg@XGRID=\@ne
  \dg@YGRID=\@ne
  \unitlength=0.03pt\relax}
\[\begin{diagram}
  \node{a}\arrow{se,t}{b}\\
  \node[2]{(ab)^2} \arrow[2]{e,t}{\cdot2} \node[2]{(ab)^4}
  \arrow[4]{e,t}{\cdot2} \node[4]{(ab)^8} \arrow[2]{e,t,-,..}{\cdot2}\\
  \node{b}\arrow{ne,b}{a}
\end{diagram}\]
\endgroup

We shall also need the following notation: let $G$ be a regular branch
group over $K$, embedded in $G\wr(\Z/m\Z)$. For all $i\in\N$ and all
$g\in G$ define the maps
\[i(g)=\ad((1,\dots,m))^i(g,1,\dots,1)=\left(g^{\binom i0},
  g^{-\binom i1},\dots,g^{(-1)^{m-1}\binom i{m-1}}\right);\]
concretely, for $m=2$ one has
\[0(g)=(g,1),\quad1(g)=(g,g)\]
and for $m=3$ one has
\[0(g)=(g,1,1),\quad1(g)=(g,g^{-1},1),\quad
2(g)=(g,g^{-2},g)\equiv(g,g,g)\mod\mho_3(G).\]
When $m$ is prime, one clearly has $i(g)=0$ for all $i\ge m$, and if
$g\in K$ then $i(g)\in K$ for all $i\in\{0,1,\dots,d-1\}$.

\subsection{The group $\Gg$}\label{subs:lcsGg}
We give an explicit description of the Lie algebra of $\Gg$, and
compute its Hilbert-Poincar\'e series. These results were obtained
in~\cite{bartholdi-g:lie}.

Set $x=(ab)^2$. Then $\Gg$ is branch over $K=\langle x\rangle^\Gg$, and
$K/(K\times K)$ is cyclic of order $4$, generated by $x$.

\newcommand\sm[2]{{\{\begin{smallmatrix}#1\\#2\end{smallmatrix}\}}}

Extend the generating set of $\Gg$ to a formal set
$S=\big\{a,b,c,d,\sm bc,\sm cd,\sm db\big\}$, whose meaning shall be
made clear later.  Define the transformation $\sigma$ on words in
$S^*$ by
\[\sigma(a)=a\sm bca,\quad\sigma(b)=d,\quad\sigma(c)=b,\quad\sigma(d)=c,\]
extended to subsets by $\sigma\sm xy=\sm{\sigma x}{\sigma y}$. Note
that for any fixed $g\in G$, all elements $h\in\Stab_\Gg(1)$ such that
$\psi(h)=(g,*)$ are obtained by picking a letter from each set in
$\sigma(g)$. This motivates the definition of $S$.

\begin{theorem}\label{theorem:G:struct}
  Consider the following Lie graph: its vertices are the symbols
  $X(x)$ and $X(x^2)$, for words $X\in\{0,1\}^*$. Their degrees are
  given by
  \begin{align*}
    \deg X_1\dots X_n(x) = 1 + \sum_{i=1}^n X_i2^{i-1} + 2^n,\\
    \deg X_1\dots X_n(x^2) = 1 + \sum_{i=1}^n X_i2^{i-1} + 2^{n+1}.\\
  \end{align*}
  There are four additional vertices: $a,b,d$ of degree $1$, and
  $[a,d]$ of degree $2$.

  Define the arrows as follows: an arrow labelled $\sm xy$ stands for
  two arrows, labelled $x$ and $y$, and the arrows labelled $c$ are
  there to expose the symmetry of the graph (indeed $c=bd$ is not in
  our chosen basis of $G/[G,G]$):
  \[\begin{diagram}
    \node a\arrow{e,t}{b,c}\node x \node a\arrow{e,t}{c,d}\node{[a,d]}\\
    \node b\arrow{e,t}a\node x \node d\arrow{e,t}a\node{[a,d]}\\
    \node x\arrow{e,t}{a,b,c}\node{x^2}\node x\arrow{e,t}{c,d}\node{0(x)}\\
    \node{[a,d]}\arrow{e,t}{b,c}\node{0(x)} \node{0*}\arrow{e,t}a\node{1*} \\
    \node{1^n(x)}\arrow{e,t}{\sigma^n\sm cd}\node{0^{n+1}(x)}
    \node{1^n(x)}\arrow{e,t}{\sigma^n\sm bd}\node{0^n(x^2)}\\
    \node{1^n0*}\arrow{e,t}{\sigma^n\sm
    cd}\node{0^n1*\makebox[0mm][l]{ if $n\ge1$.}}
  \end{diagram}\]

  Then the resulting graph is the Lie graph of $\Lie(\Gg)$. A slight
  modification gives the Lie graph of $\Lie_{\F_2}(\Gg)$: the degree
  of $X_1\dots X_n(x^2)$ is $2\deg X_1\dots X_n(x)$; and the
  square maps are given by
  \begin{align*}
    X(x)&\overset{\cdot2}\longrightarrow X(x^2),\\
    1^n(x^2)&\overset{\cdot2}\longrightarrow 1^{n+1}(x^2).\\
  \end{align*}

  The subgraph spanned by $a,t$, the $X_1\dots X_i(x)$ for $i\le n-2$
  and the $X_1\dots X_i(x^2)$ for $i\le n-4$ is the Lie graph
  associated to the finite quotient $\Gg/\Stab_\Gg(n)$.
\end{theorem}

\begin{figure}
  \begingroup\makeatletter
  \renewcommand\dggeometry{\dg@XGRID=\thr@@
    \dg@YGRID=\tw@
    \unitlength=0.012pt\relax}
  \def\dgt@neeee{\dg@DX=4 \dg@DY=\@ne \dg@SIZE=4}
  \tiny
  \iffast\[\fbox{Lie graph of $\Lie(\Gg)$}\]\else
  \[\begin{diagram}
    \node{b}\arrow{se,t}{a} \node[2]{x^2}\\
    \node[2]{x}\arrow{ne,t}{a,b,c}\arrow{se,t}{c,d}
    \node[3]{0(x^2)}\arrow{e,t}{a} \node{1(x^2)} \node[3]{00(x^2)}
    \arrow{e,t}{a} \node{10(x^2)} \arrow{e,t}{b,c}
    \node{01(x^2)\makebox[0mm][l]{ $\cdots$}}\\
    \node{a}\arrow{ne,t}{b,c}\arrow{se,b}{c,d}
    \node[2]{0(x)}\arrow{e,b}{a}
    \node{1(x)}\arrow{se,b}{b,c}\arrow{ne,t}{c,d}
    \node[3]{01(x)}\arrow{e,b}{a}
    \node{11(x)}\arrow{se,b}{b,d}\arrow{ne,t}{b,c}\\
    \node[2]{[a,d]}\arrow{ne,b}{b,c} \node[3]{00(x)}\arrow{e,b}{a}
    \node{10(x)}\arrow{ne,b}{b,c} \node[3]{000(x)} \arrow{e,b}{a}
    \node{100(x)} \arrow{e,b}{b,c} \node{010(x)\makebox[0mm][l]{ $\cdots$}}\\
    \node{d}\arrow{ne,b}{a}\\
    \node 1 \arrow{n,-,..} \arrow{s,-,..} \node 2 \node 3 \node 4
    \arrow[3]{n,-,..} \arrow[2]{s,-,..} \node 5 \node 6 \node 7 \node 8
    \arrow[3]{n,-,..} \arrow[2]{s,-,..} \node 9 \node{10} \node{11}\\
    \node{b}\arrow{se,t}{a}\\
    \node[2]{x}\arrow{se,t}{c,d} \arrow[2]{e,t}{\cdot2} \node[2]{x^2}
    \qbezier(122,-45)(184,-25)(246,-45)\put(180,-30){$\cdot2$}
    \node[2]{0(x^2)} \node[2]{1(x^2)} \node[2]{00(x^2)}\\
    \node{a}\arrow{ne,t}{b,c}\arrow{se,b}{c,d}
    \node[2]{0(x)}\arrow{e,b}{a} \arrow{neee,t}{\cdot2}
    \node{1(x)}\arrow{e,b}{b,c} \arrow{neeee,t}{\cdot2}
    \node{00(x)}\arrow{e,b}{a} 
    \put(164,-69){\lamsarrow(170,25)}\put(240,-62){$\cdot2$}
    \node{10(x)}\arrow{e,b}{b,c} \node{01(x)}\arrow{e,b}{a}
    \node{11(x)}\arrow{e,b}{b,d}
    \node{000(x)} \arrow{e,b}{a} \node{100(x)} \arrow{e,b}{b,c}
    \node{010(x)\makebox[0mm][l]{ $\cdots$}}\\
    \node[2]{[a,d]}\arrow{ne,b}{b,c}\\
    \node{d}\arrow{ne,b}{a}
  \end{diagram}\]
  \fi
  \endgroup
  \caption{The beginning of the Lie graphs of $\Lie(\Gg)$ (top) and
    $\Lie_{\F_2}(\Gg)$ (below).}
\end{figure}

\begin{proof}
  The proof proceeds by induction on length of words, or, what amounts
  to the same, on depth in the lower central series.

  First, the assertion is checked ``manually'' up to degree $3$. The
  details of the computations are the same as
  in~\cite{bartholdi-g:lie}.

  We claim that for all words $X,Y$ with $\deg Y(x)>\deg X(x)$ we have
  $Y(x)\in\langle X(x)\rangle^\Gg$, and similarly $Y(x^2)\in\langle
  X(x^2)\rangle^\Gg$. The claim is verified by induction on $\deg X$.

  We then claim that for any non-empty word $X$, either
  $\ad(a)X(*)=0$ (if $X$ starts by ``$1$'') or
  $\ad(v)X(*)=0$ for $v\in\{b,c,d\}$ (if $X$ starts by ``$0$'').
  Again this holds by induction.

  We then prove that the arrows are as described above. For instance,
  for the last one,
  \begin{align*}
    \ad(\sigma^n\sm cd)1^n0* &= \begin{cases}
      \big(\ad(\sigma^n\sm db)1^{n-1}0*,\ad(\sm a1)1^{n-1}0*\big)\\
      \hfill= 0\ad(\sigma^{n-1}\sm cd)1^{n-1}0* = 0^n1* & \text{ if }n\ge2,\\
      (\ad(\sm bc)0*,\ad(a)0*) = 01* & \text{ if }n=1.
    \end{cases}
  \end{align*}

  Finally we check that the degrees of all basis elements are as
  claimed. Fix a word $X(*)$, and consider the largest $n$ such that
  $X(*)\in\gamma_n(\Gg)$. Thus there is a sequence of $n-1$ arrows
  leading from the left of the Lie graph of $\Lie(\Gg)$ to $X(*)$, and
  no longer sequence, so $\deg X(*)=n$.

  The modification giving the Lie graph of $\Lie_{\F_2}(\Gg)$ is
  justified by the fact that in $\Lie(\Gg)$ we always have $\deg
  X(x^2)\le2\deg X(x)$, so the element $X(x^2)$ appears always last as
  the image of $X(x)$ through the square map. The degrees are modified
  accordingly. Now $X(x^2)=X1(x^2)$, and $2\deg X1(x)\ge 4\deg X(x)$,
  with equality only when $X=1^n$. This gives an additional square map
  from $1^n(x^2)$ to $1^{n+1}(x^2)$, and requires no adjustment of
  the degrees.
\end{proof}

\begin{corollary}\label{cor:g:rk}
  Define the polynomials
  \begin{align*}
    Q_2&=-1-\hbar,\\
    Q_3&=\hbar+\hbar^2+\hbar^3,\\
    Q_n(\hbar)&=(1+\hbar)Q_{n-1}(\hbar^2)+\hbar+\hbar^2\text{ for }n\ge4.
  \end{align*}
  Then $Q_n$ is a polynomial of degree $2^{n-1}-1$, and the first
  $2^{n-3}-1$ coefficients of $Q_n$ and $Q_{n+1}$ coincide.  The
  term-wise limit $Q_\infty=\lim_{n\to\infty}Q_n$ therefore exists.

  The Hilbert-Poincar\'e series of $\Lie(\Gg/\Stab_\Gg(n))$ is
  $3\hbar+\hbar^2+\hbar Q_n$, and the Hilbert-Poincar\'e series of $\Lie(\Gg)$ is
  $3\hbar+\hbar^2+\hbar Q_\infty$.

  The Hilbert-Poincar\'e series of $\Lie_{\F_2}(\Gg)$ is
  $3+\frac{2\hbar+\hbar^2}{1-\hbar^2}$.

  As a consequence, $\Gg/\Stab_\Gg(n)$ is nilpotent of class
  $2^{n-1}$, and $\Gg$ has finite width.
\end{corollary}

\subsection{The group $\FGg$}
We give here an explicit description of the Lie algebra of $\FGg$,
and compute its Hilbert-Poincar\'e series.

\begin{theorem}[\cite{bartholdi:lcs}]\label{theorem:delta:struct}
  In $\FGg$ write $c=[a,t]$ and $u=[a,c]\equiv2(at)$. For words
  $X=X_1\dots X_n$ with $X_i\in\{0,1,2\}$ define symbols
  $\overline{X_1\dots X_n}(c)$ (representing elements of
  $\FGg$) by
  \begin{align*}
    \overline{i0}(c) &= i0(c) / i(u),\\
    \overline{i2^{m+1}1^n}(c) &
    =i\big(\overline{2^{m+1}1^n}(c)\cdot01^m0^n(u)^{(-1)^n}\big),\\
    \text{and }\overline{iX}(c) &= i\overline X(c)\text{ for all
      other }X.\\
  \end{align*}

  Consider the following Lie graph: its vertices are the symbols
  $\overline{X}(c)$ and $X(u)$. Their degrees are given by
  \begin{align*}
    \deg\overline{X_1\dots X_n}(c)&=1+\sum_{i=1}^nX_i3^{i-1}+\frac12(3^n+1),\\
    \deg X_1\dots X_n(u)&=1+\sum_{i=1}^nX_i3^{i-1}+(3^n+1).
  \end{align*}
  There are two additional vertices, labelled $a$ and $t$, of degree
  $1$.

  Define the arrows as follows, for all $n\ge1$:
  \[\begin{diagram}
    \node a\tarrow{e,t}{-t}\node c \node t\aarrow{e,t}a\node c\\
    \node c\tarrow{e,t}{-t}\node{0(c)} \node c\tarrow{e,t}a\node u\\
    \node u\tarrow{e,t}{-t}\node{1(c)}
    \node{\overline{2^n}(c)}\arrow{e,t}{-t}\node{\overline{0^{n+1}}(c)}\\
   \node{0*}\aarrow{e,t}a\node{1*} \node{1*}\aarrow{e,t}a\node{2*}\\
    \node{2^n0*}\tarrow{e,t}t\node{0^n1*}
    \node{2^n1*}\tarrow{e,t}t\node{0^n2*}\\
    \node[2]{\makebox[5mm][r]{$\overline{X_1\dots X_n}(c)$}}
    \tarrow{e,t}{-(-1)^{\sum X_i}t}
    \node{\makebox[5mm][l]{$(X_1-1)\dots(X_n-1)(u)$}}
  \end{diagram}\]
  Then the resulting graph is the Lie graph of $\Lie(\FGg)$.

  The subgraph spanned by $a,t$, the $\overline{X_1\dots X_i}(c)$ for
  $i\le n-2$ and the $X_1\dots X_i(u)$ for $i\le n-3$ is the Lie graph
  associated to the finite quotient $\FGg/\Stab_\FGg(n)$.
\end{theorem}

\begin{figure}
  \begingroup\makeatletter
  \global\renewcommand\dggeometry{\dg@XGRID=\tw@
    \dg@YGRID=\thr@@
    \unitlength=0.008pt\relax}
  \tiny
  \iffast\[\fbox{Lie graph of $\Lie(\FGg)$}\]\else
  \[\begin{diagram}
    \node 1 \node 2 \node 3 \node 4 \node 5 \arrow[2]{s,-,..}
    \node 6 \node 7 \node 8 \node 9 \node{10} \arrow{s,-,..}
    \node{11} \node{12} \node{13} \node{14} \node{15} \arrow[2]{s,-,..}
    \node{16} \node{17} \node{18} \node{19} \node{20} \arrow{s,-,..}
    \node{21}\\
    \node[6]{\overline{00}(c)} \aarrow{se} \node[2]{\overline{20}(c)}
    \tarrow{se} \node[2]{\overline{11}(c)} \aarrow{se} \tarrow{ssse}
    \node[2]{\overline{02}(c)} \aarrow{se} \node[2]{\overline{22}(c)}
    \tarrow{se} \tarrow{ssse} \node[2]{\overline{100}(c)} \aarrow{se}
    \node[2]{\overline{010}(c)} \aarrow{se} \node[2]{\overline{210}(c)}
    \tarrow{se}\\
    \node{a} \tarrow{se} \node[2]{\overline{0}(c)} \aarrow{se}
    \node[2]{\overline{2}(c)} \tarrow{se} \tarrow{ne}
    \node[2]{\overline{10}(c)} \aarrow{ne} \node[2]{\overline{01}(c)}
    \aarrow{ne} \node[2]{\overline{21}(c)} \tarrow{se} \tarrow{ne}
    \node[2]{\overline{12}(c)} \aarrow{ne} \tarrow{se}
    \node[2]{\overline{000}(c)} \aarrow{ne}
    \node[2]{\overline{200}(c)} \tarrow{ne}
    \node[2]{\overline{110}(c)} \aarrow{ne} \node[2]{}\\
    \node[2]{c} \tarrow{ne} \aarrow{se} \node[2]{\overline{1}(c)}
    \aarrow{ne} \tarrow{se} \node[2]{1(u)} \aarrow{se} \node[6]{10(u)}
    \aarrow{se} \node[2]{01(u)} \aarrow{se} \node[2]{21(u)} \tarrow{se}
    \node[2]{12(u)} \aarrow{se}\\
    \node{t} \aarrow{ne} \node[2]{u} \tarrow{ne} \node[2]{0(u)}
    \aarrow{ne} \node[2]{2(u)} \node[4]{00(u)} \aarrow{ne}
    \node[2]{20(u)} \tarrow{ne} \node[2]{11(u)} \aarrow{ne}
    \node[2]{02(u)} \aarrow{ne} \node[2]{22(u)}\\
    \node 2 \node 1 \node 2 \node 1 \arrow[2]{n,-,..} \node 2 \node 2 \node
    2 \node 1 \node 1 \node 1 \arrow[4]{n,-,..} \node 2 \node 2 \node 2
    \node 2 \node 2 \node 2 \node 2 \node 2 \node 2 \node 1 \node 1
  \end{diagram}\]
  \fi
  \endgroup
  \caption{The beginning of the Lie graph of $\Lie(\FGg)$. The
    generator $\ad(t)$ is shown by plain/blue arrows, and the generator
    $\ad(a)$ is shown by dotted/red arrows.}
\end{figure}

\begin{corollary}\label{cor:delta:rk}
  Define the integers $\alpha_n=\frac12(5\cdot3^{n-2}+1)$, and the
  polynomials
  \begin{align*}
    Q_2&=1,\\
    Q_3&=1+2\hbar+\hbar^2+\hbar^3+\hbar^4+\hbar^5+\hbar^6,\\
    Q_n(\hbar)&=(1+\hbar+\hbar^2)Q_{n-1}(\hbar^3)+\hbar+\hbar^{\alpha_n-2}\text{ for }n\ge4.
  \end{align*}
  Then $Q_n$ is a polynomial of degree $\alpha_n-2$, and the first
  $3^{n-2}+1$ coefficients of $Q_n$ and $Q_{n+1}$ coincide.  The
  term-wise limit $Q_\infty=\lim_{n\to\infty}Q_n$ therefore exists.

  The Hilbert-Poincar\'e series of $\Lie(\FGg/\Stab_\FGg(n))$ is
  $2\hbar+\hbar^2Q_n$, and the Hilbert-Poincar\'e series of $\Lie(\FGg)$ is
  $2\hbar+\hbar^2Q_\infty$.

  As a consequence, $\FGg/\Stab_\FGg(n)$ is nilpotent of class
  $\alpha_n$, and $\FGg$ has finite width.
\end{corollary}

In quite the same way as for $\GSg$, we may improve the general
result $\FGg^{(k)}\le\gamma_{2^k}(\FGg)$:
\begin{theorem}\label{theorem:lcsdsFGg}
  The derived series of $\FGg$ satisfies $\FGg'=\gamma_2(\FGg)$
  and $\FGg^{(k)}=\gamma_5(\FGg)^{\times3^{k-2}}$ for $k\ge2$. We
  have for all $k\in\N$
  \[\FGg^{(k)}\le\gamma_{2+3^{k-1}}(\FGg).\]
\end{theorem}

\begin{theorem}
  Keep the notations of Theorem~\ref{theorem:delta:struct}. Define now
  furthermore symbols $\overline{X_1\dots X_n}(u)$ (representing
  elements of $\FGg$) by
  \begin{align*}
    \overline{2^n}(u)&=2^n(u)\cdot 2^{n-1}0(c)\cdot 2^{n-2}01(c)\cdots
    201^{n-2}(c),\\
    \text{and }\overline X(u) &= X(u)\text{ for all other }X.\\
  \end{align*}

  Consider the following Lie graph: its vertices are the symbols
  $\overline{X}(c)$ and $\overline X(u)$. Their degrees are given by
  \begin{align*}
    \deg\overline{X_1\dots X_n}(c)&=1+\sum_{i=1}^nX_i3^{i-1}+\frac12(3^n+1),\\
    \deg 2^n(u)&=3^{n+1},\\
    \deg X_1\dots X_n(u)&=\max\{1+\sum_{i=1}^nX_i3^{i-1}+(3^n+1),
    \frac12(9-3^n)+3\sum_{i=1}^nX_i3^{i-1}\}.
  \end{align*}
  There are two additional vertices, labelled $a$ and $t$, of degree
  $1$.

  Define the arrows as follows, for all $n\ge1$:
  \[\begin{diagram}
    \node a\tarrow{e,t}{-t}\node c \node t\aarrow{e,t}a\node c\\
    \node c\tarrow{e,t}{-t}\node{0(c)} \node c\tarrow{e,t}a\node u\\
    \node u\tarrow{e,t}{-t}\node{1(c)}
    \node{\overline{2^n}(c)}\arrow{e,t}{-t}\node{\overline{0^{n+1}}(c)}\\
    \node{0*}\aarrow{e,t}a\node{1*} \node{1*}\aarrow{e,t}a\node{2*}\\
    \node{2^n0*}\tarrow{e,t}t\node{0^n1*}
    \node{2^n1*}\tarrow{e,t}t\node{0^n2*}\\
    \node[2]{\makebox[5mm][r]{$\overline{X_1\dots X_n}(c)$}}
    \tarrow{e,t}{-(-1)^{\sum X_i}t}
    \node{\makebox[5mm][l]{$(X_1-1)\dots(X_n-1)(u)$}}\\
    \node c \arrow{e,t}{\cdot3} \node{\overline{00}(c)}
    \node{\overline{2^n}(u)} \arrow{e,t}{\cdot3}\node{\overline{2^{n+1}}(u)}\\
    \node{*0(c)} \arrow{e,t}{\cdot3} \node{*2(u)\makebox[0mm][l]{ if $3\deg*0(c)=\deg*2(u)$}}
  \end{diagram}\]
  Then the resulting graph is the Lie graph of $\Lie_{\F_3}(\FGg)$.

  The subgraph spanned by $a,t$, the $\overline{X_1\dots X_i}(c)$ for
  $i\le n-2$ and the $X_1\dots X_i(u)$ for $i\le n-3$ is the Lie graph
  of the Lie algebra $\Lie_{\F_3}(\FGg/\Stab_\FGg(n))$.

  As a consequence, the dimension series of $\FGg/\Stab_\FGg(n)$
  has length $3^{n-1}$ (the degree of $\overline{2^n}(u)$), and
  $\FGg$ has finite width.
\end{theorem}

\subsection{The group $\GSg$}\label{subs:lcsGS}
We give here an explicit description of the Lie algebra of $\GSg$,
and compute its Hilbert-Poincar\'e series.

Introduce the following sequence of integers:
\[\alpha_1=1,\quad\alpha_2=2,
\quad\alpha_n=2\alpha_{n-1}+\alpha_{n-2}\text{ for }n\ge3,\]
and $\beta_n=\sum_{i=1}^n\alpha_i$. One has
\begin{align*}
  \alpha_n &= \frac1{2\sqrt2}\left((1+\sqrt2)^n-(1-\sqrt2)^n\right),\\
  \beta_n &= \frac14\left((1+\sqrt2)^{n+1}+(1-\sqrt2)^{n+1}-2\right).
\end{align*}
The first few values are
\[\begin{array}{c|cccccccc}
  n & 1 & 2 & 3 & 4 & 5 & 6 & 7 & 8\\ \hline
  \alpha_n & 1 & 2 & 5 & 12 & 29 & 70 & 169 & 398\\
  \beta_n & 1 & 3 & 8 & 20 & 49 & 119 & 288 & 686
\end{array}\]

\begin{theorem}[\cite{bartholdi:lcs}]\label{theorem:gamma:struct}
  In $\GSg$ write $c=[a,t]$ and $u=[a,c]=2(t)$.  Consider the following
  Lie graph: its vertices are the symbols $X_1\dots X_n(x)$ with
  $X_i\in\{0,1,2\}$ and $x\in\{c,u\}$.  Their degrees are given by
  \begin{align*}
    \deg X_1\dots X_n(c)&=1+\sum_{i=1}^nX_i\alpha_i+\alpha_{n+1},\\
    \deg X_1\dots X_n(u)&=1+\sum_{i=1}^nX_i\alpha_i+2\alpha_{n+1}.
  \end{align*}
  There are two additional vertices, labelled $a$ and $t$, of degree
  $1$.

  Define the arrows as follows:
  \[\begin{diagram}
    \node a\tarrow{e,t}{-t}\node c \node c\tarrow{e,t}t\node{0(c)}\\
    \node t\aarrow{e,t}a\node c \node c\aarrow{e,t}a\node u\\
    \node[2]u\tarrow{e,t}t\node{1(c)}\\
    \node{0*}\aarrow{e,t}a\node{1*} \node{1*}\aarrow{e,t}a\node{2*}\\
    \node[2]{2*}\tarrow{e,t}{t}\node{0\#\makebox[0mm][l]{ whenever
    $*\blue{\overset t\longrightarrow}\#$}}\\
    \node{2(c)}\tarrow{e,t}t\node{1(u)}
    \node{1(c)}\tarrow{e,t}{-t}\node{0(u)}\\
    \node{10*}\tarrow{e,t}{-t}\node{01*}
    \node{11*}\tarrow{e,t}{-t}\node{02*}\\
    \node{20*}\tarrow{e,t}t\node{11*} \node{21*}\tarrow{e,t}t\node{12*}\\
  \end{diagram}\]
  (Note that these last $3$ lines can be replaced by the rules
  ${2*}\blue{\overset t\longrightarrow}1\#$ and
  ${1*}\blue{\overset{-t}\longrightarrow}0\#$ for all arrows
  $*\red{\overset a\longrightarrow}\#$.)

  Then the resulting graph is the Lie graph of $\Lie(\GSg)$. It is
  also the Lie graph of $\Lie_{\F_3}(\GSg)$, with the only
  non-trivial cube maps given by
  \[2^n(c)\overset{\cdot3}\longrightarrow 2^n00(c),\qquad
  2^n(c)\overset{\cdot3}\longrightarrow2^n1(u).\]

  The subgraph spanned by $a,t$, the $X_1\dots X_i(c)$ for $i\le n-2$
  and the $X_1\dots X_i(u)$ for $i\le n-3$ is the Lie graph associated
  to the finite quotient $\GSg/\Stab_\GSg(n)$.
\end{theorem}

\begin{figure}  \begingroup\makeatletter
  \renewcommand\dggeometry{\dg@XGRID=\@ne
    \dg@YGRID=\tw@
    \unitlength=0.016pt\relax}
  \tiny
  \iffast\[\fbox{Lie graph of $\Lie(\GSg)$}\]\else
  \[\begin{diagram}
    \node 1 \node 2 \node 3 \node 4 \node 5 \arrow[5]{s,-,..}
    \node 6 \node 7 \node 8 \node 9 \node{10} \arrow[4]{s,-,..}
    \node{11} \node{12} \node{13} \node{14} \node{15} \arrow{s,-,..}
    \node{16} \node{17} \node{18} \node{19} \node{20} \arrow[2]{s,-,..}
    \node{21}\\
    \node[15]{200(c)} \tarrow{se} \node[2]{020(c)} \aarrow{se}
    \node[2]{220(c)} \tarrow{se}\\
    \node[14]{100(c)} \tarrow{se} \aarrow{ne} \node[2]{110(c)}
    \aarrow{se}
    \tarrow{ne} \node[2]{120(c)} \aarrow{ne} \node[2]{011(c)}
    \aarrow{se} \node[2]{}\\
    \node[13]{000(c)} \aarrow{ne} \node[2]{010(c)} \aarrow{ne}
    \node[2]{210(c)} \tarrow{se} \tarrow{ne} \node[2]{101(c)}
    \tarrow{ne} \aarrow{se} \node[2]{111(c)} \aarrow{ne} \tarrow{se}\\
    \node[6]{00(c)} \aarrow{se} \node[2]{20(c)} \tarrow{se}
    \node[2]{02(c)} \aarrow{se} \node[2]{22(c)} \tarrow{ne} \tarrow{se}
    \node[6]{001(c)} \aarrow{ne} \node[2]{201(c)} \tarrow{ne} \node[2]{}\\
    \node{a} \tarrow{se} \node[2]{0(c)} \aarrow{se} \node[2]{2(c)}
    \tarrow{ne} \tarrow{se} \node[2]{10(c)} \aarrow{ne} \tarrow{se}
    \node[2]{11(c)} \tarrow{ne} \aarrow{se} \node[2]{12(c)} \aarrow{ne}
    \node[2]{01(u)} \aarrow{se} \node[2]{21(u)} \tarrow{se}
    \node[2]{22(u)} \tarrow{ne}\\
    \node[2]{c} \tarrow{ne} \aarrow{se} \node[2]{1(c)} \aarrow{ne}
    \tarrow{se} \node[2]{1(u)} \aarrow{se} \node[2]{01(c)} \aarrow{ne}
    \node[2]{21(c)} \tarrow{ne} \tarrow{se} \node[2]{10(u)} \tarrow{ne}
    \aarrow{se} \node[2]{11(u)} \tarrow{se} \aarrow{ne} \node[2]{12(u)}
    \aarrow{ne}\\
    \node{t} \aarrow{ne} \node[2]{u} \tarrow{ne} \node[2]{0(u)}
    \aarrow{ne} \node[2]{2(u)} \tarrow{ne} \node[4]{00(u)} \aarrow{ne}
    \node[2]{20(u)}
    \tarrow{ne} \node[2]{02(u)} \aarrow{ne}\\
    \node 2 \node 1 \node 2 \node 1 \arrow[2]{n,-,..} \node 2 \node 2 \node
    2 \node 2 \node 1 \arrow[3]{n,-,..} \node 2 \node 2 \node 2 \node 3
    \node 2 \node 4 \node 2 \node 3 \node 2 \node 2 \node 2 \node 1
    \arrow[5]{n,-,..}
  \end{diagram}\]
  \fi
  \endgroup  \caption{The beginning of the Lie graph of
  $\Lie(\GSg)$. The
  generator $\ad(t)$ is shown by plain/blue arrows, and the
  generator $\ad(a)$ is shown by dotted/red arrows.}
\end{figure}

\begin{corollary}\label{cor:gamma:rk}
  Define the following polynomials:
  \begin{align*}
    Q_1&=0,\\
    Q_2&=\hbar+\hbar^2,\\
    Q_3&=\hbar+\hbar^2+2\hbar^3+\hbar^4+\hbar^5,\\
    Q_n&=(1+\hbar^{\alpha_n-\alpha_{n-1}})Q_{n-1}
    +\hbar^{\alpha_{n-1}}(\hbar^{-\alpha_{n-2}}+1+\hbar^{\alpha_{n-2}})Q_{n-2}
    \text{ for }n\ge3.
  \end{align*}
  Then $Q_n$ is a polynomial of degree $\alpha_n$, and the polynomials
  $Q_n$ and $Q_{n+1}$ coincide on their first $2\alpha_{n-1}$ terms.
  The coefficient-wise limit $Q_\infty=\lim_{n\to\infty}Q_n$ therefore
  exists.

  The largest coefficient in $Q_{2n+1}$ is $2^n$, at position
  $\frac12(\alpha_{2n+1}+1)$, so the coefficients of $Q_\infty$ are
  unbounded. The integers $k$ such that $\hbar^k$ has coefficient $1$ in
  $Q_\infty$ are precisely the $\beta_n+1$.

  The Hilbert-Poincar\'e series of $\Lie(\GSg/\Stab_\GSg(n))$ is
  $\hbar+Q_n$, and the Hilbert-Poincar\'e series of $\Lie(\GSg)$ is
  $\hbar+Q_\infty$. The same holds for the Lie algebra
  $\Lie_{\F_3}(\GSg/\Stab_\GSg(n))$ and $\Lie_{\F_3}(\GSg)$.

  As a consequence, $\GSg/\Stab_\GSg(n)$ is nilpotent of class
  $\alpha_n$, and $\GSg$ does not have finite width.
\end{corollary}

We note as an immediate consequence that
\[[\GSg:\gamma_{\beta_n+1}(\GSg)]=3^{\frac12(3^n+1)},\]
so that the asymptotic growth of
$l_n=\dim(\gamma_n(\GSg)/\gamma_{n+1}(\GSg))$ is polynomial of
degree $d=\log3/\log(1+\sqrt2)-1$, meaning that $d$ is minimal such
that
\[\limsup_{n\to\infty}\frac{\sum_{i=1}^nl_i}{\sum_{i=1}^n i^d}<\infty.\]
We then have by Proposition~\ref{prop:grlgr} the
\begin{corollary}\label{cor:lcsGSlowerbd}
  The growth of $\GSg$ is at least
  $e^{n^{\frac{\log3}{\log(1+\sqrt2)+\log3}}}\cong e^{n^{0.554}}$.
\end{corollary}

We may also improve the general result
$\GSg^{(k)}\le\gamma_{2^k}(\GSg)$ to the following
\begin{theorem}\label{theorem:lcsdsGSg}
   For all $k\in\N$ we have
  \[\GSg^{(k)}\le\gamma_{\alpha_{k+1}}(\GSg).\]
\end{theorem}

\section{Subgroup growth}
For a finitely generated group $G$, its \emdef[subgroup!growth
function]{subgroup growth function}\index{growth!subgroup} is
$a_n(G)=|\setsuch{H \leq G}{[G:H]=n}|$, and its \emdef{normal
subgroup growth} is the function $b_n(G)=|\setsuch{N\triangleleft
G}{[G:N]=n}|$. Building on their earlier
papers~\cite{segal:part1,segal:part2,grunewald-s-s:subgroups,
mann-s:uniform,lubotzky-m:subgp} \ALubotzky, \AMann\ and \DSegal\
obtained in~\cite{lubotzky-m-s:subgp} a characterization of
finitely generated groups with polynomial subgroup growth. Namely,
the finitely generated groups of polynomial subgroup growth are
precisely the virtually solvable groups of finite rank. For a well
written survey refer to~\cite{lubotzky:sg}.

Since $a_n$ and $b_n$ only count finite-index subgroups, it is
especially interesting to estimate the subgroup growth of
just-infinite groups, and branch groups appear naturally in this
context. A recent result lends support to that view:

\begin{theorem}[\cite{segal:fimages}]\label{thm:segalsbgp}
  Let $f:\mathbb N\to \R_{\geq 0}$ be non-decreasing and such that
  $\log(f(n))/\log(n)$ is unbounded as $n\to\infty$. Then there
  exists a $4$-generator branch (spinal) group $G$ whose subgroup growth
  is not polynomial, but satisfies
  \[a_n(G)\precsim n^{f(n)}.\]
\end{theorem}

In order to prove the above theorem \DSegal\ uses the construction
described in Subsection~\ref{subs:segal} with $A_i=\PSL(2,p_i)$,
$i\in\N$, as rooted subgroups and the action is the natural doubly
transitive action of $\PSL(2,p_i)$ on a set of $p_i+1$ elements.
The subgroup growth function can be made slow by a choice of a
sequence of primes $(p_i)_{i\in\N}$ that grows quickly enough. In
addition, \DSegal\ shows that there exists a continuous range of
possible ``slow'' subgroup growths.

\date{October 27, 2002}
\chapter{Representation Theory of Branch Groups}\label{chapter:rep}
This chapter deals with the representation theory of branch groups
and of their finite quotients over vector spaces of finite
dimension and over Hilbert spaces. For infinite discrete groups
the theory of infinite-dimensional unitary representations is a
quite a difficult subject. Fortunately, for branch groups the
situation is easier to handle, and we produce several results, all
taken from~\cite{bartholdi-g:spectrum}, confirming this.

The study of unitary representations of profinite branch groups is of
great importance and is motivated from several directions. For
just-infinite branch groups it is equivalent to the study of irreducible
representations of finite quotients $G/\Stab_G(\LL_n)$. The first step in
this direction is to consider the quasi-regular representations
$\rho_{G/P_n}$, where $P_n$ is the stabilizer of a point of level $n$, and
this will be done below, again following~\cite{bartholdi-g:spectrum}.
Another direction is the study of representations of all $3$ types of
groups (discrete branch $p$-groups, branch pro-$p$-groups, and their
finite quotients) in vector spaces over the finite field $\F_p$. This last
study was initiated by \DPassman\ and \WTemple\ in~\cite{passman-t:reps}.

The spectral properties of the quasi-regular representations
$\rho_{G/P}$ will be described in Chapter~\ref{chapter:spectrum}.

We introduce the following notion:
\begin{definition}
  Let $G$ be a group, and $\Bbbk$ an algebraically closed field. We
  define $F_G(n)\in\N\cup\{\infty\}$ as the number of irreducible
  representations of degree at most $n$ of $G$ over
  $\Bbbk$. Similarly, $f_G(n)$ denotes the number of such
  representations of degree exactly $n$.
\end{definition}
Therefore, $f_G(n)$ is the growth function of the representation ring
of $G$ over $\Bbbk$, whose degree-$n$ component is generated by $\Bbbk
G$-modules of dimension $n$, and whose addition and multiplication are
$\oplus$ and $\otimes$.

First, we remark that if $G$ is finitely generated, $\Bbbk$ is
algebraically closed, and $G$ does not have any $\charact(\Bbbk$)-torsion,
then $F_G(n)$ is finite for all $n$. This follows from a theorem of Weil
(see~\cite{farkas:affinerings}). We assume these conditions are satisfied
by $G$.

If $H$ is a finite-index subgroup of $G$, we have $F_G\sim F_H$, as shown
in~\cite{passman-t:reps}.

The first lower bound on $F(n)$ appears in a paper by \DPassman\ and
\WTemple~\cite{passman-t:reps}, where it was stated for the Gupta-Sidki
$p$-groups. We improve slightly the result:
\begin{theorem}[\cite{passman-t:reps}]
  Let $G$ be a finitely-generated $m$-regular branch group over $K$,
  and consider its representations over any field $\Bbbk$. Then
  \[F_G\succsim n^{(m-1)\log_{[\psi(K):K^m]}[K:K']-1}.\]
\end{theorem}
The proof given in~\cite{passman-t:reps} extends easily to all regular
branch groups.  We note that this result is obtained by considering all
possible inductions of degree-$1$ representations from $K$ up to $G$; it
may well be that the function $F_G$ grows significantly faster than
claimed, and the whole representation theory of $(K)\psi/K^m$ should be
taken into account.

As a consequence, we obtain:
\begin{itemize}
\item $F_\Gg\succsim n^2$, since $[K:K']=64$ and $[(K)\psi:K^2]=4$;
\item $F_\FGg\succsim n^3$, since $[\FGg':\FGg'']=81$ and
  $[(\FGg')\psi:(\FGg')^2]=9$;
\item $F_\GSg\succsim n^3$, for the same reason;
\item for the Gupta-Sidki $p$-groups $G_p$ of
  Subsection~\ref{subs:gspg}, the general result $F_{G_p}\succsim n^{p-2}$.
\end{itemize}

Note that since these groups are just-infinite, non-faithful
representations in vector spaces must factor through a finite quotient;
and since these groups are of intermediate growth, they cannot be linear
(by Tits' alternative), so in fact all finite-dimensional representations
factor through a finite quotient of $G$, which may even be taken to be of
the form $G_n=G/\Stab_G(\LL_n)$ if $G$ has the congruence subgroup
property.

For concrete cases, like $\Gg$ and $\FGg$, we may exhibit some unitary
irreducible representations as follows:

\subsection{Quasi-regular representations}\label{subs:quasireg}
The representations we consider here are associated to parabolic
subgroups, i.e., stabilizers of an infinite ray in the tree (see
Section~\ref{sec:parabolic}).

For $G$ a group acting on a tree and $P$ a parabolic subgroup, we let
$\rho_{G/P}$ denote the quasi-regular representation of $G$, acting by
right-multiplication on the space $\ell^2(G/P)$. This representation
is infinite-dimensional, and a criterion for irreducibility, due to
George Mackey, follows:
\begin{definition}
  The \emdef{commensurator} of a subgroup $H$ of $G$ is
  \[\comm_G(H) = \{g\in G|\,H\cap H^g\text{ is of finite index in
    }H\text{ and }H^g\}.\]
  Equivalently, letting $H$ act on the left on the right cosets $\{gH\}$,
  \[\comm_G(H) = \{g\in G|\,H\cdot(gH)\text{ and }H\cdot(g^{-1}H)\text{ are finite orbits}\}.\]
\end{definition}

\begin{theorem}[Mackey~\cite{mackey:representations,burger-h:irreducible}]\label{thm:mackey}
  Let $G$ be an infinite group and let $P$ be any subgroup of $G$.
  Then the quasi-regular representation $\rho_{G/P}$ is irreducible if
  and only if $\comm_G(P)=P$.
\end{theorem}

The following results appear in~\cite{bartholdi-g:parabolic}:
\begin{theorem}\label{thm:commP}
  If $G$ is weakly branch, then $\comm_G(P)=P$,
  and therefore $\rho_{G/P}$ is irreducible.
\end{theorem}

Note that the quasi-regular representations we consider are good
approximants of the regular representation, in the sense that $\rho_G$
is a subrepresentation of $\bigotimes_{P\text{ parabolic
    }}\rho_{G/P}$.  We have a continuum of parabolic subgroups
$P_e=\Stab_G(e)$, where $e$ runs through the boundary of a tree, so
we also have a continuum of quasi-regular representations.
If $G$ is countable, there are uncountably many non-equivalent
representations, because among the uncountably many $P_e$ only
countably many are conjugate. As a consequence,
\begin{theorem}
  There are uncountably many non-equivalent representations of the
  form $\rho_{G/P}$, where $P$ is a parabolic subgroup.
\end{theorem}

We now consider the finite-dimensional representations $\rho_{G/P_n}$,
where $P_n$ is the stabilizer of the vertex at level $n$ in the ray
defining $P$. These are permutational representations on the sets
$G/P_n$ of cardinality $m_1\dots m_n$. The $\rho_{G/P_n}$ are factors
of the representation $\rho_{G/P}$.  Noting that
$P=\bigcap_{n\ge0}P_n$, it follows that
\[\rho_{G/P_n}\Rightarrow\rho_{G/P},\]
in the sense that for any non-trivial $g\in G$ there is an $n\in\N$
with $\rho_{G/P_n}(g)\neq1$.

We describe now the decomposition of the finite quasi-regular
representations $\rho_{G/P_n}$. It turns out that it is closely
related to the orbit structure of $P_n$ on $G/P_n$. We state the
result for the examples $\Gg,\Sg,\FGg,\BGg,\GSg$:
\begin{theorem}[\cite{bartholdi-g:parabolic}]\label{thm:finitedec}
  $\rho_{\Gg/P_n}$ and $\rho_{\Sg/P_n}$ decompose as a direct sum of
  $n+1$ irreducible components, one of degree $2^i$ for each
  $i\in\{1,\dots,n-1\}$ and two of degree $1$.

  $\rho_{\FGg/P_n}$, $\rho_{\BGg/P_n}$ and $\rho_{\GSg/P_n}$ decompose
  as a direct sum of $2n+1$ irreducible components, two of degree
  $2^i$ for each $i\in\{1,\dots,n-1\}$ and three of degree $1$.
\end{theorem}

\part{Geometric and Analytic Aspects}
\date{October 27, 2002}
\chapter{Growth} \label{chapter:growth}
The notion of growth in finitely generated groups was introduced
by Efremovich in~\cite{efremovich:proximity} and Shvarts
in~\cite{svarts:growth} in their study of Riemannian manifolds.
The works of \JMilnor\ in the late sixties
(\cite{milnor:curvature,milnor:solvable}) contributed to the
current reinforced interest in the topic. Before we make brief
historical remarks on the research made in connection to word
growth in finitely generated groups, let us introduce the
necessary definitions. We concentrate solely on finitely generated
infinite groups.

Let $S=\{s_1,\dots,s_k\}$ be a non-empty set of symbols. A \emdef{weight
function} on $S$ is any function $\tau:S\to \R_{>0}$. Therefore, each
symbol in $S$ is assigned a positive weight. The \emdef{weight} of any
word over $S$ is then defined by the extension of $\tau$ to a
homomorphism, still written $\tau:S^*\to\R_{\ge0}$, defined on the free
monoid $S^*$ of words over $S$. Therefore, for any word over $S$ we have
\[\tau(s_{i_1}s_{i_2}\dots s_{i_\ell}) = \sum_{j=1}^\ell \tau(s_{i_j}).\]
Note that the empty word is the only word of weight 0. For any
non-negative real number $n$ there are only finitely many words in $S^*$
of weight at most $n$.

Let $G$ be an infinite group and $\rho:S^*\twoheadrightarrow G$ a
surjective monoid homomorphism. Therefore, $G$ is finitely
generated and $\rho(S)=\{\rho(s_1),\dots,\rho(s_k)\}$ generates
$G$ as a monoid. The \emdef{weight} of an element $g$ in $G$ with
respect to the triple $(S,\tau,\rho)$ is, by definition, the
smallest weight of a word $u$ in $S^*$ that represents $g$, i.e.,
the smallest weight of a word in $\rho^{-1}(g)$. The weight of $g$
with respect to $(S,\tau,\rho)$ is denoted by
$\partial_G^{(S,\tau,\rho)}(g)$.

For $n$ a non-negative real number, the elements in $G$ that have weight
at most $n$ with respect to $(S,\tau,\rho)$ constitute the \emdef{ball} of
radius $n$ in $G$ with respect to $(S,\tau,\rho)$, denoted by
$B_G^{(S,\tau,\rho)}(n)$.

Let $G$ act transitively on a set $X$ on the right, and let $x$ be
an element of $X$. We define
$B_{x,G}^{(S,\tau,\rho)}=\setsuch{x^g}{g\in
  B_G^{(S,\tau\rho)}}$ to be the \emdef{ball} in $X$ of radius $n$ with center
  at $x$.

The number of elements in $B_{x,G}^{(S,\tau,\rho)}(n)$ is finite and we
denote it by $\gamma_{x,G}^{(S,\tau,\rho)}(n)$. The function
$\gamma_{x,G}^{(S,\tau,\rho)}$, defined on the non-negative real numbers,
is called the \emdef{growth function} of $X$ at $x$ as a $G$-set with
respect to $(S,\tau,\rho)$.

The equivalence class of $\gamma_{x,G}^{(S,\tau,\rho)}$ under $\sim$
(recall the notation from the introduction) is called the \emdef{degree of
growth} of $G$ and it does not depend on the (finite) set $S$, the weight
function $\tau$ defined on $S$, the homomorphism $\rho$, nor the choice of
$x$.

\begin{proposition}[Invariance of the growth function]
If $\gamma_{x,G}^{(S,\tau,\rho)}$ and $\gamma_{x',G}^{(S',\tau',\rho')}$
are two growth functions of the group $G$, then they are equivalent with
respect to $\sim$.
\end{proposition}

When we define a weight function on a group $G$ we usually pick a finite
generating subset of $G$ closed for inversion and not containing the
identity, assign a weight function to those generating elements and extend
the weight function to the whole group $G$ in a natural way, thus blurring
the distinction between a word over the generating set and the element in
$G$ represented by that word, and completely avoiding the discussion of
$\rho$.

A standard way to assign a weight function is to assign the weight
1 to each generator. In that case we use the standard notation and
terminology, i.e., we denote the weight of a word $u$ by $|u|$ and
call it the \emdef{length} of $u$. In this setting, the length of
the group element $g$ is the distance from $g$ to the identity in
the Cayley graph of the group with respect to the generating set
$S$.

If we let $G$ act on itself by right multiplication we see that the growth
function of $G$ just counts the number of elements in the corresponding
ball in $G$, i.e $\gamma_{1,G}^{(S,\tau,\rho)}(n) =
|B_G^{(S,\tau,\rho)}(n)|$. Since the degree of growth is an invariant of
the group we are more interested in it than in the actual growth function
for a given generating set.

In the next few examples the weight 1 is assigned to the
generating elements. If $G=\Z$ and $S=\{1,-1\}$ then $\gamma_G(n)
= 2n+1 \sim n$. More generally, if $G$ is the free abelian group
of rank $k$, we have $\gamma_G \sim n^k$. If $G$ is the free group
of rank $k \geq 2$ with the standard generating set together with
the inverses, then $\gamma_G(n)= \frac{k(2k-1)^n-1}{k-1} \sim
e^n$. It is clear that the growth functions of groups on $k$
generators are bounded above by the growth function of the free
group on $k$ generators. Therefore, the exponential degree of
growth is the largest possible degree of growth.

For any finitely generated infinite group $G$, the following trichotomy
exists: $G$ is of
\begin{itemize}
\item \emdef{polynomial growth} if $\gamma_G(n)\precsim n^d$, for some
  $d\in\N$;
\item \emdef{intermediate growth} if
  $n^d \precnsim \gamma_G(n) \precnsim e^n$, for all $d\in\N$;
\item \emdef{exponential growth} if $e^n\sim\gamma_G(n)$.
\end{itemize}
We say $G$ is of \emdef{subexponential growth} if
$\gamma_G(n)\precnsim e^n$ and of \emdef{superpolynomial growth}
if $n^d\precnsim\gamma_G(n)$ for all $d\in\N$.

By the results of \JMilnor, Joseph Wolf and Brian Hartley
(see~\cite{milnor:solvable} and~\cite{wolf:solvable}), solvable
groups have exponential growth unless they are virtually nilpotent
in which case the growth is polynomial. There is a formula giving
the degree of polynomial growth in terms of the lower central
series of $G$, due to Yves Guivarc'h and Hyman
Bass~\cite{guivarch:poly1,bass:nilpotent}. Namely, if $d_j$
represents the torsion-free rank of the $j-th$ factor in the lower
central series of the finitely generated nilpotent group $G$, then
the degree of growth of $G$ is polynomial of degree $\sum jd_j$.

By the Tits' alternative~\cite{tits:linear}, a finitely generated
linear group is either virtually solvable or it contains the free
group of rank two. Therefore, the growth of a linear group must be
exponential or polynomial. In the same spirit, Ching Chou showed
in~\cite{chou:elementary} that the growth of elementary amenable
groups must be either polynomial or exponential (recall that the
class of elementary amenable groups is the smallest class
containing all finite and all abelian groups that is closed under
subgroups, homomorphic images, extensions and directed unions).
Any non-elementary hyperbolic group must have exponential growth
(\cite{ghys-h:gromov}).

A fundamental result of \MGromov~\cite{gromov:nilpotent} states that
every group of polynomial growth is virtually nilpotent. The results
of \JMilnor, \JWolf, \YGuivarch, \HBass\ and \MGromov\ together imply
that a group has a polynomial growth if and only if it is virtually
nilpotent, and in this case the growth function is equivalent to $n^d$
where $d$ is the integer $\sum jd_j$, as stated above.

We recall the definition of amenable group:

\begin{definition}\label{def:amen}
  Let $G$ be a group acting on a set $X$. This action is
  \emdef{amenable} in the sense of von
  Neumann~\cite{vneumann:masses} if there exists a finitely additive
  measure $\mu$ on $X$, invariant under the action of $G$, with
  $\mu(X)=1$.

  A group $G$ is \emdef{amenable} if its action on itself by
  right-multiplication is amenable.
\end{definition}

The following criterion, due to F{\o}lner, can sometimes be used
to show that a certain action is amenable.

\begin{theorem}\label{theorem:folner}
Let $G$ act transitively on a set $X$ and let $S$ be a generating
set for $G$. The action is amenable if and only if, for every
positive real number $\lambda$ there exists a finite set $F$ such
that $|F \Delta Fs|< \lambda|F|$, for all $s\in S$, where $\Delta$
denotes the symmetric difference.
\end{theorem}

Using the F{\o}lner's criterion one can show that groups of
subexponential growth are always amenable.

In~\cite{milnor:5603} John Milnor asked the following question:
''Is the function $\gamma(n)$ necessarily equivalent either to a
power of $n$ or to the exponential function $e^n$?''. In other
words, Milnor asked if the growth is always polynomial or
exponential.

The first examples of groups of intermediate growth were constructed by
the second author in~\cite{grigorchuk:burnside}, and they are known as the
first and the second Grigorchuk group. Both examples are $2$-groups of
intermediate growth and they are both amenable but not elementary. These
examples show that the answer to Milnor's question is ``no''. In the same
time, these examples show that there exist amenable but not elementary
amenable groups, thus answering a question of \MDay\ from~\cite{day:amen}.

Other examples of groups of intermediate growth are $\FGg$, see
Section~\ref{sec:FGg} and the first torsion-free example
from~\cite{grigorchuk:pgps}, which is based on the first
Grigorchuk group $\Gg$.

It was shown in~\cite{grigorchuk:gdegree,grigorchuk:pgps} that the
Grigorchuk $p$-groups are of intermediate growth and that there are
uncountably large chains and antichains of growth functions associated
with these examples, which are all branch groups.

An important unanswered question in the theory of groups of intermediate
growth is the existence of a group whose degree of growth is
$e^{\sqrt{n}}$. We remarked before that residually-$p$ groups that have
degree of growth strictly below $e^{\sqrt{n}}$ must be virtually
nilpotent, and therefore of polynomial growth (see~\ref{cor:lcslowerbd}).
The same conclusion holds for residually nilpotent groups
(see~\cite{lubotzky-m:subgp}).

The historical remarks on the research on growth above are biased towards
the existence and development of examples of groups of intermediate
growth. There is plenty of great research on growth in group theory that
is not concerned with this aspect. An excellent review of significant
results and a good bibliography can be found in~\cite{harpe:ggt}
and~\cite{grigorchuk-h:groups}.

Before we move on to more specific examples, we make an easy observation.
\begin{proposition}\label{theorem:branch-NP}
No weakly branch group has polynomial growth, i.e. no weakly branch group
is virtually nilpotent.
\end{proposition}

\section{Growth of \GG\ groups with finite directed part}
We will concentrate on the case of \GG\ groups with finite
directed part $B$. Note that all branching indices, except maybe
the first one, are bounded above by $|B|$ since each homomorphic
image $A_{\sigma^r\om}= (B)\om_r$ acts transitively on a set of
$m_{r+1}$ elements. Therefore, $|B| \geq |A_{\sigma^r\om}| \geq
m_{r+1}$, for all $r$. Let us denote, once and for all, the
largest branching index by $M$ and the smallest one by $m$.

All the estimates of word growth that we give in this and the following
sections are done with respect to the canonical generating set $S_\om =
(A_\om \cup B_\om)-1$. As a shorthand, we use $\gamma_\om(n)$ instead of
$\gamma_{G_\om}(n)$.

In order to express some of the results we use the notions of complete
subsequence and $r$-homogeneous and $r$-factorable sequence (see the
remarks before Theorem~\ref{theorem:(r)pgrowth}). We recall that all
sequences in $\Omhat$, i.e. sequences that define \GG\ groups (see
Definition~\ref{defn:GG}), can be factored into finite complete
subsequences.

\begin{theorem} \label{theorem:hatgrowth}
  All \GG\ groups with finite directed part have subexponential growth.
\end{theorem}

The following lemma is a direct generalization
of~\cite[Lemma~1]{grigorchuk:pgps}. The proof is similar, but adapted to
our more general setting.

\begin{lemma} [3/4-Shortening] \label{lemma:3/4}
  Let $\omover$ be a sequence that starts with a complete sequence of
  length $r$. Then the following inequality holds for every reduced
  word $F$ representing an element in $\Stab_\om(\LL_r)$:
  \[|L_r(F)| \leq  \frac{3}{4}|F| + M^r,\]
  where $|L_r(F)|$ represents the total length of the words on the
  level $r$ of the $r$-level decomposition of $F$.
\end{lemma}
\begin{proof}
Define $\xi_i$ to be the number of $B$-letters from $K_i \setminus
(K_{i-1}\cup\dots\cup K_1)$ appearing in the words at the level $i-1$, and
$\nu_i$ to be the number of simple reductions performed to get the words
$F_{j_1\dots j_i}$ on the level $i$ from their unreduced versions
$\overline{F_{j_1\dots j_i}}$.

  A reduced word $F$ of length $n$ has at most $(n+1)/2$ $B$-letters.
  Every $B$-letter in $F$ that is in $K_1$ contributes one $B$-letter
  and no $A$-letters to the unreduced words
  $\overline{F_1},\dots,\overline{F_{m_1}}$. The $B$-letters in $F$
  that are not in $K_1$, and there are at most $(n+1)/2-\xi_1$ such
  letters, contribute one $B$-letter and one $A$-letter.  Finally, the
  $\nu_1$ simple reductions reduce the number of letters on level 1 by
  at least $\nu_1$.  Therefore,
  \[|L_1(F)| \leq 2((n+1)/2-\xi_1) + \xi_1 - \nu_1 = n+1-\xi_1-\nu_1.\]

  In the same manner, each of the $\xi_2$ $B$-letters on level $1$
  that is in $K_2 \setminus K_1$ contributes one $B$-letter to the unreduced
  words on level $2$. The other $B$-letters, and there are at most
  $(|L_1(F)|+M)/2-\xi_2$ of them, contribute at most $2$ letters, so,
  after simple reductions, we have
  \[|L_2(F)| \leq n+1 + M - \xi_1 - \xi_2 - \nu_1 - \nu_2.\]

  Proceeding in the same manner, we obtain the estimate
  \begin{equation} \label{eq:L}
    |L_r(F)| \leq n+1+M+\dots +M^{r-1}- \xi_1-\xi_2-\dots-\xi_r-\nu_1-\nu_2-\dots-\nu_r.
  \end{equation}
  If $\nu_1+\nu_2+\dots+\nu_r\geq n/4$, the claim of the lemma
  immediately follows.

  Let us therefore consider the case when
  \begin{equation} \label{eq:nusmall}
    \nu_1+\nu_2+\dots+\nu_r < n/4.
  \end{equation}

  For $i=0,\dots,r-1$, define $|L_i(F)|^+$ to be the number of
  $B$-letters from $B\setminus(K_1\cup\dots\cup K_i)$ appearing in the words
  at the level $i$.  Clearly, $|L_0(F)|^+$ is the number of
  $B$-letters in $F$ and
  \[|L_0(F)|^+ \geq \frac{n-1}{2}.\]
  Going from the level $0$ to the level $1$, each $B$-letter
  contributes one $B$ letter of the same type. Therefore, the words
  $\overline{F_1},\dots,\overline{F_{m_1}}$ from the first level
  before the reduction takes place have exactly $|L_0(F)|^+ -\xi_1$
  letters that come from $B-K_1$. Since we lose at most $2\nu_1$
  letters due to the simple reductions, we obtain
  \[|L_1(F)|^+ \geq \frac{n-1}{2} -\xi_1 - 2\nu_1.\]

  Next, we go from level $1$ to level $2$. There are $|L_1(F)|^+$
  $B$-letters on level $1$ that come from $B\setminus K_1$, so there are
  exactly $|L_1(F)|^+ -\xi_2$ $B$-letters from $B\setminus (K_1 \cup K_2)$ in
  the words $\overline{F_{11}},\dots,\overline{F_{m_1m_2}}$, and then
  we lose at most $2\nu_2$ $B$-letters due to the simple reductions.
  We get
  \[|L_2(F)|^+ \geq \frac{n-1}{2} -\xi_1 - \xi_2 - 2\nu_1 - 2\nu_2,\]
  and, by proceeding in a similar manner,
  \begin{equation} \label{eq:L+}
    |L_{r-1}(F)|^+ \geq \frac{n-1}{2} -\xi_1-\dots-\xi_{r-1} -2\nu_1-\dots -2\nu_{r-1}.
  \end{equation}
  Since $\om_1\dots\om_r$ is complete, we have $K_r \setminus (K_1 \cup \dots
  \cup K_{r-1}) = B \setminus (K_1 \cup \dots \cup K_{r-1})$ and
  $\xi_r=|L_{r-1}(F)|^+$ so that the
  inequalities~(\ref{eq:L}),~(\ref{eq:nusmall}) and~(\ref{eq:L+}) yield
  \[|L_r(F)| \leq \frac{n}{2}+\frac{1}{2}+1+M+\dots + M^{r-1} +
  \nu_1 + \dots+\nu_{r-1} -\nu_r,\]
  which implies the claim.
\end{proof}

Now we can finish the proof of Theorem~\ref{theorem:hatgrowth}
using the approach used by the second author
in~\cite{grigorchuk:pgps} (see
also~\cite[Theorem~VIII.61]{harpe:ggt}). We use the following easy
lemma, which follows from the fact every subgroup of index $L$ has
a transversal whose representatives have length at most $L-1$ (use
a \Schreier\ transversal).

\begin{lemma}
  Let $G$ be a group and $H$ be a subgroup of finite index $L$ in $G$.
  Let $\gamma(n)$ denote the growth function of $G$ with respect to
  some finite generating set $S$ and the standard length function on
  $S$, and let $\beta(n)$ denote the number of words of length at most
  $n$ that are in $H$, i.e.
  \[\beta(n) = |\setsuch{g}{g\in H,\;|g|\leq n}| = |B_G^S(n) \cap H|.\]
  Then
  \[\gamma(n) \leq L\beta(n+L-1).\]
\end{lemma}

Let
\[e_\om=\lim_{n\to\infty}\sqrt[n]{\gamma_\om(n)}\]
denote the \emdef{exponential growth rate} of $G_\om$. It is known that
this rate is 1 if and only if the group in question has  subexponential
growth. Therefore, all we need to show is that this rate is 1.

\begin{proof}[Proof of Theorem~\ref{theorem:hatgrowth}]

  For any $\epsilon>0$ there exists $n_0$ such that
  \[\gamma_\om(n) \leq (e_\om + \epsilon)^n,\]
  for all $n \geq n_0$. If we denote $C_\om=\gamma_\om(n_0)$, we
  obtain
  \[\gamma_\om(n) \leq C_\om(e_\om + \epsilon)^n,\]
  for all $n$. Note that $C_\om$ depends on $\epsilon$.

  Let $\om$ start with a complete sequence of length $r$. Denote
  \[\beta_\om(n) = |\setsuch{g}{g\in \Stab_\om(\LL_r),\;|g|\leq n}|.\]
  By Lemma~\ref{lemma:3/4} and the fact that $\psi_r$ is an embedding we
  have
  \[\beta_\om(n) \leq\sum \gamma_{\sigma^r\om}(n_1)\gamma_{\sigma^r\om}(n_2)\dots\gamma_{\sigma^r\om}(n_s),\]
  where $s=m_1 \cdot m_2\dots \cdot m_r$, and the summation is over
  all tuples $(n_1,\dots,n_s)$ of non-negative integers with
  $n_1+n_2+\dots+n_s \leq \frac{3}{4}n+M^r$. Let $L$ be the index of
  $\Stab_\om(\LL_r)$ in $G_\om$.  By the previous lemma and the above
  discussions we have
  \[\gamma_\om(n) \leq L\beta_\om(n+L-1) \leq LC_{\sigma^r\om}^s \sum (e_{\sigma^r\om}+\epsilon)^{n_1+\dots+n_s}.\]
  where the summation is over all tuples $(n_1,\dots,n_s)$ of
  non-negative integers with $n_1+\dots+n_s \leq
  \frac{3}{4}(n+L-1)+M^r$. The number of such tuples is polynomial
  $P(n)$ in $n$ (depending also on the constants $L$, $M$ and $r$).
  Therefore,
  \[\gamma_\om(n) \leq LC_{\sigma^r\om}^s P(n)(e_{\sigma^r\om}+\epsilon)^{\frac{3}{4}(n+L-1)+M^r}.\]
  Taking the $n$-th root on both sides and the limit as $n$ tends to
  infinity gives
  \[e_\om \leq (e_{\sigma^r\om}+\epsilon)^{\frac{3}{4}},\]
  and since this inequality holds for all positive $\epsilon$, we
  obtain
  \[e_\om \leq (e_{\sigma^r\om})^{\frac{3}{4}}.\]
  Since the exponential growth rate $e_{\sigma^t\om}$ is bounded above
  by $|A_{\sigma^t\om}|+|B_{\sigma^t\om}|-1 \leq 2|B|-1$, for all
  $t>0$, it follows that $e_\om=1$ for all $\om\in\Omhat$.
\end{proof}

A general lower bound, tending to $e^n$ when $m\to\infty$, exists on the
word growth, and holds for all \GG\ groups:
\begin{theorem}\label{theorem:wdlowerbd}
  All \GG\ groups with finite directed part have superpolynomial
  growth.  Moreover, the growth of $G_\om$ satisfies
  \[e^{n^\alpha}\precsim\gamma_\om(n),\]
  where $\alpha = \frac{\log(m)}{\log(m)-\log\frac12}$.
\end{theorem}

A proof can be found in~\cite{bartholdi-s:wpg} for the more special case
that is considered there. Note that $\alpha > \frac{1}{2}$, as long as $m
\neq 2$. Putting the last several results together gives
\begin{theorem}
  All \GG\ groups with finite directed part have intermediate growth.
\end{theorem}

For the special case of $\Gg$ a slightly better lower bounds exist, due to
\YLeonov~\cite{leonov:pont} who obtained $e^{n^{0.5041}} \preceq
\gamma(n)$, and to the first author~\cite{bartholdi:lowerbd} who obtained
$e^{n^{0.5157}} \preceq \gamma(n)$.

\section{Growth of \GG\ groups defined by homogeneous sequences}\label{sec:homo}
Implicitly, all defining triples $\om$ have $r$-homogeneous defining
sequence $\omover$ for some fixed $r$.  The following estimate holds:

\begin{theorem} [$\eta$-Estimate] \label{theorem:(r)growth}
  If $\omover$ is an $r$-homogeneous sequence, then the growth
  function of the group $G_\om$ satisfies
  \[\gamma_\om(n) \precsim e^{n^{\alpha}}\]
  where $\alpha= \frac{\log(M)}{\log(M)-log(\eta_r)}<1$ and $\eta_r$ is
  the positive root of the polynomial $x^r + x^{r-1} + x^{r-2} -2$.
\end{theorem}

We mentioned already that the weight assignment is irrelevant as far as
the degree of growth is concerned. But, appropriately chosen weight
assignments can make calculations easier, and this is precisely how the
above estimates are obtained.

Let $G$ be a group that is generated as a monoid by the set of generators
$S$ that does not contain the identity. A weight function $\tau$ on $S$ is
called triangular if $\tau(s_1)+\tau(s_2) \geq \tau(s_3)$ whenever
$s_1,s_2,s_3 \in S$ and $s_1s_2=s_3$.

In the case of spinal groups, in order to define a triangular weight we
must have
\[\tau(a_1) + \tau(a_2) \geq \tau(a_1a_2) \qquad \text{and} \qquad
\tau(b_1) + \tau(b_2) \geq \tau(b_1b_2),\] for all $a_1,a_2\in A$ and
$b_1,b_2 \in B$ such that $a_1a_2 \neq 1$ and $b_1b_2 \neq 1$.

Every $g$ in $G_\omega$ admits a minimal form with respect to a triangular
weight $\tau$
\[[a_0]b_1a_1b_2a_2 \dots a_{k-1}b_k[a_k]\]
where all $a_i$ are in $A-1$, all $b_i$ are in $B-1$, and the appearances
of $a_0$ and $a_k$ are optional.

The following weight assignment generalizes the approach taken
in~\cite{bartholdi:upperbd} by the first author in order to estimate the
growth of $\Gg$ (see also~\cite{bartholdi-s:wpg}).

The linear system of equations in the variables $\tau_0,\dots,\tau_r$:
\begin{equation} \label{eq:system}
  \begin{cases}
    \eta_r(\tau_0+\tau_i) = \tau_0+\tau_{i-1} & \text{ for }i=r,\dots,2,\\
    \eta_r(\tau_0+\tau_1) = \tau_r.
  \end{cases}
\end{equation}
has a solution, up to a constant multiple, given by
\begin{equation} \label{eq:solution}
  \begin{cases}
    \tau_i = \eta_r^r + \eta_r^{r-i} -1 & \text{ for }i=r,\dots,1,\\
    \tau_0 = 1-\eta_r^r.
  \end{cases}
\end{equation}
We also require $\tau_1+\tau_2=\tau_r$ and we get that $\eta_r$ must be a
root of the polynomial $x^r+x^{r-1}+x^{r-2}-2$. We choose $\eta_r$ to be
the root of this polynomial that is between $0$ and $1$ obtain that the
solution~(\ref{eq:solution}) of the system~(\ref{eq:system}) satisfies the
additional properties
\begin{align}
  0 < \tau_1 < \dots < \tau_r < 1, \qquad 0 < \tau_0 < 1, \\
  \tau_i+\tau_j \geq \tau_k \text{ for all } 1 \leq i,j,k \leq r\text
  { with }i\neq j.
  \label{eq:triangular}
\end{align}
The index $r$ in $\eta_r$ will be omitted without warning.

Now, given $\om \in \Omr$, we define the weight of the generating
elements in $S_\om$ as follows: $\tau(a)=\tau_0$, for $a$ in $A-1$
and $\tau(b_\om)=\tau_i$, where $i$ is the smallest index with
$(b)\om_{i}=1$, i.e., the smallest index with $b\in
K_i=\Ker(\om_i)$.

Clearly, $\tau$ is a triangular weight function. The only point worth
mentioning is that if $b$ and $c$ are two $B$-letters of the same weight
and $bc=d \not = 1$ then $d$ has no greater weight than $b$ or $c$ (this
holds because $b_\om,c_\om \in K_i$ implies $d_\om \in K_i$).

For obvious reasons, the weight $\partial_{G_\om}^{(S_\om,\tau,\rho)}(g)$,
for $g \in G_\om$, is denoted by $\partial^{\tau}(g)$ and, more often,
just by $\partial(g)$.

We define the \emdef{portrait of an element of $G_\om$ of size $n$} with
respect to the weight $\tau$ as the portrait with respect to the sequence
of profile sets $B_{G_\om}^\tau(n)$, $B_{G_{\sigma^\om}}^\tau(n)$,
$\dots$, which is the sequence of balls of radius $n$ in the corresponding
companion groups. Therefore, in the process of building a portrait only
those vertices that would be decorated by an element that has weight
larger than $n$ are decomposed further at least one more level and become
interior vertices decorated by vertex permutations.

\begin{figure}[!ht]
  \begin{center}
    \includegraphics{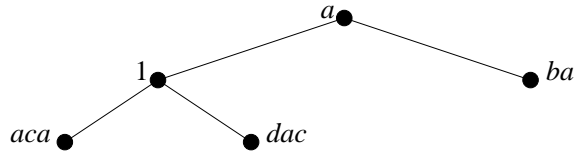}
  \end{center}
  \caption{Portrait of $g$ of size 3 (same as the portrait of size 4)}\label{figure:size3}
\end{figure}

\begin{figure}[!ht]
  \begin{center}
    \includegraphics{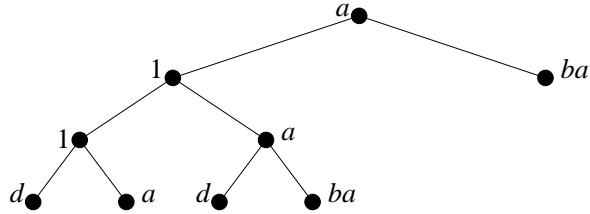}
  \end{center}
  \caption{Portrait of $g$ of size 2}\label{figure:size2}
\end{figure}

Just as an easy example, we give the portraits of size 3 and size 2 with
respect to the standard word length of the element $g=abacadacabadac$ in
$\Gg$ in Figure~\ref{figure:size3} and Figure~\ref{figure:size2}.
Other portraits of the same element are given in
Subsection~\ref{subs:portraits}.

The following lemma says that the total sum of the weights of the sections
of an element $f$ is significantly shorter (by a factor less than 1) than
the weight of $f$. This observation leads to upper bounds on the word
growth in \GG\ groups.

\begin{lemma} [$\eta$-Shortening] \label{lemma:<eta}
  Let $f \in G_\om$ and $\omover$ be an $r$-homogeneous sequence. Then
  \[\sum_{i=1}^{m_1} \partial^{\tau}(f_i) \leq
  \eta_r\big(\partial^{\tau}(f)+ \tau_0\big).\]
\end{lemma}
\begin{proof}
  Let a minimal form of $f$ be
  \[f = [a_0]b_1a_1 \dots b_{k-1}a_{k-1}b_k[a_k].\]
  Further let $f=hg$ where $g \in A_\om$ and $h \in \Stab_\om(\LL_1)$. Then $h$
  can be written in the form $h = [a_0]b_1a_1 \dots
  a_{k-1}b_k[a_k]g^{-1}$ and rewritten in the form
  \begin{equation} \label{eq:h}
    h = b_1^{g_1} \dots b_k^{g_k},
  \end{equation}
  where $g_i=([a_0]a_1 \dots a_{i-1})^{-1} \in A_\om$. Clearly, $\partial(f)
  \geq (k-1)\tau_0 +\sum_{j=1}^k \tau(b_j)$, which yields
  \begin{equation} \label{eq:eta-interm}
    \sum_{i=1}^k \eta(\tau_0 + \tau(b_j)) \leq \eta(\partial(f)+\tau_0).
  \end{equation}
  Now, observe that if the $B$-generator $b$ is of weight $\tau_i$
  with $i>1$ then $(b^g)\psi$ has as components one $B$-generator of
  weight $\tau_{i-1}$ and one $A$-generator (of weight $\tau_0$ of
  course) with the rest of the components trivial. Therefore, such a $b^g$
  (from~(\ref{eq:h})) contributes at most
  $\tau_0+\tau_{i-1}=\eta(\tau_0+\tau(b))$ to the sum $\sum
  \partial(f_i)$. On the other hand, if $b$ is a $B$-generator of
  weight $\tau_1$ then $(b^g)\psi$ has as components one $B$-generator
  of weight at most $\tau_r$, and the rest of the components are
  trivial. Such a $b^g$ contributes at most
  $\tau_r=\eta(\tau_0+\tau(b))$ to the sum $\sum \partial(f_i)$.
  Therefore
  \begin{equation} \label{eq:eta-interm2}
    \sum_{i=1}^{m_1} \partial(f_i) \leq \sum_{j=1}^k \eta(\tau_0+\tau(b_j))
  \end{equation}
  and the claim of the lemma follows by combining~(\ref{eq:eta-interm})
  and~(\ref{eq:eta-interm2}).
\end{proof}

For a chosen $C$ we can now construct the portraits of size $C$ of the
elements of $G_\om$. In case $C$ is large enough, the previous lemma
guarantees that these portraits are finite.

\begin{lemma} \label{lemma:Ln}
  There exists a positive constant $C$ such that
  \[L(n) \precsim n^\alpha,\]
  with $\alpha=\frac{\log{(M)}}{\log(M)-\log(\eta_r)}$, where $L(n)$
  is the maximal possible number of leaves in the portrait of size $C$
  of an element of weight at most $n$.
\end{lemma}

We give a proof that does not work exactly, but gives the right idea. One
can find a complete proof in~\cite{bartholdi-s:wpg}.

\begin{proof}[Sketchy proof]
  We present how the proof would work if
  \begin{equation} \label{eq:falseeta}
    \sum_{i=1}^{m_1} \partial(f_i) \leq \eta \partial(f)
  \end{equation}
  holds rather than the inequality in Lemma~\ref{lemma:<eta}.

  We choose $C$ big enough so that the portraits of size $C$ are
  finite.  Define a function $L'(n)$ on $\R_{\geq 0}$ by
  \[L'(n) = \begin{cases} 1 & \text{ if }n \leq C,\\
    n^\alpha & \text{ if }n > C.
  \end{cases}\]
  We prove, by induction on $n$, that $L(n)\leq L'(n)$. If the weight
  $n$ of $g$ is $\leq C$, the portrait has $1$ leaf and $L'(n)=1$.
  Otherwise, the portrait of $g$ is made up of those of
  $g_1,\dots,g_{m_1}$. Let the weights of these $m_1$ elements be
  $n_1,\dots,n_{m_1}$. By induction, the number of leaves in the
  portrait of $g_i$ is at most $L'(n_i)$, $i=1,\dots,m_1$, and the
  number of leaves in the portrait of $g$ is, therefore, at most
  $\sum_{i=1}^{m_1} L'(n_i)$.

  There are several cases, but let us just consider the case when all
  of the numbers $n_1,\dots,n_{m_1}$ are greater than $C$. Using
  Jensen's inequality, the inequality~\ref{eq:falseeta}, the facts that
  $\eta^\alpha=M^{\alpha-1}$, $0 < \alpha < 1$, and direct
  calculation, we see that
  \[\begin{aligned}
    \sum_{i=1}^{m_1} L'(n_i) &=
        \sum_{i=1}^{m_1} n_i^\alpha \leq
        m_1\left(\frac{1}{m_1}\sum_{i=1}^{m_1} n_i\right)^\alpha \leq \\
    &\leq \frac{1}{(m_1)^{\alpha-1}}(\eta n)^\alpha =
        \left(\frac{M}{m_1}\right)^{\alpha-1}n^\alpha
        \leq n^\alpha = L'(n).
  \end{aligned}\]
\end{proof}

\begin{proof}[Proof of Theorem~\ref{theorem:(r)growth}]
The number of labelled, rooted trees with at most $L(n)$ leaves,
whose branching indices do not exceed $M$, is $\precsim
e^{n^\alpha}$. A tree with $L(n)$ leaves has $\sim L(n)$ interior
vertices, so there are $\precsim e^{n^\alpha}$ ways to decorate
the interior vertices. The decoration of the leaves can also be
chosen in $\precsim e^{n^\alpha}$ ways. Therefore $\gamma_\om(n)
\precsim e^{n^\alpha}$.
\end{proof}

\subsection{Growth in the case of factorable sequences}
An upper bound on the degree of word growth in case of $r$-factorable
sequences can thus be obtained from Theorem~\ref{theorem:(r)growth}, since
every $r$-factorable sequence is $(2r-1)$-homogeneous, but we can do
slightly better if we combine Lemma~\ref{lemma:3/4} with the idea of
portrait of an element. We omit the proof because of its similarity to the
other proofs in this chapter.

\begin{theorem} [3/4-Estimate] \label{theorem:[r]growth3/4}
  If $\omover$ is an $r$-factorable sequence, then the growth function
  of the group $G_\om$ satisfies
  \[\gamma_\om(n) \precsim e^{n^{\alpha}}\]
  where
  $\alpha=\frac{\log(M^r)}{\log(M^r)-log(3/4)}=\frac{\log(M)}{\log(M)-\log\left(\sqrt[\uproot{2}r]{3/4}\right)}<1$.
\end{theorem}

The 3/4-Estimate was obtained only for the class of Grigorchuk $p$-groups
defined by $r$-homogeneous (not $r$-factorable as above) sequences by
Roman Muchnik and Igor Pak in~\cite{muchnik-p:growth} by different means.
The same article contains some sharper considerations in case $p=2$.

We can provide a small improvement in a special case that includes all
Grigorchuk 2-groups. Namely, we are going to assume that $\omover$ is an
$r$-factorable sequence such that each factor contains three homomorphisms
whose kernels cover $B$.

\begin{lemma} [2/3-Shortening] \label{lemma:2/3}
  Let $\om \in \Omhat$ defines a group acting on the rooted binary
  tree. In addition, let the defining sequence $\omover$ be such that
  there exist $3$ terms $\om_k$, $\om_\ell$ and $\om_s$, $1\leq
  k<\ell<m \leq r$, with the property that $K_k \cup K_\ell \cup
  K_s=B$. Then the following inequality holds for every reduced word
  $F$ representing an element in $\Stab_\om(\LL_r)$:
  \[|L_r(F)| < \frac{2}{3}|F| + 3\cdot M^r.\]
\end{lemma}

We note that the above shortening lemma cannot be improved, unless
one starts paying attention to reductions beyond the simple ones.
Indeed, in the first Grigorchuk group $\Gg$ the word
$F=(abadac)^{4k}$ has length $24k$, while
$|L_3(F)|=16k=\frac{2}{3}|F|$. On the other hand
$F=(abadac)^{16}=1$ in $\Gg$, so by taking into account other
relations the multiplicative constant of Lemma~\ref{lemma:2/3}
could possibly be sharpened.

As a corollary to the shortening lemma above, we obtain:
\begin{theorem} [2/3-Estimate] \label{theorem:[r]growth2/3}
  If $\omover$ is an $r$-factorable sequence such that each factor
  contains three letters whose kernels cover $B$, then the growth
  function of the group $G_\om$ satisfies
  \[\gamma_\om(n) \precsim e^{n^{\alpha}}\]
  where
  $\alpha=\frac{\log(M^r)}{\log(M^r)-log(2/3)}=\frac{\log(M)}{\log(M)-\log\left(\sqrt[\uproot{2}r]{2/3}\right)}<1$.
\end{theorem}

\section{Parabolic space and Schreier graphs}\label{sec:pspace}
We describe here the aspects of the parabolic subgroups (introduced in
Section~\ref{sec:parabolic}) related to growth. The $G$-space we study
is defined as follows:
\begin{definition}
  Let $G$ be a branch group, and let $P=\Stab_G(e)$ be a parabolic
  subgroup of $G$. The associated \emdef{parabolic space} is the $G$-set
  $G/P$.
\end{definition}

We consider in this section only finitely-generated, contracting groups.
We assume a branch group $G$, with fixed generating set $S$, has been
chosen.

The proof of the following result appears
in~\cite{bartholdi-g:spectrum}.
\begin{proposition}\label{prop:polygrowth}
  Let $G$, a group of automorphisms of the regular tree of branching
  index $m$, satisfy the conditions of
  Proposition~\ref{prop:confinal}, and suppose it is contracting (see
  Definition~\ref{defn:contract}). Let $P$ be a parabolic subgroup.
  Then $G/P$, as a $G$-set, is of polynomial growth of degree at most
  $\log_{1/\lambda}(m)$, where $\lambda$ is the contracting constant.
  If moreover $G$ is spherically transitive, then $G/P$'s asymptotical
  growth is polynomial of degree $\log_{1/\lambda'}(m)$, with
  $\lambda'$ the infimum of the $\lambda$ as above.
\end{proposition}

The growth of $G/P$ is also connected to the growth of the associated
Lie algebra (see Chapter~\ref{chapter:lie}). The following appears
in~\cite{bartholdi:lcs}:

\begin{theorem}\label{theorem:growth}
  Let $G$ be a branch group, with parabolic space $G/P$. Then there
  exists a constant $C$ such that
  \[\frac{C\gr(G/P)}{1-\hbar}\ge \frac{\gr\Lie(G)}{1-\hbar},\]
  where $\gr(X)$ denotes the growth formal power series of $X$.
\end{theorem}

Anticipating, we note that in the above theorem equality holds for
$\Gg$, but does not hold for $\FGg$ nor $\GSg$ (see
Corollaries~\ref{cor:lcslowerbd} and~\ref{cor:lcsGSlowerbd}).

The most convenient way to describe the parabolic space $G/P$ by
giving it a graph structure:
\begin{definition}\label{defn:schreier}
  Let $G$ be a group generated by a set $S$ and $H$ a subgroup of
  $G$. The \emdef{Schreier graph} $\sch(G,H,S)$ of $G/H$ is the
  directed graph on the vertex set $G/H$, with for every $s\in S$ and
  every $gH\in G/H$ an edge from $gH$ to $sgH$. The base point of
  $\sch(G,H,S)$ is the coset $H$.
\end{definition}
Note that $\sch(G,1,S)$ is the Cayley graph of $G$ relative to
$S$. It may happen that $\sch(G,P,S)$ have loops and multiple
edges even if $S$ is disjoint from $H$. Schreier graphs are
$|S|$-regular graphs, and any degree-regular graph ${\mathcal G}$
containing a $1$-factor (i.e., a regular subgraph of degree $1$;
there is always one if ${\mathcal G}$ has even degree) is a
Schreier graph~\cite[Theorem~5.4]{lubotzky:cayley}.

For all $n\in\N$ consider the finite quotient $G_n=G/\Stab_G(n)$, acting
on the $n$-th level of the tree.  We first consider the finite graphs
${\mathcal G}_n=\sch(G,P_n\Stab_G(n),S)=\sch(G_n,P/(P\cap\Stab_G(n),S)$.
Due to the fractal (or recursive) nature of branch groups, there are
simple local rules producing ${\mathcal G}_{n+1}$ from ${\mathcal G}_n$,
the limit of this process being the Schreier graph of $G/P$.  Before
stating a general result, we start by describing these rules for two of
our examples: $\Gg$ and $\FGg$.

\subsection{$\sch(\Gg,P,S)$}\label{subs:schG}
Assume the notation of Section~\ref{sec:Gg}. The graphs ${\mathcal
G}_n=\sch(\Gg_n,P_n,S)$ will have edges labelled by
$S=\{a,b,c,d\}$ (and not oriented, because all $s\in S$ are
involutions) and vertices labelled by $Y^n$, where $Y=\{1,2\}$.

First, it is clear that ${\mathcal G}_0$ is a graph on one vertex,
labelled by the empty sequence $\emptyset$, and four loops at this
vertex, labelled by $a,b,c,d$. Next, ${\mathcal G}_1$ has two
vertices, labelled by $1$ and $2$; an edge labelled $a$ between
them; and three loops at $1$ and $2$ labelled by $b,c,d$.

Now given ${\mathcal G}_n$, for some $n\ge1$, perform on it the
following transformation: replace the edge-labels $b$ by $d$, $d$
by $c$, $c$ by $b$; replace the vertex-labels $\sigma$ by
$2\sigma$; and replace all edges labelled by $a$ connecting
$\sigma$ and $\tau$ by: edges from $2\sigma$ to $1\sigma$ and from
$2\tau$ to $1\tau$, labelled $a$; two edges from $1\sigma$ to
$1\tau$ labelled $b$ and $c$; and loops at $1\sigma$ and $1\tau$
labelled $d$. We claim the resulting graph is ${\mathcal
G}_{n+1}$.

To prove the claim, it suffices to check that the letters on the
edge-labels act as described on the vertex-labels. If
$b(\sigma)=\tau$, then $d(2\sigma)=2\tau$, and similarly for $c$ and
$d$; this explains the cyclic permutation of the labels $b,c,d$. The
other substitutions are verified similarly.

As an illustration, here are ${\mathcal G}_2$ and ${\mathcal G}_3$
for $\Gg$. Note that the sequences in $Y^*$ that appear correspond
to ``Gray enumeration'', i.e., enumeration of integers in base $2$
where only one bit is changed from a number to the next:
\begin{center}
  \begin{picture}(120,40)
    \put(0,20){\line(1,0){40}}\put(20,30){\msmash{a}}
    \put(40,20){\curve(0,0,20,3,40,0)\curve(0,0,20,-3,40,0)}
    \put(60,30){\msmash{b,c}}
    \put(80,20){\line(1,0){40}}\put(100,30){\msmash{a}}
    \put(0,15){\msmash{21}}
    \put(0,20){\spline(0,0)(18,27)(-18,27)(0,0)}\put(-25,35){\msmash{b,c,d}}
    \put(40,15){\msmash{11}}
    \put(40,20){\spline(0,0)(15,25)(-15,25)(0,0)}\put(40,35){\msmash{d}}
    \put(80,15){\msmash{12}}
    \put(80,20){\spline(0,0)(15,25)(-15,25)(0,0)}\put(80,35){\msmash{d}}
    \put(120,15){\msmash{22}}
    \put(120,20){\spline(0,0)(15,25)(-15,25)(0,0)}\put(142,35){\msmash{b,c,d}}
  \end{picture}\\
  \begin{picture}(280,40)
    \put(0,20){\line(1,0){40}}\put(20,30){\msmash{a}}
    \put(40,20){\curve(0,0,20,3,40,0)\curve(0,0,20,-3,40,0)}\put(60,30){\msmash{b,c}}
    \put(80,20){\line(1,0){40}}\put(100,30){\msmash{a}}
    \put(120,20){\curve(0,0,20,3,40,0)\curve(0,0,20,-3,40,0)}\put(140,30){\msmash{b,d}}
    \put(160,20){\line(1,0){40}}\put(180,30){\msmash{a}}
    \put(200,20){\curve(0,0,20,3,40,0)\curve(0,0,20,-3,40,0)}\put(220,30){\msmash{b,c}}
    \put(240,20){\line(1,0){40}}\put(260,30){\msmash{a}}
    \put(0,10){\msmash{221}}
    \put(0,20){\spline(0,0)(15,25)(-15,25)(0,0)}\put(-25,35){\msmash{b,c,d}}
    \put(40,15){\msmash{121}}
    \put(40,20){\spline(0,0)(15,25)(-15,25)(0,0)}\put(40,35){\msmash{d}}
    \put(80,15){\msmash{111}}
    \put(80,20){\spline(0,0)(15,25)(-15,25)(0,0)}\put(80,35){\msmash{d}}
    \put(120,15){\msmash{211}}
    \put(120,20){\spline(0,0)(15,25)(-15,25)(0,0)}\put(120,35){\msmash{c}}
    \put(160,15){\msmash{212}}
    \put(160,20){\spline(0,0)(15,25)(-15,25)(0,0)}\put(160,35){\msmash{c}}
    \put(200,15){\msmash{112}}
    \put(200,20){\spline(0,0)(15,25)(-15,25)(0,0)}\put(200,35){\msmash{d}}
    \put(240,15){\msmash{122}}
    \put(240,20){\spline(0,0)(15,25)(-15,25)(0,0)}\put(240,35){\msmash{d}}
    \put(280,15){\msmash{222}}
    \put(280,20){\spline(0,0)(15,25)(-15,25)(0,0)}\put(302,35){\msmash{b,c,d}}
  \end{picture}
\end{center}

\subsection{$\sch(\FGg,P,S)$}
Assume the notation of Section~\ref{sec:FGg}. First, ${\mathcal
G}_0$ has one vertex, labelled by the empty sequence $\emptyset$,
and four loops, labelled $a,a^{-1},t,t^{-1}$.  Next, ${\mathcal
G}_1$ has three vertices, labelled $1,2,3$, cyclically connected
by a triangle labelled $a,a^{-1}$, and with two loops at each
vertex, labelled $t,t^{-1}$. In the pictures only geometrical
edges, in pairs $\{a,a^{-1}\}$ and $\{t,t^{-1}\}$, are
represented.

Now given ${\mathcal G}_n$, for some $n\ge1$, perform on it the
following transformation: replace the vertex-labels $\sigma$ by
$2\sigma$; replace all triangles labelled by $a,a^{-1}$ connecting
$\rho,\sigma,\tau$ by: three triangles labelled by $a,a^{-1}$
connecting respectively $1\rho,2\rho,3\rho$ and
$1\sigma,2\sigma,3\sigma$ and $1\tau,2\tau,3\tau$; a triangle
labelled by $t,t^{-1}$ connecting $1\rho,1\sigma,1\tau$; and loops
labelled by $t,t^{-1}$ at $2\rho,2\sigma,2\tau$. We claim the
resulting graph is ${\mathcal G}_{n+1}$.

As above, it suffices to check that the letters on the edge-labels act
as described on the vertex-labels. If $a(\rho)=\sigma$ and
$t(\rho)=\tau$, then $t(1\rho)=1\sigma$, $t(2\rho)=2\rho$ and
$t(3\sigma)=3\tau$. The verification for $a$ edges is even simpler.

\begin{figure}
  \includegraphics{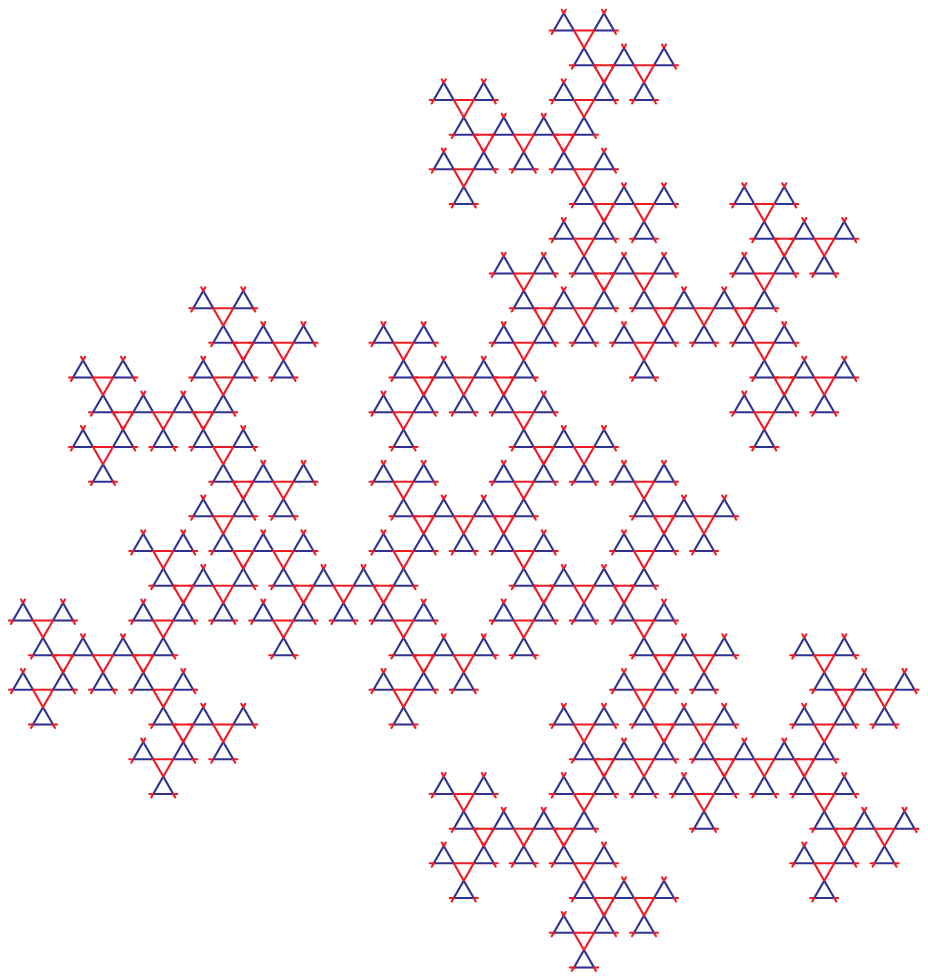}
  {\makeatletter\ifmath@color
    \caption{The Schreier Graph of $\Gamma_6$. The red and blue edges
      represent the generators \red{$s$} and \blue{$a$}.}
    \else
    \caption{The Schreier Graph of $\Gamma_6$. The edges
      represent the generators $s$ and $a$.}
    \fi}
  \label{figure:schreier}
\end{figure}

\subsection{Substitutional graphs}
The two Schreier graphs presented in the previous subsection are
special cases of \emdef{substitutional graphs}, which we define below.

Substitutional graphs were introduced in the late 70's to describe
growth of multicellular organisms. They bear a strong similarity to
$L$-systems~\cite{rozenberg-s:l}, as was noted by Mikhael
Gromov~\cite{gromov:iggo}. Another notion of graph substitution has
been studied by~\cite{previte:subs}, where he had the same convergence
preoccupations as us.

Let us make a convention for this section: all graphs ${\mathcal
G}=(V({\mathcal G}),E({\mathcal G}))$ shall be \emdef{labelled},
i.e., endowed with a map $E({\mathcal G})\to C$ for a fixed set
$C$ of colors, and \emdef{pointed}, i.e., shall have a
distinguished vertex $*\in V({\mathcal G})$. A \emdef{graph
embedding} ${\mathcal G}'\hookrightarrow{\mathcal G}$ is just an
injective map $E({\mathcal G}')\to E({\mathcal G})$ preserving the
adjacency operations.
\begin{definition}
  A \emdef{substitutional rule} is a tuple $(U,R_1,\dots,R_n)$, where
  $U$ is a finite $m$-regular edge-labelled graph, called the
  \emdef{axiom}, and each $R_i$, $i=1,\dots,n$, is a rule of the
  form $X_i\to Y_i$, where $X_i$ and $Y_i$ are finite edge-labelled graphs.
  The graphs $X_i$ are required to have no common edge. Furthermore,
  for each $i=1,\dots,n$, there is an inclusion, written $\iota_i$, of
  the vertices of $X_i$ in the vertices of $Y_i$; the degree of $\iota_i(x)$
  is the same as the degree of $x$ for all $x\in V(X_i)$, and all vertices
  of $Y_i$ not in the image of $\iota_i$ have degree $m$.
\end{definition}

Given a substitutional rule, one sets ${\mathcal G}_0=U$ and constructs
iteratively ${\mathcal G}_{n+1}$ from ${\mathcal G}_n$ by listing all
embeddings of all $X_i$ in ${\mathcal G}_n$ (they are disjoint), and
replacing them by the corresponding $Y_i$. If the base point $*$ of
${\mathcal G}_n$ is in a graph $X_i$, the base point of ${\mathcal
G}_{n+1}$ will be $\iota_i(*)$.

Note that this expansion operation preserves the degree, so
${\mathcal G}_n$ is a $m$-regular finite graph for all $n$. We are
interested in fixed points of this iterative process.

For any $R\in\N$, consider the balls $B_{*,n}(R)$ of radius $R$ at the
base point $*$ in ${\mathcal G}_n$. Since there is only a finite number of
rules, there is an infinite sequence $n_0<n_1<\dots$ such that the balls
$B_{*,n_i}(R)\subseteq{\mathcal G}_{n_i}$ are all isomorphic. We consider
${\mathcal G}$, a limit graph in the sense of~\cite{grigorchuk-z:infinite}
(the limit exists), and call it a \emph{substitutional graph}.

\begin{theorem}[\cite{bartholdi-g:spectrum}]
  The following four substitutional rules describe the Schreier graph of $\Gg$:
  \begin{center}
    \begin{picture}(360,50)
      \put(0,40){\line(1,0){60}}\put(30,50){\msmash a}
      \put(0,45){\footnotesize\msmash\sigma}
      \put(60,45){\footnotesize\msmash\tau}
      \put(30,35){\vector(0,-1){20}}
      \put(80,40){\line(1,0){60}}\put(110,50){\msmash b}
      \put(80,45){\footnotesize\msmash\sigma}
      \put(140,45){\footnotesize\msmash\tau}
      \put(110,35){\vector(0,-1){20}}
      \put(160,40){\line(1,0){60}}\put(190,50){\msmash c}
      \put(160,45){\footnotesize\msmash\sigma}
      \put(220,45){\footnotesize\msmash\tau}
      \put(190,35){\vector(0,-1){20}}
      \put(240,40){\line(1,0){60}}\put(270,50){\msmash d}
      \put(240,45){\footnotesize\msmash\sigma}
      \put(300,45){\footnotesize\msmash\tau}
      \put(270,35){\vector(0,-1){20}}
      \put(0,10){\line(1,0){20}}\put(10,6){\msmash a}
      \put(20,10){\curve(0,0,10,2,20,0)\curve(0,0,10,-2,20,0)}
      \put(30,4){\msmash{b,c}}
      \put(40,10){\line(1,0){20}}\put(50,6){\msmash a}
      \put(20,10){\spline(0,0)(12,20)(-12,20)(0,0)}\put(10,20){\msmash d}
      \put(40,10){\spline(0,0)(12,20)(-12,20)(0,0)}\put(50,20){\msmash d}
      \put(0,7){\footnotesize\msmash{2\sigma}}
      \put(20,7){\footnotesize\msmash{1\sigma}}
      \put(40,7){\footnotesize\msmash{1\tau}}
      \put(60,7){\footnotesize\msmash{2\tau}}
      \put(80,10){\line(1,0){60}}\put(110,6){\msmash d}
      \put(80,7){\footnotesize\msmash{2\sigma}}
      \put(140,7){\footnotesize\msmash{2\tau}}
      \put(160,10){\line(1,0){60}}\put(190,6){\msmash b}
      \put(160,7){\footnotesize\msmash{2\sigma}}
      \put(220,7){\footnotesize\msmash{2\tau}}
      \put(240,10){\line(1,0){60}}\put(270,6){\msmash c}
      \put(240,7){\footnotesize\msmash{2\sigma}}
      \put(300,7){\footnotesize\msmash{2\tau}}
      \put(340,40){\msmash{\text{axiom}}}
      \put(320,10){\spline(0,0)(12,20)(-12,20)(0,0)}\put(340,22){\msmash{b,c,d}}
      \put(320,10){\line(1,0){40}}\put(340,6){\msmash a}
      \put(360,10){\spline(0,0)(12,20)(-12,20)(0,0)}
      \put(320,7){\footnotesize\msmash{1}}
      \put(360,7){\footnotesize\msmash{2}}
    \end{picture}
  \end{center}
  where the vertex inclusions are given by $\sigma\mapsto 2\sigma$ and
  $\tau\mapsto 2\tau$. The base point is the vertex $222\dots$.
\end{theorem}

\begin{theorem}[\cite{bartholdi-g:spectrum}]
  The substitutional rules producing the Schreier graphs of
  $\FGg$, $\BGg$ and $\GSg$ are given below. Note that the Schreier graphs of
  $\BGg$ and $\GSg$ are isomorphic:
  \begin{center}
    \begin{picture}(300,100)(0,-10)
      \put(0,40){\msmash{\FGg:}}
      \put(0,0){\blue{\line(1,0){100}}}\put(50,-4){\blue{\msmash a}}
      \put(0,0){\blue{\line(15,26){50}}}\put(15,43){\blue{\msmash a}}
      \put(100,0){\blue{\line(-15,26){50}}}\put(85,43){\blue{\msmash a}}
      \put(0,-3){\footnotesize\msmash\rho}
      \put(100,-3){\footnotesize\msmash\sigma}
      \put(50,92){\footnotesize\msmash\tau}
      \put(100,40){\vector(1,0){40}}
      \put(130,0){\blue{\line(1,0){40}}}\put(150,7){\blue{\msmash a}}
      \put(130,0){\blue{\line(15,26){20}}}\put(143,17){\blue{\msmash a}}
      \put(170,0){\blue{\line(-15,26){20}}}\put(157,17){\blue{\msmash a}}
      \put(130,-3){\footnotesize\msmash{3\rho}}
      \put(144,37){\footnotesize\msmash{2\rho}}
      \put(168,-3){\footnotesize\msmash{1\rho}}
      \put(190,0){\blue{\line(1,0){40}}}\put(210,7){\blue{\msmash a}}
      \put(190,0){\blue{\line(15,26){20}}}\put(203,17){\blue{\msmash a}}
      \put(230,0){\blue{\line(-15,26){20}}}\put(217,17){\blue{\msmash a}}
      \put(195,-3){\footnotesize\msmash{2\sigma}}
      \put(217,37){\footnotesize\msmash{1\sigma}}
      \put(230,-3){\footnotesize\msmash{3\sigma}}
      \put(160,52){\blue{\line(1,0){40}}}\put(180,59){\blue{\msmash a}}
      \put(160,52){\blue{\line(15,26){20}}}\put(173,69){\blue{\msmash a}}
      \put(200,52){\blue{\line(-15,26){20}}}\put(187,69){\blue{\msmash a}}
      \put(158,57){\footnotesize\msmash{1\tau}}
      \put(180,92){\footnotesize\msmash{3\tau}}
      \put(205,57){\footnotesize\msmash{2\tau}}
      \put(170,0){\red{\line(-10,52){10}}}\put(157,37){\red{\msmash s}}
      \put(170,0){\red{\line(40,35){40}}}\put(185,7){\red{\msmash s}}
      \put(160,52){\red{\line(50,-17){50}}}\put(195,48){\red{\msmash s}}
      \put(150,35){\red{\spline(0,0)(10,17)(-10,17)(0,0)}}
      \put(190,0){\red{\spline(0,0)(-20,0)(-10,-17)(0,0)}}
      \put(200,52){\red{\spline(0,0)(20,0)(10,-17)(0,0)}}
      \put(300,70){\msmash{\text{axiom}}}
      \put(280,0){\blue{\line(1,0){40}}}\put(300,7){\blue{\msmash a}}
      \put(280,0){\blue{\line(15,26){20}}}\put(293,17){\blue{\msmash a}}
      \put(320,0){\blue{\line(-15,26){20}}}\put(307,17){\blue{\msmash a}}
      \put(283,-3){\footnotesize\msmash{3}}\put(272,7){\red{\msmash s}}
      \put(294,37){\footnotesize\msmash{2}}\put(328,7){\red{\msmash s}}
      \put(318,-3){\footnotesize\msmash{1}}\put(310,45){\red{\msmash s}}
      \put(300,35){\red{\spline(0,0)(10,17)(-10,17)(0,0)}}
      \put(280,0){\red{\spline(0,0)(-20,0)(-10,-17)(0,0)}}
      \put(320,0){\red{\spline(0,0)(20,0)(10,-17)(0,0)}}
    \end{picture}\\[2mm]
    \begin{picture}(300,100)(0,-10)
      \put(0,40){\msmash{\BGg,\GSg:}}
      \put(0,0){\blue{\line(1,0){100}}}\put(50,-4){\blue{\msmash a}}
      \put(0,0){\blue{\line(15,26){50}}}\put(15,43){\blue{\msmash a}}
      \put(100,0){\blue{\line(-15,26){50}}}\put(85,43){\blue{\msmash a}}
      \put(0,-3){\footnotesize\msmash\rho}
      \put(100,-3){\footnotesize\msmash\sigma}
      \put(55,87){\footnotesize\msmash\tau}
      \put(100,40){\vector(1,0){40}}
      \put(130,0){\blue{\line(1,0){40}}}\put(150,7){\blue{\msmash a}}
      \put(130,0){\blue{\line(15,26){20}}}\put(143,17){\blue{\msmash a}}
      \put(170,0){\blue{\line(-15,26){20}}}\put(157,17){\blue{\msmash a}}
      \put(130,-3){\footnotesize\msmash{3\rho}}
      \put(144,37){\footnotesize\msmash{2\rho}}
      \put(168,-3){\footnotesize\msmash{1\rho}}
      \put(190,0){\blue{\line(1,0){40}}}\put(210,7){\blue{\msmash a}}
      \put(190,0){\blue{\line(15,26){20}}}\put(203,17){\blue{\msmash a}}
      \put(230,0){\blue{\line(-15,26){20}}}\put(217,17){\blue{\msmash a}}
      \put(195,-3){\footnotesize\msmash{2\sigma}}
      \put(217,37){\footnotesize\msmash{1\sigma}}
      \put(230,-3){\footnotesize\msmash{3\sigma}}
      \put(160,52){\blue{\line(1,0){40}}}\put(180,59){\blue{\msmash a}}
      \put(160,52){\blue{\line(15,26){20}}}\put(173,69){\blue{\msmash a}}
      \put(200,52){\blue{\line(-15,26){20}}}\put(187,69){\blue{\msmash a}}
      \put(158,57){\footnotesize\msmash{1\tau}}
      \put(186,87){\footnotesize\msmash{3\tau}}
      \put(205,57){\footnotesize\msmash{2\tau}}
      \put(170,0){\red{\line(-10,52){10}}}\put(159,35){\red{\msmash t}}
      \put(190,0){\red{\line(10,52){10}}}
      \put(170,0){\red{\line(40,35){40}}}\put(188,10){\red{\msmash t}}
      \put(190,0){\red{\line(-40,35){40}}}
      \put(160,52){\red{\line(50,-17){50}}}\put(194,46){\red{\msmash t}}
      \put(200,52){\red{\line(-50,-17){50}}}
      \put(300,70){\msmash{\text{axiom}}}
      \put(280,0){\blue{\line(1,0){40}}}\put(300,7){\blue{\msmash a}}
      \put(280,0){\blue{\line(15,26){20}}}\put(293,17){\blue{\msmash a}}
      \put(320,0){\blue{\line(-15,26){20}}}\put(307,17){\blue{\msmash a}}
      \put(283,-3){\footnotesize\msmash{3}}\put(272,7){\red{\msmash t}}
      \put(294,37){\footnotesize\msmash{2}}\put(328,7){\red{\msmash t}}
      \put(318,-3){\footnotesize\msmash{1}}\put(310,45){\red{\msmash t}}
      \put(300,35){\red{\spline(0,0)(10,17)(-10,17)(0,0)}}
      \put(280,0){\red{\spline(0,0)(-20,0)(-10,-17)(0,0)}}
      \put(320,0){\red{\spline(0,0)(20,0)(10,-17)(0,0)}}
    \end{picture}
  \end{center}
  where the vertex inclusions are given by $\rho\mapsto 3\rho$,
  $\sigma\mapsto 3\sigma$ and $\tau\mapsto 3\tau$. The base point is the
  vertex $333\dots$.
\end{theorem}

By Proposition~\ref{prop:polygrowth}, these two limit graphs have
asymptotically polynomial growth of degree no higher than $\log_2(3)$.

Note that there are maps $\pi_n:V({\mathcal G}_{n+1})\to
V({\mathcal G}_n)$ that locally (i.e., in each copy of some
right-hand rule $Y_i$) are the inverse of the embedding $\iota_i$.
In case these $\pi_n$ are graph morphisms one can consider the
projective system $\{{\mathcal G}_n,\pi_n\}$ and its inverse limit
$\widehat{\mathcal G}=\varprojlim{\mathcal G}_n$, which is a
profinite graph~\cite{ribes-z:profinite}. We devote our attention
to the discrete graph ${\mathcal G}=\varinjlim{\mathcal G}_n$.

The growth series of ${\mathcal G}$ can often be described as an infinite
product. We give such an expression for the graph in
Figure~\ref{figure:schreier}, making use of the fact that ${\mathcal G}$
``looks like a tree'' (even though it is amenable).

Consider the finite graphs ${\mathcal G}_n$; recall that ${\mathcal G}_n$
has $3^n$ vertices. Let $D_n$ be the diameter of ${\mathcal G}_n$ (maximal
distance between two vertices), and let
$\gamma_n=\sum_{i\in\N}\gamma_n(i)X^i$ be the growth series of ${\mathcal
G}_n$ (here $\gamma_n(i)$ denotes the number of vertices in ${\mathcal
G}_n$ at distance $i$ from the base point $*$).

The construction rule for ${\mathcal G}$ implies that ${\mathcal
G}_{n+1}$ can be constructed as follows: take three copies of
${\mathcal G}_n$, and in each of them mark a vertex $V$ at
distance $D_n$ from $*$. At each $V$ delete the loop labelled $s$,
and connect the three copies by a triangle labelled $s$ at the
three $V$'s. It then follows that $D_{n+1}=2D_n+1$, and
$\gamma_{n+1}=(1+2X^{D_n+1})\gamma_n$. Using the initial values
$\gamma_0=1$ and $D_0=0$, we obtain by induction
\[D_n=2^n-1,\qquad\gamma_n=\prod_{i=0}^{n-1}(1+2X^{2^i}).\]
We have also shown that the ball of radius $2^n$ around $*$
contains $3^n$ points, so the growth of ${\mathcal G}$ is at least
$n^{log_2(3)}$. But Proposition~\ref{prop:polygrowth} shows that
it is also an upper bound, and we conclude:
\begin{proposition}
  $\Gamma$ is an amenable $4$-regular graph whose growth function is
  transcendental, and admits the product decomposition
  \[\gamma(X)=\prod_{i\in\N}(1+2X^{2^i}).\]
  It is planar, and has polynomial growth of degree $\log_2(3)$.
\end{proposition}

Any graph is a metric space when one identifies each edge with a disjoint
copy of an interval $[0,L]$ for some $L>0$. We turn ${\mathcal G}_n$ in a
diameter-$1$ metric space by giving to each edge in ${\mathcal G}_n$ the
length $L=\operatorname{diam}({\mathcal G}_n)^{-1}$. The family
$\{{\mathcal G}_n\}$ then converges, in the following sense:

Let $A,B$ be closed subsets of the metric space $(X,d)$. For any
$\epsilon$, let $A_\epsilon=\{x\in X|\,d(x,A)\le\epsilon\}$, and
define the \emdef{Hausdorff distance}
\[d_X(A,B) = \inf\{\epsilon|\,A\subseteq B_\epsilon,B\subseteq
A_\epsilon\}.\]
This defines a metric on closed subsets of $X$. For
general metric spaces $(A,d)$ and $(B,d)$, define their
\emdef{Gromov-Hausdorff distance}
\[d^{GH}(A,B) = \inf_{X,i,j}d_X(i(A),j(B)),\]
where $i$ and $j$ are isometric embeddings of $A$ and $B$ in a metric
space $X$.

We may now rephrase the considerations above as follows: the sequence
$\{{\mathcal G}_n\}$ is convergent in the Gromov-Hausdorff metric. The
limit set ${\mathcal G}_\infty$ is a compact metric space.

The limit spaces are then: for $\Gg$ and $\Sg$, the limit ${\mathcal
G}_\infty$ is the interval $[0,1]$ (in accordance with its linear growth,
see Proposition~\ref{prop:polygrowth}). The limit spaces for $\FGg$,
$\BGg$ and $\GSg$ are fractal sets of dimension $\log_2(3)$.

\date{October 27, 2002}
\chapter{Spectral Properties of Unitary Representations}\label{chapter:spectrum}
We describe here some explicit computations of the spectrum of the
Schreier graphs defined in Section~\ref{sec:pspace}. The natural viewpoint
is that of spectra of representations, which we define now:
\begin{definition}
  Let $G$ be a group generated by a finite symmetric set $S$. The
  \emdef{spectrum} $\spec(\tau)$ of a representation $\tau:G\to\uhilb$
  with respect to the given set of generators is the spectrum of
  $\Delta_\tau=\sum_{s\in S}\tau(s)$ seen as an bounded operator on
  $\hilb$.
\end{definition}
(The condition that $S$ be symmetric ensures that $\spec(\tau)$ is a
subset of $\R$.)

Let $H$ be a subgroup of $G$. Then the spectrum of the ($\ell^2$-adjacency
operator of the) Schreier graph $\sch(G,H,S)$ is the spectrum of the
quasi-regular representation $\rho_{G/H}$ of $G$ in $\ell^2(G/H)$. This
establishes the connection with the previous chapter.

\section{Unitary representations}
Let $G$ act on the rooted tree $\tau$. Then $G$ also acts on the boundary
$\partial\tree$ of the tree. Since $G$ preserves the uniform measure on
this boundary, we have a unitary representation $\pi$ of $G$ in
$L^2(\partial\tree,\nu)$, where $\nu$ is the Bernoulli measure;
equivalently, we have a representation in $L^2([0,1],\text{Lebesgue})$.
Let $\hilb_n$ be the subspace of $L^2(\partial\tree,\nu)$ spanned by the
characteristic functions $\chi_\sigma$ of the rays $e$ starting by
$\sigma$, for all $\sigma\in Y^n$. It is of dimension $k_n$, and can
equivalently be seen as spanned by the characteristic functions in
$L^2([0,1],\text{Lebesgue})$ of intervals of the form $[(i-1)/k_n,i/k_n]$,
$1\le i\le k_n$. These $\hilb_n$ are invariant subspaces, and afford
representations $\pi_n=\pi_{|\hilb_n}$. As clearly $\pi_{n-1}$ is a
subrepresentation of $\pi_n$, we set $\pi^\perp_n=\pi_n\ominus\pi_{n-1}$,
so that $\pi=\oplus_{n=0}^\infty\pi^\perp_n$.

Let $P$ be a parabolic subgroup (see Section~\ref{sec:parabolic}), and set
$P_n=P\Stab_G(n)$.  Denote by $\rho_{G/P}$ the quasi-regular
representation of $G$ in $\ell^2(G/P)$ and by $\rho_{G/P_n}$ the
finite-dimensional representations of $G$ in $\ell^2(G/P_n)$, of degree
$k_n$. Since $G$ is level-transitive, the representations $\pi_n$ and
$\rho_{G/P_n}$ are unitary equivalent.

The $G$-spaces $G/P$ are of polynomial growth when the conditions of
Proposition~\ref{prop:polygrowth} are fulfilled, and therefore, according
to the criterion of F{\o}lner given in~Theorem~\ref{theorem:folner}, they
are amenable.

The following result belongs to the common lore:
\begin{proposition}\label{prop:burger}
  Let $H$ be a subgroup of $G$. Then the quasi-regular representation
  $\rho_{G/H}$ is weakly contained in $\rho_G$ if and only if
  $H$ is amenable.
\end{proposition}
(``Weakly contained'' is a topological extension of ``contained''; it
implies for instance inclusion of spectra.)

\begin{theorem}
  Let $G$ be a group acting on a regular rooted tree, and let $\pi$,
  $\pi_n$ and~$\pi^\perp_n$ be as above.
  \begin{enumerate}
  \item If $G$ is weakly branch, then $\rho_{G/P}$ is an irreducible
    representation of infinite dimension.
  \item $\pi$ is a reducible representation of infinite dimension
    whose irreducible components are precisely those of the $\pi^\perp_n$
    (and thus are all finite-dimensional). Moreover
    \[\spec(\pi)=\overline{\bigcup_{n\ge0}\spec(\pi_n)}=\overline{\bigcup_{n\ge0}\spec(\pi^\perp_n)}.\]
  \item The spectrum of $\rho_{G/P}$ is contained in
    $\overline{\cup_{n\ge0}\spec(\rho_{G/P_n})}=\overline{\cup_{n\ge0}\spec(\pi_n)}$,
    and thus is contained in the spectrum of $\pi$. If moreover either
    $P$ or $G/P$ are amenable, these spectra coincide, and if $P$ is
    amenable, they are contained in the spectrum of $\rho_G$:
    \[\spec(\rho_{G/P})=\spec(\pi)=\overline{\bigcup_{n\ge0}\spec(\pi_n)}\subseteq\spec(\rho_G).\]
  \item $\Delta_\pi$ has a pure-point spectrum, and its spectral radius
    $r(\Delta_\pi)=s\in\R$ is an eigenvalue, while the spectral
    radius $r(\Delta_{\rho_{G/P}})$ is not an eigenvalue of
    $\Delta_{\rho_{G/P}}$. Thus $\Delta_{\rho_{G/P}}$ and $\Delta_\pi$
    are different operators having the same spectrum.
  \end{enumerate}
\end{theorem}

\subsection{Can one hear a representation?}\label{subs:kac}
We end by turning to a question of Mark Ka\'c~\cite{kac:drum}:
``Can one hear the shape of a drum?'' This question was answered
in the negative in~\cite{gordon-w-w:drum}, and we here answer by
the negative to a related question: ``Can one hear a
representation?'' Indeed $\rho_{G/P}$ and $\pi$ have same spectrum
(i.e., cannot be distinguished by hearing), but are not
equivalent.  Furthermore, if $G$ is a branch group, there are
uncountably many nonequivalent representations within
$\{\rho_{G/\Stab_G(e)}|\,e\in\partial\tree\}$, as is shown
in~\cite{bartholdi-g:parabolic}.

The same question may be asked for graphs: ``are there two non-isomorphic
graphs with same spectrum?'' There are finite examples, obtained through
the notion of \emph{Sunada pair}~\cite{lubotzky:cayley}. C\'edric
B\'eguin, Alain Valette and Andrzej \.Zuk produced the following example
in~\cite{beguin-:spectrum}: let $\Gamma$ be the integer Heisenberg group
(free $2$-step nilpotent on $2$ generators $x,y$). Then
$\Delta=x+x^{-1}+y+y^{-1}$ has spectrum $[-2,2]$, which is also the
spectrum of $\Z^2$ for an independent generating set. As a consequence,
their Cayley graphs have same spectrum, but are not quasi-isometric (they
do not have the same growth).

Using the result of Nigel Higson and Gennadi Kasparov~\cite{higson-k:bc}
(giving a partial positive answer to the Baum-Connes conjecture), we may
infer the following

\begin{proposition}\label{prop:bc}
  Let $\Gamma$ be a torsion-free amenable group with finite generating
  set $S=S^{-1}$ such that there is a map $\phi:\Gamma\to\Z/2\Z$ with
  $\phi(S)=\{1\}$.  Then
  \[\spec(\sum_{s\in S}\rho(s))=[-|S|,|S|].\]
\end{proposition}

In particular, there are countably many non-quasi-isometric graphs with
the same spectrum, including the graphs of $\Z^d$, of free nilpotent
groups and of suitable torsion-free groups of intermediate growth (for the
first examples, see~\cite{grigorchuk:pgps}).

In contrast to Proposition~\ref{prop:bc} stating that the spectrum of the
regular representation is an interval, the spectra of the representation
$\pi$ may be totally disconnected. The first two authors prove
in~\cite{bartholdi-g:spectrum} the following

\begin{theorem}\label{thm:spectra}
  For $\lambda\in\R$, define
  \[J(\lambda)=\overline{\pm\sqrt{\lambda\pm\sqrt{\lambda\pm\sqrt{\lambda\pm\sqrt{\dots}}}}}\]
(note the closure operator written in bar notation on the top).
  Then we have the following spectra:
  \[\begin{array}{r|l}
    \mathbf{G} & \mathbf{\spec(\pi)}\\ \hline
    \Gg & [-2,0]\cup[2,4]\\
    \Sg & [0,4]\\
    \FGg & \{4,1\}\cup 1+J(6)\\
    \BGg,\GSg & \{4,-2,1\}\cup 1\pm\sqrt{\frac92\pm 2J(\frac{45}{16})}
  \end{array}\]
\end{theorem}

In particular, the spectrum of the graph in Figure~\ref{figure:schreier}
is the totally disconnected set $\{4,1\}\cup1+J(6)$, which is the union of
a Cantor set of null Lebesgue measure and a countable collection of
isolated points.

\section{Operator recursions}\label{subs:rec}
We now describe the computations of the spectra for the examples in
Theorem~\ref{thm:spectra}, which share the property of acting on an
$m$-regular tree $\tree$.  Let $\hilb$ be an infinite-dimensional Hilbert
space, and suppose $\Phi:\hilb\to\hilb\oplus\dots\oplus\hilb$ is an
isomorphism, where the domain of $\Phi$ is a sum of $m\ge2$ copies of
$\hilb$. Let $S$ be a finite subset of $\uhilb$, and suppose that for all
$s\in S$, if we write $\Phi^{-1}s\Phi$ as an operator matrix
$(s_{i,j})_{i,j\in\{1,\dots,d\}}$ where the $s_{i,j}$ are operators in
$\hilb$, then $s_{i,j}\in S\cup\{0,1\}$.

This is precisely the setting in which we will compute the spectra of our
five example groups: for $\Gg$, we have $m=2$ and $S=\{a,b,c,d\}$ with
\begin{gather*}
  a = \begin{pmatrix}0&1\\1&0\end{pmatrix},\qquad b = \begin{pmatrix}a&0\\0&c\end{pmatrix},\\
  c = \begin{pmatrix}a&0\\0&d\end{pmatrix},\qquad d = \begin{pmatrix}1&0\\0&b\end{pmatrix}.
\end{gather*}

For $\Sg$, we also have $m=2$, and $S=\{a,\tilde b,\tilde c,\tilde d\}$
given by
\begin{gather*}
  b = \begin{pmatrix}a&0\\0&c\end{pmatrix},\qquad c = \begin{pmatrix}1&0\\0&d\end{pmatrix},\qquad d = \begin{pmatrix}1&0\\0&b\end{pmatrix}.
\end{gather*}

For $\FGg=\langle a,s\rangle$, $\BGg=\langle a,t\rangle$ and $\GSg=\langle
a,r\rangle$, we have $m=3$ and

\begin{gather*}
  a = \begin{pmatrix}0&1&0\\0&0&1\\1&0&0\end{pmatrix},\qquad
  s = \begin{pmatrix}a&0&0\\0&1&0\\0&0&s\end{pmatrix},\\
  t = \begin{pmatrix}a&0&0\\0&a&0\\0&0&t\end{pmatrix},\qquad
  r = \begin{pmatrix}a&0&0\\0&a^2&0\\0&0&r\end{pmatrix}.
\end{gather*}

Each of these operators is unitary. The families $S=\{a,b,c,d\},\dots$
generate subgroups $G(S)$ of $\uhilb$.  The choice of the isomorphism
$\Phi$ defines a unitary representation of $\langle S\rangle$.

We note, however, that the expression of each operator as a matrix of
operators does not uniquely determine the operator, in the sense that
different isomorphisms $\Phi$ can yield non-isomorphic operators
satisfying the same recursions. Even if $\Phi$ is fixed, it may happen
that different operators satisfy the same recursions. We considered two
types of isomorphisms: $\hilb=\ell^2(G/P)$, where $\Phi$ is derived from
the decomposition map $\psi$; and $\hilb=L^2(\partial\tree)$, where
$\Phi:\hilb\to\hilb^Y$ is defined by
$\Phi(f)(\sigma)=(f(0\sigma),\dots,f((m-1)\sigma))$, for $f\in
L^2(\partial\tree)$ and $\sigma\in\partial\tree$. There are actually
uncountably many non-equivalent isomorphisms, giving uncountably many
non-equivalent representations of the same group, as indicated in
Subsection~\ref{subs:kac}.

\subsection{The spectrum of $\Gg$}
For brevity we shall only describe the spectrum of $\pi$ for the first
Grigorchuk group. Details and computations for the other examples appear
in~\cite{bartholdi-g:spectrum}.

Denote by $a_n, b_n, c_n, d_n$ the permutation matrices of the
representation $\pi_n=\rho_{G/P_n}$. We have
\begin{gather*}
  a_0 = b_0 = c_0 = d_0 = (1),\\
  a_n = \begin{pmatrix}0&1\\1&0\end{pmatrix},\qquad b_n = \begin{pmatrix}a_{n-1}&0\\0&c_{n-1}\end{pmatrix},\\
  c_n = \begin{pmatrix}a_{n-1}&0\\0&d_{n-1}\end{pmatrix},\qquad d_n = \begin{pmatrix}1&0\\0&b_{n-1}\end{pmatrix}.
\end{gather*}

The Hecke-Laplace operator of $\pi_n$ is
\[\Delta_n = a_n + b_n + c_n + d_n = \begin{pmatrix}2a_{n-1}+1&1\\1&\Delta_{n-1}-a_{n-1}\end{pmatrix},\]
and we wish to compute its spectrum.  We start by proving a slightly
stronger result: define
\[Q_n(\lambda,\mu) = \Delta_n - (\lambda+1)a_n - (\mu+1)\]
and
\begin{align*}
  \Phi_0 &= 2-\mu-\lambda,\\
  \Phi_1 &= 2-\mu+\lambda,\\
  \Phi_2 &= \mu^2-4-\lambda^2,\\
  \Phi_n &= \Phi_{n-1}^2 - 2(2\lambda)^{2^{n-2}}\qquad(n\ge3).
\end{align*}

Then the following steps compute the spectrum of $\pi_n$: first, for
$n\ge2$, we have
\[|Q_n(\lambda,\mu)| = (4-\mu^2)^{2^{n-2}}\left|Q_{n-1}\left(\frac{2\lambda^2}{4-\mu^2},\mu+\frac{\mu\lambda^2}{4-\mu^2}\right)\right|\qquad(n\ge2).\]
Therefore, for all $n$ we have
\[|Q_n| = \Phi_0\Phi_1\cdots\Phi_n.\]
Then, for all $n$ we have
\begin{multline*}
  \{(\lambda,\mu):\,Q_n(\lambda,\mu)\text{ non invertible}\} = \{(\lambda,\mu):\,\Phi_0(\lambda,\mu)=0\}\cup \{(\lambda,\mu):\,\Phi_1(\lambda,\mu)=0\}\\
  {}\cup\{(\lambda,\mu):\,4-\mu^2+\lambda^2+4\lambda\cos\left(\frac{2\pi j}{2^n}\right)=0\text{ for some }j = 1,\dots,2^{n-1}-1\}.
\end{multline*}

In the $(\lambda,\mu)$ system, the spectrum of $Q_n$ is thus a collection
of $2$ lines and $2^{n-1}-1$ hyperbol\ae.

\begin{figure}
\begin{center}
\setlength\unitlength{1pt}
\begin{picture}(360,300)
\put(180,270){\makebox(0,0)[lb]{\smash{\normalsize $\mu$}}}
\put(310,130){\makebox(0,0)[lb]{\smash{\normalsize $\lambda$}}}
\put(0,-40){\includegraphics{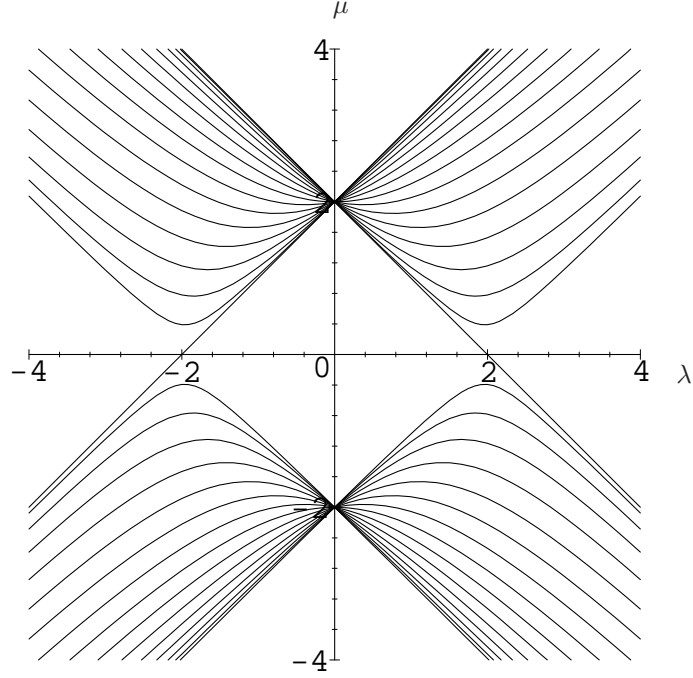}}
\end{picture}
\end{center}
\caption{The spectrum of $Q_n(\lambda,\mu)$ for $\Gg$}
\end{figure}

The spectrum of $\Delta_n$ is precisely the set of $\theta$ such that
$|Q_n(-1,\theta-1)|=0$. From the computations above we obtain

\begin{proposition}
\[\spec(\Delta_n) = \{1\pm\sqrt{5-4\cos\phi}:\,\phi\in2\pi\Z/2^n\}\setminus\{0,-2\}.\]
  Therefore the spectrum of $\pi$, for the group $\Gg$, is
  \[\spec(\Delta) = [-2,0]\cup[2,4].\]
\end{proposition}

The first eigenvalues of $\Delta_n$ are $4$; $2$; $1\pm\sqrt5$;
$1\pm\sqrt{5\pm2\sqrt2}$; $1\pm\sqrt{5\pm2\sqrt{2\pm\sqrt2}}$; etc.

The numbers of the form
$\pm\sqrt{\lambda\pm\sqrt{\lambda\pm\sqrt{\dots}}}$ appear as preimages of
the quadratic map $z^2-\lambda$, and after closure produce a Julia set for
this map (see~\cite{barnsley:fe}).  In the given example this Julia set is
just an interval. For the groups $\FGg,\BGg,\GSg$, however, the spectra
are simple transformations of Julia sets which are totally disconnected,
as similar computations show.

\date{October 27, 2002}
\chapter{Open Problems}
This chapter collects questions for which no answer is yet known. Here we
use the geometric definition of branch group given in
Definition~\ref{defn:gbranch}, except in
Question~\ref{question:kernel-action}.

The theory of branch group is a recent development and the
questions arise and die almost every day, so there are no
longstanding problems in the area and there is no easy way to say
which questions are difficult and which are not. The list below
just serves as a list of problems that one should naturally ask at
this given moment. We kindly invite the interested readers to
solve as many of the problems below.

\textbf{Added in the proof.} The authors had a chance to proofread
the article 18 months after the initial submission. Some of the
proposed problems have been solved in the meantime and we include
appropriate comments to that effect.

\begin{question}
Is there a finitely generated fractal regular branch group $G$,
branched over $K$, acting on the binary tree, such that the index
of the geometric embedding of $K \times K$ into $K$ is two?
\end{question}

\begin{question}\label{question:kernel-action}
Every branch group from Definition~\ref{defn:branch} acts canonically on
the tree determined by its branch structure as a group of tree
automorphisms. Is the kernel of this action necessarily central?
\end{question}

\begin{question}
Does every finitely generated branch $p$-group, where $p$ is a prime,
satisfy the congruence subgroup property?
\end{question}

\begin{question}\label{question:neumann-spinal}
Is every finitely generated branch group isomorphic to a spinal group?
\end{question}

\begin{question}\label{question:word-conj}
Is the conjugacy problem solvable in all branch groups with solvable word
problem?
\end{question}

\begin{question}\label{question:non-isomorphic}
When do the defining triples $\om$ and $\om'$ define non-isomorphic
examples of branch groups in \cite{grigorchuk:gdegree,grigorchuk:pgps} and
more generally in \GG\ and \GGS\ groups?
\end{question}

\begin{question}\label{question:L-presentations}
Which branch groups have finite $L$-presentations? Finite
ascending $L$-presentations? In particular, what is the status of
the Gupta-Sidki 3-group?
\end{question}

All that is known at present is that the Gupta-Sidki 3-group has a
finite endomorphic presentation.

\begin{question}
Do there exist finitely presented branch groups?
\end{question}

\begin{question}
Is it correct that there are no finitely generated hereditary
just-infinite torsion groups?
\end{question}

\begin{question}
Is every finitely generated just-infinite group of intermediate growth
necessarily a branch group?
\end{question}

\begin{question}
Is every finitely generated hereditarily just-infinite group necessarily
linear?
\end{question}

\begin{question}
Recall that $G$ has bounded generation if there exist elements
$g_1\dots,g_k$ in $G$ such that every element in $G$ can be
written as $g_1^{n_1}\dots g_k^{n_k}$ for some
$n_1,\dots,n_k\in\Z$. Can a just-infinite branch group have
bounded generation? Can infinite simple group have bounded
generation?
\end{question}

\begin{question}\label{question:height}
  What is the height (in the sense of~\cite{pride:large}, see
  Chapter~\ref{chapter:just-infinite}) of a Grigorchuk $2$-group
  $G_\om$ when the defining sequence $\omover$ is not periodic?
  Same question for arbitrary \GG\ groups. In particular, can the
  height be infinite?
\end{question}

\begin{question}
Is every maximal subgroup in a finitely generated branch group necessarily
of finite index?
\end{question}

\begin{question}
Is there is a finitely generated branch group containing the free
group $F_2$ on two generators?
\end{question}

Positive answer is provided by Said Sidki and John Wilson
in~\cite{sidki-w:freebranch}

\begin{question}
Are there finitely generated branch groups with exponential growth
that do not contain the free group $F_2$?
\end{question}

\begin{question}
Is there a finitely generated branch group whose degree of growth
is $e^{\sqrt{n}}$? Is there such a group in the whole class of
finitely generated groups?
\end{question}

\begin{question}
What is the exact degree of growth of any of the basic examples of regular
branch groups (for example $\Gg$, $\FGg$, $\BGg$, $\dots$)?
\end{question}

\begin{question}
What is the growth of the Brunner-Sidki-Vieira group
(see~\cite{brunner-s-v:nonsolvable} and
Section~\ref{sec:l-examples})?
\end{question}

It is known that the Brunner-Sidki-Vieira group does not contain
any non-abelian free groups, but it is not known whether it
contains a non-abelian free monoid. Note that this group is not a
branch group, but it is a weakly branch group.

\begin{question}
Are there finitely generated non-amenable branch groups not
containing the free group $F_2$?
\end{question}

\begin{question}
Is it correct that in each finitely generated fractal branch group
$G$ every finitely generated subgroup is either finite or Pride
equivalent with $G$?
\end{question}

A stronger property holds for $\Gg$, namely, \JWilson\ and the
second author have proved in~\cite{grigorchuk-w:rigidity} that
every finitely generated subgroup of $\Gg$ is either finite or
commensurable with $\Gg$. \CRover\ has announced that the answer
is also positive for the Gupta-Sidki group $\GSg$. However, the
answer is negative in general.

\begin{question}
Do there exist branch groups with the property (T)?
\end{question}

\bibliography{mrabbrev,people,sunik,bartholdi,grigorchuk,math}
\nocite{grigorchuk:onbranch}
\end{document}